\documentclass[10pt,leqno]{article}

\usepackage[francais,english]{babel}
\usepackage{amsmath,amssymb,amscd}
\usepackage{a4wide}
\usepackage[latin1]{inputenc}
\usepackage{entete-hanoi}
\usepackage[all]{xy} 
\usepackage{fancyhdr}
\usepackage[T1]{fontenc}
\usepackage{hyperref}
\usepackage{xcolor}

\usepackage{mathdots}

\title{The global Gan-Gross-Prasad conjecture for unitary groups: the endoscopic case} 
\author{Raphaël Beuzart-Plessis, Pierre-Henri Chaudouard and Micha\l \ Zydor}
\date{}

\begin{document}
\selectlanguage{english}
\maketitle

\begin{abstract}
  In this paper, we prove the Gan-Gross-Prasad conjecture and the Ichino-Ikeda conjecture for unitary groups $U_n\times U_{n+1}$ in all the endoscopic cases. Our main technical innovation  is the computation of the contributions of certain cuspidal data, called $*$-generic, to the Jacquet-Rallis trace formula for linear groups. We offer two different computations of these contributions: one, based on  truncation, is expressed in terms of regularized Rankin-Selberg periods of Eisenstein series and Flicker-Rallis intertwining periods. The other, built upon Zeta integrals, is expressed in terms of functionals on the Whittaker model. A direct proof of the equality between the two expressions is also given. Finally several useful auxiliary results about the spectral expansion of the Jacquet-Rallis trace formula are provided.
\end{abstract}

\tableofcontents

\section{Introduction}\label{sec:Introduction}

\subsection{The endoscopic cases of the Gan-Gross-Prasad conjecture}

\begin{paragr}
  One of the main motivation of the paper is the obtention of the remaining cases, the so-called ``endoscopic cases'', of the Gan-Gross-Prasad and the Ichino-Ikeda conjectures for unitary groups. To begin with, we shall give  the main statements we prove.
\end{paragr}

\begin{paragr}
  Let  $E/F$ be a quadratic extension of number fields and $c$ be the non-trivial element of the Galois group $\Gal(E/F)$. Let $\AAA$ be the ring of adèles of $F$. Let $n\geq 1$ be an integer. Let $\hc_n$ be the set of isomorphism classes of non-degenerate $c$-Hermitian spaces $h$ over $E$ of rank $n$. For any  $h_n\in \hc_n$, we  identify $h_n$ with a representative and we shall denote by $U(h_n)$ its automorphisms group. Let $h_0\in \hc_1$ be the element of rank $1$ given by the norm $N_{E/F}$. 

We attach to any $h\in \hc_n$ the following algebraic groups over $F$:
\begin{itemize}
\item the unitary group $U'_h$ of automorphisms of $h$;
\item the product of unitary groups $U_h=U(h)\times U(h\oplus h_0)$ where $h\oplus h_0$ denoted the orthogonal sum.
\end{itemize}
We have an obvious diagonal embedding $U'_h\hookrightarrow U_h$.
\end{paragr}

\begin{paragr}[Arthur parameter.] --- \label{S:Aparam}Let $G_n$ be the group of automorphims of the $E$-vector space $E^n$. We view $G_n$ as an $F$-group by Weil restriction. By a {\em Hermitian Arthur parameter}\footnote{Strictly speaking, it is a \emph{discrete} Arthur parameter. By simplicity, we shall omit the adjective \emph{discrete}.} of $G_n$, we mean an irreducible automorphic representation $\Pi$ for which there exists a partition $n_1+\ldots+n_r=n$ of $n$ and for any $1\leq i \leq r$ a cuspidal automorphic representation $\Pi_i$ of $G_{n_i}(\AAA)$ such that
  \begin{enumerate}
  \item each $\Pi_i$ is conjugate self-dual and the Asai $L$-function $L(s,\Pi_i,As^{(-1)^{n+1}})$ has a pole at $s=1$;
  \item the representations $\Pi_i$ are mutually non-isomorphic for $1\leq i\leq r$;
\item the representation $\Pi$ is isomorphic to the full induced representation $\Ind_P^{G_n}(\Pi_1\boxtimes \ldots \boxtimes \Pi_r)$ where $P$ is a parabolic subgroup of $G_n$ of Levi factor $G_{n_1}\times\ldots \times G_{n_r}$.
  \end{enumerate}

  \begin{remarque}
    It is well-known (see  \cite{Flicker}) that condition 1 above is equivalent to the fact that $\Pi_i$ is  $(GL_{n_i,F},\eta^{n+1})$-distinguished in the sense of §\ref{S:disting} below.
  \end{remarque}

The integer $r$ and the representations $(\Pi_i)_{1\leq i \leq r}$ are unique (up to a permutation). We set $S_\Pi=(\mathbb{Z}/2\mathbb{Z})^r$.

Let $G=G_n\times G_{n+1}$. By a {\em Hermitian Arthur parameter} of $G$, we mean an automorphic representation of the form $\Pi=\Pi_n\boxtimes \Pi_{n+1}$ where $\Pi_i$ is a Hermitian Arthur parameter of $G_i$ for $i=n,n+1$. For such a Hermitian Arthur parameter, we set $S_\Pi=S_{\Pi_n}\times S_{\Pi_{n+1}}$.
\end{paragr}

\begin{paragr}
  Let $h\in \hc_n$ and $\sigma$ be a cuspidal automorphic representation of $U_h(\AAA)$. We say that a Hermitian Arthur parameter $\Pi$ of $G_n$ is a {\em weak base-change} of $\sigma$ if for almost all places of $F$ that split in $E$, the local component $\Pi_v$ is the split local base change of $\sigma_v$. If this is the case, we write $\Pi=BC(\sigma)$.
  
  \begin{remarque}
  By the work of Mok \cite{Mok} and Kaletha-Minguez-Shin-White \cite{KMSW}, we know that if $\sigma$ admits a weak base-change then it admits a {\em strong base-change} that is a Hermitian Arthur parameter $\Pi$ of $G_n$ such that $\Pi_v$ is the base-change of $\sigma_v$ for every place $v$ of $F$ (where the local base-change in the ramified case is also constructed in \textit{loc. cit.} and characterized by certain local character relations). Moreover, this is the case if and only if $\sigma$ has a {\em generic Arthur parameter} in the sense of \textit{loc. cit.} Besides, a result of Ramakrishnan \cite{Ram2} implies that a weak base-change is automatically a strong base-change. Therefore, we could have used the notion of strong base-change instead. However, we prefer to stick with the terminology of weak base-change in order to keep the statement of the next theorem independent of \cite{Mok} and \cite{KMSW}.
  \end{remarque}
  \end{paragr}

  \begin{paragr}[Gan-Gross-Prasad conjecture.] --- \label{S-intro:GGP}
  Our first main result is the global Gan-Gross-Prasad conjecture \cite[Conjecture 24.1]{GGP} in the case of $U(n)\times U(n+1)$ and can be stated as follows.

    \begin{theoreme}\label{thm:GGP}
      Let $\Pi$ be a Hermitian Arthur parameter of $G$. The following two statements are equivalent:
      \begin{enumerate}
      \item The complete Rankin-Selberg $L$-function of $\Pi$ (including Archimedean places) satisfies
      $$L(\frac12,\Pi)\not=0;$$
      \item There exists $h\in \hc_n$ and an irreducible cuspidal automorphic subrepresentation $\sigma$ of $U_h$ such that $\Pi$ is a weak base change of $\sigma$ and the period integral $\pc_h$ defined by 
        \begin{align*}
          \pc_h(\varphi)= \int_{[U'_h]} \varphi(h)\, dh
        \end{align*}
induces a non-zero linear form on the space of $\sigma$.
      \end{enumerate}
    \end{theoreme}

    \begin{remarque}
      If the Arthur parameter is moreover simple (that is if $\Pi$ is cuspidal), the theorem is proved by Beuzart-Plessis-Liu-Zhang-Zhu (cf. \cite[Theorem 1.7]{BPLZZ}). Previous works had to assume extra local hypothesis on $\Pi$, which implied that $\Pi$ was also simple (see \cite{Z1}, \cite{Xue}, \cite{RBP} and  \cite{BPPlanch}) or only proved the direction $2.\Rightarrow 1.$ of the theorem (\cite{GJR}, \cite{IYunit}, \cite{JZ}).
    \end{remarque}

As observed in \cite[Theorem 1.2]{Z1} and \cite[Theorem 1.8]{BPLZZ} we can deduce from Theorem \ref{thm:GGP} the following statement (whose proof is word for word that of \cite{Z1}):

\begin{theoreme}
  Let $\Pi_{n+1}$ be a Hermitian Arthur parameter of $G_{n+1}$. Then there exists a simple Hermitian Arthur parameter $\Pi_n$ of $G_{n}$ such that the Rankin-Selberg $L$-function satisfies: 
$$L(\frac12,\Pi_n\times \Pi_{n+1})\not=0.$$
\end{theoreme}
  \end{paragr}

  \begin{paragr}[Ichino-Ikeda conjecture.] --- Let $\sigma=\bigotimes_v'\sigma_v$ be an irreducible cuspidal automorphic representation of $U_h$ that is tempered everywhere in the following sense: for every place $v$, the local representation $\sigma_v$ is tempered. By \cite{Mok} and \cite{KMSW}, $\sigma$ admits a weak (hence a strong) base-change $\Pi$ to $G$. Set
  $$\displaystyle \mathcal{L}(s,\sigma)=\prod_{i=1}^{n+1}L(s+i-1/2,\eta^i)\frac{L(s,\Pi)}{L(s+1/2,\sigma,\Ad)}$$
  where $\eta$ denotes the quadratic idele class character associated to the extension $E/F$, $L(s,\eta^i)$ is the completed Hecke $L$-function associated to $\eta^i$ and $L(s,\sigma,\Ad)$ is the completed adjoint $L$-function of $\sigma$ (defined using the local Langlands correspondence for $G$ from \cite{Mok}, \cite{KMSW}). We denote by $\mathcal{L}(s,\sigma_v)$ the corresponding quotient of local $L$-factors. For each place $v$ of $F$, we define a {\em local normalized period} $\pc^\natural_{h,\sigma_v}: \sigma_v\times \sigma_v\to \CC$ as follows. It depends on the choice of a Haar measure on $U'_h(F_v)$ as well as an invariant inner product $(.,.)_v$ on $\sigma_v$ and is given by
  $$\displaystyle \pc^\natural_{h,\sigma_v}(\varphi_v,\varphi'_v)=\mathcal{L}(\frac12,\sigma_v)^{-1}\int_{U'_h(F_v)} (\sigma_v(h_v)\varphi_v,\varphi'_v)_v dh_v,\;\; \varphi_v,\varphi'_v\in \sigma_v,$$
  where, thanks to the temperedness assumption, the integral is absolutely convergent \cite[Proposition 2.1]{NHar} and the local factor $\mathcal{L}(s,\sigma_v)$ has no zero (nor pole) at $s=\frac12$. Moreover, by \cite[Theorem 2.12]{NHar}, if $\varphi=\otimes_v' \varphi_v\in \sigma$, then for almost all places $v$ we have
  \begin{align}\label{intro:localunrcomp}
  \pc^\natural_{h,\sigma_v}(\varphi_v,\varphi_v)=\vol(U'_h(\oc_v))(\varphi_v,\varphi_v)_v.
  \end{align}
  We also recall that the global representation $\sigma$ has a natural invariant inner product given by
  $$\displaystyle (\varphi,\varphi)_{\Pet}=\int_{[U_h]} \lvert \varphi(g)\rvert^2 dg,\;\; \varphi\in \sigma.$$
  
  Our second main result is the global Ichino-Ikeda conjecture for unitary groups formulated in \cite[Conjecture 1.3]{NHar} and can be stated as follows (this result can be seen as a refinement of Theorem \ref{thm:GGP}, the precise relation requiring the local Gan-Gross-Prasad conjecture and Arthur's multiplicity formula for unitary groups will not be discussed here).
  
  \begin{theoreme}\label{thm:II}
  Assume that $\sigma$ is a cuspidal automorphic representation of $U_h$ that is tempered everywhere and let $\Pi=\Pi_n\boxtimes \Pi_{n+1}$ be the weak (hence the strong) base-change of $\sigma$ to $G$. Suppose that we normalize the period integral $\pc_h$ and the Peterssen inner product $(.,.)_{\Pet}$ by choosing the invariant Tamagawa measures\footnote{We warn the reader that our convention is to include the global normalizing $L$-values in the definition of Tamagawa measures, cf. Section \ref{S:Haar-measures} for precise definitions.} $d_{\Tam}h$ and $d_{\Tam}g$ on $U'_h(\AAA)$ and $U_h(\AAA)$ respectively. Assume also that the local Haar measures $dh_v$ on $U'_h(F_v)$ factorize the Tamagawa measure: $d_{\Tam}h=\prod_v dh_v$. Then, for every nonzero factorizable vector $\varphi=\otimes_v' \varphi_v\in \sigma$, we have
  \begin{align*}
  \frac{\lvert \pc_h(\varphi)\rvert^2}{(\varphi, \varphi)_{\Pet}}=\lvert S_\Pi\rvert^{-1} \mathcal{L}(\frac12,\sigma)\prod_v \frac{\pc^\natural_{h,\sigma_v}(\varphi_v,\varphi_v)}{(\varphi_v,\varphi_v)_v}
  \end{align*}
  where we recall that $S_\Pi$ denotes the finite group $S_{\Pi_n}\times S_{\Pi_{n+1}}$.
  \end{theoreme}
  
  Note that the product over all places in the theorem is well-defined by \eqref{intro:localunrcomp}. Moreover, once again, this theorem is proved in \cite{BPLZZ} under the extra assumption that $\Pi$ is cuspidal (in which case $\lvert S_\Pi\rvert=4$). Previous results in that direction includes \cite{Zhang2}, \cite{RBP}, \cite{BPPlanch} where some varying local assumptions on $\sigma$ entailing the cuspidality of $\Pi$ were imposed. In a slightly different direction, the paper \cite{GL} establishes the above identity up to an unspecified algebraic number under some arithmetic assumptions on $\sigma$.
      
  \end{paragr}

\subsection{The spectral expansion of the Jacquet-Rallis trace formula for the linear groups}
\label{ssec:intro-spec-exp}
\begin{paragr}[Motivations.] --- 
As in \cite{Z1}, \cite{Zhang2}, \cite{Xue}, \cite{RBP}, \cite{BPPlanch} and \cite{BPLZZ}, our proofs of Theorems \ref{thm:GGP} and \ref{thm:II} follow the strategy of Jacquet and Rallis \cite{JR} and are thus based on a comparison of {\em relative trace formulas} on unitary groups $U_h$ for $h\in \hc_n$ and the group $G$. Let's recall that these trace formulas have two different expansions:  one, called the geometric side, in terms of distributions indexed by geometric classes and the other, called the spectral side, in terms  of distributions indexed by cuspidal  data. As usual, the point is to get enough test functions to first compare the geometric sides which gives a comparison of spectral sides.

For specific test functions, the trace formula boils down to a simple and quite easy  equality between a sum of relative regular orbital integrals and a sum of  relative characters attached to cuspidal representations. This is the simple trace formula used by Zhang in \cite{Z1} and \cite{Zhang2} to  prove special cases of Theorems \ref{thm:GGP} and \ref{thm:II}.  In return one has to impose restrictive local conditions on  the representations one considers. 

In \cite{Z0}, \cite{Z2}, \cite{Z3}, Zydor established  general Jacquet-Rallis  trace formulas. Besides, in \cite{CZ}, Chaudouard-Zydor proved the comparison of all the geometric terms for matching test functions, that is functions  with matching local orbital integrals. Using these results, Beuzart-Plessis-Liu-Zhang-Zhu in \cite{BPLZZ} proved \ref{thm:GGP} and \ref{thm:II} when $\Pi$ is cuspidal. Their main innovation  is a construction of  Schwartz test functions only detecting certain cuspidal data. In this way, they were able to construct matching test functions for which the spectral expansions reduce to some relative characters attached to cuspidal representations.
\end{paragr}

\begin{paragr} In this paper, we also want to use the construction of Beuzart-Plessis-Liu-Zhang-Zhu. But for this,  we need  two extra ingredients. First we need the slight extension of  Zydor's work to the    space of  Schwartz test functions.  For the geometric sides, this was done in \cite{CZ}.  For the test functions we need, the spectral side of the trace formulas for  unitary groups still reduces   to relative characters attached to cuspidal representations and we need nothing more. But,  for the group $G$, we  shall extend the spectral side of  the trace formula  to the space of Schwartz functions.  Second there is an even more serious question: since the representation $\Pi$ is no longer assumed to be  cuspidal, the spectral contribution  associated to $\Pi$ is much more involved. In this section, we shall explain alternative and somewhat  more tractable expressions for the spectral contributions in the trace formula for $G$. For the specific cuspidal datum attached to $\Pi$, we get a precise result as we shall see in section \ref{ssec:intro-generic} below.
\end{paragr}

\begin{paragr}[The spectral expansion for the Schwartz space.] --- Let $\Xgo(G)$ be the set of cuspidal data of $G$ (see § \ref{S:cuspidal-data}).  To $\chi$ is associated a direct invariant factor $L^2_\chi([G])$ of $L^2([G])$ (see \cite[Chap. II]{MW} or section \ref{ssec:cuspidal-data} for a review). Let $f$ be a function in the  Schwartz space $\Sc(G(\AAA))$ (cf. \S \ref{S:Schw} for a definition). Let  $K_f$, resp. $K_{f,\chi}$, be  the kernel associated to the action by right convolution of $f$ on $L^2([G])$, resp. $L^2_\chi([G])$. 

Following \cite{Z3} (see §\ref{S:KTchi}), we introduce the modified kernel $K_{f, \chi}^T$ depending on a parameter $T$ in a certain real vector space. Set $H=G_n$ and $G'=\GL_{n,F}\times \GL_{n+1,F}$ both seen as subgroups of $G$ (the embedding $H\hookrightarrow G$ being the ``diagonal'' one where the inclusion $G_n \hookrightarrow G_{n+1}$ is induced by the identification of  $E^n$  with the hyperplane of $E^{n+1}$ of vanishing last coordinate). The following theorem is an extension to Schwartz functions of \cite[théorème 0.1]{Z3}.

\begin{theoreme}\label{intro:jfDef}(see theorem \ref{thm:jfDef}) 
  \begin{enumerate}
  \item For any $T$ in a certain positive  Weyl chamber, we have
\[
\sum_{\chi \in \Xgo(G)} \int_{[H]}\int_{[G']}| K^T_{f,\chi}(h,g') \eta_{G'}(g')| \, dg' dh < \infty
\]
\item Let  $\eta_{G'}$  be the quadratic  character  of $G'(\AAA)$  defined in  \S \ref{S:character etaG'}. For each $\chi \in \Xgo(G)$, the integral
\begin{align}\label{intro:polynexp}
\int_{[H]} \int_{[G']} K^T_{f,\chi}(h,g') \eta_{G'}(g') \, dg' dh
\end{align}
coincides with  a polynomial-exponential function in $T$ whose purely polynomial part is contant and denoted by $I_\chi(f)$ .
\item The distributions $I_\chi$ are continuous, left $H(\AAA)$-equivariant and right $(G'(\AAA),\eta_{G'})$-equivariant. Moreover the sum
  \begin{align}\label{eq:intro-spec}
    I (f)= \sum_{\chi} I_{\chi}(f)
  \end{align}
is absolutely convergent and defines a continuous distribution.
  \end{enumerate}
\end{theoreme}
The (coarse) spectral expansion of the trace formula for $G$ is precisely the expression \eqref{eq:intro-spec}.
\end{paragr}

\begin{paragr}\label{S-intro:Ichi}
  The definition of $I_\chi$ given in theorem \ref{intro:jfDef} is convenient to relate the spectral expansion to the geometric expansion. However, to get more explicit forms of the distributions $I_\chi$, we shall use  the following three  expressions:
\begin{align}
\label{eq-intro:Lah}  \int_{[H]}\int_{[G']}(\La_{r}^{T} K_{f,\chi})(h,g')\,\eta_{G'}(g')\,dg'dh\\
   \label{eq-intro:Fh} \int_{[H]} F^{G_{n+1}}(h,T) \int_{[G']}K_{f,\chi}(h,g')\,\eta_{G'}(g')\,dg' dh\\
  \label{eq-intro:Fgn} \int_{[G']}F^{G_{n+1}}(g_n',T) \int_{[H]}K_{f,\chi}(h,g')\,dh \,\eta_{G'}(g')dg'
  \end{align}
Essentially they are given by integration of the kernel $K_{f,\chi}$ along $[H]\times [G']$. However, to have a convergent expression for a general $\chi$, one needs to use some truncation depending on the same parameter $T$ as above. We introduce the  Ichino-Yamana truncation operator, denoted by $\La_{r}^{T}$, whose  definition is recalled in §\ref{S:IY-trunc}. In \eqref{eq-intro:Lah}, we apply it to the left-variable of $K_{f,\chi}$. But one can also use the Arthur characteristic function $F^{G_{n+1}}(\cdot,T)$ whose definition is recalled in \ref{S:Arthur-F}. In \eqref{eq-intro:Fh}, this  function  is evaluated at $h\in H(\AAA)$ through the embedding $H=G_n\hookrightarrow G_{n+1}$. In \eqref{eq-intro:Fgn}, it is evaluated at the component $g'_n$ of the variable $g'=(g'_n,g'_{n+1})\in G'(\AAA)=GL_n(\AAA)\times GL_{n+1}(\AAA)$.

 The link with the distribution $I_\chi$ is provided by the following theorem (which is a combination of Propositions  \ref{prop:cv-LaTrKchi} and \ref{prop:FonHG} and Theorem \ref{thm:asym-trio}). Note that we shall not need the full strength of the theorem in this paper. However it will be used in a greater generality in a subsequent paper.

\begin{theoreme} \label{intro:jgood} Let $f \in \Sc(G(\AAA))$ and $\chi \in \Xgo(G)$. 
  \begin{enumerate}
  \item For any $T$ in some positive Weyl chamber, the expressions \eqref{eq-intro:Lah}, \eqref{eq-intro:Fh} and  \eqref{eq-intro:Fgn} are absolutely convergent.
\item Each of the three expressions  is asymptotically equal (in the technical sense of  §\ref{S:asym-equal}) to a polynomial-exponential function of $T$ whose purely polynomial term is constant and equal to $I_\chi(f)$.
\end{enumerate}
\end{theoreme}
\end{paragr}

\subsection{On the $*$-generic contribution for the Jacquet-Rallis trace formula for the linear groups}\label{ssec:intro-generic}

\begin{paragr}
From now on we assume that the cuspidal datum $\chi$ is relevant $*$-generic  that is $\chi$ is the class of a pair $(M,\pi)$ with the property that the normalized induction $\Pi:=\Ind_{P(\AAA)}^{G(\AAA)}(\pi)$, where we have fixed a  parabolic subgroup $P$ with Levi component $M$, is a Hermitian Arthur parameter of $G$. To $\Pi$ we associate, following \cite[\S 3.4]{Zhang2}, a {\em relative character} $I_\Pi$. The precise definition of this object is recalled in \S \ref{S:IPI}. Let us just say here that it is associated to two functionals $\lambda$ and $\beta_\eta$ on the Whittaker model $\wc(\Pi,\psi_N)$ of $\Pi$, where $\psi_N$ is a certain generic automorphic character of the standard maximal unipotent subgroup $N$ of $G$, that naturally show up in integrals of Rankin-Selberg type. More precisely, $\lambda$ is the value at $s=\frac12$ of a family of Zeta integrals, studied by Jacquet-Piatetski-Shapiro-Shalika \cite{JPSS}, representing the Rankin-Selberg $L$-function $L(s,\Pi)$ whereas $\beta_\eta$ is essentially the pole at $s=1$ of another family of Zeta integrals, first introduced by Flicker \cite{Flicker}, representing the (product of) Asai $L$-functions $L(s,\Pi,\As_G):=L(s,\Pi_n,\As^{(-1)^{n+1}})L(s,\Pi_{n+1},\As^{(-1)^n})$. The relative character $I_\Pi$ is then given in terms of these functionals by
\begin{align}\label{eq-intro:Ipi}
I_\Pi(f)=\sum_{\varphi\in \Pi} \lambda(\Pi(f)W_\varphi) \overline{\beta_\eta(W_\varphi)},\;\; f\in \Sc(G(\AAA)),
\end{align}
where the sum runs over an orthonormal basis of $\Pi$ (for the Petersson inner product) and $W_\varphi$ denotes the Whittaker function associated to the Eisenstein series $E(\varphi)$ (obtained, as usual, by integrating $E(\varphi)$ against $\psi_N^{-1}$ over $[N]$).

The following is our main technical result whose proof occupies most part of the paper.

\begin{theoreme}\label{thm:specJR}
Let $\chi$ be a cuspidal datum associated to a Hermitian Arthur parameter $\Pi$ as above. Then, for every function $f\in \Sc(G(\AAA))$ we have
$$\displaystyle I_\chi(f)=2^{-\dim(A_M)} I_\Pi(f)$$
where $A_M$ denotes the maximal central split torus of $M$.
\end{theoreme}

\begin{remarque}
  It is perhaps worth emphasizing that the contribution of $\chi$ is \emph{purely discrete} in the Jacquet-Rallis trace formula. Such a phenomenon happens in Jacquet relative trace formula, see \cite{LapFRTF}. By contrast, the contribution of the same kind of cuspidal datum $\chi$ to the Arthur-Selberg trace formula is purely continuous (unless, of course, if $\Pi$ is cuspidal).
\end{remarque}

We shall provide two different proofs of theorem \ref{thm:specJR}, one based on truncations, the other using integral representations of Asai and Rankin-Selberg $L$-functions. Let's explain separately the main steps of each approach. 
\end{paragr}

\begin{paragr}[A journey through truncations.] --- We first begin with the approach based on truncations. The first step  is to get a spectral decomposition of the function
  \begin{align}\label{intro:integ-G'}
 \int_{[G']} K_{f,\chi}(g,g') \eta_{G'}(g') dg'.
  \end{align}
of    the variable $ g\in [G]$.

  The kernel itself    $K_{f,\chi}$ has a well-known spectral decomposition based on the Langlands decomposition. Then the  problem is basically to invert an adelic integral and a complex integral. It is solved by Lapid in \cite{LapFRTF} (up to some non-explicit constants) but we will use a slightly different method avoiding  delicate Lapid's contour moving. Instead we replace the integral \eqref{intro:integ-G'} by its truncated version
\begin{align}\label{intro:integ-G'-tronque}
 \int_{[G']} (K_{f,\chi}\Lambda^T_m)(g,g') \eta_{G'}(g') dg'.
  \end{align}
where the {\em mixed truncation operator} $\Lambda^T_m$ defined by Jacquet-Lapid-Rogawski \cite{JLR} is applied to the right variable of the kernel. We can recover \eqref{intro:integ-G'} by taking the limit when $T\to +\infty$. It is easy to get the spectral decomposition of \eqref{intro:integ-G'-tronque} (see Proposition \ref{prop:spec-exp}). Using an analog of the famous Maa\ss-Selberg relations due to Jacquet-Lapid-Rogawski (see \cite{JLR} and Lemma \ref{lem:Maass} below), we get in   proposition \ref{prop:expanLmT}) that \eqref{intro:integ-G'-tronque} is equal to a finite sum of contributions (up to an explicit constant) of the following type
\begin{align}\label{eq-intro:integ-GM}
   \int_{i\ago_{P}^{G,*}}  \sum_{Q\in \pc(M)}  J_{Q,\chi}(g,\la,f) \frac{\exp(-\bg \la,T_Q\bd) }{\theta_Q(-\la) } \,d\la.
\end{align}
 Here it suffices to say that $i\ago_{P}^{G,*}$ is some space of unramified unitary characters and that    $J_{Q,\chi}(g,\la,f) $ is a a certain relative character built upon Flicker-Rallis intertwining periods (introduced by  Jacquet-Lapid-Rogawski). The integrand is a familiar expression of Arthur's theory of $(G,M)$-families with quite standard notations.  It turns out that the family $(J_{Q,\chi}(g,\la,f))_{Q\in \pc(M)}$ is indeed an Arthur  $(G,M)$-family of Schwartz functions in the parameter $\la$. Let's emphasize that this Schwartz property relies in fact on deep estimates introduced by Lapid in \cite{LapFRTF} and \cite{LapHC}. By a standard argument, it is then easy to get the limit of \eqref{eq-intro:integ-GM} when $T\to+\infty$ which gives the spectral decomposition of \eqref{intro:integ-G'}  (see Theorem \ref{thm:spec-exp-FR-kernel}). Note that the spectral decomposition we get is already discrete at this stage.

 From this result, one gets the equality 
\begin{align}\label{eq-intro:LaTr-Ipi}
   \int_{[H]}\int_{[G']}(\La_{r}^{T} K_{f,\chi})(h,g')\,\eta_{G'}(g')\,dhdg'=2^{-\dim(A_M)} I_{P,\pi}(f).
\end{align}
The left-hand side has been defined in §\ref{S-intro:Ichi} and the  relative character $I_{P,\pi}$ is defined as follows:
\begin{align*}
  \sum_{\varphi\in \Pi}  I_{RS}(\Pi(f)\varphi)\cdot \overline{J_\eta(\varphi)}
\end{align*}
where the sum  is over an orthonormal basis, $I_{RS}(\varphi)$ is the regularized Rankin-Selberg period of the Eisenstein series $E(\varphi)$ defined by Ichino-Yamana and $J_\eta(\varphi)$ is a Flicker-Rallis intertwining period (for more detail we refer to §\ref{S:RelcharJPpi}).

In particular, the  left-hand side of \eqref{eq-intro:LaTr-Ipi} does not depend on $T$. So Theorem \ref{intro:jgood} implies

\begin{theoreme}\label{thm-intro:IchiIpi}(see theorem \ref{thm:intLaT} for a slightly more precise statement)
  \begin{align*}
  I_\chi(f)=2^{-\dim(A_M)} I_{P,\pi}(f).
\end{align*}
\end{theoreme}

\begin{remarque}
As the reading of Section \ref{ssec:pf-GG} should make it  clear, this statement suffices to prove the Gan-Gross-Prasad conjecture namely Theorem \ref{thm:GGP}. However, to get the Ichino-Ikeda conjecture, namely Theorem \ref{thm:II}, we will want to use statements about comparison of local relative characters  written in terms of Whittaker functions. For this purpose, Theorem \ref{thm:specJR} will be more convenient.
\end{remarque}

The link between regularized Rankin-Selberg period of  Eisenstein series and Whittaker functionals has been investigated by Ichino-Yamana (see \cite{IY}). The following theorem  relates the Flicker-Rallis intertwining periods to the functional $\beta_\eta(W_\varphi)$ in \eqref{eq-intro:Ipi}. It uses a local unfolding method inspired from \cite[Appendix A]{FLO}  (see Chapter \ref{Chap:FR-functional-computation}). 

\begin{theoreme} \label{thm-intrp:J-beta}For all $\varphi\in \Pi$, we have 
  \begin{align*}
      J_{\eta}(\varphi)=\beta_\eta(W_\varphi)
  \end{align*}
\end{theoreme}

In this way, one proves the following theorem which implies Theorem \ref{thm:specJR}.

\begin{theoreme}\label{thm-intrp:deuxI}
  \begin{align*}
    I_{P,\pi}=I_\Pi.
  \end{align*}
\end{theoreme}

\end{paragr}

\begin{paragr}[Second proof: the use of Zeta integrals.] ---  The spectral decomposition of \eqref{intro:integ-G'}  essentially boils down to a spectral expansion of the period integral
$$\displaystyle P_{G',\eta}(\varphi):=\int_{[G']} \varphi(g') \eta_{G'}(g') dg'$$
for test functions $\varphi\in \Sc_\chi([G])$, where $\Sc_\chi([G])$ denotes the {\em Schwartz space} of $[G]$ consisting of smooth functions rapidly decaying with all their derivatives  that are ``supported on $\chi$'' (see Section \ref{S:space-of-functions} for a precise definition). 
Choose a parabolic subgroup $P=MN_P$ with Levi component $M$. By Langlands $L^2$ spectral decomposition and of the special form of $\chi$,  any $\varphi\in \Sc_\chi([G])$ admits a spectral decomposition
\begin{align}\label{intro:spectral identity}
\varphi=\int_{i\ago_M^*} E(\varphi_\lambda) d\lambda
\end{align}
where $i\ago_M^*$ denotes the real vector space of unramified unitary characters of $M(\AAA)$ and $\varphi_\lambda$ belongs to the normalized induction space $\Ind_{P(\AAA)}^{G(\AAA)}(\pi\otimes \lambda)$ and $E(\varphi_\lambda)$ is the associated Eisenstein series.

\begin{theoreme}\label{thm:specFR2}
For every $\varphi\in \Sc_\chi([G])$, we have
\begin{align*}
P_{G',\eta}(\varphi)=2^{-\dim(A_M)}\beta_\eta(W_{\varphi_0}).
\end{align*}
where $W_{\varphi_0}$ stands for the Whittaker function of the Eisenstein series $E(\varphi_0)$.
\end{theoreme}

The proof of Theorem \ref{thm:specFR2}  is close to the computation by Flicker \cite{Flicker} of the Flicker-Rallis period of cusp forms in terms of an Asai $L$-function and local Zeta integrals. More precisely, we first realize $P_{G',\eta}(\varphi)$ as the residue at $s=1$ of the inner product of the restriction $\varphi_{\mid [G']}$ with some Eisenstein series $E(s,\phi)$ (where $\phi$ is an auxiliary Schwartz function on $\AAA^n\oplus \AAA^{n+1}$). Mimicking the unfolding of {\em loc. cit.} we connect this inner product with an Eulerian Zeta integral $Z^{\FR}(s,W_\varphi,\phi)$ involving the Whittaker function $W_\varphi$ of $\varphi$ (obtained as before by integration against $\psi_N^{-1}$). We should emphasize here that, since $\varphi$ is not a cusp form, the unfolding gives us more terms but using the special nature of the cuspidal datum $\chi$ we are able to show that these extra terms do not contribute to the residue at $s=1$. The formation of $Z^{\FR}(s,W_{\varphi},\phi)$ commutes with the spectral expansion \eqref{intro:spectral identity} when $\Re(s)\gg 1$ and, as follows from the local theory, the Zeta integrals $Z^{\FR}(s,W_{\varphi_\lambda},\phi)$ for $\lambda\in i\ago_M^*$ are essentially Asai $L$-functions whose meromorphic continuations, poles and growths in vertical strips are known. Combining this with an application of the Phragmen-Lindel\"of principle, we are then able to deduce Theorem \ref{thm:specFR2}.

Let us mention here that, as in the proof of \eqref{eq-intro:LaTr-Ipi}, a key point is the fact (due to Lapid (\cite{LapHC} or \cite{LapFRTF})) that the spectral transform $\lambda\mapsto \varphi_\lambda$ is, in a suitable technical sense, ``Schwartz'' that is rapidly decreasing together with all its derivatives.

The second step is to integrate \eqref{intro:integ-G'} over $g\in [H]$. To do so, we define  a regularization of the integral over $[H]$ that doesn't require truncation. More precisely, denoting by $\tc([G])$ the space of functions of {\em uniform moderate growth} on $[G]$, we can define the ``$\chi$-part'' $\tc_\chi([G])$ of $\tc([G])$ (see Section \ref{S:space-of-functions}) of which $\Sc_\chi([G])$ is a dense subspace. Moreover, starting with $\varphi\in \tc([G])$ we can also form its Whittaker function $W_\varphi$ and consider the usual Rankin-Selberg integral $Z^{\RS}(s,W_\varphi)$ that converges for $\Re(s)\gg 1$ and represents, when $\varphi$ is an automorphic form, the Rankin-Selberg $L$-function for $G_n\times G_{n+1}$.

\begin{theoreme} (see theorems \ref{theo Zeta integral RS} and \ref{theo RS period})\label{thm:regularized-RS-period}
The functional
$$\displaystyle \varphi\in \Sc_\chi([G])\mapsto \int_{[H]} \varphi(h) dh$$
extends by continuity to a functional on $\tc_\chi([G])$ denoted by $\displaystyle \varphi\in \tc_\chi([G])\mapsto \int^*_{[H]} \varphi(h) dh.$ Moreover, for every $\varphi\in \tc_\chi([G])$, the Zeta function $s\mapsto Z^{\RS}(s,W_\varphi)$ extends to an entire function on $\CC$ and we have
$$\displaystyle \int^*_{[H]} \varphi(h) dh=Z^{\RS}(\frac12,W_\varphi).$$
\end{theoreme}

The proof of this theorem is similar to that of Theorem \ref{thm:specFR2}: we first show that, for $\varphi\in \Sc_\chi([G])$ and $\Re(s)\gg 1$, we have
$$\displaystyle \int_{[H]} \varphi(h) \lvert \det h\rvert^sdh=Z^{\RS}(s+\frac12,W_\varphi)$$
by mimicking the usual unfolding for the Rankin-Selberg integral. Once again, as $\varphi$ is not necessarily a cusp form, we get extra terms in the course of the unfolding but, thanks to the special nature of the cuspidal datum $\chi$, we are able to show that they all vanish. At this point, we use the spectral decomposition \eqref{intro:spectral identity} to express $Z^{\RS}(s,W_\varphi)$ as the integral of $Z^{\RS}(s,W_{\varphi_\lambda})$ when $\Re(s)\gg 1$. By Rankin-Selberg theory, $Z^{\RS}(s,W_{\varphi_\lambda})$ is essentially a Rankin-Selberg $L$-function whose meromorphic continuation, location of the poles, control in vertical strips and functional equation are known. Combining this with another application of the Phragmen-Lindel\"of principle, we are able to bound $\displaystyle Z^{\RS}(\frac12,W_\varphi)= \int_{[H]}\varphi(h) dh$ in terms of $Z^{\RS}(s,W_\varphi)$ for $\Re(s)\gg 1$ and this readily gives the theorem.

One direct consequence of Theorem \ref{thm:regularized-RS-period} is that the regularized period $\int^*_{[H]} E(h,\Pi(f)\varphi) dh$ coincides with $Z^{\RS}(\frac12,\Pi(f)W_\varphi)=\lambda(\Pi(f)W_\varphi)$. Thus by a combination of Theorems \ref{thm:specFR2} and \ref{thm:regularized-RS-period}, we get

\begin{align}\label{intro:eqIchiregRSperiod}
\int_{[H]}^*\int_{[G']} K_{f,\chi}(h,g') \eta_{G'}(g') dg' \,dh=2^{-\dim(A_M)} I_\Pi(f).
\end{align}

Finally we have to show that the left-hand side is equal to $I_\chi(f)$. In fact, we show (see Theorem \ref{theo JR TF} and §\ref{S:RRS-period-cv}) that we have
\begin{align*}
  \int_{[H]}^*\int_{[G']} K_{f,\chi}(h,g') \eta_{G'}(g') dg' \,dh=\int_{[G']}\int_{[H]}K_{f,\chi}(h,g') \,dh \, \eta_{G'}(g') dg' 
\end{align*}
where the right-hand side is (conditionally) convergent. We can conclude that it is equal to $I_\chi(f)$ by applying Theorem \ref{intro:jgood} to the expression \eqref{eq-intro:Fgn}.

\end{paragr}

\subsection{Outline of the paper}

We now give a quick outline of the content of the paper. Chapter \ref{Chapter:preliminaries} contains preliminary material. Notably, we fix most notation to be used in the paper, we explain our convention on normalization of measures, we introduce the various spaces of functions we need and we discuss several properties of Langlands decomposition along cuspidal data as well as kernel functions that are important for us. Chapter \ref{Chap:Ichilimit} contains the statements and proofs concerning the spectral expansion of the Jacquet-Rallis trace formula for $G$ that were discussed in Section \ref{ssec:intro-spec-exp} above.

In Chapter \ref{sec:sym}, we introduce the Flicker-Rallis intertwining periods and prove the spectral expansion of the Flicker-Rallis period of the kernel associated to $*$-generic cuspidal datum. In Chapter \ref{sec:generic contribution} we deduce from it Theorem \ref{thm-intro:IchiIpi} namely the spectral expansion for $I_\chi(f)$. Chapters \ref{chapter FR} and \ref{chap canonical RS} are devoted to the proofs of Theorems \ref{thm:specFR2} and \ref{thm:regularized-RS-period} respectively. These two theorems are combined in Chapter \ref{Chap:generic-JRTF} to give another proof of the spectral expansion of $I_\chi(f)$ (Theorem \ref{thm:specJR}).  In Chapter \ref{Chap:FR-functional-computation}, we relate the Flicker-Rallis intertwining periods to the functional $\beta_\eta$. From this, we deduce Theorem \ref{thm-intrp:J-beta} and Theorem \ref{thm-intrp:deuxI}. The final Chapter \ref{Chap:proofmaintheorems} explain the deduction of Theorems \ref{thm:GGP} and \ref{thm:II} from Theorem \ref{thm:specJR}. Finally, we have gathered in Appendix \ref{appendice tvs} some useful facts on topological vector spaces and holomorphic functions valued in them. It contains in particular some variations on the theme of the Phragmen-Lindel\"of principle for such functions that will be crucial for the proofs of Theorems \ref{thm:specFR2} and \ref{thm:regularized-RS-period}.

\subsection{Acknowledgement}

The project leading to this publication of R. B.-P. has received funding from Excellence Initiative of Aix-Marseille University-A*MIDEX, a French ``Investissements d'Avenir'' programme.


\section{Preliminaries}\label{Chapter:preliminaries}

\subsection{General notation}

\begin{paragr}
For $f$ and $g$ two positive functions on a set $X$, we way that $f$ is {\em essentially bounded} by $g$ and we write
$$\displaystyle f(x)\ll g(x), \; x\in X,$$
if there exists a constant $C>0$ such that $f(x)\leqslant Cg(x)$ for every $x\in X$. If we want to emphasize that the constant $C$ depends on auxilliary parameters $y_1,\ldots,y_k$, we will write $f(x)\ll_{y_1,\ldots,y_k}g(x)$. We say that the functions $f$ and $g$ are {\em equivalent} and we write
$$\displaystyle f(x)\sim g(x), \; x\in X,$$
if $f(x)\ll g(x)$ and $g(x)\ll f(x)$. If moreover $f$ and $g$ take values in the set of real numbers greater than $1$, we say that $f$ and $g$ are {\em weakly equivalent}, and we write
$$\displaystyle f(x)\approx g(x), \; x\in X,$$
if there exists $N>0$ such that
$$\displaystyle g(x)^{1/N}\ll f(x)\ll g(x)^N,\;\; x\in X.$$
\end{paragr}

\begin{paragr}
For every $C,D\in \RR\cup \{-\infty\}$ with $D>C$, we set $\cH_{>C}=\{z\in \CC\mid \Re(z)>C \}$ and $\cH_{]C,D[}=\{z\in \CC\mid C<\Re(z)<D \}$. A {\em vertical strip} is a subset of $\CC$ which is the closure of $\cH_{]C,D[}$ for some $C,D\in \RR$ with $D>C$.

When $f$ is a meromorphic function on some open subset $U$ of $\CC$ and $s_0\in U$, we denote by $f^*(s_0)$ the leading term in the Laurent expansion of $f$ at $s_0$.
\end{paragr}

\begin{paragr}
When $G$ is a group and we have a space of functions on it invariant by right translation, we denote by $R$ the corresponding representation of $G$. If $G$ is a Lie group and the representation is differentiable, we will also denote by the same letter the induced action of the Lie algebra or of its associated enveloping algebra. If $G$ is a topological group equipped with a bi-invariant Haar measure, we denote by $\ast$ the convolution product (whenever it is well-defined).
\end{paragr}

\subsection{Algebraic groups and adelic points}\label{ssec:alggps}

\begin{paragr}
Let $ F $ be a number field and $\AAA $ its adele ring. We write $ \AAA_f $ for the ring of finite adeles and $F_\infty=F\otimes_{\QQ} \RR$ for the product of Archimedean completions of $F$ so that $\AAA=F_\infty\times \AAA_f$. Let $ V_F $ be the set of places of $F$ and $V_{F,\infty}\subset V_F$ be the subset of Archimedean places.  For every $ v \in V_F $, we let $ F_v $ be the local field obtained by completion of  $ F $ at $ v $. We denote by  $ | \cdot | $  the morphism $ \AAA^\times \to \RR_+^\times $ given by the product of normalized absolute values $|\cdot|_v$ on each $ F_v $. For any finite subset $S\subset V_F\setminus V_{F,\infty}$, we denote by $\oc_F^S$ the ring of $S$-integers in $F$.
\end{paragr}

\begin{paragr}
Let $ G $ be an algebraic group defined over $ F $. We denote by $N_G$ the unipotent radical of $G$. Let $ X^*(G) $ be the group of  characters of $ G $ defined over $ F $. Let $\ago_G^*=X^*(G)\otimes_\ZZ\RR$ and $\ago_G=\Hom_\ZZ(X^*(G),\RR)$. We have a canonical pairing
\begin{align}\label{eq:pairing}
\bg \cdot,\cdot\bd : \ago_G^*\times \ago_G\to \RR.
\end{align}
We have also a canonical homomorphism
\begin{align}\label{eq:HG}
H_G: G (\AAA) \to \ago_ {G} 
\end{align}
such that $ \bg \chi, H_G(g) \bd = \log | \chi (g) | $ for any $g\in G(\AAA)$. The kernel of $H_G$ is denoted by  $ G(\AAA) ^1 $. We define $[G]=G(F)\back G(\AAA)$ and $[G]^1=G(F)\back G(\AAA)^1$.

We let $\ggo_\infty$ be the Lie algebra of $G(F_\infty)$, $\uc(\ggo_\infty)$ be the enveloping algebra of its complexification and $\zc(\ggo_\infty)\subset \uc(\ggo_\infty)$ be its center.
\end{paragr}

\begin{paragr} From now on we assume that $ G $ is also reductive. We will mainly use the notations of Arthur's works. For the convenience of the reader, we recall some of them. Let $ P_0 $ be a parabolic subgroup of $ G $ defined over $ F $  and minimal for these properties. Let $ M_0 $ be a Levi factor of $ P_0$ defined over  $ F $.
  
We call a parabolic  (resp. and semi-standard, resp. and standard) subgroup of $G$ a parabolic subgroup of $G$ defined over $ F $ (resp.  which contains $ M_0 $, resp. which contains $ P_0 $). For any semi-standard parabolic subgroup $P$, we have a Levi decomposition $ P = M_P N_P $ where $ M_P $ contains $ M_0 $ and we define $[G]_P=M_P(F)N_P(\AAA)\back G(\AAA)$.  We call a Levi subgroup of $ G $ (resp. semi-standard,  resp. standard) a  Levi factor defined over $ F $ of a parabolic subgroup of $G$ (resp. semi-standard, resp.  standard).
\end{paragr}

\begin{paragr}
Let $ K=\prod_{v\in V_F} K_v\subset G (\AAA) $ be a ``good'' maximal compact subgroup in good position relative to $M_0$. We write
\begin{align*}
K=K_\infty K^\infty
\end{align*}
where  $K_\infty=\prod_{v\in V_{F,\infty}} K_v$ and $K^\infty=\prod_{v\in V_F\setminus V_{F,\infty}} K_v$. We let $\kgo_\infty$ be the Lie algebra of $K_\infty$, $\uc(\kgo_\infty)$ be the enveloping algebra of its complexification and $\zc(\kgo_\infty)\subset \uc(\kgo_\infty)$ be its center.
\end{paragr}

\begin{paragr}   Let $ P $ be a semi-standard parabolic subgroup. We extend the homomorphism  $H_P:P(\AAA)\to \ago_P$ (see \eqref{eq:HG}) into the Harish-Chandra map 
$$ H_P: G (\AAA) \to \ago_ {P} $$
in such a way that for every $g \in G(\AAA)$ we have $H_P(g)=H_P(p)$ where $ p \in P (\AAA) $ is given by the Iwasawa decomposition namely  $ g \in pK $. 
\end{paragr}

\begin{paragr}\label{S:AP}
Let $A$ be a split torus over $F$. Then, $A$ admits an unique split model over $\QQ$ (which is also the maximal split subtorus of $\Res_{F/\QQ}(A)$) and by abuse of notation we denote by $A(\RR)$ the group of $\RR$-points of this model. In particular, this gives an embedding $\RR^\times\subset F_\infty^\times\subset \AAA^\times$. We also write $A^\infty$ for the neutral component of $A(\RR)$. Let  $A_G$ be the maximal central $F$-split torus of $G$. We define  $[G]_0=A_G^\infty G(F)\back G(\AAA)$.

Let $P$ be a semi-standard parabolic subgroup of $G$. We define $A_P=A_{M_P}$, $A_P^\infty=A_{M_P}^\infty$ and $[G]_{P,0}=A_P^\infty M_P(F) N_P(\AAA)\back G(\AAA)$. The restrictions maps  $X^*(P)\to X^*(M_P)\to X^*(A_P)$ induce isomorphisms $\ago_P^*\simeq \ago_{M_P}^* \simeq \ago_{A_P}^*$. Let $\ago_0^*=\ago_{P_0}^* $,  $\ago_0=\ago_{P_0}$, $A_0=A_{P_0}$ and $A_0^\infty=A_{P_0}^\infty$.
\end{paragr}

\begin{paragr}\label{S:proj} For any semi-standard parabolic subgroups $P\subset Q$ of $G$, the restriction map $X^*(Q)\to X^*(P)$ induces maps  $\ago_Q^*\to \ago_P^*$ and $\ago_P\to \ago_Q$. The first one is injective whereas the kernel of the second one is denoted by $\ago_P^Q$. The restriction map $X^*(A_P)\to X^*(A_Q)$ gives a surjective map  $\ago_P^*\to \ago_Q^*$ whose kernel is denoted by  $\ago_P^{Q,*}$. We get also an  injective map $\ago_Q\to \ago_P$. In this way, we get dual decompositions $\ago_P=\ago_P^Q\oplus \ago_Q$ and $\ago_P^*=\ago_P^{Q,*}\oplus \ago_Q^*$. Thus we have  projections $\ago_0\to \ago_P^Q$  and $\ago_0^*\to \ago_P^{Q,*}$ which we will denote by  $X\mapsto X_P^Q$.

We denote by $\ago_{P,\CC}^{Q,*}$ and $\ago_{P,\CC}^Q$ the $\CC$-vector spaces obtained by extension of scalars from $\ago_{P}^{Q,*}$ and  $\ago_{P}^{Q}$. We still denote by $\bg\cdot,\cdot\bd$ the pairing \eqref{eq:pairing} we get by extension of the scalars to $\CC$. We have a decomposition
$$\ago_{P,\CC}^{Q,*}=\ago_{P}^{Q,*}\oplus i\ago_{P}^{Q,*}$$
where $i^2=-1$.  We shall denote by $\Re$ and $\Im$ the associated projections  and call them real and imaginary  parts. The same holds for the dual spaces $\ago_{P,\CC}^Q$. In the obvious way, we define the complex conjugate denoted by $ \bar{\la}$ of $\la\in \ago_{P,\CC}^{Q,*}$.
\end{paragr}

\begin{paragr}
Let  $\Ad_{P}^Q$ be the adjoint action of $M_P$ on the Lie algebra of $M_Q\cap N_P$.  Let  $\rho_P^Q$ be the unique element of $\ago_{P}^{Q,*}$ such that for every $m\in M_P(\AAA)$
$$|\det(\Ad_P^Q(m))|=\exp(\bg 2\rho_P^Q,H_P(m)\bd).
$$
For $Q=G$, the exponent $G$ is omitted. For every $g\in G(\AAA)$, we set
$$\displaystyle \delta_P(g)=\exp(\bg 2\rho_P,H_P(g)\bd)$$
so that, in particular, the restriction of $\delta_P$ to $P(\AAA)$ coincides with the modular character of the latter.
\end{paragr}

\begin{paragr} \label{S:root-coroot}Let $P_0'=M_0N_{P_0'}$ be a minimal semi-standard parabolic subgroup such that $P_0'\subset P$. Let $\Delta_{P_0'}^P$ be the set of simple roots of  $A_{0}$ in  $M_P\cap P_0'$. We denote this set by $\Delta_0^P$ if $P_0'=P_0$. Let $\Delta_P$ be the image of $\Delta_{P_0'}\setminus \Delta_{P_0'}^P$ (viewed as a subset of $\ago_0^*$) by the projection $\ago_{0}^*\to \ago_P^*$. It does not depend on the choice of $P_0'$. More generally  one defines $\Delta_P^Q$. We have also the set of coroots $\Delta_P^{Q,\vee}\subset \ago_P^{Q}$. By duality, we get a set of simple weights $\hat{\Delta}_P^Q$.  The sets $\Delta_P^Q$ and  $\hat{\Delta}_P^Q$ determine open  cones  in  $\ago_{0}$ whose characteristic functions are denoted respectively by $\tau_P^Q$ and  $\hat{\tau}_P^Q$. If $Q=G$, the exponent $G$ is omitted. We set
$$\displaystyle A_P^{\infty,+}=\left\{a\in A_P^\infty\mid \langle \alpha,H_P(a)\rangle\geqslant 0,\; \forall \alpha\in \Delta_P \right\}.$$
In the same way, we define $\ago_P^+$ and $\ago_P^{*,+}$ (this time using coroots).
\end{paragr}

\begin{paragr}[Weyl group.] --- \label{S:Weyl}Let $W$ be the Weyl group of $(G,A_0)$ that is the  quotient by $M_0$ of the normalizer of  $A_0$ in $G(F)$. For  $P=M_P N_P$ and  $Q=M_QN_Q$ two standard parabolic subgroups of $ G $, we denote by $W(P,Q)$ the set of $w\in W$ such that  $w\Delta^P_0=\Delta_0^Q$. For $w\in W(P,Q)$, we have  $wM_Pw^{-1}=M_Q$. When $P=Q$, the group $W(P,P)$ is simply denoted by $W(P)$. Sometimes, we shall also denote  $W(P,Q)$ by $W(M_P,M_Q)$  if we want to emphasize the Levi components (and $W(M_P)=W(P)$).
\end{paragr}

\begin{paragr}\label{S:pcM}
Let $M$ be a standard Levi subgroup of $G$. We denote by $\pc(M)$ the set of semi-standard parabolic subgroups $P$ of $G$ such that  $M_P=M$. There is an unique element $P\in \pc(M)$ which is standard and the map
\begin{align}\label{eq:wQ}
(Q,w)\mapsto w^{-1}Qw
\end{align}
induces a  bijection from the disjoint union $\bigcup_Q W(P,Q)$ where $Q$ runs over the set of standard parabolic subgroups of $G$ onto $\pc(M)$. 
\end{paragr}

\begin{paragr}[Truncation parameter.] ---\label{S:trunc-param}
We shall denote by $T$ a point of $\ago_{0}$ such that $\bg \al,T\bd$ is large enough for every $\al\in \Delta_{0}$. We do not want to be precise here. We just need that Arthur's formulas about truncation functions hold for the $T's$ we consider (see \cite{ar1} \S\S 5,6). The point $T$ plays the role of a truncation parameter.

For any semi-standard parabolic subgroup $P$, we define a point $T_P\in \ago_{P}$ such that for any  $w\in W$ such that $wP_0w^{-1}\subset P$, the point $T_P$ is the projection of $w\cdot T$ on $ \ago_{P}$ (this does not depend on the choice of $w$). The reader should be warned that it is not consistent with the notation of \S\ref{S:proj} since there $T_P$ denotes instead the projection of $T$ onto $\ago_{P}$ (of course, the two conventions coincide when $P$ is standard).
\end{paragr}

\begin{paragr}\label{Siegel domain}
Let $P$ be a standard parabolic subgroup. By a {\em Siegel domain} for $[G]_P$ we mean a subset of $G(\AAA)$ of the form
$$\displaystyle \sgo_P=\omega_0 \left\{ a\in A_0^\infty\mid \langle \alpha, H_0(a)+T\rangle\geqslant 0, \; \forall \alpha\in \Delta^P_0 \right\}K$$
where $T\in \ago_0$ and $\omega_0\subset P_0(\AAA)^1$ is a compact such that $G(\AAA)=M_P(F)N_P(\AAA)\sgo_P$.
\end{paragr}

\subsection{Haar measures}\label{S:Haar-measures}

\begin{paragr}\label{S:Haar-norm}
We equip $\ago_P$ with the Haar measure that gives a covolume $1$ to the lattice $\Hom(X^*(P),\ZZ)$. The space  $i\ago_{P}^{*}$ is then equipped with  the dual Haar measure so that we have 
\begin{align*}
\int_{i\ago_{P}^*} \int_{\ago_{P}}  \phi(H) \exp (-\bg \la ,H\bd)\,dHd\la=\phi(0) 
\end{align*}
for all $\phi\in \Cc(\ago_{P})$. Note that this implies that the covolume of $iX^*(P)$ in $i\ago_P^*$ is given by
\begin{align}\label{eq measure iaP*}
\vol(i\ago_P^*/iX^*(P))=(2\pi)^{-\dim(\ago_P)}.
\end{align}

The group $ A_P ^ \infty $  is equipped with the Haar measure compatible with the isomorphism $ A_P ^ \infty \simeq  \ago_ {P} $ induced by the map $H_P$. The groups $\ago_P^G\simeq \ago_P/\ago_G$ and  $i\ago_P^{G,*}\simeq i\ago_P^*/i\ago_G^*$ are provided with the quotient Haar measures.
For any basis $B$ of $\ago_P^G$ we denote by  $\ZZ(B)$ the lattice generated by $B$ and by  $\vol(\ago_P^G/\ZZ(B))$ the  covolume of this lattice. We have on $\ago_0^*$ the  polynomial function:
$$\theta_P(\la)= \vol(\ago_P^G/\ZZ(\Delta_P^{\vee}))^{-1} \prod_{\al \in \Delta_P} \bg \la,\al^\vee\bd.
$$
\end{paragr}

\begin{paragr}
Let $H$ be a linear algebraic group over $F$. In this paper, we will always equip $H(\AAA)$ with its right-invariant Tamagawa measure simply denoted $dh$. Let us recall how it is defined in order to fix some notation. We choose a right-invariant rational volume form $\omega_H$ on $H$ as well as a non-trivial continuous additive character $\psi':\AAA/F\to \CC^\times$. For each place $v\in V_F$, the local component $\psi'_v$ of $\psi'$ induces an additive measure on $F_v$ which is the unique Haar measure autodual with respect to $\psi'_v$. Then using local $F$-analytic charts, we associate to $\omega_H$ a right Haar measure $dh_v=\lvert \omega_H\rvert_{\psi'_v}$ on $H(F_v)$ as in \cite[\S 2.2]{Wei}. By \cite{Gros}, there exists an Artin-Tate $L$-function $L_H(s)$ such that, denoting by $L_{H,v}(s)$ the corresponding local $L$-factor and setting $\Delta_{H,v}=L_{H,v}(0)$, for any model of $H$ over $\oc_F^S$ for some finite set $S\subseteq V_F\setminus V_{F,\infty}$, we have
\begin{equation}\label{eq0 Tamagawa measures}
\displaystyle \vol(H(\oc_v))=\Delta_{H,v}^{-1}
\end{equation}
for almost all $v\in V_F$. Setting $\Delta^*_H=L^*_H(0)$, the Tamagawa measure on $H(\AAA)$ is defined as the product
\begin{equation}\label{eq Tamagawa measures}
\displaystyle dh=(\Delta_H^*)^{-1}\prod_v \Delta_{H,v}dh_v.
\end{equation}
Although the local measures $dh_v$ depend on choices, the global measure $dh$ doesn't (by the product formula).
\end{paragr}

\begin{paragr}\label{S:dgS}
For $S\subseteq V_F$ a finite subset, we put $\Delta_H^{S,*}=L^{S,*}_H(0)$ where $L^S_H(s)$ stands for the corresponding partial $L$-function and we equip $H(F_S)$, $H(\AAA^S)$ with the right Haar measures $dh_S=\prod_{v\in S} dh_v$ and $dh^S=(\Delta_H^{S,*})^{-1}\prod_{v\notin S} \Delta_{H,v}dh_v$ respectively. Note that we have the decomposition
\begin{equation}\label{eq decomposition Tamagawa measure}
\displaystyle dh=dh_S\times dh^S.
\end{equation}
In particular, this means that $H(F_v)$ is equipped with the right Haar measure $dh_v$ for every $v\in V_F$.
\end{paragr}

\begin{paragr}
We have $L_{H}(s)=L_{H_{\red}}(s)$ where $H_{\red}=H/H_u$ denotes the quotient of $H$ by its unipotent radical $H_u$. When $H=N$ is unipotent we have $\vol([N])=1$. For $H=\GL_n$, the $L$-function $L_H(s)$ is given by
$$\displaystyle L_H(s)=\zeta_F(s+1)\ldots\zeta_F(s+n)$$
where $\zeta_F$ stands for the (completed) zeta function of the number field $F$. In this case, we will take
$$\displaystyle \omega_H=(\det h)^{-1}\bigwedge_{1\leqslant i,j\leqslant n}dh_{i,j}$$
so that \eqref{eq0 Tamagawa measures} is satisfied for every non-Archimedean place $v$ where $\psi_v$ is unramified.
\end{paragr}

\begin{paragr}\label{S:Haar} The homogeneous space $[G]$ (resp. $[G]^1\simeq [G]_0$) is equipped with the quotient of the Tamagawa measure on $G(\AAA)$ by the counting measure on $G(F)$ (resp. by the product of the counting measure on $G(F)$ with the Haar measure we fixed on $A_G^\infty$). For $P$ a standard parabolic subgroup, we equip similarly $[G]_{P}$ with the quotient of the Tamagawa measure on $G(\AAA)$ by the product of the counting measure on $M_P(F)$ with the Tamagawa measure on $N_P(\AAA)$. Since the action by left translation of $a\in A_P^\infty$ on $[G]_P$ multiplies the measure by $\delta_P(a)^{-1}$, taking the quotient by the Haar measure on $A_P^\infty$ induces a ``semi-invariant'' measure on $[G]_{P,0}=A_P^\infty\backslash [G]_P$ that is a positive linear form on the space of continuous functions $\varphi:[G]_P\to \CC$ satisfying $\varphi(ag)=\delta_P(a)\varphi(g)$ for $a\in A_P^\infty$ and compactly supported modulo $A_P^\infty$. 
\end{paragr}

\subsection{Norms and Harish-Chandra $\Xi$ function}\label{ssec:norms}

\begin{paragr}
Let $X$ be an algebraic variety over $F$. We define following \cite[\S A.1]{RBP} a weak equivalence class of norms $\lVert .\rVert_{X(\AAA)}:X(\AAA)\to \RR_{\geqslant 1}$. For all the algebraic varieties that we encounter in this paper, we implicitely fix a norm in this equivalence class and we set $\sigma_{X(\AAA)}(x)=1+\log \lVert x\rVert_{X(\AAA)}$, $x\in X(\AAA)$. If the context is clear, we simply write $\lVert .\rVert$ and $\sigma$ for $\lVert .\rVert_{G(\AAA)}$ and $\sigma_{G(\AAA)}$ respectively.

Let $P\subset G$ be a semi-standard parabolic subgroup. We set
$$\lVert g\rVert_{[G]_P}=\inf_{\gamma\in M_P(F)N_P(\AAA)} \lVert \gamma g\rVert \;\; \mbox{ and } \;\;\sigma_{[G]_P}(g)=\inf_{\gamma\in M_P(F)N_P(\AAA)} \sigma(\gamma g)$$
for every $g\in [G]_P$. We have
\begin{equation}\label{eq0bis Function spaces}
\displaystyle \lVert mk\rVert_{[G]_P}\approx\lVert m\rVert_{[M_P]},\;\; m\in [M_P], k\in K.
\end{equation}

Let $\sgo_P$ be a Siegel domain as in Section \ref{Siegel domain}. By Proposition \cite[A.1.1 (viii)]{RBP} and \eqref{eq0bis Function spaces}, we have
\begin{equation}\label{eq 0bis Function spaces}
\displaystyle \lVert g\rVert_{[G]_P}\approx \lVert g\rVert_{G(\AAA)},\;\; g\in \sgo_P.
\end{equation}
\end{paragr}

\begin{paragr}
Let $\omega_G\subseteq G(\AAA)$ be a compact subset with nonempty interior and set
$$\displaystyle \Xi^{[G]_P}(g)=\vol_{[G]_P}(g\omega_G)^{-1/2},\;\; g\in [G]_P.$$
For another choice of compact subset $\omega_G'\subseteq G(\AAA)$ with nonempty interior, the resulting functions are equivalent which is why we dropped the subset $\omega_G$ from the notation. We have
\begin{equation}\label{eq1 Function spaces}
\displaystyle \Xi^{[G]_P}(g)\sim \exp(\langle \rho_0,H_0(g)\rangle),\;\; g\in \sgo_P.
\end{equation}
Indeed, we readily check that $\Xi^{[G]_P}(mk)\sim \exp(\langle \rho_P, H_P(m)\rangle)\Xi^{[M_P]}(m)$ for every $m\in [M_P]$ and $k\in K$ so that we are reduced to the case $P=G$. By invariance we have $\vol_{[G]}(g\omega_G)=\vol_{[G]}(g\omega_G g^{-1})$ and it is easy to see that, if $\omega_G$ is chosen sufficiently small, there exists a compact subset $\omega_G'\subseteq G(\AAA)$ such that $g\omega_G g^{-1}\subseteq N_0(F)\omega_G'$ for every $g\in \sgo_G$. As $G(F)$ is discrete inside $G(\AAA)$, the set
$$\displaystyle \left\{\gamma\in P_0(F)\back G(F)\mid \gamma \omega'_G\cap P_0(F)\omega_G'\neq \emptyset \right\}$$ 
is finite. Therefore, we have
$$\displaystyle \vol_{[G]}(g\omega_G g^{-1})\sim \vol_{[G]_{P_0}}(g\omega_Gg^{-1})\sim \exp(-\langle 2\rho_0,H_0(p)\rangle)$$
for $g\in \sgo_G$ and \eqref{eq1 Function spaces} follows.

By \cite[\S 2, (9)]{LapHC}, we also have\footnote{Note that the definition of the $\Xi$ function in {\it loc. cit.} coincides, up to equivalence, with ours by \eqref{eq1 Function spaces}.}
\begin{num}
\item\label{eq 0 Function spaces} There exists $d>0$ such that $\displaystyle \int_{[G]_P} \Xi^{[G]_P}(g)^2 \sigma_{[G]_P}(g)^{-d}dg<\infty$;
\end{num}

From \eqref{eq 0bis Function spaces} and \eqref{eq1 Function spaces}, we deduce the existence of $N\geqslant 1$ such that
\begin{equation}\label{eq2 Function spaces}
\displaystyle \Xi^{[G]_P}(g)\ll \lVert g\rVert^N_{[G]_P},\;\; g\in [G]_P.
\end{equation}
\end{paragr}

\subsection{Spaces of functions}\label{S:space-of-functions}

\begin{paragr}
We say of a function $f:G(\AAA)\to \CC$ that it is {\em smooth} if it is right invariant by a compact-open subgroup $J$ of $G(\AAA_{f})$ and for every $g_f\in G(\AAA_{f})$ the function $g_\infty\in G(F_\infty)\mapsto f(g_fg_\infty)$ is $C^\infty$.
\end{paragr}

\begin{paragr}\label{S:Schw}
Let $C$ be a compact subset of $G(\AAA_f)$ and let $J\subset K^\infty$ be a compact open subgroup.  The left and right  actions of the envelopping algebra $\uc(\ggo_\infty)$ are denoted by $L$ and $R$ respectively. Let $\Sc(G(\AAA),C,J)$ be  the space of smooth functions $f:G(\AAA)\to \CC$ which are biinvariant by $J$, supported in the subset $G(F_\infty)\times C$ and such that the semi-norms
\begin{align*}
\|f\|_{r,X,Y}=\sup_{g\in G(\AAA)} \|g\|_{G(\AAA)}^r  |(R(X)L(Y)f)(g)|
\end{align*}
are finite for every integer $r\geq 1$ and $X,Y\in \uc(\ggo_\CC)$. This family of semi-norms define a topology on $\Sc(G(\AAA),C,K_0)$ making it into a Fr\'echet space. The global \emph{Schwartz space} $\Sc(G(\AAA))$ is the topological direct limit over all pairs $(C,J)$ of the spaces  $\Sc(G(\AAA),C,J)$. The Schwartz space is an algebra for the convolution product denoted by $*$. It contains the dense subspace $\Cc(G(\AAA))$ of smooth and compactly supported functions. For an integer $r\geq 0$, we will also consider the space $C^r_c(G(\AAA))$ generated by products $f_\infty f^\infty$ where $f_\infty$ is a compactly supported function on $G(F\otimes_\QQ \RR)$ which admits derivatives up to the order $r$ and $f^\infty$ is a smooth compactly supported  function on $G(\AAA_f)$.

For every integer $n\geqslant 1$, we define similarly the global Schwartz space $\Sc(\AAA^n)$ (the definition is the same as above up to replacing $G$ by the additive group $\GG_a^n$).
\end{paragr}

\begin{paragr}
Let $\cH$ be a Hilbert space carrying a continuous representation of $G(\AAA)$ (not necessarily unitary). A vector $v\in \cH$ is {\em smooth} if it is invariant by a compact-open subgroup of $G(\AAA_f)$ and the function $g_\infty\in G(F_\infty)\mapsto g_\infty v\in \cH$ is $C^\infty$. We denote by $\cH^\infty$ the subspace of smooth vectors. It is a $G(\AAA)$-invariant subspace carrying its own locally convex topology making it into a strict LF space: for every compact-open subgroup $J\subseteq G(\AAA_f)$, the subspace $(\cH^\infty)^J$ of vectors fixed by $J$ is equipped with the Fr\'echet topology associated to the semi-norms $v\mapsto \lvert Xv\rvert_{\cH}$ where $\lvert .\rvert_{\cH}$ is the norm on $\cH$, $X\in \uc(\ggo_\infty)$ and $Xv:=R(X)(g_\infty\mapsto g_\infty v)_{\mid g_\infty=1}$. 
\end{paragr}

\begin{paragr}\label{S:L2spaces}
Let $P$ be a semi-standard parabolic subgroup of $G$. We denote by $L^2([G]_P)$ the space of $L^2$-measurable functions on $[G]_P$. It is a Hilbert space when equipped with the scalar product
$$\displaystyle \langle \varphi_1,\varphi_2\rangle_{[G]_P}=\int_{[G]_P} \varphi_1(g)\overline{\varphi_2(g)} dg$$
associated to the Tamagawa invariant measure on $[G]_P$. We denote similarly by $L^2([G]_{P,0})$ the Hilbert space of measurable functions $\varphi$ on $[G]_P$ satisfying $\varphi(ag)=\delta_P(a)^{1/2}\varphi(g)$ for almost all $a\in A_P^\infty$ and such that $\displaystyle \int_{[G]_{P,0}} \lvert \varphi(g)\rvert^2 dg$ is convergent.

We can define more generally weighted $L^2$ spaces as follows. Let $w$ be a weight on $[G]_P$ \cite[\S 3.1]{Ber} that is a positive measurable function on $[G]_P$ such that for every compact subset $\omega_G\subset G(\AAA)$ we have $w(xg)\ll w(x)$ for $x\in [G]_P$ and $g\in \omega_G$. The Hilbert space $L^2_w([G]_P)=L^2([G]_P,w(g)dg)$ is then equipped with a continuous (non-unitary) representation of $G(\AAA)$ by right-translation. In particular, following the previous paragraph, we denote by $L^2_w([G]_P)^\infty$ its subspace of smooth vectors: it consists of smooth functions $\varphi:[G]_P \to \CC$ such that $R(X)\varphi\in L_w^2([G]_P)$ for every $X\in \uc(\ggo_\infty)$. By the Sobolev inequality, see \cite[\S 3.4, Key Lemma]{Ber}, for every $\varphi\in L^2_w([G]_P)^\infty$ we have
\begin{equation}\label{basic inequality smooth L2 spaces}
\displaystyle \lvert \varphi(g)\rvert \ll \Xi^{[G]_P}(g) w(g)^{-1/2},\;\; g\in [G]_P.
\end{equation}

In this paper, we will only use three kind of weights. First, for $N\in \RR$, $\lVert .\rVert_{[G]_P}^N$ is a weight and we set $L^2_N([G]_P)=L^2_{\lVert .\rVert_{[G]_P}^N}([G]_P)$. Secondly, for $d>0$, $\sigma_{[G]_P}^d$ is a weight and we set $L^2_{\sigma,d}([G]_P)=L^2_{\sigma_{[G]_P}^d}([G]_P)$. Thirdly, if $\lambda$ is a weight on $\ago_P$, i.e. a real measurable function such that for every compact $\omega\subset \ago_P$ we have $\lvert \lambda(X+Y)-\lambda(X)\rvert\ll 1$ for $X\in \ago_P$ and $Y\in \omega$ (e.g. $\lambda\in \ago_P^*$), then $\exp(\lambda\circ H_P)$ is a weight on $[G]_P$ and we set $L^2_\lambda([G]_P)=L^2_{\exp(\lambda\circ H_P)}([G]_P)$.
\end{paragr}

\begin{paragr}
The {\em Schwartz space} $\Sc([G]_P)$ of $[G]_P$ is defined as the space of smooth functions $\varphi:[G]_P\to \CC$ such that for every $N>0$ and $X\in \uc(\ggo_\infty)$ we have
$$\displaystyle \lvert (R(X)\varphi)(g)\rvert\ll \lVert g\rVert_{[G]_P}^{-N},\;\; g\in [G]_P.$$
Then, $\Sc([G]_P)$ is naturally equipped with a locally convex topology making it into a strict LF space (it is the inductive limit of the Fr\'echet spaces $\Sc([G]_P)^J$ for $J$ a compact-open subgroup of $G(\AAA_f)$). From \eqref{basic inequality smooth L2 spaces}, we have the alternative description
\begin{equation}\label{eq3 Function spaces}
\displaystyle \Sc([G]_P)=\bigcap_{N>0} L^2_N([G]_P)^\infty.
\end{equation}
\end{paragr}

\begin{paragr}
The {\em Harish-Chandra Schwartz space} $\cC([G]_P)$ of $[G]_P$ is defined as the space of smooth functions $\varphi:[G]_P\to \CC$ such that for every $d>0$ and $X\in \uc(\ggo_\infty)$ we have
$$\displaystyle \lvert (R(X)\varphi)(g)\rvert\ll \Xi^{[G]_P}(g)\sigma_{[G]_P}(g)^{-d},\;\; g\in [G]_P.$$
Once again, $\cC([G]_P)$ is naturally equipped with a locally convex topology making it into a strict LF space. Alternatively, we have
\begin{equation}\label{eq4 Function spaces}
\displaystyle \cC([G]_P)=\bigcap_{d>0} L^2_{\sigma,d}([G]_P)^\infty.
\end{equation}
\end{paragr}

\begin{paragr}
The {\em space of functions of uniform moderate growth} $\tc([G]_P)$ of $[G]_P$ is defined as the space of smooth functions $\varphi:[G]_P\to \CC$ for which there exists $N>0$ such that for every $X\in \uc(\ggo_\infty)$ we have
$$\displaystyle \lvert (R(X)\varphi)(g)\rvert\ll \lVert g\rVert_{[G]_P}^N,\;\; g\in [G]_P.$$
For $N>0$, we denote by $\tc_N([G]_P)$ the subspace of functions $\varphi\in \tc([G]_P)$ satisfying the above inequality for every $X\in \uc(\ggo_\infty)$. Then, $\tc_N([G]_P)$ is naturally equipped with a locally convex topology making it into a strict LF space and $\tc([G]_P)=\bigcup_{N>0} \tc_N([G]_P)$ is a (non-strict) LF space. We also have the alternative description
\begin{equation}\label{eq5 Function spaces}
\displaystyle \tc([G]_P)=\bigcup_{N>0} L_{-N}^2([G]_P)^\infty.
\end{equation}
\end{paragr}

\begin{paragr}
The spaces $\Sc([G]_P)$, $\cC([G]_P)$ and $\tc([G]_P)$ are all topological representations of $G(\AAA)$ for the action by right translation $R$. For $J\subseteq G(\AAA_{f})$ a compact-open subgroup and $N>0$, $\Sc([G]_P)^J$, $\cC([G]_P)^J$ and $\tc_N([G]_P)^J$ are even SF representations of $G(F_\infty)$ in the sense of \cite{BK}. It follows that the action of $G(\AAA)$ on $\Sc([G]_P)$, $\cC([G]_P)$ and $\tc([G]_P)$ integrates to an action of the algebra $(\Sc(G(\AAA)),\ast)$ (by right convolution). Moreover, by \cite[end of Section 3.5]{Ber} (see also \cite[Corollary 2.6]{CasSchwartz}) we have
\begin{num}
\item\label{eq: Schwartz spaces nuclear} The Fr\'echet spaces $\Sc([G]_P)^J$ and $\cC([G]_P)^J$ are nuclear.
\end{num}

Let $\Sc([G]_P)'$ be the topological dual of $\Sc([G]_P)$. By duality, it is also equipped with an action by convolution of $\Sc(G(\AAA))$. By the alternative descriptions \eqref{eq3 Function spaces} and \eqref{eq5 Function spaces}, we have
\begin{num}
\item\label{eq6 Function spaces} For every distribution $D\in \Sc([G]_P)'$ and $f\in \Sc(G(\AAA))$, the distribution $R(f)D$ is representable by a function in $\tc([G]_P)$.
\end{num}
\end{paragr}

\begin{paragr}
Assume that $G=G_1\times G_2$ where $G_1$ and $G_2$ are two connected reductive groups over $F$. Let $J_1\subset G_1(\AAA_{f})$ , $J_2\subseteq G_2(\AAA_{f})$ be two compact open subgroups and set $J=J_1\times J_2$. By \eqref{eq: Schwartz spaces nuclear}, \eqref{eq 5 tvs} and a reasoning similar to (the proof of) \cite[Proposition 4.4.1 (v)]{BP3} we obtain:

\begin{num}
\item\label{eq7 Function spaces} There are topological isomorphisms
$$\displaystyle \cC([G_1])^{J_1}\widehat{\otimes} \cC([G_2])^{J_2}\simeq \cC([G])^{J},\; \Sc([G_1])^{J_1}\widehat{\otimes} \Sc([G_2])^{J_2}\simeq \Sc([G])^{J}$$
sending a pure tensor $\varphi_1\otimes \varphi_2$ to the function $(g_1,g_2)\mapsto \varphi_1(g_1)\varphi_2(g_2)$.
\end{num}
By the above, given two continuous linear forms $L_1$, $L_2$ on $\cC([G_1])$, $\cC([G_2])$ respectively, the linear form $L_1\otimes L_2$ on $\cC([G_1])\otimes \cC([G_2])$ extends by continuity to a linear form on $\cC([G])$ that we shall denote by $L_1\widehat{\otimes} L_2$.
\end{paragr}

\begin{paragr}[Constant terms and pseudo-Eisenstein series.] --- Let $Q\subset P$ be another standard parabolic subgroup. We have two continuous $G(\AAA)$-equivariant linear maps
$$\displaystyle \tc([G]_P)\to \tc([G]_Q),\;\; \varphi\mapsto \varphi_Q\;\; \mbox{ and } \;\; \Sc([G]_Q)\to \Sc([G]_P),\;\; \varphi\mapsto E_Q^P(\varphi)$$
defined by
$$\displaystyle \varphi_Q(g)=\int_{[N_Q]} \varphi(ug) du \;\; \mbox{ and } \;\; E_Q^P(\varphi,g)=\sum_{\gamma\in Q(F)\backslash P(F)} \varphi(\gamma g)$$
respectively.

\begin{lemme}\label{lemma estimates constant term}
There is a constant $c>0$ such that for every $N\geqslant 0$,
$$\displaystyle f\mapsto \sup_{g\in [G]_P} \lVert g\rVert_{[G]_P}^{N}\delta_P(g)^{cN}\lvert f_P(g)\rvert$$
is a continuous semi-norm on $\Sc([G])$.
\end{lemme}

\begin{preuve}
As $[N_P]$ is compact, it suffices to show the existence of $N\geqslant 1$ and $c>0$ such that
$$\displaystyle \delta_P(g)^c \lVert g\rVert_{[G]_P}\ll \lVert g\rVert_{[G]}^N$$
for every $g\in G(\AAA)$. By \eqref{eq 0bis Function spaces}, we are easily reduced to show the existence of $N\geqslant 1$ and $c>0$ such that
$$\displaystyle \delta_P(a)^c \lVert a\rVert\ll \lVert ua\rVert_{[G]}^N$$
for every $(u,a)\in N_P(\AAA)\times A_0^\infty$. It is then equivalent to establish the following:
\begin{num}
\item\label{eq1 majorations} For every algebraic character $\chi$ of $A_0$, there exist $N\geqslant 1$ and $c>0$ such that
$$\displaystyle \delta_P(a)^c\lvert \chi(a)\rvert\ll \lVert ua\rVert_{[G]}^N,\;\mbox{for } (a,u)\in A_0^\infty\times [N_P].$$
\end{num}
Let $\chi$ be such a character and let $\delta_P^{alg}$ be the algebraic character of $M_P$ given by the determinant of the adjoint action on the Lie algebra of $N_P$. We have $\lvert \delta_P^{alg}(m)\rvert=\delta_P(m)$ for every $m\in M_P(\AAA)$ and we can find an integer $c\geqslant 1$, a rational representation $(V,\rho)$ of $G$ and a nonzero vector $v_0\in V$ such that $\rho(ua)v_0=\chi(a)\delta_P^{alg}(a)^{c} v_0$ for every $a\in A_0$ and $u\in N_P$. Fixing a basis $v_1,\ldots,v_d$ of $V$, for every $v=\lambda_1v_1+\ldots+\lambda_dv_d\in V_{\AAA}:=V\otimes_{F} \AAA$, we set
$$\displaystyle \lvert v\rvert_V=\prod_{v} \max(\lvert \lambda_1\rvert_v,\ldots,\lvert \lambda_d\rvert_v).$$
Then, there exists $N\geqslant 1$ such that $\lvert \rho(g)v_0\rvert_V\ll \lVert g\rVert_{G(\AAA)}^N$ for all $g\in G(\AAA)$ and moreover we have $\lvert v\rvert_V\geqslant 1$ for every $v\in V\setminus \{0 \}$. Therefore, we obtain
$$\displaystyle \lVert \gamma ua\rVert^N_{G(\AAA)} \gg \lvert \rho(\gamma ua)v_0\rvert_V=\lvert \chi(a)\rvert \delta_P(a)^c\lvert \rho(\gamma) v_0\rvert_V\geqslant \lvert \chi(a)\rvert \delta_P(a)^c$$
for $\gamma\in G(F)$, $a\in A_0^\infty$ and $u\in N_P(\AAA)$. The estimates \eqref{eq1 majorations} follows by replacing the left hand side in the inequality above by its the infimum over $\gamma\in G(F)$.
\end{preuve}
\end{paragr}

\begin{paragr}
For $Q\subset P$ two standard parabolic subgroup, we let $(\ago_Q^*)^{P++}$ be the set of $\lambda\in \ago_Q^*$ such that $\langle \alpha^\vee, \lambda\rangle>0$ for every $\alpha\in \Delta_0^P\setminus \Delta_0^Q$.

\begin{proposition}\label{prop:caract Schwartz}
Let $\varphi\in \tc([G]_P)$. Then, $\varphi\in \Sc([G]_P)$ if and only if for every standard parabolic subgroup $Q\subset P$ and every $\lambda\in (\ago_Q^*)^{P++}$ we have $\varphi_Q\in L^2_{2\rho_Q+\lambda}([G]_Q)^\infty$. Moreover, for each compact-open subgroup $J\subset G(\AAA_f)$, the semi-norms
$$\displaystyle \varphi\mapsto \lVert R(X)\varphi_Q\rVert_{L^2_{2\rho_Q+\lambda}}$$
for $Q\subset P$ a standard parabolic subgroup, $\lambda\in (\ago_Q^*)^{P++}$ and $X\in \uc(\ggo_\infty)$ generate the topology on $\Sc([G]_P)^J$.
\end{proposition}

\begin{preuve}
Let $\varphi\in \Sc([G]_P)$, $Q\subset P$ be a standard parabolic subgroup, $\lambda\in (\ago_Q^*)^{P++}$ and $X\in \uc(\ggo_\infty)$. By Cauchy-Schwarz we have
$$\displaystyle \int_{[G]_Q} \lvert R(X)\varphi_Q(x)\rvert^2 e^{\langle 2\rho_Q+\lambda,H_Q(x)\rangle }dx\leqslant \int_{[G]_P} \lvert R(X)\varphi(x)\rvert^2 \sum_{\gamma\in Q(F)\backslash P(F)} e^{\langle 2\rho_Q+\lambda,H_Q(\gamma x)\rangle } dx.$$
and by the convergence of Eisenstein series the inner sum above converges to a function which is essentially bounded by $\lVert x\rVert^N$ for some $N$. Therefore, $\varphi_Q\in L^2_{2\rho_Q+\lambda}([G]_Q)^\infty$ and this shows the direct implication.

We now prove the converse: let $\varphi\in \tc([G]_P)$ be such that for every standard parabolic subgroup $Q\subset P$ and every $\lambda\in (\ago_Q^*)^{P++}$ we have $\varphi_Q\in L^2_{2\rho_Q+\lambda}([G]_Q)^\infty$. Applying this assumption to $Q=P$ and $\lambda$ varying over fixed basis of $\ago_P^*$ and its opposite, from \eqref{basic inequality smooth L2 spaces} and \eqref{eq1 Function spaces} we get for every $N>0$ the estimates
\begin{equation}\label{eq1:L2 caract Schwartz}
\displaystyle \lvert \varphi(g)\rvert\ll \exp(\langle \rho_0, H_0(g)\rangle- N\lVert H_P(g)\rVert),\;\;\; g\in \sgo_P.
\end{equation}
Similarly, applying the hypothesis to the maximal parabolic subgroup $Q_\alpha\subset P$ associated to a simple root $\alpha\in \Delta_0^P$ and $\lambda$ a multiple of $\rho_Q^P$, from \eqref{basic inequality smooth L2 spaces} and \eqref{eq1 Function spaces} for every $N>0$ we obtain
\begin{equation}\label{eq2:L2 caract Schwartz}
\displaystyle \lvert \varphi_{Q_\alpha}(g)\rvert\ll \exp(\langle \rho_0, H_0(g)\rangle-N\langle \alpha, H_0(g)\rangle),\;\;\; g\in \sgo_Q.
\end{equation}
Moreover, by the approximation property of the constant term \cite[Lemma I.2.10]{MWlivre}, there exists $M>0$ such that for every $N>0$ we have
\begin{equation}\label{eq3:L2 caract Schwartz}
\displaystyle \lvert \varphi(g)-\varphi_{Q_\alpha}(g)\rvert\ll \exp(-N\langle \alpha, H_0(g)\rangle+M\lVert H_0(g)\rVert),\;\;\; g\in \sgo_P.
\end{equation}
Finally, there exist $\epsilon>0$ and $C>0$ such that $\max \{\langle \alpha, H_0(g)\rangle\mid\alpha\in \Delta_0^P\}  \cup \{ \lVert H_P(g)\rVert\}\geqslant \epsilon \lVert H_0(g)\rVert-C$ for all $g\in \sgo_P$. Combining this with \eqref{eq1:L2 caract Schwartz}, \eqref{eq2:L2 caract Schwartz} and \eqref{eq3:L2 caract Schwartz}, gives the estimates
$$\displaystyle \lvert \varphi(g)\rvert\ll \lVert g\rVert^{-N}, \; g\in \sgo_P,$$
for every $N>0$. Applying the same reasoning to derivatives of $\varphi$, this shows that $\varphi\in \Sc([G]_P)$.

For the last part of the statement, it suffices to notice that the linear map
$$\displaystyle \Sc(G(\AAA))\to \prod_{Q\subset P, \lambda\in (\ago_Q^*)^{P++}} L^2_{2\rho_Q+\lambda}([G]_Q)^\infty$$
$$\displaystyle \varphi\mapsto (\varphi_Q)_{Q}$$
is injective with a closed image (by the previous characterization) hence is a topological embedding by the open mapping theorem.
\end{preuve}
\end{paragr}

\subsection{Estimates on Fourier coefficients}

\begin{paragr}
Let $P$ be a standard parabolic subgroup of $G$, $\psi:\AAA/F\to \CC^\times$ be a non-trivial additive character and $\ell:N_P\to \mathbb{G}_a$ be an algebraic character. We set $\psi_\ell:=\psi\circ \ell_{\AAA}: [N_P]\to \CC^\times$. For $f\in C^\infty([G])$, we set
$$\displaystyle f_{N_P,\psi_\ell}(g)=\int_{[N_P]} f(ug) \psi_\ell(u)^{-1} du,\;\;\; g\in G(\AAA).$$
Let $N_{P,\der}$ denote the derived subgroup of $N_P$ and set $N_{P,\ab}=N_P/N_{P,\der}$ (a vector space over $F$). Then, $\ell$ can be seen as an element in the dual space $N_{P,\ab}^*$ that we also consider as an algebraic variety over $F$. The adjoint action of $M_P$ on $N_P$ induces one on $N_{P,\ab}^*$ that we denote by $\Ad^*$.

\begin{lemme}\label{lem1 majorations Fourier}
	\begin{enumerate}
\item There exists $c>0$ such that for every $N_1,N_2\geqslant 0$,
		$$\displaystyle f\mapsto \sup_{m\in M_P(\AAA)}\lVert \Ad^*(m^{-1})\ell\rVert_{N_{P,\ab}^*(\AAA)}^{N_1} \lVert m\rVert_{[M_P]}^{N_2} \delta_P(m)^{cN_2} \lvert f_{N_P,\psi_\ell}(m)\rvert$$
		is a continuous semi-norm on $\Sc([G])$.
\item For every $N_1,N_2\geqslant 0$,
		$$\displaystyle f\mapsto \sup_{m\in M_P(\AAA)} \lVert \Ad^*(m^{-1})\ell\rVert_{N_{P,\ab}^*(\AAA)}^{N_1} \lVert m\rVert_{[M_P]}^{-N_2} \lvert f_{N_P,\psi_\ell}(m)\rvert$$
		is a continuous semi-norm on $\tc_N([G])$.
	\end{enumerate}
\end{lemme}

\begin{preuve}
Bounding brutally under the integral sign, we have
$$\displaystyle \lvert f_{N_P,\psi_\ell}(g)\rvert\leqslant \lvert f\rvert_P(g)$$
for $f\in C^\infty([G])$ and $g\in G(\AAA)$. Let $N\geqslant 0$ and $J\subseteq G(\AAA_f)$ be a compact-open subgroup. By Lemma \ref{lemma estimates constant term} and \eqref{eq0bis Function spaces}, it suffices to show the existence of elements $X_1,\ldots,X_M\in \uc(\ggo_\infty)$ such that
\begin{align}\label{eq1 majorations FC}
\displaystyle \lvert f_{N_P,\psi_\ell}(m)\rvert\leqslant \lVert \Ad^*(m^{-1})\ell\rVert_{N_{P,\ab}^*(\AAA)}^{-N} \sum_{i=1}^M \lvert (R(X_i)f)_{N_P,\psi_\ell}(m)\rvert
\end{align}
for every $f\in C^\infty([G])^J$ and $m\in M_P(\AAA)$. Let $u\in N_P(\AAA)$. By definition of $\lVert .\rVert_{N_{P,\ab}^*(\AAA)}$, we are readily reduced to show the existence of $X_1,\ldots,X_M\in \uc(\ggo_\infty)$ such that
\begin{align}
\displaystyle \lvert f_{N_P,\psi_\ell}(m)\rvert\leqslant \lVert \ell(\Ad(m)u)\rVert_{\AAA}^{-1} \sum_{i=1}^M \lvert (R(X_i)f)_{N_P,\psi_\ell}(m)\rvert
\end{align}
for every $f\in C^\infty([G])^J$ and $m\in M_P(\AAA)$. This last claim is a consequence of the two following facts whose proofs are elementary and left to the reader.

\begin{num}
	\item For every non-Archimedean place $v$, there exists a constant $C_v\geqslant 1$ with $C_v=1$ for almost all $v$ such that $\lvert \ell(\Ad(m_v)u_v)\rvert_v>C_v$ implies $f_{N_P,\psi_\ell}(m)=0$ for every $f\in C^\infty([G])^J$ and $m\in M_P(\AAA)$.
	
	\item Let $v$ be an Archimedean place and let $X\in \ggo_v$ be such that $u_v=e^X$ . Then, we have $(R(X)f)_{N_P,\psi_\ell}(m)=d\psi_v(\ell(\Ad(m)u)_v) f_{N_P,\psi_\ell}(m)$ for all $f\in C^\infty([G])$ and $m\in M_P(\AAA)$ where $d\psi_v: F_v\to i\RR$ is the differential of $\psi_{v}$ at the origin.
\end{num}
\end{preuve}
\end{paragr}

\begin{paragr}
Let $n\geqslant 1$ be a positive integer. We let $\GL_n$ acts on $F^n$ by right multiplication and we denote by $e_n=(0,\ldots,0,1)$ the last element of the standard basis of $F^n$. We also denote by $\pc_n$ the {\em mirabolic} subgroup of $\GL_n$, that is the stabilizer of $e_n$ in $\GL_n$. We identify $A_{\GL_n}$ with $\GG_m$, and thus $A_{\GL_n}^\infty$ with $\RR_{>0}$, in the usual way. The next lemma will be used in conjunction with Lemma \ref{lem1 majorations Fourier} to show the convergence of various Zeta integrals.

\begin{lemme}\label{lem conv GLn}
Let $C>1$. Then, for $N_1\gg_C 1$ and $N_2\gg_C 1$ the integral
	$$\displaystyle \int_{\pc_n(F)\backslash \GL_n(\AAA)\times \RR_{>0}} \lVert ag\rVert^{-N_1}_{[\GL_n]} \lVert e_n g\rVert_{\AAA^n}^{-N_2} \lvert \det g\rvert^s dadg$$
	converges for $s\in \cH_{]1,C[}$ uniformly on every (closed) vertical strip.
\end{lemme}

\begin{preuve}
The integral of the lemma can be rewritten as
\begin{equation}\label{eq2 majorations FC}
\displaystyle \int_{[\GL_n]} \lVert g\rVert^{-N_1}_{[\GL_n]} \int_{\RR_{>0}} \sum_{\xi\in F^n\setminus \{ 0\}} \lVert \xi ag\rVert_{\AAA^n}^{-N_2} \lvert \det ag\rvert^s da dg.
\end{equation}
There exists $N_3>0$ such that $\lVert v\rVert_{\AAA^n}\ll \lVert vg\rVert_{\AAA^n}^{N_3} \lVert g\rVert_{\GL_n(\AAA)}^{N_3}$ for $(v,g)\in \AAA^n\times \GL_n(\AAA)$. Therefore, the inner integral above is essentially bounded by
$$\displaystyle \lvert \det g\rvert^{\Re(s)} \lVert g\rVert_{\GL_n(\AAA)}^{N_2} \int_{\RR_{>0}} \sum_{\xi\in F^n\setminus \{ 0\}} \lVert a\xi\rVert_{\AAA^n}^{-N_2/N_3} \lvert a\rvert^{ns} da$$
hence, for $1<\Re(s)<C$, by
$$\displaystyle \lVert g\rVert_{\GL_n(\AAA)}^{N_2+N_4} \int_{\RR_{>0}} \sum_{\xi\in F^n\setminus \{ 0\}} \lVert a\xi\rVert_{\AAA^n}^{-N_2/N_3} \lvert a\rvert^{ns} da$$
for some $N_4>0$. However, since the inner integral in \eqref{eq2 majorations FC} is left invariant by $\GL_n(F)$, as a function of $g$, we may replace $\lVert g\rVert_{\GL_n(\AAA)}$ in the estimate above by $\lVert g\rVert_{[\GL_n]}$. As for $N\gg 1$ we have $\displaystyle \int_{[\GL_n]} \lVert g\rVert_{[\GL_n]}^{-N}dg<\infty$ \cite[Proposition A.1.1 (vi)]{RBP}, it only remains to show that for $N\gg 1$ the integral
$$\displaystyle \int_{\RR_{>0}} \sum_{\xi\in F^n\setminus \{ 0\}} \lVert a\xi\rVert_{\AAA^n}^{-N} \lvert a\rvert^{ns} da$$
converges for $1<\Re(s)<C$ uniformly in vertical strips. This is an easy consequence of the following claim:

\begin{num}
	\item\label{eq3 majorations FC} For every $k\geqslant n$, if $N$ is sufficiently large we have
	$$\displaystyle \sum_{\xi\in F^n\setminus \{ 0\}} \lVert a\xi\rVert_{\AAA^n}^{-N}\ll \lvert a\rvert^{-k},\; a\in \RR_{>0}.$$
\end{num}

There exists $M_0\geqslant 1$ such that
$$\displaystyle \lvert a\rvert=\max_{1\leqslant i\leqslant n}(\lvert a\xi_i\rvert)\ll \lVert a\xi\rVert_{\AAA^n}^{M_0}$$
for $(a,\xi)\in \RR_{>0} \times (F^n\setminus \{ 0\})$. Therefore, we just need to prove \eqref{eq3 majorations FC} when $k=n$. Let $C\subset \AAA^n$ be a compact subset which surjects onto $\AAA^n/F^n$. There exists $M_1\geqslant 0$ such that
$$\displaystyle \lVert a\xi+av\rVert_{\AAA^n}\ll \lVert a\xi\rVert_{\AAA^n}^{M_1} \max(1,\lvert a\rvert)^{M_1}\ll \lVert a\xi\rVert_{\AAA^n}^{M_1+M_0M_1}$$
for $(a,\xi,v)\in \RR_{>0} \times (F^n\setminus \{ 0\})\times C$. Hence, for every $N'>0$ if $N$ is sufficiently large, we have
\[\begin{aligned}
\displaystyle \sum_{\xi\in F^n\setminus \{ 0\}} \lVert a\xi\rVert_{\AAA^n}^{-N} & \ll \int_C \sum_{\xi\in F^n\setminus \{ 0\}} \lVert a\xi+av\rVert_{\AAA^n}^{-N'}dv \\
& \ll \int_{\AAA^n/F^n} \sum_{\xi\in F^n} \lVert a\xi+av\rVert_{\AAA^n}^{-N'}dv=\lvert a\rvert^{-n}\int_{\AAA^n} \lVert v\rVert_{\AAA^n}^{-N'} dv
\end{aligned}\]
for $a\in \RR_{>0}$. The last integral above is absolutely convergent when $N'\gg 1$ \cite[Proposition A.1.1 (vi)]{RBP} and the claim \eqref{eq3 majorations FC} follows.
\end{preuve}
\end{paragr}

\subsection{Automorphic forms and representations}

\begin{paragr}
Let $P$ be a standard parabolic subgroup of $G$. The space $\Ac_P(G)$ of automorphic forms on $[G]_P$ is defined as the subspace of $\zc(\ggo_\infty)$-finite functions in $\tc([G]_P)$. The subspace $\Ac_{P,\cusp}(G)$ of cuspidal automorphic forms consists of the $\varphi\in \Ac_P(G)$ such that for every proper standard parabolic subgroup $Q\subseteq P$ we have $\varphi_Q=0$. 

For $\Jc\subset \zc(\ggo_\infty)$ an ideal of finite codimension, we denote by $\Ac_{P,\Jc}(G)$ the subspace of automorphic forms $\varphi\in \Ac_P(G)$ such that $R(z)\varphi=0$ for every $z\in \Jc$ and we set $\Ac_{P,\cusp,\Jc}(G)=\Ac_{P,\Jc}(G)\cap \Ac_{P,\cusp}(G)$. Then, there exists $N\geqslant 1$ such that $\Ac_{P,\Jc}(G)$ is a closed subspace of $\tc_N([G]_P)$ and we equip $\Ac_{P,\Jc}(G)$ with the induced topology from $\tc_N([G]_P)$ (this topology does not depend on the choice of $N$ by the open mapping theorem). Similarly, $\Ac_{P,\cusp,\Jc}(G)$ is a closed subspace of $\Sc([G]_P)$ and we equip $\Ac_{P,\cusp,\Jc}(G)$ with the induced topology (which also coincides with the topology induced from $\Ac_{P,\Jc}(G)$). We have $\Ac_{P}(G)=\bigcup_{\Jc} \Ac_{P,\Jc}(G)$ and $\Ac_{P,\cusp}(G)=\bigcup_{\Jc} \Ac_{P,\cusp,\Jc}(G)$ where $\Jc$ runs over all ideals of finite codimensions in $\zc(\ggo_\infty)$ and we equip these spaces with the convex inductive limit topology (these are strict LF spaces). We shall also consider the closed subspace  $\Ac_{P,\disc}(G)$ of $\Ac_P(G)$ of automorphic forms on $[G]_{P,0}$ that are square-integrable. Its topology is also induced from the strict LF space $L^2([G]_P)^\infty$.

For $P=G$, we simply set $\Ac(G)=\Ac_G(G)$, $\Ac_{\disc}(G)=\Ac_{G,\disc}(G)$ and $\Ac_{\cusp}(G)=\Ac_{G,\cusp}(G)$. 
\end{paragr}

\begin{paragr}\label{S:automorphic repn}
By a {\em cuspidal automorphic representation} (resp. discrete) $\sigma$ of $M_P(\AAA)$ we mean a topologically irreducible subrepresentation of $\Ac_{\cusp}(M_P)$ (resp. $\Ac_{\disc}(M_P)$). Let $\sigma$ be a cuspidal or discrete automorphic representation of $M_P(\AAA)$. For every $\lambda\in \ago_{P,\CC}^*$, the twist $\sigma_\lambda=\sigma\otimes \lambda$ is defined as the space of automorphic forms
$$m\in [M_P]\mapsto \exp(\langle \lambda, H_P(m)\rangle) \varphi(m)$$
for $\varphi\in \sigma$. If $\sigma$ is cuspidal, $\sigma_\lambda$ is again a cuspidal automorphic representation. We denote by $\Ac_{\sigma}(M_P)$ the $\sigma$-isotypic component of $\Ac_{\cusp}(M_P)$ (resp. $\Ac_{\disc}(M_P)$) i.e. the sum of all cuspidal (resp. discrete) automorphic representations of $M_P(\AAA)$ that are isomorphic to $\sigma$. We let $\Pi=I_{P(\AAA)}^{G(\AAA)}(\sigma)$ (resp. $\Ac_{P,\sigma}(G)=I_{P(\AAA)}^{G(\AAA)}(\Ac_\sigma(M_P))$) be the {\em normalized smooth induction} of $\sigma$ (resp. $\Ac_\sigma(M_P)$) that we identify with the space of forms $\varphi\in \Ac_{P}(G)$ such that
$$m\in [M_P]\mapsto \exp(-\langle \rho_P, H_P(m)\rangle) \varphi(mg)$$
belongs to $\sigma$ (resp. $\Ac_\sigma(M_P)$) for every $g\in G(\AAA)$. We have $\Ac_{P,\sigma}(G)\subset \Ac_{P,\cusp}(G)$ if $\sigma$ is cuspidal and $\Ac_{P,\sigma}(G)\subset \Ac_{P,\disc}(G)$ if $\sigma$ is discrete.
The algebra $\Sc(G(\AAA))$ acts on $\Ac_{P,\sigma}(G)$ by right convolution. For every $\lambda\in \ago_{P,\CC}^*$, we denote by $I(\la)$ the action on $\Ac_{P,\sigma}(G)$ we get by transport from the action of $\Sc(G(\AAA))$ on $\Ac_{P,\sigma_\la}$ and the identification $\Ac_{P,\sigma}\to \Ac_{P,\sigma_\la}$.

The spaces of $\sigma$ and $\Pi$ naturally carry topologies making them into strict LF spaces. More precisely, for every compact-open subgroups $J_M\subset M_P(\AAA_f)$ and $J\subset G(\AAA_{f})$, $\sigma^{J_M}$ and $\Pi^J$ are SF representations of $M_P(F_\infty)$ and $G(F_\infty)$ in the sense of \cite{BK} respectively and when $J'_M\subseteq J_M$, $J'\subseteq J$ are smaller compact-open subgroups the inclusions $\sigma^{J'_M}\subseteq \sigma^{J_M}$, $\Pi^{J'}\subseteq \Pi^J$ are closed embeddings.

If the central character of $\sigma$ is unitary, we equip $\Pi$ with the {\em Petersson inner product}
$$\| \varphi\|^2_{\Pet}= \displaystyle \langle \varphi,\varphi\rangle_{\Pet}=\int_{[G]_{P,0}} \lvert \varphi(g)\rvert^2 dg,\; \varphi\in \Pi.$$
\end{paragr}

\begin{paragr}[Eisenstein series.] ---\label{S:Eis-series}
  Let $P$ be  a standard parabolic subgroup of $G$. For any $\varphi\in \Ac_P(G)$, $g\in G(\AAA)$ and $\la\in \ago_{P,\CC}^{*}$, we denote by
  \begin{align*}
    E(g,\varphi,\la)=\sum_{\delta\in P(F)\back G(F)} \exp(\bg \la,H_P(\delta g)\bd \varphi(\delta g)
  \end{align*}
the Eisenstein series where the left-hand side is obtained from the analytic continuation of the right-hand side which is only defined for $\Re(\la)$ in a suitable cone.
\end{paragr}

\begin{paragr}
  Let $P$ and $Q$ be  standard parabolic subgroups of $G$. For any $w\in W(P,Q)$ and $\la\in  \ago_{P,\CC}^{*}$, we have the intertwining operator
  \begin{align*}
    M(w,\la):\Ac_P(G)\to \Ac_{Q}(G)
  \end{align*}
defined by analytic continuation from the integral
\begin{align*}
   (M(w,\la)\varphi)(g)=   \exp(-\bg w\la,H_P(g)\bd) \int_{  (N_Q\cap w N_Pw^{-1})(\AAA) \back N_Q(\AAA) } \exp(\bg \la,H_P(w^{-1}ng)\bd) \varphi(w^{-1}ng)\, dn.
\end{align*}
\end{paragr}

\begin{paragr}\label{autom repns product}
Assume that $G=G_1\times G_2$ where $G_1$ and $G_2$ are connected reductive groups over $F$. This induces decompositions $P=P_1\times P_2$, $M_P=M_{P_1}\times M_{P_2}$ and there exist two, uniquely determined, cuspidal automorphic representations $\sigma_1$, $\sigma_2$ of $M_{P_1}(\AAA)$ and $M_{P_2}(\AAA)$ respectively such that, setting $\Pi_1=I_{P_1(\AAA)}^{G_1(\AAA)}(\sigma_1)$ and $\Pi_2=I_{P_2(\AAA)}^{G_2(\AAA)}(\sigma_2)$, for every compact-open subgroups $J_1\subseteq G_1(\AAA_{f})$, $J_2\subseteq G_2(\AAA_{f})$ (resp. $J_1\subseteq M_{P_1}(\AAA_{f})$, $J_2\subseteq M_{P_2}(\AAA_{f})$), setting $J=J_1\times J_2$, there is a topological isomorphism
\begin{align}\label{isom autom repns product}
\displaystyle \Pi_1^{J_1}\widehat{\otimes} \Pi_2^{J_2}\simeq \Pi^{J}\;\; (\mbox{resp. } \sigma_1^{J_1}\widehat{\otimes} \sigma_2^{J_2}\simeq \sigma^{J})
\end{align}
sending $\varphi_1\otimes \varphi_2\in \Pi_1^{J_1}\otimes \Pi_2^{J_2}$ (resp. $\varphi_1\otimes \varphi_2\in \sigma_1^{J_1}\otimes \sigma_2^{J_2}$) to the function $(g_1,g_2)\mapsto \varphi_1(g_1)\varphi_2(g_2)$. We then write $\Pi=\Pi_1\boxtimes \Pi_2$ and $\sigma=\sigma_1\boxtimes \sigma_2$ respectively.
\end{paragr}

\begin{paragr}
Assume now that $G$ is quasi-split. Let $\psi_N:N_0(\AAA)\to \CC^\times$ be a continuous non-degenerate character which is trivial on $N_0(F)$. If the representation $\Pi$ is $\psi_N$-generic, i.e. if it admits a continuous nonzero linear form $\ell:\Pi\to \CC$ such that $\ell\circ \Pi(u)=\psi_N(u)\ell$ for every $u\in N(\AAA)$, it is (abstractly) isomorphic to its {\em Whittaker model}
$$\displaystyle \wc(\Pi,\psi_N)=\{g\in G(\AAA)\mapsto \ell(\Pi(g)\varphi)\mid\; \varphi\in \Pi \}.$$
We equip this last space with the topology coming from $\Pi$ (thus it is a strict LF space).

If we are moreover in the situation of \S \ref{autom repns product}, there are decompositions $N_0=N_{0,1}\times N_{0,2}$, $\psi_N=\psi_1\boxtimes \psi_2$ and the isomorphism \eqref{isom autom repns product} induces one between Whittaker models
$$\displaystyle \wc(\Pi_1,\psi_1)^{J_1}\widehat{\otimes} \wc(\Pi_2,\psi_2)^{J_2}\simeq \wc(\Pi,\psi_N)^{J}.$$
\end{paragr}

\subsection{Relative characters}

\begin{paragr}
    Let $B$ a $G(F_\infty)$-invariant nondegenerate symmetric bilinear form on $\ggo_\infty$. We assume that the restriction of $B$ to  $\mathfrak{k}_\infty$ is negative and the  restriction of $B$ to  the orthogonal complement of   $\mathfrak{k}_\infty$  is positive. Let $(X_i)_{i\in I}$ be an orthonormal basis of   $\mathfrak{k}_\infty$  relative to $-B$. Let $C_{K}=-\sum_{i\in I}X_i^2$ : this is a ``Casimir element'' of $\uc(\mathfrak{k}_\infty)$.
  \end{paragr}

  \begin{paragr}
    Let $\hat K_\infty$ and $\hat K$ be respectively the sets of isomorphism classes of irreducible unitary representations of  $K_\infty$ and of $K$. 
  \end{paragr}

\begin{paragr}\label{S:K0ON}
Let $\pi$ be a discrete automorphic representation of $M_P$.   For any $\tau\in \hat K$, let $\Ac_{P,\pi}(G,\tau)$ be the  (finite dimensional) subspace of  functions in $\Ac_{P,\pi}(G)$  which transform under $K$ according to $\tau$. A $K$-basis $\bc_{P,\pi}$ of  $\Ac_{P,\pi}(G)$ is by definition the union over of $\tau\in \hat K$ of orthonormal bases $\bc_{P,\pi,\tau}$ of $\Ac_{P,\pi}(G,\tau)$ for the Pertersson inner product.
\end{paragr}

\begin{paragr}
  Let
  $$B: \Ac_{P,\pi}(G) \times \Ac_{P,\pi}(G) \to \CC$$
  be a continuous sesquilinear form.

  \begin{proposition}\label{prop:car-relJB}
    Let $\om$ be a compact subset of $\ago_{P}^{G,*}$.
    \begin{enumerate}
    \item Let $f\in \Sc(G(\AAA))$ and  $\bc_{P,\pi}$  be a $K$-basis of $\Ac_{P,\pi}(G)$. The sum
    \begin{align}\label{eq: sumJB proj tensor product}
    \sum_{\varphi\in \bc_{P,\pi}  } I_P(\la,f) \varphi\otimes \overline{\varphi}
    \end{align}
    converges absolutely in the completed projective tensor product $\Ac_{P,\pi}(G)\widehat{\otimes}\overline{\Ac_{P,\pi}(G)}$ uniformly for $\la\in \ago_{P,\CC}^{G,*}$ such that $\Re(\la)\in \om$. In particular, the sum 
    \begin{align}\label{eq:sumJB}
      J_B(\la,f)=\sum_{\varphi\in \bc_{P,\pi}  } B(I_P(\la,f) \varphi,\varphi)
    \end{align}
    is absolutely convergent uniformly for $\la\in \ago_{P,\CC}^{G,*}$ such that $\Re(\la)\in \om$. Moreover these sums do not depend on the choice of $\bc_{P,\pi}$.
  \item The map
    \begin{align*}
      f\mapsto  J_B(\la, f)
    \end{align*}
is a continuous linear form on  $\Sc(G(\AAA))$. More precisely for $C\subset G(\AAA_f)$ a compact subset and $K_0$ as above,  there exist $c>0$ and a continuous semi-norm $\|\cdot\|$ on $\Sc(G(\AAA),C,K_0)$  such that for all  $\la\in \ago_{P,\CC}^{G,*}$ such that $\Re(\la)\in \om$ and $f\in \Sc(G(\AAA),C,K_0)$ we have
\begin{align*}
  |J_B(\la, f)|\leq c\|f\|.
\end{align*}
\end{enumerate}
\end{proposition}

\begin{remarque}\label{rq:classeCr}
  An examination of the proof below show that the assertion 2 also holds \emph{mutatis mutandis} if $f\in C_c^r(G(\AAA))$ with $r$ large enough. The semi-norm is then taken among the norms $\|\cdot\|_{r,X,Y}$ for which the sum of the degrees of $X $ and $Y$ is less than $r$.
\end{remarque}

\begin{preuve} By definition of the projective tensor product topology, it suffices to show the following: for every continuous semi-norm $p$ on $\Ac_{P,\pi}(G)$, the series
$$\displaystyle \sum_{\varphi\in \bc_{P,\pi}  } p(I_P(\la,f) \varphi) p(\varphi)$$
is absolutely convergent uniformly for $\la\in \ago_{P,\CC}^{G,*}$ such that $\Re(\la)\in \om$. Let $K_0\subset K^\infty$ be a normal open compact subgroup by which $f$ is biinvariant. The series above can be rewritten as
  \begin{align} \label{eq:sumJB2}
      \sum_{\tau\in \hat K}\sum_{\varphi\in \bc_{P,\pi,\tau}  } p(I_P(\la,f) \varphi) p(\varphi).
  \end{align}
  and the only representations that contribute to \eqref{eq:sumJB2} are the representations $\tau$ that admits $K_0$-invariant vectors. Then the elements $\varphi\in \bc_{P,\pi,\tau} $ are automatically $K_0$-fixed.  For any $\tau\in \hat K_\infty$, let  $\Ac_{P,\pi}(G,K_0,\tau)\subset  \Ac_{P,\pi}(G)$ be the subspace of functions that are right $K_0$-invariant and transform under $K_\infty$ according to $\tau$. By  \cite{WRRG2} §10.1, there exist $c>0$ and  an integer $r$ such that for all $\varphi\in \Ac_{P,\pi}(G)^{K_0}$ we have
    \begin{align*}
       p(\varphi)\leq c \|R(1+C_K)^r\varphi\|_{\Pet}.
    \end{align*}
    For any $\tau\in \hat K_\infty$ or $\hat K$, let $\la_\tau\geq 0$ be the eigenvalue of $C_K$ acting on $\tau$. We have 
    \begin{align*}
  \| R(1+C_K)^r I_P(\la,f)\varphi\|_{\Pet}&=\| I_P(\la, L( (1+C_K)^r)f)\varphi \|_{\Pet}\\
&=(1+\la_\tau)^{-N} \|I_P(\la,f_{r,N})\varphi\|_{\Pet}
    \end{align*}
    where $f_{r,N}=R((1+C_K)^N)L((1+C_{K}))^rf $.
    Let $C\subset G(\AAA_f)$ be a compact subset.  There exists $c_1>0$ and a semi-norm $\|\cdot\|$ on $\Sc(G(\AAA),C,K_0)$ (among those of  §\ref{S:Schw}) such that for any $f\in \Sc(G(\AAA),C,K_0)$  and $\la\in \ago_{P,\CC}^{G,*}$ such that $\Re(\la)\in \om$ we have 
\begin{align*}
  \|I_P(\la,f_{r,N})\varphi\|_{P,\pi}\leq c_1 \|f\|\|\varphi\|_{\Pet}.
\end{align*}
Thereby we are reduced to prove for large enough $N$ the convergence of
\begin{align}\label{eq:cv-partielle}
 \sum_{\tau\in \hat K_\infty}   (1+\la_\tau)^{r-N} \dim(\Ac^\infty_{P,\pi}(G,K_0,\tau)).
    \end{align}
However there exist $c_2>0$ and $m\geq 1$ such that $\dim(\Ac^\infty_{P,\pi}(G,K_0,\tau))\leq c_2(1+\la_\tau)^m$ (see the proof of \cite{Muller-sing} lemma 6.1). So the convergence of \eqref{eq:cv-partielle}is  reduced to that  of $   \sum_{\tau\in \hat K_\infty}   (1+\la_\tau)^{-N} $ which is well-known.

Finally it is easy to show that $J_B(f)$ does not depend on the choice of the basis $\bc_{P,\pi,\tau}$.
\end{preuve}

\begin{proposition}\label{prop:dirac}
  Let $K_0\subset K^\infty$ be a normal   open compact subgroup.  For any integer $m\geq 1$ there exist $Z\in \uc(\ggo_\CC)$, $g_1\in \Cc(G(\AAA))$ and $g_2\in C_c^m(G(\AAA))$ such that
  \begin{itemize}
  \item $Z$, $g_1$ and $g_2$ are  invariant under $K_\infty$-conjugation ;
  \item $g_1$ and $g_2$ are $K_0$-biinvariant ;
  \item  for any $f\in \Sc(G(\AAA))$ that is  $K_0$-biinvariant we have:
    \begin{align*}
      f=f*g_1+(f*Z)*g_2.
    \end{align*}
  \end{itemize}
For large enough $m$, we have
\begin{align*}
   J_B(\la,f)=\sum_{\varphi\in \bc_{P,\pi}  } B(I_P(\la,f) \varphi,I_P(\bar{\la},g_1^\vee)\varphi)+ \sum_{\varphi\in \bc_{P,\pi}  } B(I_P(\la,f*Z) \varphi,I_P(\bar{\la},g_2^\vee)\varphi)
\end{align*}
where the sums are absolutely convergent and $g_i^\vee(x)=\overline{g_i(x^{-1})}$.
  \end{proposition}

  \begin{preuve}
    The first part of the proposition is lemma 4.1 and corollary 4.2 of \cite{ar1}. Once we have noticed that the operators $I_P(\la,g_i)$ preserve the spaces  $\Ac_{P,\pi}(G,\tau)$, the second part results  from an easy computation in a finite dimensional space. 
  \end{preuve}
\end{paragr}

\subsection{Decomposition according to cuspidal data and automorphic kernels}\label{ssec:cuspidal-data} 

\begin{paragr}[Cuspidal data.] ---\label{S:cuspidal-data} 
Let $\underline{\Xgo}(G)$ be the set of pairs $(M_P,\sigma)$ where $P$ is a standard parabolic subgroup of $G$ and $\sigma$ is an isomorphism class of cuspidal automorphic representations of $M_P(\AAA)$ with central character trivial on $A_P^\infty$. We let $\Xgo(G)$ be the quotient of $\underline{\Xgo}(G)$ by the equivalence relation defined as follows: $(M_P,\sigma)\sim (M_Q,\tau)$ if there exists $w\in W(P,Q)$ such that $w\sigma w^{-1}\simeq \tau$. We call $\Xgo(G)$ the set of {\em cuspidal data} for $G$. For every standard parabolic subgroup $P$ of $G$, the natural inclusion $\underline{\Xgo}(M_P)\subset \underline{\Xgo}(G)$ descends to a finite map $\Xgo(M_P)\to \Xgo(G)$. For $\chi\in \Xgo(G)$ represented by a pair $(M_P,\sigma)$, we denote by $\chi^\vee$ the cuspidal datum associated to $(M_P,\sigma^\vee)$ where $\sigma^\vee$ stands for the complex conjugate of $\sigma$.
\end{paragr}

\begin{paragr}[Langlands decomposition]
For $(M_P,\sigma)\in \underline{\Xgo}(G)$, we let $\Sc_{\sigma}([G]_P)$ be the subspace of $\varphi\in \Sc([G]_P)$ such that the function
$$\displaystyle \varphi_\lambda(g):=\int_{A_{P}^\infty} \exp(-\langle \rho_P+\lambda, H_{P}(a)\rangle)\varphi(ag) da, \; g\in [G]_P,$$
belongs to $\Ac_{P, \sigma_{\lambda}}(G)$ for every $\lambda\in \ago_{P,\CC}^*$.
	
Let $P\subset G$ a standard parabolic subgroup, $\chi\in \Xgo(G)$ be a cuspidal datum and $\{(M_{Q_i},\sigma_i)\mid i\in I \}$ be the (possiby empty but finite) inverse image of $\chi$ in $\underline{\Xgo}(M_P)$. We define $L^2_\chi([G]_P)$ as the closure of the subspace $\displaystyle \sum_{i\in I} E_{Q_i}^P(\Sc_{\sigma_i}([G]_{Q_i}))$ in $L^2([G]_P)$. More generally, for $w$ be a weight on $\ago_P$ (see \S \ref{S:L2spaces}), we let $L^2_{w,\chi}([G]_{P})$ be the closure of $\displaystyle \sum_{i\in I} E_{Q_i}^P(\Sc_{\sigma_i}([G]_{Q_i}))$ in $L^2_w([G]_{P})$. We define similarly a subspace $L^2_\chi([G]_{P,0})\subset L^2([G]_{P,0})$. By Langlands (see e.g. \cite[Proposition II.2.4]{MWlivre}), we have Hilbert decompositions
\begin{align}
\label{eq:langlands}
\displaystyle L^2([G]_P)=\widehat{\bigoplus_{\chi\in \Xgo(G)}} L^2_\chi([G]_P) \text{ and } L^2([G]_{P,0})=\widehat{\bigoplus_{\chi\in \Xgo(G)}} L^2_\chi([G]_{P,0})  .
\end{align}
More generally, for every weight $w$ on $\ago_P$ we have a Hilbert decomposition
\begin{equation}\label{eq:weighted langlands}
\displaystyle L^2_w([G]_P)=\widehat{\bigoplus_{\chi\in \Xgo(G)}} L^2_{w,\chi}([G]_P).
\end{equation}
Let $\Xgo\subseteq \Xgo(G)$ be a subset and $w$ be a weight on $\ago_P$. We set
$$\displaystyle L^2_{w,\Xgo}([G]_P):=\widehat{\bigoplus_{\chi\in \Xgo}} L^2_\chi([G]_P),\; L^{2,\Xgo}_w([G]_P):=\widehat{\bigoplus_{\chi\in \Xgo(G)\setminus \Xgo}} L^2_{w,\chi}([G]_P)$$
and when $w=0$, we just drop the index $w$. We have
\begin{num}
\item\label{compatibility projections weightedL2} Let $w$ and $w'$ be two weights on $\ago_P$. The orthogonal projections $L^2_w([G]_P)\to L^2_{w,\Xgo}([G]_P)$ and $L^2_{w'}([G]_P)\to L^2_{w',\Xgo}([G]_P)$ coincide on the intersection $L^2_w([G]_P)\cap L^2_{w'}([G]_P)$.
\end{num}
Indeed, we have $L^2_w([G]_P)\cap L^2_{w'}([G]_P)=L^2_{w''}([G]_P)$ where $w''=\max(w,w')$ is a weight on $\ago_P$ and this allows to restrict to the case where $w\leqslant w'$. The claim then follows from the fact, easy to see from the definition, that the natural inclusion $L^2_{w'}([G]_P)\subset L^2_{w}([G]_P)$ sends $L^2_{w', \Xgo}([G]_P)$ (resp. $L^{2,\Xgo}_{w'}([G]_P)$) into $L^2_{w,\Xgo}([G]_P)$ (resp. $L^{2,\Xgo}_{w}([G]_P)$).

We will denote by $\varphi\mapsto \varphi_{\Xgo}$ the orthogonal projection $L^2_w([G]_P)\to L^2_{w,\Xgo}([G]_P)$ (by \eqref{compatibility projections weightedL2} such a notation shouldn't lead to any confusion). These projections are $G(\AAA)$-equivariant and so preserve the subspaces of smooth vectors.
\end{paragr}

\begin{paragr}
Let again $\Xgo\subseteq \Xgo(G)$ be a subset. We set
$$\displaystyle \Sc_{\Xgo}([G]_P):=\Sc([G]_P)\cap L^2_{\Xgo}([G]_P),\;\; \Sc^{\Xgo}([G]_P):=\Sc([G]_P)\cap L^{2,\Xgo}([G]_P),\;\; \cC_{\Xgo}([G]_P):=\cC([G]_P)\cap L^2_{\Xgo}([G]_P).$$
We also define $\tc_{\Xgo}([G]_P)$ (resp. $\tc_{N,\Xgo}([G]_P)$ for $N>0$) as the orthogonal of $\Sc^{\Xgo}([G]_P)$ in $\tc([G]_P)$ (resp. in $\tc_{N}([G]_P)$). Then, $\tc_{N,\Xgo}([G]_P)$ is a closed subspace of $\tc_{N}([G]_P)$ hence is a strict LF space. We equip $\tc_{\Xgo}([G]_P)=\bigcup_{N>0} \tc_{N,\Xgo}([G]_P)$ with the inductive limit locally convex topology (it is a LF space).

	\begin{proposition}\label{prop decomposition cuspidal data}
		Let $\Xgo\subseteq \Xgo(G)$ be a subset and $P\subseteq G$ be a standard parabolic subgroup.
		\begin{enumerate}
			\item For every standard parabolic subgroup $Q\subset P$, we have $E_Q^P(\Sc_{\Xgo}([G]_Q))\subset \Sc_{\Xgo}([G]_P)$ and $E_Q^P(\Sc^{\Xgo}([G]_P))\subset \Sc^{\Xgo}([G]_P)$.
			
			\item For every $\varphi\in \Sc([G]_P)$, we have $\varphi_{\Xgo}\in \Sc_{\Xgo}([G]_P)$. Moreover, the linear map $\varphi\in \Sc([G]_P)\mapsto \varphi_{\Xgo}\in \Sc_{\Xgo}([G]_P)$ is continuous and the series
			\begin{equation}\label{eq:series cuspidal projections}
			\displaystyle \sum_{\chi\in \Xgo(G)} \varphi_\chi
			\end{equation}
			converges absolutely (see \S \ref{S:absolute convergence}) to $\varphi$ in $\Sc([G]_P)$.
			
			\item For every $\varphi\in \tc([G]_P)$, there exists an unique function $\varphi_{\Xgo}\in \tc_{\Xgo}([G]_P)$ such that $\langle \varphi,\psi_{\Xgo}\rangle_{[G]_P}=\langle \varphi_{\Xgo},\psi\rangle_{[G]_P}$ for every $\psi\in \Sc([G]_P)$. Moreover, the linear map $\varphi\in \tc([G]_P)\mapsto \varphi_{\Xgo}\in \tc_{\Xgo}([G]_P)$ is $A_P^\infty\times G(\AAA)$-equivariant and continuous and for every parabolic subgroup $Q\subseteq P$ and $\varphi\in \tc([G]_P)$ we have
			\begin{equation}\label{eq:cuspidal projection and constant term}
			(\varphi_{\Xgo})_Q=(\varphi_Q)_{\Xgo}.
			\end{equation}
			\item For every compact-open subgroup $J\subset G(\AAA_f)$ and $N>0$, there exists $M>N$ such that for $\varphi\in \tc_N([G]_P)^J$ the series \eqref{eq:series cuspidal projections} converges absolutely to $\varphi$ in $\tc_M([G]_P)$.
		\end{enumerate}
	\end{proposition}
	
	\begin{preuve}
		\begin{enumerate}
			\item Up to replacing $\Xgo$ by its complement in $\Xgo(G)$, the two inclusions are equivalent and so we just need to prove the first. By definition of $\Sc_{\Xgo}([G]_P)$, we need to establish that $E_Q^P(\Sc_{\Xgo}([G]_Q))$ is orthogonal to $\Sc^{\Xgo}([G]_P)$ or, by adjunction, that for every $\varphi\in \Sc^{\Xgo}([G]_P)$, $\varphi_Q$ is orthogonal to $\Sc_{\Xgo}([G]_Q)$. Denote by $\{(M_{Q_i},\sigma_i)\mid i\in I\}$ the (possiby infinite) inverse image of $\Xgo$ in $\underline{\Xgo}(M_Q)$ and let $\varphi\in \Sc^{\Xgo}([G]_P)$. By adjunction again, and the definition of $\Sc^{\Xgo}([G]_P)$, $\varphi_Q$ is orthogonal to $\displaystyle \sum_{i\in I} E_{Q_i}^Q(\Sc_{\sigma_i}([G]_{Q_i}))$. Let $\kappa\in C_c^\infty(\ago_Q)$. Then, $(\kappa\circ H_Q) E_{Q_i}^Q(\psi)=E_{Q_i}^Q((\kappa\circ H_Q) \psi)$ and $(\kappa\circ H_Q) \psi\in \Sc_{\sigma_i}([G]_{Q_i})$ for every $i\in I$ and $\psi\in \Sc_{\sigma_i}([G]_{Q_i})$. It follows that $(\kappa\circ H_Q) \varphi_Q$ is also orthogonal to $\displaystyle \sum_{i\in I} E_{Q_i}^Q(\Sc_{\sigma_i}([G]_{Q_i}))$. Besides, by Lemma \ref{lemma estimates constant term}, the function $(\kappa\circ H_Q) \varphi_Q$ belongs to $\Sc([G]_Q)$. Therefore, by definition of $L^2_{\Xgo}([G]_Q)$, $(\kappa\circ H_Q) \varphi_Q$ is orthogonal to $L^2_{\Xgo}([G]_Q)$ and in particular to $\Sc_{\Xgo}([G]_Q)$. Finally, there certainly exists a sequence $\kappa_n\in C_c^\infty(\ago_Q)$ such that $(\kappa_n\circ H_Q) \varphi_Q$ converges to $\varphi_Q$ in $\tc([G]_Q)$ and we conclude that $\varphi_Q$ is indeed orthogonal to $\Sc_{\Xgo}([G]_Q)$.
			
			\item We prove the first part by induction on $\dim(\ago_0)-\dim(\ago_P)$. Let $\varphi\in \Sc([G]_P)$, $Q\subset P$ be a standard parabolic subgroup, $\lambda\in (\ago_Q^*)^{P++}$. By Proposition \ref{prop:caract Schwartz}, in order to show that $\varphi_{\Xgo}\in \Sc([G]_P)$ it suffices to check that $(\varphi_{\Xgo})_Q\in L^2_{2\rho_Q+\lambda}([G]_Q)^\infty$. Again by Proposition \ref{prop:caract Schwartz} and \eqref{compatibility projections weightedL2}, this will follow from the equality
			\begin{equation}\label{eq:cuspidal projection and constant termbis}
			(\varphi_{\Xgo})_Q=(\varphi_Q)_{\Xgo}.
			\end{equation}
			If $Q=P$ the above identity is tautological and this already settles the case $P=P_0$. If $Q\neq P$, by the induction hypothesis we have $\Sc([G]_Q)=\Sc_{\Xgo}([G]_Q)\oplus \Sc^{\Xgo}([G]_Q)$ and \eqref{eq:cuspidal projection and constant termbis} is equivalent to
			$$\displaystyle \langle (\varphi_{\Xgo})_Q,\psi\rangle_{[G]_Q}=\langle (\varphi_Q)_{\Xgo},\psi\rangle_{[G]_Q}$$
			for every $\psi\in \Sc_{\Xgo}([G]_Q)\cup \Sc^{\Xgo}([G]_Q)$. This last equality follows directly, by adjunction, from point 1.
			
			The continuity of $\varphi\in \Sc([G]_P)\mapsto \varphi_{\Xgo}\in \Sc_{\Xgo}([G]_P)$ is easy to obtain from the closed graph theorem. We now show the absolute convergence of $\sum_{\chi\in \Xgo(G)} \varphi_\chi$. Fix a compact-open subgroup $J\subset G(\AAA_f)$. By Proposition \ref{prop:caract Schwartz} and \eqref{eq:cuspidal projection and constant termbis}, there exists an increasing family $(\lVert . \rVert_n)_n$ of Hilbertian norms defining the topology of $\Sc([G]_P)^J$ such that for each $n$, denoting by $\cH_n$ the Hilbert completion of $\Sc([G]_P)^J$ with respect to $\lVert .\rVert_n$, the linear maps $\varphi\mapsto \varphi_\chi$, $\chi\in \Xgo(G)$, extend to orthogonal projections onto two by two orthogonal subspaces of $\cH_n$. In particular, for each $n$ and $\varphi\in \cH_n$ we have
			$$\displaystyle \sum_{\chi\in \Xgo(G)} \lVert \varphi_\chi\rVert_n^2<\infty.$$
			On the other hand, the space $\Sc([G]_P)^J$ being nuclear \eqref{eq: Schwartz spaces nuclear}, for each $n$ there exists $m\geqslant n$ such that the induced linear map $\cH_m\to \cH_n$ is Hilbert-Schmidt. By the Cauchy-Schwarz inequality, this implies
			\[\begin{aligned}
			\displaystyle \sum_{\chi\in \Xgo(G)} \lVert \varphi_\chi\rVert_n\leqslant \left(\sum_{\chi\in \Xgo(G)} \lVert \varphi_\chi\rVert_m^2 \right)^{1/2} \left(\sum_{\chi\in \Xgo(G)} \lVert \frac{\varphi_\chi}{\lVert \varphi_\chi\rVert_m}\rVert_n^2 \right)^{1/2}<\infty
			\end{aligned}\]
			for $\varphi\in \cH_m$. Since this holds for every $n$, the absolute convergence of $\sum_{\chi\in \Xgo(G)} \varphi_\chi$ follows. That the sum necessarily converges to $\varphi$ is obvious (e.g. because it converges to $\varphi$ in $L^2([G]_P)$).
			
			\item The first part follows from 2., \eqref{eq6 Function spaces} and the Dixmier-Malliavin theorem. The continuity of $\varphi\mapsto \varphi_{\Xgo}$ follows from the closed graph theorem and the equivariance with respect to the $A_P^\infty\times G(\AAA)$-action follows from uniqueness and equivariance of the projection $L^2([G]_P)\to L^2_{\Xgo}([G]_P)$. Also, \eqref{eq:cuspidal projection and constant term} can be proven the same way as \eqref{eq:cuspidal projection and constant termbis}, using 1. and adjunction. 
			
			\item Let $(\lVert . \rVert_n)_n$ be an increasing family of Hilbertian norms defining the topology of $\Sc([G]_P)^J$ satisfying the same property as before: denoting by $\cH_n$ the completion of $\Sc([G]_P)^J$ with respect to $\lVert .\rVert_n$, for each $n$ the maps $\varphi\mapsto \varphi_\chi$, $\chi\in \Xgo(G)$, extend to orthogonal projections onto two by two orthogonal subspaces of $\cH_n$. We may also assume that the norms $\lVert .\rVert_n$ are $G(F_\infty)$-continuous (that is, the $G(F_\infty)$-action extends to $\cH_n$). Let $\cH_{-n}$ be the dual Hilbert space equipped with the dual Hilbertian norm $\lVert .\rVert_{-n}$. The adjoint of the natural linear map $\cH_m\to \cH_n$, $m\geqslant n$, is a continuous $G(F_\infty)$-equivariant inclusion $\cH_{-n}\hookrightarrow \cH_{-m}$. Moreover, by \eqref{eq6 Function spaces} and the Dixmier-Malliavin theorem we have a natural inclusion $\cH_{-n}^\infty\hookrightarrow \tc([G]_P)^J$ and an equality $\tc([G]_P)^J=\bigcup_{n} \cH_{-n}^\infty$ of LF spaces (that this identification is a topological isomorphism follows from the open mapping theorem for LF spaces). Finally, there exists $n$ such that $\tc_N([G]_P)^J\subset \cH_{-n}^\infty$ and for each $m$ there exists $M$ such that $\cH_{-m}^\infty\subset \tc_M([G]_P)^J$. Therefore, it suffices to show that for every $n$ there exists $m\geqslant n$ such that the series \eqref{eq:series cuspidal projections} converges absolutely in $\cH_{-m}$ for every $\varphi\in \cH_{-n}$. This follows from the same argument as before since for $m$ large enough the inclusion $\cH_{-n}\hookrightarrow \cH_{-m}$ is Hilbert-Schmidt.
		\end{enumerate}
\end{preuve}
\end{paragr}

\begin{paragr}
	Let $\Xgo\subseteq \Xgo(G)$ be a subset. By the previous proposition, the projection $\varphi\mapsto \varphi_{\Xgo}$ maps $\tc([G]_P)$ and $\Sc([G]_P)$ continuously onto $\tc_{\Xgo}([G]_P)$ and $\Sc_{\Xgo}([G]_P)$ respectively. As $\Sc([G]_P)$ is dense in $\tc([G]_P)$ this entails that
	\begin{equation}\label{density projection cuspidal component}
	\Sc_{\Xgo}([G]_P) \mbox{ is dense in } \tc_{\Xgo}([G]_P).
	\end{equation}
	
	More precisely, let $N>0$ and $J\subset G(\AAA_f)$ be a compact-open subgroup. There exists $M>L>N$ such that the closure of $\Sc([G]_P)^J$ in $\tc_L([G]_P)^J$ contains $\tc_N([G]_P)^J$ and the projection $\varphi\mapsto \varphi_{\Xgo}$ restricts to a continuous linear mapping $\tc_L([G]_P)^J\to \tc_{M,\Xgo}([G]_P)^J$. Therefore:
	\begin{equation}\label{density projection cuspidal componentbis}
	\mbox{There exists } M>0 \mbox{ such that the closure of } \Sc_{\Xgo}([G])^J \mbox{ in } \tc_{M,\Xgo}([G])^J \mbox{ contains } \tc_{N,\Xgo}([G])^J.
	\end{equation}
\end{paragr}

\begin{paragr}
	Assume that $G=G_1\times G_2$ where $G_1$ and $G_2$ are connected reductive groups over $F$. We have then a natural identification $\Xgo(G)=\Xgo(G_1)\times \Xgo(G_2)$. For subsets $\Xgo_i\subseteq \Xgo(G_i)$, $i=1,2$, setting $\Xgo=\Xgo_1\times \Xgo_2$, the space $L^2_{\Xgo_1}([G_1])\otimes L^2_{\Xgo_2}([G_2])$ is dense in $L^2_{\Xgo}([G])$ (this follows from the fact that for $(M_{P_i},\sigma_i)\in \Xgo_i$, $i=1,2$, $\Sc_{\sigma_1}([G_1]_{P_1})\otimes \Sc_{\sigma_2}([G_2]_{P_2})$ is dense in $\Sc_{\sigma_1\boxtimes \sigma_2}([G]_{P_1\times P_2})$). Using again \eqref{eq: Schwartz spaces nuclear} and \eqref{eq 5 tvs}, we see that for any compact-open subgroups $J_i\subset G_i(\AAA_f)$, $i=1,2$, setting $J=J_1\times J_2$, the isomorphisms \eqref{eq7 Function spaces} restrict to isomorphisms
	\begin{equation}\label{eq isom HCS product with cuspidal datum}
	\displaystyle \cC_{\Xgo_1}([G_1])^{J_1}\widehat{\otimes} \cC_{\Xgo_2}([G_2])^{J_2}\simeq \cC_{\Xgo}([G])^{J},\;\; \Sc_{\Xgo_1}([G_1])^{J_1}\widehat{\otimes} \Sc_{\Xgo_2}([G_2])^{J_2}\simeq \Sc_{\Xgo}([G])^{J}.
	\end{equation}
\end{paragr}

\begin{paragr}[Generic cuspidal data.]\label{generic cuspidal data} ---
	Let $n\geqslant 1$. We say that a cuspidal datum $\chi\in \Xgo(\GL_n)$ is {\em generic} if it is represented by a pair $(M_P,\sigma)$ with
	$$\displaystyle M_P=\GL_{n_1}\times \ldots\times \GL_{n_k}$$
	a standard Levi subgroup of $\GL_n$ and
	$$\displaystyle \sigma=\sigma_1\boxtimes \ldots\boxtimes \sigma_k$$
	a cuspidal automorphic representation of $M_P(\AAA)$ whose central character is trivial on $A_P^\infty$ and such that $\sigma_i\neq \sigma_j$ for $1\leqslant i<j\leqslant k$.
	
	More generally, assume that $G$ is a product of the form $\Res_{K_1/F} \GL_{n_1}\times \ldots\times  \Res_{K_r/F} \GL_{n_r}$, where $K_1,\ldots,K_r$ are finite extensions of $F$. Then a cuspidal datum $\chi\in \Xgo(G)$ can be seen as a $r$-uple $(\chi_1,\ldots,\chi_r)$ where $\chi_i\in \Xgo(\GL_{n_i,K_i})$, $i=1,\ldots,r$, and we say that $\chi$ is {\em generic} if each of the $\chi_i$ is so.
	
	The following proposition is a consequence of Langlands's spectral decomposition \cite[VI.2.2]{MWlivre}, the description of the discrete spectrum of $\GL_n$ by M\oe{}glin-Waldspurger \cite{MW} (which, implies that for $\chi$ a generic cuspidal datum and $Q$ a parabolic subgroup, the discrete spectrum of $L^2_\chi([G]_Q)$ is entirely cuspidal), the computation of the constant terms of cuspidal Eisenstein series \cite[Proposition II.1.7]{MWlivre} and the fact that intertwining operators are unitary for purely imaginary arguments.
	
	\begin{proposition}\label{prop smooth L2 spaces constant term}
		Let $\chi\in \Xgo(G)$ be a generic cuspidal datum and let $P$ be a standard parabolic subgroup of $G$. Then, for every $\varphi\in L^2_\chi([G])$ we have $\varphi_P\in L^2_\chi([G]_P)$.
	\end{proposition}
	
	\begin{corollaire}\label{cor constant term generic cuspidal datum}
		Let $\chi\in \Xgo(G)$ be a generic cuspidal datum, $P$ be a standard parabolic subgroup of $G$ and $\chi_M$ be the inverse image of $\chi$ in $\Xgo(M_P)$. Then, for every $\varphi\in \Sc_\chi([G])$ and $s\in \cH_{>0}$, the function
		$$\displaystyle \varphi_{P,s}: m\in [M_P]\mapsto \delta_P(m)^{s-1/2} \varphi_P(m)$$
		belongs to $\cC_{\chi_M}([M_P])$. Moreover, the family of linear maps
		$$\displaystyle \Sc_\chi([G])\to \cC_{\chi_M}([M_P]), \varphi \mapsto \varphi_{P,s}$$
		for $s\in \cH_{>0}$ is holomorphic.
	\end{corollaire}
	
	\begin{preuve}
		Let $\varphi\in \Sc_\chi([G])$. Note that by Proposition \ref{prop decomposition cuspidal data}.1, $\varphi_{P,s}$ is orthogonal to $\Sc^{\chi_M}([M_P])$ for every $s\in \CC$. Hence, we just need to show that for every $s\in \cH_{>0}$ the function $\varphi_{P,s}$ belongs to $\cC([M_P])$ and the map $s\in \cH_{>0}\mapsto \varphi_{P,s}\in \cC([M_P])$ is holomorphic. As for every $X\in \mgo_\infty$ we have $R(X)\varphi_{P,s}=(2s-1)\langle \rho_P, X\rangle \varphi_{P,s}+(R(X)\varphi)_{P,s}$ (where we consider $\rho_P$ as an element of the dual space $\mgo_\infty^*$), by the equality \eqref{eq4 Function spaces}, it suffices to show that for every $d>0$ and $s\in \cH_{>0}$ we have $\varphi_{P,s}\in L^2_{\sigma,d}([M_P])$ and that the map $s\in \cH_{>0}\mapsto \varphi_{P,s}\in L^2_{\sigma,d}([M_P])$ is holomorphic. By Lemma \ref{lemma estimates constant term} and \eqref{eq2 Function spaces}, there exists $c>0$ such that $\varphi_{P,s}\in L^2_{\sigma,d}([M_P])$ for every $s\in \cH_{>c}$ and $d>0$. On the other hand, by Proposition \ref{prop smooth L2 spaces constant term}, we have $\varphi_{P,0}\in L^2([M_P])$. By H\"older inequality, for every $s\in \cH_{>0}$, $t>\Re(s)$ and $d>0$, we have
		$$\displaystyle \lVert \varphi_{P,s}\rVert_{L^2_{\sigma,d}}\leqslant \lVert \varphi_{P,0}\rVert_{L^2}^{1-\Re(s)/t} \lVert \varphi_{P,t}\rVert_{L^2_{\sigma,td/\Re(s)}}^{\Re(s)/t}$$
		and it follows that $\varphi_{P,s}\in L^2_{\sigma,d}([M_P])$ for every $s\in \cH_{>0}$ and $d>0$. The holomorphy of the map $s\in \cH_{>0}\mapsto \varphi_{P,s}\in L^2_{\sigma,d}([M_P])$ is equivalent to the holomorphy of $s\in \cH_{>0}\mapsto \langle \varphi_{P,s},\psi\rangle_{[M_P]}$ for every $\psi\in L^2_{\sigma,-d}([M_P])$ but this follows from the inequality
		$$\displaystyle \lvert \varphi_{P,s}\rvert\leqslant \lvert \varphi_{P,t_1}\rvert + \lvert \varphi_{P,t_2}\rvert$$
		for every $s\in \CC$ and $t_2>\Re(s)>t_1$. 
	\end{preuve}
\end{paragr}

\begin{paragr}\label{section PWLapid}
Let $\chi\in \Xgo(G)$ be a generic cuspidal datum represented by a pair $(M_P,\pi)\in \underline{\Xgo}(G)$. Set $\Pi=I_{P(\AAA)}^{G(\AAA)}(\pi)=\Ac_{P,\pi}(G)$ for the normalized smooth induction of $\pi$. Let $\bc_{P,\pi}$ be a $K$-basis of $\Pi$ as in \S \ref{S:K0ON}. For $\varphi\in \Sc([G])$ and $\lambda\in i\ago_P^*$ the series
\begin{equation}\label{spec 1}
\displaystyle \varphi_{\Pi_\lambda}=\sum_{\psi \in \bc_{P,\pi}} \langle \varphi, E(\psi,\lambda)\rangle_{[G]} E(\psi,\lambda)
\end{equation}
converges absolutely in $\tc_N([G])$ for some $N$ (that may a priori depend on $\lambda$). Indeed by \cite[Corollary 6.5]{BL} and \eqref{eq 2 tvs}, there exists $N$ such that the linear map $\psi\in \Ac_{P,\sigma}(G)\mapsto E(\psi,\lambda)\in \tc_N([G])$ is continuous so that the claim follows from Proposition \ref{prop:car-relJB} and the Dixmier-Malliavin theorem. The next theorem is a slight restatement of (part of) the main result of \cite{LapHC}\footnote{Note that in {\em loc. cit.} the Harish-Chandra Schwartz space $\cC([G])$ is denoted by $\Sc(G(F)\backslash G(\AAA))$}. We refer the reader to \S \ref{Schwartz functions in TVS} for the notion of Schwartz function valued in a TVS.
	
\begin{theoreme}[(Lapid)]\label{theo Lapid}
There exists $N>0$ such that for $\varphi\in \cC([G])$, the series \eqref{spec 1} still makes sense (that is the scalar products $\langle \varphi, E(\psi,\lambda)\rangle_{[G]}$ are convergent) and converges in $\tc_N([G])$ for every $\lambda\in i\ago_P^*$. Moreover, the function $\lambda\in i\ago_P^*\mapsto \varphi_{\Pi_\lambda}\in \tc_N([G])$ is Schwartz and if $\varphi\in \cC_\chi([G])$ we have the equality
$$\displaystyle \varphi=\int_{i\ago_P^*} \varphi_{\Pi_\lambda} d\lambda$$
(the right hand-side being absolutely convergent in $\tc_N([G])$).
\end{theoreme}

\begin{preuve}
Note that $G$ satisfies condition (HP) of \cite{LapHC}: it is proven in {\em loc. cit.} that general linear groups satisfy (HP) and it is straightforward to check that products of groups satisfying (HP) again satisfy (HP). The first part of the theorem is then a consequence of \cite[Proposition 5.1]{LapHC}. Indeed, by Dixmier-Malliavin we may assume that $\varphi=R(f)\varphi'$ where $\varphi'\in \cC([G])$ and $f\in C_c^\infty(G(\AAA))$. By {\em loc. cit.} the scalar product $\langle \varphi,E(\psi,\lambda)\rangle_{[G]}$ converges for every $\psi\in \bc_{P,\pi}$ and there exists $N>0$ such that $\psi\mapsto E(\psi,\lambda)$ factorizes through a continuous linear mapping $\Pi\to \tc_N([G])$ for every $\lambda\in i\ago_P^*$. As
$$\displaystyle \langle \varphi,E(\psi,\lambda)\rangle_{[G]}=\langle \varphi',E(R(f^*)\psi,\lambda)\rangle_{[G]}, \; \varphi\in \bc_{P,\pi},$$
we deduce by Proposition \ref{prop:car-relJB} that the series \eqref{spec 1} converges absolutely in $\tc_N([G])$ for every $\lambda\in i\ago_P^*$. That the function $\lambda\in i\ago_P^* \mapsto \varphi_{\Pi_\lambda}\in \tc_N([G])$ is Schwartz follows similarly from \cite[Corollary 5.7]{LapHC}. The last part of the theorem is a consequence of \cite[Theorem 4.5]{LapHC} since $\chi$ is generic and therefore $L^2_\chi([G])$ is included in the ``induced from cuspidal part'' $L^2_c([G])$ of $L^2([G])$, with the notation of {\em loc. cit.}, and moreover the stabilizer of the pair $(M_P,\pi)$ in $W$ is trivial.
\end{preuve}
\end{paragr}

\begin{paragr}[Automorphic kernels.] ---\label{S:Autom-kernel}
The right convolution by $f\in \Sc(G(\AAA))$ on each space of the decompositions \eqref{eq:langlands} gives integral operators whose kernels are respectively  denoted by $K_{f}(x,y)$, $K_{f,\chi}(x,y)$,  $K_{f}^0(x,y)$  and $K_{f,\chi}^0(x,y)$ where $x,y\in G(\AAA)$. If the context is clear, we shall omit the subscript $f$ in the notation. The kernels  are related  by the following equality for all $x,y\in G(\AAA)$
  \begin{align*}
    K^0_\chi(x,y)=\int_{A_G^\infty}   K_\chi(x,ay)\,da.
  \end{align*}

\begin{lemme}\label{lem:maj-noy}
For every $N>0$, there exists $N'>0$ such that
\begin{align}\label{ineq noy}
\displaystyle \sum_{\chi\in \Xgo(G)}\lvert K_{\chi}(x,y)\rvert \leqslant \lVert x\rVert_{[G]}^{N'} \lVert y\rVert_{[G]}^{-N},\; x,y\in [G],
\end{align}

\begin{align}\label{ineq noy0}
\displaystyle \sum_{\chi\in \Xgo(G)} \lvert K^0_{\chi}(x,y)\rvert \leqslant \lVert x\rVert_{[G]}^{N'} \lVert y\rVert_{[G]}^{-N},\; x,y\in [G]^1.
\end{align}
More generally, for every $N>0$, there exists $N'>0$ such that for each continuous semi-norm $\lVert . \rVert_{N'}$ on $\tc_{N'}([G])$
\begin{align}\label{ineq noychiseminorm}
\displaystyle \sum_{\chi\in \Xgo(G)}\lVert K_{\chi}(.,y)\rVert_{N'}\ll \lVert y\rVert_{[G]}^{-N}, \; y\in [G].
\end{align}
\end{lemme}

\begin{preuve}
Obviously, \eqref{ineq noychiseminorm} implies \eqref{ineq noy} and \eqref{ineq noy} implies \eqref{ineq noy0}. Let $N>0$ and choose $L>N$. By Proposition \ref{prop decomposition cuspidal data}.4 and the uniform boundedness principle, there exists $N'>L$ such that for every continuous semi-norm $\lVert . \rVert_{N'}$ on $\tc_{N'}([G])$ there exists a continuous semi-norm $\lVert . \rVert_L$ on $\tc_L([G])$ satisfying
$$\displaystyle \sum_{\chi\in \Xgo(G)}\lVert K_{\chi}(.,y)\rVert_{N'}\leqslant \lVert K_f(.,y)\rVert_L$$
for every $y\in [G]$. Therefore, it suffices to show that if $L$ is large enough then
$$\displaystyle \lVert K_f(.,y)\rVert_L \ll \lVert y\rVert_{[G]}^{-N}, \; y\in [G].$$
As $R(X)K_f(.,y)=K_{L(X)f}(.,y)$ and by definition of the topology on $\tc_L([G])$, it even suffices to prove that for $L$ large enough
\begin{equation}\label{ineq base noyau}
\displaystyle \sum_{\gamma\in G(F)} \lVert x^{-1}\gamma y\rVert^{-L}\ll \lVert x\rVert_{[G]}^{L} \lVert y\rVert_{[G]}^{-N},\;\; x,y\in [G].
\end{equation}
There exists $N_0>0$ such that $\sum_{\gamma\in G(F)}\lVert \gamma\rVert^{-N_0}<\infty$. Fix such a $N_0$. Then, as $\lVert y\rVert_{[G]}\leqslant \lVert \gamma y\rVert$ for every $\gamma\in G(F)$ and $y\in G(\AAA)$, we have
\[\begin{aligned}
\displaystyle \sum_{\gamma\in G(F)} \lVert x^{-1}\gamma y\rVert^{-2N_0-N}\ll \lVert x\rVert^{2N_0+N} \lVert y\rVert^{N_0}\lVert y\rVert_{[G]}^{-N-N_0}\sum_{\gamma\in G(F)} \lVert \gamma\rVert^{-N_0}\ll \lVert x\rVert^{2N_0+N} \lVert y\rVert^{N_0}\lVert y\rVert_{[G]}^{-N-N_0}
\end{aligned}\]
for $x,y\in G(\AAA)$. Since the left hand side of the above inequality is invariant by left translations of both $x$ and $y$ by $G(F)$, we may replace $\lVert x\rVert^{2N_0+N}$ and $\lVert y\rVert^{N_0}$ in the right hand side by $\lVert x\rVert_{[G]}^{2N_0+N}$ and $\lVert y\rVert_{[G]}^{N_0}$ respectively. This gives \eqref{ineq base noyau} for $L=2N_0+N$ and this ends the proof of the lemma.
\end{preuve}

\end{paragr}

\begin{paragr}\label{S:exp-kernel}
  Let $P$ be standard parabolic subgroup of $G$ and let $M=M_P$. Let $\chi\in \Xgo(G)$ and
  
  $$ \Ac_{P,\chi}(G)=\oplus_{\pi}\Ac_{P,\pi}(G)$$
  where the sum is over cuspidal representations $\pi$ of $M$ such that the image of $(M,\pi)$ by the map $\Xgo(M)\to \Xgo(G)$ is $\chi$. Let $\bc_{P,\chi}$ be a $K$-basis of $ \Ac_{P,\chi}(G)$ that is the union $\cup_{\pi}\bc_{P,\pi} $ over $\pi$ as above of $K$-bases of of  $\Ac_{P,\pi}(G)$ (see § \ref{S:K0ON}). In the same way we define $\bc_{P,\chi,\tau}=\cup_{\pi}\bc_{P,\pi,\tau} $ for any $\tau\in \hat K$.

  By a slight variant of \cite{ar1} §4 and \cite{ar2}  section 3, we have the following lemma.

  \begin{lemme}(Arthur)\label{lem:ar2} Let $C\subset G(\AAA_f)$ be a compact subset and let $K_0\subset K^\infty$ be a normal open compact subgroup.
     There exists a continuous semi-norm $\|\cdot\|$ on $\Sc(G(\AAA),C,K_0)$ and an integer $N$ such that for all $X,Y\in \uc(\ggo_\CC)$, all $x,y\in G(\AAA)^1$ and all $f\in \Sc(G(\AAA),C,K_0)$ we have
 \begin{align*}
     \sum_{\chi \in \Xgo(G)}\sum_{P_0\subset P} |\pc(M_P)|^{-1} \int_{i\ago_P^{G,*}}  \sum_{\tau\in \hat K} | \sum_{\varphi \in \bc_{P,\chi,\tau}} (R(X)E)(x,I_P(\la,f)\varphi,\la)\overline{R(Y)E(y,\varphi,\la)}|\, d\la \\
\leq \|f\| \|x\|^N_{G(\AAA)} \|y\|_{G(\AAA)}^N.
 \end{align*}
 Moreover for all $x,y\in G(\AAA)$ and all $\chi\in \Xgo(G)$ we have
    \begin{align*}
      K_{\chi}^0(x,y)=\sum_{P_0\subset P} |\pc(M_P)|^{-1} \int_{i\ago_P^{G,*}}  \sum_{\varphi \in \bc_{P,\chi}} E(x,I_P(\la,f)\varphi,\la)\overline{E(y,\varphi,\la)}\, d\la.
    \end{align*}
  \end{lemme}
  
\begin{preuve}
      One point is to remove the $K_\infty$-finiteness assumption in lemma 4.4 of  of \cite{ar1}. This can be done by approximation by  $K_\infty$-finite functions and this also enables us to put the sum over $\tau$ outside the absolute value. The other point is to remove the hypothesis about the compactness of the support of $f$. However the key point is in fact lemma 4.3 of \cite{ar1} which can be replaced by Lemma \ref{lem:maj-noy}.
\end{preuve}
\end{paragr}

\section{The spectral expansion of the Jacquet-Rallis trace formula for general linear groups}\label{Chap:Ichilimit}
\newcommand{\brtau}{\bar \tau}
\newcommand{\sbs}{\subset}
\newcommand{\smin}{\smallsetminus}

This chapter has two goals. The first, accomplished in Theorem \ref{thm:kappak}, is to extend the coarse spectral expansion $I=\sum_{\chi\in \Xgo(G)} I_\chi$ of the Jacquet-Rallis trace formula for linear groups $G$ (as proved in \cite{Z3}) to the Schwartz space. The second, given in Theorem \ref{thm:asym-trio},  is to  provide spectral expressions more suitable for explicit calculations.

\subsection{Notations}\label{ssec:sym-not}

\begin{paragr}\label{S:Gn}
  Let  $E/F$ be a quadratic extension of number fields. For convenience, we will fix $\tau\in F^\times$ such that $E=F[\sqrt{\tau}]$. Let $\eta$ be the quadratic character of $\AAA_F^\times$ attached to $E/F$. Let  $n\geq 1$ be an integer. Let $G'_n=GL_{n,F}$ be the algebraic group of $F$-linear automorphisms of $F^n$. Let  $G_n=\Res_{E/F} (G_n' \times_FE)$ be the $F$-group obtained by restriction of scalars from the algebraic group $GL_{n,E}$ of  $E$-linear automorphisms of $E^n$.  We denote by $c$ the Galois involution. We have a natural inclusion $G'_n\subset G_n$ which induces an inclusion $A_{G'_n}\subset A_{G_n}$ which is in fact an equality. The  restriction map $X^*(G_n)\to X^*(G'_n)$ gives an isomorphism $\ago_{G_n}^*\simeq \ago_{G'_n}^*$. 
\end{paragr}

\begin{paragr}\label{S:min-par-G}
  Let $(B_n',T_n')$ be a pair where $B_n'$ is the Borel subgroup $G_n'$ of  upper triangular matrices and $T_n'$ is the maximal torus of $ G' _n$ of diagonal matrices. Let $(B_n,T_n)$ be the pair deduced from $(B_n',T_n')$ by extension of scalars to $E$ and restriction to $F$: it is a pair of a minimal parabolic subgroup of $ G_n $ and its Levi factor.

 Let $ K_n \subset G_n (\AAA) $ and $ K '_n= K_n \cap G'_n (\AAA) \subset G '_n(\AAA) $ be the ``standard'' maximal compact subgroups. Notice  that we have $ K '_n\subset K_n $.
\end{paragr}

\begin{paragr}\label{S:PtoP'}
  The map $ P '\mapsto P = \Res_{E / F} (P' \times_F E) $ induces a bijection between the sets of standard parabolic subgroups of $ G '_n$ and  $ G _n$ whose inverse bijection is given by
$$ P \mapsto P '= P \cap G'_n. $$
Let $ P $ be a standard parabolic subgroup of $ G_n $. The restriction map $X^*(P)\to X^*(P')$ identifies $X^*(P)$ with a subgroup of $X(P')$ of index $2^{\dim(\ago_P)}$. It also induces an isomorphism $ \ago_{P'} \to   \ago_ {P }$ which fits into the commutative diagram:

\begin{align*}
  \xymatrix{ G'_n (\AAA) \ar[r]^{H_{P'}}\ar[d]  & \ago_{P'} \ar[d]\\
  G'_n(\AAA) \ar[r]^{H_{P}} & \ago_P }
\end{align*}

For any standard parabolic subgroups $P\subset Q$, the restriction of the function $ \tau_P ^ Q $ to $\ago_{P'}$ coincides with the function $ \tau_ {P '} ^ {Q'} $. However we have for all $x\in G'(\AAA)$
$$ \bg  \rho_P ^ Q, H_P (x) \bd = 2 \bg  \rho_ {P '} ^ {Q'}, H_ {P '} (x) \bd. 
$$

\begin{remarque}\label{rq:mesure-AP}
The map $ \ago_{P'} \to   \ago_ {P }$ does not preserve Haar measures. In fact, the pull-back on $\ago_{P'}$ of the Haar measure on $\ago_P$ is $2^{\dim(\ago_P)}$ times the  Haar measure on $\ago_{P'}$. In particular,  although the groups $A_P^\infty$ and $A_{P'}^\infty$ can be canonically identified,  the Haar  measure on $A_P^\infty$ is $2^{\dim(\ago_P)}$ times the Haar measure on $A_{P'}^\infty$.
\end{remarque}
\end{paragr}

\begin{paragr}
   We shall use the  natural embeddings $G_n'\subset G'_{n+1}$ and $G_n\subset G_{n+1}$ where the smaller group is identified with the subgroup of the bigger one  that fixes $e_{n+1}$ and  preserves the space generated by $(e_1,\ldots,e_n)$ where $(e_1,\ldots,e_{n+1})$ denotes the canonical basis of $F^{n+1}$.
\end{paragr}

\begin{paragr}\label{S:G=GnxGn+1}
  Let  $G=G_n\times G_{n+1}$ and $G'=G_n'\times G_{n+1}'$. Thus $G'$ is an $F$-subgroup of $G$.   Let
  $$\iota: G_n \hookrightarrow G_n\times G_{n+1}$$
  be the diagonal embedding. Let $H$ be the image of $\iota$ (so $H$ is isomorphic to $G_n$).
\end{paragr}

\begin{paragr}
  Let $K=K_n\times K_{n+1}$: it is a maximal compact subgroup of $G(\AAA)$. We define pairs $(P_0,M_0)=(B_{n}\times B_{n+1}, T_{n}\times T_{n+1})$ and  $(P_0',M_0')=(B_{n}'\times B_{n+1}', T_{n}'\times T_{n+1}')$ of minimal parabolic $F$-subgroups of $G$ and $G'$ with their Levi components. As in §\ref{S:PtoP'}, we have a bijection denoted $P\mapsto P'$ between the sets of standard parabolic subgroups of $G$ and $G'$.
\end{paragr}

\bpar 
In general, for a subgroup $P$ (usually a parabolic subgroup) of $G_{n}$, $G_{n+1}$ or $G$, 
we write $P'$ for the intersection of $P$ with $G_{n}'$, $G_{n+1}'$ or $G'$ respectively. 
\epar

\begin{paragr}\label{S:character etaG'} Let $\det_n$ (resp. $\det_{n+1}$) be the morphism $G'\mapsto \mathbb{G}_{m,F}$ given by the determinant on the first (resp. second) component. Let $\eta_{G'}$ be the character $G'(\AAA)\to \{\pm 1\}$ given by 
$$\eta_{G'}(h)=\eta({\det}_{n}(h))^{n+1}\eta({\det}_{n+1}(h))^{n}.$$
  \end{paragr}
  
\begin{paragr}
  We set $\ago_{n+1}=\ago_{B_{n+1}}$ and $\ago_{n+1}^+=\ago_{B_{n+1}}^+$ (see §\ref{S:root-coroot}).
\end{paragr}

\bpar 
We let $\fc_{RS}$ to be the set of $F$-parabolic subgroups of $G$ of the form $P=P_n\times P_{n+1}$ where $P_n$ is a standard parabolic subgroup of $G_n$  and $P_{n+1}$ is a semi-standard parabolic subgroup of $G_{n+1}$ such that $P_{n+1} \cap G_{n} = P_{n}$ (using the embedding $G_n \hookrightarrow G_{n+1}$).
\epar

\bpar 
For $P, Q \in \fc_{RS}$ such that $P \sbs Q$ we let $\epsilon_P^Q =(-1)^{\dim(a_{P_{n+1}}^{Q_{n+1}})}$.
\epar

\subsection{The coarse spectral expansion for Schwartz functions}

\begin{paragr} Let  $f \in \Sc(G(\AAA))$ be a Schwartz test function (see §\ref{S:Schw}). 
\end{paragr}

\begin{paragr}
  Let $P$ be  a parabolic subgroup of $G$. The right convolution by $f$ on $L^2([G]_{P})$  gives an integral operator whose kernel is denoted by $K_{P,f}$. Let $\chi \in \Xgo(G)$. Replacing  $L^2([G]_{P})$ by its closed subspace $L_{\chi}^2([G]_{P})$  (see \eqref{eq:langlands}), we get a kernel  denoted by $K_{P,\chi,f}$. If $P=G$, we omit the subscript $P$. Most of the time, we will also omit the subscript $f$. We have $K_P=\sum_{\chi\in \Xgo(G)} K_{P,\chi}$.
\end{paragr}

\begin{paragr}\label{S:KTchi} For $x \in H(\AAA)$, $y = (y_n,y_{n+1}) \in G'_n(\AAA) \times G'_{n+1}(\AAA)$, $\chi \in \Xgo(G)$
and $T \in \ago_{n+1}$ we set
\begin{align}
\label{eq:KchiT}&K_{\chi}^{T}(x,y) = \\
\nonumber &\sum_{ P \in \fc_{RS}} \epsilon_P^G 
\sum_{\gamma \in (P\cap H)(F)\bsl H(F)} 
\sum_{\delta  \in P' (F) \bsl G'(F)}
\htau_{P_{n+1}}(H_{P_{n+1}}(\delta_n y_n) - T_{P_{n+1}})K_{P, \chi}(\gamma x, \delta  y).
\end{align}
where
\begin{itemize}
\item we write $\delta=(\delta_{n},\delta_{n+1})$ and $y=(y_n,y_{n+1})$ according to the decomposition $G'=G'_n\times G'_{n+1}$;
\item in the notation $H_{P_{n+1}}(\delta_n y_n)$,   we consider $\delta_{n} y_{n}$  as an element of $G_{n+1}'(\AAA)$ (via the embedding $G_n'\hookrightarrow G'_{n+1}$);
\item $T_{P_{n+1}}$ is defined as in §\ref{S:trunc-param};
\end{itemize}
and the rest of notation is explained in \S \ref{ssec:sym-not}.

\begin{remarque}
  \label{rq:KTchi}
This is the kernel used in \cite{Z3} for compactly supported functions. Since we are considering a Schwartz function $f$, the sums over $\gamma$ and $\delta$ are not finite. However, the component $\delta_n$ may be taken in a \emph{finite} set depending on $y_n$ (see \cite{ar1} Lemma 5.1). We can then easily  show that the sums are absolutely convergent using majorization of Lemma \ref{lem:maj-noy}.
\end{remarque}
\end{paragr}

\begin{paragr}

\begin{theoreme}\label{thm:jfDef} Let $T \in \ago_{n+1}^+$.
  \begin{enumerate}
  \item We have
\[
\sum_{\chi \in \Xgo(G)} \int_{[H]}\int_{[G']}| K^T_{f,\chi}(h,g') | \, dg' dh < \infty
\]
\item As a function of $T$, the integral
\begin{align}\label{eq:polynexp}
I_\chi^T(f)=\int_{[H]} \int_{[G']} K^T_{f,\chi}(h,g') \eta_{G'}(g') \, dg' dh
\end{align}
coincides with  a polynomial-exponential function in $T$ whose purely polynomial part is constant and denoted by $I_\chi(f)$.
\item The distributions $I_\chi$ are continuous, left $H(\AAA)$-invariant and right $(G'(\AAA),\eta_{G'})$-equivariant. 
\item The sum
  \begin{align}\label{eq:expansion-spec}
    I (f)= \sum_{\chi} I_{\chi}(f)
  \end{align}
is absolutely convergent and defines a continuous distribution $I$.
  \end{enumerate}
\end{theoreme}

\begin{remarque}
  The last statement is the coarse spectral expansion of the Jacquet-Rallis trace formula for $G$.
\end{remarque}

 \begin{preuve} 
All the statements but the continuity and the extension to Schwartz functions  are proved in  \cite[Theorems 3.1 and 3.9]{Z3}  for compactly supported functions. 
Assuming extension to Schwartz case, continuity is the result of the explicit formula of \cite{Z3}, Theorem 3.7 (which also holds for Schwartz functions). 
As for absolute convergence in the Schwartz case, we state and prove a twin theorem below (see Theorem \ref{thm:kappak}) whose proof can easily be adapted to 
the current theorem, so we will not repeat the arguments here.
 \end{preuve}
 \end{paragr}

\subsection{Auxiliary expressions for $I_\chi$}

\begin{paragr}
  The goal of this section is to provide new expressions for the distribution $I_\chi$ defined in Theorem \ref{thm:jfDef}. In this paper, we will use these  expressions  to explicitly compute $I_\chi$. The main results are subsumed in Theorem \ref{thm:jgood}. Before giving its statement, we have to explain the main objects. For this, we fix $f\in \Sc(G(\A))$ and we simply denote by $K_\chi$ the kernel $K_{f,\chi}$.
\end{paragr}

\begin{paragr}[The Ichino-Yamana truncation operator.] --- \label{S:IY-trunc} Let $T\in \ago_{n+1}$.  In \cite{IY}, Ichino-Yamana defined a truncation operator which transforms functions of moderate growth on $[G_{n+1}]$ into rapidly decreasing  functions on $[G_n]$. By applying it to the \emph{right} component of $[G]=[G_n]\times [G_{n+1}]$, we get a truncation operator which we denote by $\Lambda^T_r$ (the subscript $r$ is for right). It associates to any function $\varphi$ on $[G]$ the function on $[H]$ defined by the following formula: for any  $h\in [H]$:
\begin{align}  \label{eq:LaTd}
    (\Lambda^T_r\varphi)(h)= \sum_{ P  \in \fc_{RS} } \epsilon_{P}^{G} \sum_{\delta \in (P\cap H)(F)\back H(F) }   \hat\tau_{P_{n+1}}(H_{P_{n+1}}(\delta h)-T_{P_{n+1}}) \varphi_{G_n\times P_{n+1}}(  \delta  h)
\end{align}
where we follow  notations of §\ref{S:KTchi}.  Note that in the expression $H_{P_{n+1}}(\delta h)$, we view $\delta h$ as an element of $G_{n+1}(\AAA)$ by the composition $H\hookrightarrow G \to G_{n+1}$ where the second map is the second projection.
We denote by $\varphi_{G_n\times P_{n+1}}$  the constant term of $\varphi$ along $G_n\times P_{n+1}$.

For properties of $\Lambda_r^T$ we shall refer to \cite{IY}. However for our purposes it is convenient to state the following proposition.

\begin{proposition}\label{prop:decayLaT}
  For any integers $N$ and $N'$, any open compact subset $K_0\subset G(\AAA_f)$, there is an integer $r\geq 0$ and a finite family  $(X_i)_{i\in I}$ of elements of $\uc(\ggo_\CC)$ of degree $\leq r$ such  that for any $\varphi\in C^r(G(F)\back G(\AAA)/K_0)$ we have for all $h\in [H]$
  \begin{align*}
    (\Lambda_r^T\varphi)(h) \leq  \|h\|_{[H]}^{-N}  \sum_{i\in I} \big(\sup_{x\in G(\AAA)} \|x\|_{[G]}^{-N'} |(R(X_i)\varphi)(x)|\big)
  \end{align*}
\end{proposition}

\begin{preuve}
The result, a  variant of Arthur's Lemma 1.4 of \cite{ar2}, is proven in \cite{IY}, Lemma 2.4. 
\end{preuve}

\end{paragr}

\begin{paragr}[Convergence of a first integral.] --- It is given by the following proposition.

  \begin{proposition}\label{prop:cv-LaTrKchi}
    Let $\chi\in \Xgo(G)$. The integral
    \begin{align}\label{eq:LaTrKchi}
      \int_{[H]} \int_{[G'] }   \Lambda_r^TK_\chi (x,y)\,\eta_{G'}(y)dxdy
    \end{align}
is absolutely convergent.
  \end{proposition}

  \begin{preuve}  We can easily deduce from Lemma \ref{lem:maj-noy} that for all $r_1\geq 0$, there exists a continuous semi-norm $\|\cdot\|$ on $\Sc(G(\AAA))$ and an integer $N\in \NN$  such that for all  $x\in [G]$, $y\in G(\AAA)^1$, $a\in A_G^\infty$, $f\in \Sc(G(\AAA),C,K_0)$ we have
\begin{align}\label{eq:KchiG}
    |K_{f,\chi}( x, ay)|\leq \|f\|\|a\|_{G(\AAA)}^{-r_1}\|x\|^N_{[G]}. 
  \end{align}
In particular, if we restrict ourselves to  $x\in [H]$,  we have with the same hypothesis the existence of $N$ and $\|\cdot\|$ such that
\begin{align}\label{eq:KchiG'}
    \|x\|^{-N}_{[H]}|K_{f,\chi}( x, ay)|\leq \|f\|\|a\|_{G(\AAA)}^{-r_1}. 
  \end{align}
for all  $x\in [H]$, $y\in G(\AAA)^1$, $a\in A_G^\infty$ and $f\in \Sc(G(\AAA),C,K_0)$.

The right derivatives in the first variable of the kernel $K_\chi(x,y)$ can be expressed in terms of the kernel  $K_\chi$ associated to left derivatives of $f$. Thus, taking into account Proposition \ref{prop:decayLaT}, we see that for any $r_2\geq 0$ there exists a continuous semi-norm $\|\cdot\|$ on $\Sc(G(\AAA),C,K_0)$ such that for $f\in \Sc(G(\AAA),C,K_0)$ we have
\begin{align*}
|\Lambda_r^TK_\chi (x,ay)| \leq \|f\| \|a\|^{-r_1}_{G(\AAA)} \|x\|^{-r_2}_{[H]}
\end{align*}
for all $a\in  A_G^\infty$, $x\in [H]$ and $y\in G(\AAA)^1$. The convergence is then obvious.
\end{preuve}
\end{paragr}

\begin{paragr}[Arthur function $F^{G_{n+1}}(\cdot,T)$.] \label{S:Arthur-F} --- For $T \in \ago_{n+1}$ we shall use Arthur function $F^{G_{n+1}}(\cdot,T)$ (see \cite{ar1} §6). Recall that this is the characteristic function of the set of $x \in G_{n+1}(\A)$  for which there exists a $\delta \in G_{n+1}(F)$ such that $\delta x \in \sgo_{G_{n+1}}$ (see § \ref{Siegel domain}) and $\langle \varpi, H_0(\delta x) - T \rangle \le 0$ for all $\varpi \in \widehat \Delta_{B_{n+1}}$. Recall also that $F^{G_{n+1}}(\cdot,T)$ descends to characteristic function of a compact subset of $Z_{n+1}(\A)G_{n+1}(F)\back G_{n+1}(\A)$.
\end{paragr}

\begin{paragr}[Two other convergent integrals.] ---

  \begin{proposition}\label{prop:FonHG}
    The following integrals are absolutely convergent:
    \begin{align}
      \label{eq:FonH}\int_{[H]}\int_{[G']}F^{G_{n+1}}(x,T)K_{\chi}(x,y)\,\eta_{G'}(y)\,dydx\\
 \label{eq:FonG}\int_{[H]}\int_{[G']}K_{\chi}(x,y) F^{G_{n+1}}(y_n,T) \,\eta_{G'}(y)\,dydx
    \end{align}
 where $y = (y_n,y_{n+1}) \in G'_{n}(\A) \times G'_{n+1}(\A)$.
   \end{proposition}
  
   \begin{preuve}
     The convergence of the integral \eqref{eq:FonH} is proved in the same way as in the proof of Proposition \ref{prop:cv-LaTrKchi}. The only point to observe in that the restriction of $F^{G_{n+1}}(\cdot,T)$ to $[H]$ is compactly supported.

The convergence of the integral \eqref{eq:FonG} is a consequence of two facts: first the  restriction of $F^{G_{n+1}}(\cdot,T)$ to $[G_n]$ is compactly supported; second for every $N'$, $N''$ there exist $N>0$ and  a continuous semi-norm on $\Sc(G(\AAA))$ such that 
\begin{align*}
| K_{\chi}(x,y)|\leq \|f\| \|y_n\|^{N}_{[G_n]}\|y_{n+1}\|^{-N'}_{[G_{n+1}]}  \| x\|_{[H]}^{-N''}.
\end{align*}
This majorization can be proved as in the proof of Lemma \ref{lem:maj-noy}. It suffices to prove the same result for the whole kernel (without subscript $\chi$). We can assume that $f=f_n\otimes f_{n+1}$ is a product. The kernel itself is then a product We can use Lemma \ref{lem:maj-noy} assertion 1 to bound $K_{f_{n+1}}(x,y_{n+1})$ and get the negative power of  $\|y_{n+1}\|_{[G_{n+1}]} $. Then we can bound $K_{f_{n+1}}(x,y_{n})$ to get the negative power of  $\|x\|_{[H]}$. 
   \end{preuve}
\end{paragr}

\begin{paragr}\label{S:asym-equal}
We say a functions $p$ and $q$ on $\ago_{n+1}$ are asymptotically equal if for all $\varepsilon > 0$ and $m \geq 0$ there exists  $c\geq 0$ such that for all $T \in \ago_{n+1}$ such that $\langle \al, T \rangle > \varepsilon \|T \|$ for all $\al \in \Delta_{B_{n+1}}$
we have
\[
|p(T) - q(T)| \leq ce^{-m \|T\|}.
\]
\end{paragr}

\begin{paragr} We can now state the main theorem of the section.
  
    \begin{theoreme}\label{thm:asym-trio} Let  $\chi \in \Xgo(G)$. 
Each of the three expressions  \eqref{eq:LaTrKchi} (see Proposition \ref{prop:cv-LaTrKchi}),  \eqref{eq:FonH} and \eqref{eq:FonG} (see Proposition \ref{prop:FonHG}) is asymptotically equal (in the  sense of  §\ref{S:asym-equal}) to a polynomial-exponential function of $T$ whose purely polynomial term is constant and equal to $I_\chi(f)$.
\end{theoreme}

\begin{preuve}
  The theorem is a simple combination of Theorem \ref{thm:jfDef} above and Theorems \ref{thm:kappak} and \ref{thm:jgood} below.
\end{preuve}
\end{paragr}

\subsection{Convergence of a truncated kernel}

\begin{paragr}
In this section, we give a first step in the proof of Theorem \ref{thm:asym-trio}: we define a new truncated kernel and we prove that the integral over $[H]\times [G']$  of this kernel is absolutely  convergent. 
\end{paragr}

\begin{paragr}
  We fix $f\in \Sc(G(\A))$, $\chi \in \Xgo(G)$ and $T \in \ago_{n+1}^+$.
\end{paragr}

\begin{paragr}[A new truncated kernel.] ---
For  $x \in [H]$ and $y \in [G']$ we set
\[
\kappa_{\chi}^{T}(x,y) = 
\sum_{P\in \fc_{RS}} \epsilon_{P}^{G} \sum_{\gamma \in (P \cap H) (F)\bsl H(F)} 
\sum_{\delta \in P' (F) \bsl G'(F)} \htau_{P_{n+1}}(H_{P_{n+1}}(\gamma x) - T_{P_{n+1}})K_{P, \chi}(\gamma x, \delta y).
\]
The notations are those of \S \ref{ssec:sym-not}. The expression $H_{P_{n+1}}(\gamma x)$ is interpreted as in the comments following \eqref{eq:LaTd}.

\begin{remarque}
  This is a version of the truncated kernel $K_{\chi}^{T}$ defined in §\ref{S:KTchi}. We will consider the connection between the two in Theorem \ref{thm:kappak}. The sum over $\gamma$ is finite whereas the sum over $\delta$ is convergent but not finite  (see remark \ref{rq:KTchi}).
\end{remarque}

The following theorem is the main result of the section.

\begin{theoreme}\label{thm:kappak}
  \begin{enumerate}
  \item We have 
\[
\sum_{\chi \in \Xgo(G)}\int_{[H]}\int_{[G']}|\kappa_{\chi}^T(x,y)|\,dydx < \infty.
\]
\item The integral 
  \begin{align}\label{eq:ichiT}
  i_{\chi}^{T}(f) :=  \int_{[H]}\int_{[G']}\kappa_{\chi}^{T}(x,y)\,\eta_{G'}(y) dydx
  \end{align}
coincides with a polynomial-exponential function in $T$ whose purely polynomial part is constant and denoted by $i_\chi(f)$.
\item The distribution $i_\chi$ is continuous. Moreover we have  $i_\chi=I_{\chi}$ where the right-hand side is defined in Theorem \ref{thm:jfDef}. 
\end{enumerate}
\end{theoreme}

The proof of assertion 1  of Theorem \ref{thm:kappak} will be  given in §\ref{S:proof-schwartzCVG} whereas the proof of assertions 2 and 3 will be given in §\ref{S:proof-kappak}. Before that, we must introduce additional  notations and lemmas. 
\end{paragr}

\begin{paragr}
Until   §\ref{S:proof-schwartzCVG},  we  assume that  $G$ is  a general reductive group with the notations of section \ref{ssec:alggps}. In particular, a maximal $F$-split torus $A_0$ is fixed as well as a minimal parabolic subgroup $P_0$ containing it. 
\end{paragr}

\begin{paragr}
   Recall that we have defined $\ago_{P}^{*, +}$  in §\ref{S:root-coroot}.
\end{paragr}

\begin{paragr}
  For semi-standard parabolic subgroups $P, Q$ of $G$, such that $P \subset Q$ 
we have the function  $\sigma_{P}^{Q}$ as defined in §6 of \cite{ar1}.  It is a characteristic function of a region in $\ago_{P}$. 
We note the following formula satisfied by it (see proof of Theorem 7.1 in \cite{ar1})
\begin{equation}\label{eq:sigTauId}
\tau_{P}^Q\htau_Q = \sum_{R \supset Q}\sigma_{P}^R.
\end{equation}

More generally, if $P \sbs Q \sbs R$ are semi-standard parabolic subgroups of $G$, 
we note $\sigma_{P}^{Q, R} := \sigma_{P \cap M_{R}}^{Q \cap M_{R}}$ the $\sigma$ 
function with respect to the group $M_{R}$ and its parabolic subgroups $P \cap M_{R}$ 
and $Q \cap M_{R}$.

\begin{lemme}\label{lem:obtuse} Let $P, Q, R$ be standard parabolic subgroups of $G$ with $P \sbs Q \sbs R$. 
Let $\al \in \Delta_{P} \smin \Delta_{P}^R$. Suppose $H \in \ago_{P}$ 
satisfies $\sigma_{P}^{Q, R}(H) = 1$. Then, if we denote $H_{R}$ the projection of $H$ 
onto $\ago_{R}$, we have $\langle \al, H_R \rangle \geq \langle \al, H \rangle$.
\end{lemme}

\begin{preuve}
Using \eqref{eq:sigTauId} we have $\htau_{P}^R(H) = \htau_{P}^R(H- H_R) = 1$. 
This means that $H - H_{R}$ has positive coefficients in the coroot basis $\Delta_{P}^{R, \vee}$ 
of $\ago_{P}^R$. The result follows from the known fact that distinct coroots in $\Delta_{0}^{\vee}$ form obtuse angles and so 
do distinct elements of $\Delta_{P}^{\vee}$ (c.f. \cite{labWal} Lemme 1.2.4 ).
\end{preuve}
\end{paragr}

\begin{paragr}
For a standard parabolic subgroup $P$ of $G$ we fix representatives for the cosets 
$W^P \bsl W$ and $W / W^P$ as follows:
\[
W^P \bsl W := \{ w \in W \ | \ w^{-1}\al >0 \ \forall \al \in \Delta_0^P \}, \quad 
W / W^P := \{ w \in W \ | \ w\al> 0 \ \forall \al \in \Delta_0^P \},
\]
where $\beta > 0$, for a root $\beta$ of $(G, A_0)$, means that it's a sum of elements of $\Delta_0$.
\end{paragr}

\begin{paragr}[Norms] 

We use the definition of \ref{ssec:norms} to define the norm on $[G]_{P}$. As we will only use norms 
on automorphic quotients, we denote this norm simply by $\|\cdot\|$. 
We use the same symbol for a fixed $W^G$-invariant norm on $\ago_{0}$. 
\end{paragr}

\begin{paragr}[Bound on Eisenstein Series.] ---

\begin{lemme}\label{lem:EisBound}
Let $w \in W/W^P$ and $\la \in \ago_{P}^{*, +}$. For all $\varepsilon > 0$, there is an $N \geq 0$, independent of $\la$, such that we have
\[
\sum_{\delta \in (P \cap wP_0 w^{-1})(F) \bsl P_0(F)} e^{\langle \la + (2+\varepsilon)\rho_P, H_{P}(w^{-1} \delta x) \rangle } \ll 
\|x\|^{N}   e^{\langle w\la, H_{0}(x) \rangle}. 
\]
 for any  $x$ in the Siegel domain $\sgo_{P}$ (see §\ref{Siegel domain}).
\end{lemme}

\begin{preuve}
 We have for $\delta \in P_{0}(F)$
\[
\langle \la,  H_{P}(w^{-1} \delta x) \rangle = \langle \la,  H_{0}(w^{-1} \delta  x) \rangle = 
\langle w \la, H_{0}(x) \rangle + \langle \la,  H_{0}(w^{-1} n) \rangle
\]
for some $n \in N_{0}(\A)$ depending on $x$ and $\delta$. Since $\la \in \ago_{P}^{+,*}$ we have, by \cite{labWal} Lemma 3.3.1, $\langle \la, H_{0}(w^{-1} n) \rangle \le c$ for some constant $c$ depending on $\la$ but independent of $n$. 
We have then 
\[
\sum_{\delta \in (P \cap wP_0 w^{-1})(F) \bsl P_0(F)} e^{\langle \la + (2+\varepsilon)\rho_P, H_{P}(w^{-1} \delta x) \rangle } \ll
\sum_{\delta \in (P \cap wP_0 w^{-1})(F) \bsl P_0(F)} e^{\langle  (2+\varepsilon)\rho_P, H_{P}(w^{-1} \delta x) \rangle}
 e^{\langle w\la, H_{0}(x) \rangle}
\]
and the result follows by moderate growth of Eisenstein series. 
\end{preuve}
\end{paragr}

\begin{paragr}[Auxiliary characteristic functions.] ---
We introduce the function $\brtau_P$ as the characteristic function of $H \in \ago_P$ 
such that $\langle \al, H \rangle \leq 0$ for all $\al \in \Delta_P$. By Langlands Combinatorial Lemma (\cite{labWal}, Proposition 1.7.2)
we have 
\begin{equation}\label{eq:LanCom}
\sum_{Q \supset P} \htau_P^{Q} \brtau_{Q} = 1.
\end{equation}
\end{paragr}

\begin{paragr}[Proof of  assertion 1 of Theorem \ref{thm:kappak}.] --- \label{S:proof-schwartzCVG} 
Unless otherwise stated, all sums of the type $\Sigma_{P}$ or $\Sigma_{P \sbs Q}$ are over elements of $\fc_{RS}$. 

In \S 2.4 of \cite{Z3} the operator $\La^{T,P}_{d}$ 
is defined for all $P \in \fc_{RS}$. 
For a function $\phi : P(F) \bsl G(\AAA) \to \CC$,  we have
\begin{equation}\label{eq:ladDef}
\La^{T, P}_{d}\phi(x) = \sum_{P \supset Q}\epsilon_{Q}^{P}
\sum_{\delta \in (Q \cap H) (F) \bsl (P \cap H) (F)} 
\htau_{Q_{n+1}}^{P_{n+1}}(H_{Q_{n+1}}(\delta x) - T_{Q_{n+1}}) \phi_{Q}(\delta x), \quad x \in  (P \cap H) (F) \bsl H(\AAA).
\end{equation}
Note that if $P=G$, the operator $\La^{T, G}_{d}$ is close but not exactly equal to the operator $\La^T_r$ defined  in \eqref{eq:LaTd}. Indeed in the former we take the constant term along $P$ whereas in the latter we take it along $G_n\times P_{n+1}$.

Using the inversion formula of Lemme 2.7 \cite{Z3} 
together with the formula \eqref{eq:sigTauId}
we obtain that the integral is bounded by the sum over $P,Q \in \fc_{RS}$ of 
\begin{equation}\label{eq:altSpec}
\sum_{\chi}\int_{(H \cap P)(F) \bsl H(\A)}\int_{[G']}\sigma_{P_{n+1}}^{Q_{n+1}}(H_{P_{n+1}}(x) - T_{P_{n+1}})
|\La_{d}^{T, P} \dsl \sum_{P \subset \tlP \subset Q} \epsilon_{P}^{P_1}	\sum_{\delta' \in \tlP' (F) \bsl G'(F)} K_{\tlP, \chi}(x,\delta' y) \rb |\,dydx
\end{equation}
where $\tlP \in \fc_{RS}$ in the alternating sum 
and $\La_{d}^{T, P}$ is applied with respect to $x$.
We consider $P$ and $Q$ fixed from now on. Without loss of generality we can assume they are standard in $G$.

For any standard parabolic subgroup $S$ of $G$ we have the mixed truncation operator $\La_{m}^{T,S}$ of \cite{JLR}
\[
\La_{m}^{T,S}\phi(x) = 
\sum_{R \sbs S} (-1)^{\dim (\ago_{R}^S)} \sum_{\delta \in R'(F) \bsl S'(F)}
\htau^S_{R}(H_{R}(\delta x) - T)\phi_R(\delta x), \quad x \in R'(F) \bsl S'(\A).
\]
Using \eqref{eq:LanCom}, the inversion formula (\S 6 (19) in \cite{JLR}) and \eqref{eq:sigTauId}, we bound the expression \eqref{eq:altSpec}
by a sum over $R \sbs S \sbs \tlS$, all standard parabolic subgroups of $G$, 
of 
\begin{multline}\label{eq:altSpec2}
\sum_{\chi}\int_{(H \cap P)(F) \bsl H(\A)}\int_{(G'\cap R) (F) \bsl G'(\A)}\sigma_{P_{n+1}}^{Q_{n+1}}(H_{P_{n+1}}(x) - T)
\brtau_{\tlS}(H_{\tlS}(y) - T')\sigma_{R}^{S, \tlS}(H_{R}(y) - T') \\
|\La_{d}^{T, P}\La_{m}^{T', R} \dsl \sum_{P^+ \sbs \tlP \sbs Q^+} \epsilon_{P}^{\tlP} K_{\tlP,\chi}(x,y)\rb  |\,dydx
\end{multline}
where $\La_{m}^{T', R}$ is applied with respect to $y$, $P^+$ is the smallest element of $\fc_{RS}$
containing $R$ and $P$, $Q^+$ is the largest element of $\fc_{RS}$ contained in $S$ and $Q$, the sum runs over $\tlP \in \fc_{RS}$, 
$T' \in \ago_{n+1}^+$ is any parameter, and we write $T$ instead of $T_{P_{n+1}}$ because $P$ is assumed standard
so there is no difference. 

Let $P_H = H \cap P$, it is a standard parabolic subgroup of $H$. 
Let $Z_P^{\infty} = A_{P}^{\infty} \cap A_{P_H}^{\infty}$. 
We also note $\zgo_{P} = \ago_{P} \cap \ago_{P_H}$. 
Let $M(\A)^{P,1}$ be the kernel of the composition of 
$H_{P_H} : M_{P_H}(\A) \to \ago_{P_H}$ 
with the orthogonal projection $\ago_{P_H} \to \zgo_{P}$.  

We assume that $f$ is bi-$K$-invariant - it makes the subsequent computations and notation clearer. 
It is not a serious restriction, one should deal with the general case as in \cite{ar2} or \cite{Z3}.
Using the Iwasawa decomposition
and Propositions 2.3 and 2.8 of \cite{Z3} we bound \eqref{eq:altSpec3} by
\begin{multline*}
\sum_{\chi}\int_{Z_{P}^{\infty}}\int_{A_R^{\infty}}\int_{M_{P_H}(F) \bsl M(\A)^{P,1}}
\int_{[M_{R'}]^{1}}
e^{\langle -2\rho_{P_H}, H_{P_{H}}(zm_1) \rangle}e^{\langle -2\rho_{R'}, H_{R'}(am_2) \rangle}
\sigma_{P_{n+1}}^{Q_{n+1}}(H_{P_{n+1}}(zm_1) - T) \\
\brtau_{\tlS}(H_{\tlS}(a) - T')\sigma_{R}^{S,\tlS}(H_{R}(a) - T')  
\|m_1\|^{-r_1}\|m_2\|^{-r_2}
\sup_{m_1' \in M_{P_{H}}(F) \bsl M(\A)^{P,1}}
\sup_{m_2' \in M_{R'}(F) \bsl M_{R'}(\A)^{1}} \\
|\int_{[N_{P}]}\int_{[N_{R}]}   \dsl \sum_{P^+ \sbs \tlP \sbs Q^+} \epsilon_{P}^{\tlP} K_{\tlP,\chi}(n_1zm_1',n_2am_2') \rb dn_2dn_1 
\|m_1'\|^{-r_1'}\|m_2'\|^{-r_2'}
|\,dm_2dm_1dadz
\end{multline*}
for any $r_1, r_1', r_2, r_2'$. Here, we should replace the kernel $K$ by a finite sum of kernels $K_i$, 
but to make the notation simpler we will ignore this detail. 
Making a few changes of variables we bound the above expression by
\begin{multline}\label{eq:altSpec3}
\sum_{\chi}\int_{Z_{P}^{\infty}}\int_{A_R^{\infty}}\int_{M_{P_H}(F) \bsl M(\A)^{P,1}}
\int_{[M_{R'}]^{1}}
\sigma_{P_{n+1}}^{Q_{n+1}}(H_{P_{n+1}}(z) - T)
\brtau_{\tlS}(H_{\tlS}(a) - T')\sigma_{R}^{S,\tlS}(H_{R}(a) - T')  \\
\|z\|^{r_0} \|a\|^{r_0}\|m_1\|^{-r_1}\|m_2\|^{-r_2}
\sup_{m_1' \in M_{P_{H}}(F) \bsl M(\A)^{P,1}}
\sup_{m_2' \in M_{R'}(F) \bsl M_{R'}(\A)^{1}} \\
\|m_1'\|^{-r_1'}\|m_2'\|^{-r_2'}
|\int_{[N_{P}]}\int_{[N_{R}]}   \dsl \sum_{P^+ \sbs \tlP \sbs Q^+} \epsilon_{P}^{\tlP} K_{\tlP,\chi}(n_1zm_1',n_2am_2') \rb dn_2dn_1
|\,dm_2dm_1dadz
\end{multline}
where $r_0$ is a fixed number and $r_1$, $r_2$, $r_1'$, $r_2'$ are arbitrarily large.

We set
\[
\Psi(z,a,m_1, m_2) = \sum_{\chi} |\int_{[N_{P}]}\int_{[N_{R}]}   \dsl \sum_{P^+ \sbs \tlP \sbs Q^+} \epsilon_{P}^{\tlP} K_{\tlP,\chi}(n_1zm_1,n_2am_2) \rb dn_2dn_1|.
\]

Given that the proof of Theorem \ref{thm:kappak} assertion 1 has been reduced to proving convergence of the integral \eqref{eq:altSpec3}, the following lemma will conclude the proof of the Theorem. 

\begin{lemme}\label{lem:lemCle} For all $N \ge 0$, there exists an $M$ such that
\[
\sigma_{P_{n+1}}^{Q_{n+1}}(H_{P_{n+1}}(z) - T) 
\brtau_{\tlS}(H_{\tlS}(a) - T')\sigma_{R}^{S,\tlS}(H_{R}(a) - T')\Psi(z,a,m_1, m_2) \le 
\|z\|^{-N} \|a\|^{-N} \|m_1\|^{M}\|m_2\|^{M}.
\]
\end{lemme}

To prove \ref{lem:lemCle} we first need a bound on $\Psi(z,a,m_1, m_2)$. 
It will be more natural to bound $\Psi(z,a,m_1, m_2)^2$.

Before we proceed let us make some remarks and establish some notation.
Recall that $P = P_n \times P_{n+1}$ and $Q = Q_n \times Q_{n+1}$ with $P_{n+1} \cap G_{n} = P_{n}$ 
and $Q_{n+1} \cap G_{n} = Q_{n}$. We write $R = R_{n} \times R_{n+1}$, $S = S_{n} \times S_{n+1}$, $\tlS = \tlS_{n} \times \tlS_{n+1}$. 
All parabolic subgroups are standard in their respective ambient groups. 

Recall that $B_n$ and $B_{n+1}$ are fixed Borel subgroups of $G_n$ and $G_{n+1}$ with $B_{n+1} \cap G_{n} = B_{n}$. 
The inclusion $G_{n} \hookrightarrow G_{n+1}$ induces the inclusion $\Delta_{B_{n}}^{\vee} \hookrightarrow \Delta_{B_{n+1}}^{\vee}$. 
The latter inclusion induces therefore a natural inclusion $\iota : \hDelta_{B_{n}} \hookrightarrow \hDelta_{B_{n+1}}$. 
We let $\varpi^{n+1}$ be the unique element of $\hDelta_{B_{n+1}} \smin \iota( \hDelta_{B_{n}})$. 
Its restriction to $\ago_{B_{n}}$ equals the determinant divided by $n+1$. 
Note that, since $P_{n+1} \cap G_{n} = P_{n}$, the set $\hDelta_{P_{n+1}} \smin \iota(\hDelta_{P_{n}})$ is either empty or 
consists solely of $\varpi^{n+1}$.

Having introduced these, we can write $P^+ = P^+_{n} \times P^+_{n+1}$ and $Q^+ = Q^+_{n} \times Q^+_{n+1}$, both standard parabolic subgroups 
and elements of $\fc_{RS}$. We have then
\[
\hDelta_{P^+_{n+1}} = \dsl \iota(\hDelta_{R_{n}}) \cup \{\varpi^{n+1}\} \rb \cap \hDelta_{R_{n+1}} \cap \hDelta_{P_{n+1}}, \quad 
\hDelta_{Q^+_{n+1}} = \iota(\hDelta_{S_{n}}) \cup \hDelta_{S_{n+1}} \cup \hDelta_{Q_{n+1}}. 
\]

For $w \in W$ we write $w = (w_n, w_{n+1}) \in W^{G_{n}} \times W^{G_{n+1}}$. Since $W^{G_{n}}$ embeds 
naturally into $W^{G_{n+1}}$ we view $w_n$ as element of $W^{G_{n+1}}$. 
For any $s \in W^{G_{n+1}}$ let $\hDelta_{s}$ be the set of $\varpi \in \hDelta_{B_{n+1}}$  stabilized by $s$. 

Using an easy extension of Lemma 2.3 of \cite{ar2} to Schwartz functions, for $\tlP \supset P$, we have
\begin{equation}\label{eq:kpchiConstTerm}
\int_{[N_P]}K_{\tlP, \chi}(nx, y)\,dn = 
\sum_{w \in W^{\tlP}/W^{ P}}
\sum_{\delta \in (P \cap w P_0 w^{-1})(F) \bsl P_0(F)}
K_{P, \chi}(x,w^{-1} \delta y).
\end{equation}

Let $\Omega'$ be the subset of $(w,w') \in W^{Q^+} / W^P \times W^{Q^+} / W^R$ satisfying
\begin{equation}\label{eq:twistCond}
\hDelta_{P^+_{n+1}} \cap \hDelta_{w_n}  \cap  \hDelta_{w_{n+1}}  = \hDelta_{Q^+_{n+1}}, \quad 
\hDelta_{P^+_{n+1}} \cap \hDelta_{w'_n}  \cap  \hDelta_{w'_{n+1}}  = \hDelta_{Q^+_{n+1}}.
\end{equation}

Using the above notation, applying \eqref{eq:kpchiConstTerm}, and its analogue for the group $R$, 
taking into consideration cancellations in alternating sums,
we see that 
$\Psi(z,a,m_1, m_2)^2$ is bounded by a 
sum over $(w,w') \in \Omega'$ of 
\begin{equation}\label{eq:afterAlt}
\dsl 
\sum_{\chi}
\sum_{\delta \in (P \cap w P_0 w^{-1})(F) \bsl P_0(F)}
|K_{P, \chi}(z m_1,w^{-1}\delta a m_2)| \rb 
\dsl \sum_{\chi} \sum_{\delta \in (R \cap w P_0 w^{-1})(F) \bsl P_0(F)}
|K_{R, \chi}((w')^{-1}\delta z m_1, a m_2)| \rb
\end{equation}
here, we drop the remaining compact unipotent integration as it won't affect the bounds. 

Fix $(w,w') \in \Omega'$, it's enough to focus on one pair.
We have the following natural variant of Lemma \ref{lem:maj-noy}, equation \eqref{ineq noy}: 
for all $N > 0 $ there exists $N' > 0 $ such that for all $m \ge 0$, , $x,y \in [M_P]^1$, $z \in A_P^{\infty}$ 
and $\la \in \ago_{P}^*$ we have 
\begin{equation}\label{eq:sumchi2}
\sum_{\chi}|K_{P,\chi}(x,zy)| \ll \|x\|^{-N}\|y\|^{N'}e^{-m\|H_G(z)\|} e^{\langle \la, H_{P}(z) \rangle}.
\end{equation}

Using the bound \eqref{eq:sumchi2} above and Lemma \ref{lem:EisBound} we bound \eqref{eq:afterAlt} 
by
\begin{equation}\label{eq:boundExp}
\|m_1\|^{N_2'-N_1}\|m_2\|^{N_1'-N_2}\|a\|^{r_3} \|z\|^{r_3} e^{-r_4\|H_{G}(a)\|} e^{\langle w\la_1, H_{0}(a) \rangle-\langle \la_1, H_{0}(z) \rangle} \cdot 
e^{\langle w'\la_2, H_{0}(z) \rangle-\langle \la_2, H_{0}(a) \rangle}
\end{equation}
times a power of , where $\la_1 \in \ago_{P}^{*, +}$ 
and $\la_2 \in \ago_{R}^{*, +}$ are for us to be chosen appropriately and $r_3$ and $r_4$ are independent of $\la_1$ and $\la_2$. 
Moreover, choosing $r_1'$ and $r_2'$ in \eqref{eq:altSpec3} sufficiently large we don't have to worry about 
the constants $N_1$, $N_1'$, $N_2$, $N_2'$ (they play a role in the proof of Theorem \ref{thm:jgood}). 
Note that for $a \in A_{R}^{\infty}$ the term $H_{G}(a)$  is not affected by truncation 
but thanks to the factor $e^{-r_4\|H_{G}(a)\|}$ the integral over $A_{G}^{\infty}$ is convergent, we can assume then 
that $a \in A_{R}^{G,\infty}$.

Recall the space $\zgo_{P} \sbs \ago_{0}$ which equals the diagonally embedded $\ago_{P_n} \cap \ago_{P_{n+1}}$. 
Our goal is then, for all $N \ge 0$, to choose $\la_1 \in \ago_{P}^{*, +}$ 
and $\la_2 \in  \ago_{R}^{*, +}$ so that
for $Z \in \zgo_{P}$ and $H = (H_1, H_2) \in \ago_{R_n}^{G_n} \times \ago_{R_{n+1}}^{G_{n+1}}$ 
such that
\begin{equation}\label{eq:coneCond}
\sigma_{P_{n+1}}^{Q_{n+1}}(Z- T) 
\sigma_{R_n}^{S_n,\tlS_{n}}(H_{n} - T_n')
\sigma_{R_{n+1}}^{S_{n+1},\tlS_{n+1}}(H_{n+1} - T_{n+1}')
\brtau_{\tlS_{n}}(H_{n} - T_{n}')
\brtau_{\tlS_{n+1}}(H_{n+1} - T_{n+1}') = 1
\end{equation}
we have that 
\[
\langle \la_1 - w' \la_2, Z \rangle + 
\langle \la_2 - w \la_1, H \rangle  \gg \|Z\|^{N} + \|H\|^N.
\]
More specifically, 
since the natural projection of $\zgo_{P}$ onto $\ago_{P_{n+1}}^{G_{n+1}}$ is an isomorphism, 
using Corollary 6.2 of \cite{ar1}, it is enough to show the following lemma. 

\begin{lemme}\label{lem:ineqCle}
There exist
$\la_1 \in \ago_{P}^{*, +}$ 
and $\la_2 \in  \ago_{R}^{*, +}$ such that under \eqref{eq:coneCond} we have
\begin{multline*}
\langle \la_1 - w' \la_2, Z \rangle + 
\langle \la_2 - w \la_1, H \rangle \\
\gg_{T,T'}
\sum_{\al \in \Delta_{P_{n+1}}^{Q_{n+1}}} \langle\al, Z \rangle + 
 \sum_{\al \in \Delta_{R_{n}}^{S_{n}}} \langle \al, H_n \rangle + 
\sum_{\al \in \Delta_{R_{n+1}}^{S_{n+1}}} \langle \al, H_{n+1} \rangle - 
 \sum_{\al \in \Delta_{\tlS_{n}}} \langle \al, H_n \rangle - 
  \sum_{\al \in \Delta_{\tlS_{n+1}}} \langle \al, H_{n+1} \rangle
\end{multline*} 
where by $\gg_{T,T'}$ we mean an inequality 
up to an additive constant that depends on $T$ and $T'$. 
\end{lemme}

Indeed, given $\la_1, \la_2$ as in Lemma \ref{lem:ineqCle} above, 
we can then take their arbitrary multiples ensuring, through a reasoning explained 
at the end of Theorem 2.2 of \cite{Z0}, the desired property of Lemma \ref{lem:lemCle}. 

We focus on proving Lemma \ref{lem:ineqCle} from now on. 
We introduce the following notation that will save some space in what follows. 
Let 
\[
\underline{a} =  \{ a_{\varpi}\}_{\varpi \in \hDelta_{P_{n}} \sqcup \hDelta_{P_{n+1}} \sqcup \hDelta_{R_{n}} \sqcup \hDelta_{R_{n+1}} }
\]
be a set of numbers, where the sets of weights are treated as disjoint for indexing purposes. 
Define then for $Z \in \zgo_{P}$ and $H = (H_1, H_2) \in \ago_{R_n}^{G_n} \times \ago_{R_{n+1}}^{G_{n+1}}$ 
\begin{multline*}
\la(\underline{a}, Z, H_n, H_{n+1}) := 
\sum_{\varpi \in \hDelta_{P_{n}}}    a_{\varpi} (\langle \varpi, Z \rangle - \langle w_{n}\varpi, H_{n} \rangle) + 
\sum_{\varpi \in \hDelta_{P_{n+1}}}a_{\varpi} (\langle \varpi, Z \rangle - \langle w_{n+1}\varpi, H_{n+1} \rangle)  + \\
    \sum_{\varpi \in \hDelta_{R_{n}}}a_{\varpi} (\langle \varpi, H_n \rangle - \langle w'_{n}\varpi, Z \rangle)+ 
     \sum_{\varpi \in \hDelta_{R_{n+1}}}a_{\varpi} (\langle \varpi, H_{n+1} \rangle - \langle w'_{n+1}\varpi, Z \rangle).
\end{multline*}

\begin{lemme}\label{lem:lem2Cles}
Let $\al_0 \in \Delta_{P_{n+1}}^{Q_{n+1}}$ and let $\varpi_0 \in \hDelta_{P_{n+1}}$ be the corresponding weight. 
There exists a set of positive constants 
\[
\underline{a}_{\al_0} =  \{ a_{\al_0, \varpi}\}_{\varpi \in \hDelta_{P_{n}} \sqcup \hDelta_{P_{n+1}} \sqcup \hDelta_{R_{n}} \sqcup \hDelta_{R_{n+1}} }
\]
 such that, assuming \eqref{eq:coneCond}, we have
\[
\la(\underline{a}_{\al_0}, Z, H_n, H_{n+1})
     \gg_{T,T'} 
     \begin{cases}
     \langle \al_0, Z \rangle - \langle \al_0, H_{n, \tlS_{n}} \rangle \quad 
     &\text{if }\varpi_0 \in (\hDelta_{P_{n+1}} \smallsetminus \hDelta_{Q_{n+1}}) \cap \iota(\hDelta_{\tlS_{n}}) \\
     \langle \al_0, Z \rangle - \langle \al_0, H_{n+1, \tlS_{n+1}} \rangle \quad 
     &\text{if }\varpi_0 \in (\hDelta_{P_{n+1}} \smallsetminus \hDelta_{Q_{n+1}}) \cap \hDelta_{\tlS_{n+1}} \\
     \langle \al_0, Z \rangle \quad &\text{else}
     \end{cases}
\]
where $H_{n, \tlS_{n}}$ and $H_{n+1, \tlS_{n+1}}$ 
are projections of $H_{n}$ and $H_{n+1}$ onto $\ago_{\tlS_{n}}$ and $\ago_{\tlS_{n+1}}$ respectively. 
\end{lemme}

\begin{preuve}
Let $\hDelta = \hDelta_{P_{n+1}} \smallsetminus \hDelta_{Q_{n+1}}$.
Define
\begin{gather*}
\hDelta_{1} = \hDelta \cap \iota(\hDelta_{\tlS_{n}}), \quad 
\hDelta_{2} = \hDelta \cap \hDelta_{\tlS_{n+1}}, \quad
\hDelta_{3} = \hDelta \cap \iota((\hDelta_{S_{n}} \smin \hDelta_{\tlS_{n}})), \quad 
\hDelta_{4} = \hDelta \cap (\hDelta_{S_{n+1}} \smin \hDelta_{\tlS_{n+1}}), \\
\hDelta_{5} = \hDelta \cap \iota (\hDelta_{B_{n}} \smin \hDelta_{R_{n}}), \quad 
\hDelta_{6} = \hDelta \cap (\hDelta_{B_{n+1}} \smin \hDelta_{R_{n+1}}), \quad 
\hDelta_{7} = \hDelta \cap \iota(\hDelta_{R_{n}} \cup \{\varpi^{n+1}\} \smin \hDelta_{S_{n}}) \cap (\hDelta_{R_{n+1}} \smin \hDelta_{S_{n+1}}).
\end{gather*}
The union of $\hDelta_{i}$ is $\hDelta$.
Suppose $\varpi_0 \in \hDelta_{i}$ for $i = 2,4, 6$. Using Lemma 6.1 and Corollary 6.2 of \cite{ar1}, 
We obtain constants   $a_{\al_0, \varpi} \ge 0$ for $\varpi \in \hDelta_{P_{n+1}} \sqcup \hDelta_{R_{n+1}}$
directly from the proof of Lemma A.2 of \cite{Z3} (equation (A.4) specifically), 
so that 
\[
   \sum_{\varpi \in \hDelta_{P_{n+1}}}a_{\al_0, \varpi} (\langle \varpi, Z \rangle - \langle w_{n+1}\varpi, H_{n+1} \rangle) +
     \sum_{\varpi \in \hDelta_{R_{n+1}}}a_{\al_0, \varpi} (\langle \varpi, H_{n+1} \rangle - \langle w'_{n+1}\varpi, Z \rangle) 
     \gg_{T, T'}  \langle \al_0, Z \rangle - \langle \al_0, H_{n+1} \rangle.
\]
We set the rest of the constants $a_{\al_0, \varpi}$ to zero. 
If $\varpi_0 \in \hDelta_{4} \cup \hDelta_{6}$  we conclude directly from Lemma 6.1 of \cite{ar1}. 
Otherwise, if $\varpi_0 \in \hDelta_{2}$, 
the condition \eqref{eq:coneCond} gives $\sigma_{R_{n+1}}^{S_{n+1},\tlS_{n+1}}(H_{n+1} - T_{n+1}') = 1$
and Lemma \ref{lem:obtuse} implies $\langle \al_0, H_{n+1, \tlS_{n+1}} \rangle \gg_{T'} \langle \al_0, H_{n+1} \rangle$.

If $\varpi_0 \in \hDelta_{i}$ for $i = 1,3, 5$ we obtain constants $a_{\al_0, \varpi}$ in the same fashion. 
One just has to observe that elements of $\Delta_{P_{n}}^{Q_{n}}$ 
act identically on $\zgo_{P}$ as their counterparts in $\Delta_{P_{n+1}}^{Q_{n+1}}$. 

Suppose now that $\varpi_0 \in \hDelta_{7}$. By the condition \eqref{eq:twistCond}, we have  $\varpi_{0} \notin \hDelta_{w_{n}'}$ 
or $\varpi_{0} \notin \hDelta_{w_{n+1}'}$. Suppose the latter holds. 
Using the results recalled in the beginning of Appendix A of \cite{Z3}, as well as Lemma 6.1 of \cite{ar1}, 
we have $\langle \varpi_0 - w_{n+1}'\varpi_{0}, Z \rangle \gg_{T} c_{0} \langle \al_{0}, Z \rangle$ for some $c_0 > 0$ 
and $\langle \varpi_0 - w_{n+1}\varpi_{0}, H_{n+1} \rangle \gg_{T'} 0$. 
This identifies the the desired constants in this case. 
We deal with the other case in the same fashion, which completes the proof of Lemma \ref{lem:lem2Cles}.
\end{preuve}

\begin{remarque}\label{rem:lem2Cles}
In the above proof, we do not use the fact that $T \in \ago_{0, n+1}^+$. Assuming this, 
we simply have $\langle \varpi_0 - w_{n+1}'\varpi_{0}, Z\rangle \gg c_{0} \langle \al_{0}, Z \rangle$ when $\varpi_0 \in \hDelta_7$, 
and similarly, Lemma A.2 in \cite{Z3} 
invoked several times yields inequalities independent of $T$, as is clear by inspection of its proof in loc. cit. 
Therefore, if $T \in \ago_{0, n+1}^+$, the implied additive constant in Lemma \ref{lem:lem2Cles} is independent of $T$.
\end{remarque}

The proof of the following lemma is similar and easier to that of the preceding one. 

\begin{lemme}\label{lem:lem3Cles}
Let $\al_n \in \Delta_{R_{n}}^{S_{n}}$, 
$\al_{n+1} \in \Delta_{R_{n+1}}^{S_{n+1}}$
 and let $\varpi_n \in \hDelta_{R_{n}}$ 
 and $\varpi_{n+1} \in \hDelta_{R_{n+1}}$  be their corresponding weights. 
There exists positive constants 
\[
\underline{a}_{\al_n} =  \{ a_{\al_n, \varpi}\}_{\varpi \in \hDelta_{P_{n}} \sqcup \hDelta_{P_{n+1}} \sqcup \hDelta_{R_{n}} \sqcup \hDelta_{R_{n+1}} }
\]
and 
\[
\underline{a}_{\al_{n+1}} =  \{ a_{\al_{n+1}, \varpi}\}_{\varpi \in \hDelta_{P_{n}} \sqcup \hDelta_{P_{n+1}} \sqcup \hDelta_{R_{n}} \sqcup \hDelta_{R_{n+1}} }
\]
such that, assuming \eqref{eq:coneCond}, we  have
\[
\la(\underline{a}_{\al_n}, Z, H_n, H_{n+1})
     \gg_{T,T'} 
     \begin{cases}
     \langle \al_n, H_n \rangle - \langle \al_n, Z \rangle \quad 
     &\text{if }\varpi_n \in (\hDelta_{R_{n}} \smallsetminus \hDelta_{S_{n}}) \cap  (\hDelta_{P_{n}} \smallsetminus \hDelta_{Q_{n}}) \\
     \langle \al_n, H_n \rangle \quad &\text{else}
     \end{cases}
\]
and 
\[
\la(\underline{a}_{\al_{n+1}}, Z, H_n, H_{n+1})
     \gg_{T,T'} 
     \begin{cases}
     \langle \al_{n+1}, H_{n+1} \rangle - \langle \al_{n+1}, Z \rangle \ 
     &\text{if }\varpi_{n+1} \in (\hDelta_{R_{n+1}} \smallsetminus \hDelta_{S_{n+1}}) \cap  
     (\hDelta_{P_{n+1}} \smallsetminus \hDelta_{Q_{n+1}}) \\
     \langle \al_{n+1}, H_{n+1} \rangle \ &\text{else}.
     \end{cases}
\]
\end{lemme}

Finally, we have the following result

\begin{lemme}\label{lem:lemDer}
 Let $\al_{n} \in \Delta_{\tlS_n}$, $\al_{n+1} \in \Delta_{\tlS_{n+1}}$ 
 and let $\varpi_n \in \hDelta_{\tlS_{n}}$ 
 and $\varpi_{n+1} \in \hDelta_{\tlS_{n+1}}$  be their corresponding weights. 
Then, there exists positive constants 
\[
\underline{a}_{\al_n} =  \{ a_{\al_n, \varpi}\}_{\varpi \in \hDelta_{P_{n}} \sqcup \hDelta_{P_{n+1}} \sqcup \hDelta_{R_{n}} \sqcup \hDelta_{R_{n+1}} }
\]
and 
\[
\underline{a}_{\al_{n+1}} =  \{ a_{\al_{n+1}, \varpi}\}_{\varpi \in \hDelta_{P_{n}} \sqcup \hDelta_{P_{n+1}} \sqcup \hDelta_{R_{n}} \sqcup \hDelta_{R_{n+1}} }
\]
such that, assuming \eqref{eq:coneCond}, $\varpi_n \notin (\hDelta_{P_{n}} \smallsetminus \hDelta_{Q_{n}})$
 and $\varpi_{n+1} \notin (\hDelta_{P_{n+1}} \smallsetminus \hDelta_{Q_{n+1}})$ 
 we have
\[
\la(\underline{a}_{\al_n}, Z, H_n, H_{n+1})
     \gg_{T,T'} - \langle \al_n ,H_n \rangle
\]
and 
\[
\la(\underline{a}_{\al_{n+1}}, Z, H_n, H_{n+1})
     \gg_{T,T'} - \langle \al_{n+1}, H_{n+1} \rangle.
\]
Otherwise, we define $\underline{a}_{\al_n}$ and $\underline{a}_{\al_{n+1}}$ as in  \ref{lem:lem2Cles} by identification of roots. 
\end{lemme}

\begin{preuve}
Let's prove the statement concerning $\al_n$, the other case being analogous. 
Let $\gamma_{n} \in \Delta_{B_{n}} \smin \Delta_{B_{n}}^{\tlS_{n}}$ be its lift. 
We must have  $\varpi_{n} \in (\hDelta_{B_{n}} \smallsetminus \hDelta_{P_{n}}) \cup \hDelta_{Q_{n}}$. Then, as explained 
in the proof of Lemma \ref{lem:lem2Cles}, Lemmas 
 A.2 of \cite{Z3} 
and 6.1 of \cite{ar1} provide constants $\underline{a}_{\al_n}$ 
such that
\[
\la(\underline{a}_{\al_n}, Z, H_n, H_{n+1})
     \gg_{T,T'} - \langle \gamma_n, H_n \rangle
\]
Condition \eqref{eq:coneCond} and Lemma \ref{lem:obtuse} 
prove that $\langle \al_n, H_n \rangle \gg_{T'} \langle \gamma_n, H_n \rangle$ which allows to conclude. 
\end{preuve}

We are ready to define the desired $\la_1 = (\la_{1,n},\la_{1,n+1})  \in \ago_{P_{n}}^* \times \ago_{P_{n+1}}^*$ 
and $\la_2 = (\la_{2,n}, \la_{2,n+1}) \in \ago_{R_n}^* \times \ago_{R_{n+1}}^*$ 
of Lemma \ref{lem:ineqCle}. 
In Lemmas \ref{lem:lem2Cles}, \ref{lem:lem3Cles} and \ref{lem:lemDer} 
we have defined sets $\underline{a}_{\al}$ 
for all $\al \in \Delta_{P_{n+1}} \sqcup \Delta_{R_{n}}^{S_n} \sqcup \Delta_{R_{n+1}}^{S_{n+1}} 
\sqcup \Delta_{\tlS_{n}} \sqcup \Delta_{\tlS_{n+1}}$. 
Let $k \ge 0$. We set
\begin{align*}
\la_{1,n} &=& \sum_{\varpi \in \hDelta_{P_n}} \dsl \sum_{\al \in \Delta_{P_{n+1}} }ka_{\al, \varpi} + 
\sum_{\al \in \Delta_{R_{n}}^{S_n} \sqcup \Delta_{R_{n+1}}^{S_{n+1}} 
\sqcup \Delta_{\tlS_{n}} \sqcup \Delta_{\tlS_{n+1}} }a_{\al, \varpi}\rb \varpi, \\
\la_{1,n+1} &=&  \sum_{\varpi \in \hDelta_{P_{n+1}}} \dsl \sum_{\al \in \Delta_{P_{n+1}} }ka_{\al, \varpi} + 
\sum_{\al \in \Delta_{R_{n}}^{S_n} \sqcup \Delta_{R_{n+1}}^{S_{n+1}} 
\sqcup \Delta_{\tlS_{n}} \sqcup \Delta_{\tlS_{n+1}} }a_{\al, \varpi}\rb  \varpi, \\
\la_{2,n} &=& \sum_{\varpi \in \hDelta_{R_{n}}} \dsl \sum_{\al \in \Delta_{P_{n+1}} }ka_{\al, \varpi} + 
\sum_{\al \in \Delta_{R_{n}}^{S_n} \sqcup \Delta_{R_{n+1}}^{S_{n+1}} 
\sqcup \Delta_{\tlS_{n}} \sqcup \Delta_{\tlS_{n+1}} }a_{\al, \varpi}\rb  \varpi, \\
\la_{2,n+1} &=& \sum_{\varpi \in \hDelta_{R_{n+1}}} \dsl \sum_{\al \in \Delta_{P_{n+1}} }ka_{\al, \varpi} + 
\sum_{\al \in \Delta_{R_{n}}^{S_n} \sqcup \Delta_{R_{n+1}}^{S_{n+1}} 
\sqcup \Delta_{\tlS_{n}} \sqcup \Delta_{\tlS_{n+1}} }a_{\al, \varpi}\rb  \varpi.
\end{align*}

We have then for $Z \in \zgo_{P}$ and $H = (H_1, H_2) \in \ago_{R_n}^{G_n} \times \ago_{R_{n+1}}^{G_{n+1}}$ 
\begin{multline}\label{eq:alltogether}
\langle \la_1 - w' \la_2, Z \rangle + 
\langle \la_2 - w \la_1, H \rangle = \\
\sum_{\al \in \Delta_{P_{n+1}}}k\la(\underline{a}_{\al}, Z, H_n, H_{n+1}) + 
 \sum_{\al \in \Delta_{R_{n}}^{S_n} \sqcup \Delta_{R_{n+1}}^{S_{n+1}} 
\sqcup \Delta_{\tlS_{n}} \sqcup \Delta_{\tlS_{n+1}}}
\la(\underline{a}_{\al}, Z, H_n, H_{n+1}) 
\end{multline}
which, by Lemmas \ref{lem:lem2Cles}, \ref{lem:lem3Cles} and \ref{lem:lemDer} 
proves the desired inequality of Lemma \ref{lem:ineqCle} for $k$ sufficiently large (the constant $k$ 
is there is to offset the $-\langle \al, Z \rangle$ terms of 
Lemma \ref{lem:lem3Cles}).

We have proven Lemma \ref{lem:ineqCle} which entails Lemma \ref{lem:lemCle} which 
accomplishes the proof of assertion 1 of Theorem \ref{thm:kappak}. 
\end{paragr}

\begin{paragr}[Proof of assertions 2 and 3 of Theorem \ref{thm:kappak}.] --- \label{S:proof-kappak}
If we replace $\kappa_{\chi}^{T}$ by $K_\chi^T$ and assume that $f$ is compactly supported, statement 2 is proved in   \cite{Z3}. 
Practically the same method applies to $\kappa_{\chi}^{T}$ as well as Schwartz functions, since we have proven the convergence assertion of Theorem \ref{thm:kappak}, 
which gives rise to the distribution $i_{\chi}$. Its continuity follows immediately from the examination of the proof of assertion 1 and the analogue 
of Theorem 3.7 of \cite{Z3}. 
What  is not obvious however is that $i_{\chi}(f) = I_{\chi}(f)$. We will prove this for compactly supported $f$ and the general case will follow by continuity. 

As in \S\S \ref{S:proof-schwartzCVG}, 
unless otherwise stated, all sums of the type $\Sigma_{P}$ or $\Sigma_{P \sbs Q}$ are over elements of $\fc_{RS}$. 

Let $Q \in \fc_{RS}$. Define
\begin{multline*}
K_{\chi}^{T, Q}(x,y) = 
\sum_{Q \supset P}
\epsilon_{P}^{Q}
\sum_{\gamma \in (P \cap H) (F)\bsl (Q \cap H)(F)} 
\sum_{\delta \in P' (F) \bsl Q'(F)}
\htau_{P_{n+1}}^{Q_{n+1}}(H_{P_{n+1}}(\delta_n y_n) - T_{P_{n+1}})K_{P, \chi}(\gamma x, \delta y).
\end{multline*}
For $P \subset Q$, elements of $\fc_{RS}$, let
\[
\ago_{P, Q}^{*, - } = \{ \la \in \ago_{P_{n+1}}^{*} \ | \  \langle \la + 2\upla_{Q}, \varpi^{\vee} \rangle < 0, \ \forall \varpi \in \hDelta_{Q}\}
\]
where for $Q \in \fc_{RS}$ we set
\[
\upla_{Q} = \rho_{Q_{n+1}} - \rho_{Q_n} \in \ago_{0, n+1}.
\]
Define for $\la \in \ago_{Q_{n+1}, \CC}^{*}$
\[
I_{\chi}^{T}(f, Q, \la) := 
\int\limits_{(Q \cap H)(F) \bsl H(\AAA)}
\int\limits_{ Q' (F) \bsl G'(\AAA)}
e^{\langle \la, H_{Q_{n+1}}(x) \rangle}
\tau_{Q_{n+1}}(H_{Q_{n+1}}(y_n) - T_{Q_{n+1}})
K_{\chi}^{T, Q}(x,y)\upeta_{G'}(y)\,dydx.
\]
The integral converges absolutely 
for $\Re (\la) \in \ago_{Q, Q}^{*, - } $
and admits a meromorphic continuation, holomorphic at $\la = 0$ so that we have
\begin{equation}\label{eq:sumT}
I_{\chi}(f) = \sum_{Q } I_{\chi}^{T}(f, Q, 0)
\end{equation}
for any $T$. The equality above is a formal consequence of Theorem 3.7 of \cite{Z3}.

Let $P \in \fc_{RS}$. For a function $\phi : P_n (F) \bsl G_n(\AAA) \to \CC$,  we introduce the operator $\La^{T, P}_{m'}$: it is  a 
variant  of the mixed operator $\La^{T, P}_{m}$ of \cite{JLR}. It  is defined in  \cite{Z3} §2.3 (where it is denoted by  $\La^{T, P}_{m}$). We have 
 \[
\La^{T, P}_{m'}\phi(y_n) = \sum_{P \supset Q}\epsilon_{Q}^{P}
\sum_{\delta \in Q_n' (F) \bsl P_n' (F)} 
\htau_{Q_{n+1}}^{P_{n+1}}(H_{Q_{n+1}}(\delta y_n) - T_{Q_{n+1}}) \phi_{Q_n}(\delta y_n), \quad y_n \in  P_n' (F) \bsl G_n'(\AAA).
\]
We will also use the operator $\La_{d}^{T,P}$ as in \eqref{eq:ladDef}.

Fix $R \subset R^{1} \subset Q$ and $S \subset S^{1} \subset Q$ all elements of $\fc_{RS}$. 
Let $\la \in \ago_{R_{n+1},\CC}^{*}$.
We define $I_{\chi}^{T}(f, R, R^{1}, S, S^{1}, Q, \la) $ as
\begin{multline*}
\int\limits_{(R \cap H)(F) \bsl H(\AAA)}
\int\limits_{(S_n' \times G_{n+1}') (F) \bsl G'(\AAA)}
e^{\langle \la, H_{R_{n+1}}(x) \rangle}
\tau_{Q_{n+1}}(H_{Q_{n+1}}(x) - T'_{Q_{n+1}})
\tau_{Q_{n+1}}(H_{Q_{n+1}}(y_n) - T_{Q_{n+1}}) \\
\sigma_{R_{n+1}}^{R^{1}_{n+1}, Q_{n+1}}(H_{R_{n+1}}(x) - T'_{R_{n+1}})
\sigma_{S_{n+1}}^{S_{n+1}^{1}, Q_{n+1}}(H_{S_{n+1}}(y_n) - T_{S_{n+1}})
\La_{d}^{T', R}\La_{m'}^{T, S} 
\\
\dsl 
\sum_{
\begin{subarray}{c}
R \cup S \subset P \subset R^1 \cap S^1
\end{subarray}
} \epsilon_{P_{n+1}}^{Q_{n+1}}
\sum_{\delta_{n+1} \in P_{n+1}'(F) \bsl G_{n+1}'(F)}K_{P, \chi}(x, y_n, \delta_{n+1} y_{n+1} )
\rb  \upeta_{G'}(y)\,dxdy
\end{multline*}
where $T'$ is a translate of $T$ by an element of $\ago_{n+1}$ that depends only on the support of $f$, 
$\La_{m'}^{T, S}$ is applied with respect to $y$ and $\La^{T', R}_{d}$ with respect to $x$. 
The proof of Theorem 3.1 in \cite{Z3} proves that the above integral
converges absolutely for any $\Re(\la) \in \ago_{R, Q}^{*,-}$ and extends everywhere to a meromorphic function, holomorphic at $0$. 
Moreover, we have the equality of meromorphic functions on $\ago_{Q_{n+1}, \CC}^{*}$. 
\[
I_{\chi}^{T}(f, Q, \la)  = \sum_{R \subset R^{1} \subset Q}\sum_{S \subset S^{1} \subset Q}I_{\chi}^{T}(f, R, R^{1}, S, S^{1}, Q, \la).
\]

Fix $R \in \fc_{RS}$ and a $\la \in \ago_{R_{n+1}, \CC}^{*}$  such that $\Re(\la) \in \ago_{R,R}^{*,-}$. 
Note that $ \ago_{R, R}^{*,-} \subset  \ago_{R, Q}^{*,-}$ for any $Q \supset R$.
Fix  also $S \in \fc_{RS}$. 
We look at the sum
\begin{equation*}\label{eq:theGoodsum0}
\sum_{Q}
\sum_{R \subset R^{1} \subset Q}
\sum_{S \subset S^{1} \subset Q}
I_{\chi}^{T}(f, R, R^{1}, S, S^{1}, Q, \la).
\end{equation*}
We claim it equals $I_{\chi}^{T}(f, R, S, \la)$
defined by meromorphic continuation of 
\begin{multline*}\label{eq:theGoodsum}
\int\limits_{(R \cap H)(F) \bsl H(\AAA)}
\int\limits_{(S_n' \times G_{n+1}') (F) \bsl G'(\AAA)}
e^{\langle \la, H_{R_{n+1}}(x) \rangle}
\tau_{R_{n+1}}(H_{R_{n+1}}(x) - T'_{R_{n+1}})
\tau_{S_{n+1}}(H_{S_{n+1}}(y_n) - T_{S_{n+1}}) \\
\La_{d}^{T', R}\La_{m'}^{T, S} K_{\chi}(x,y)\upeta_{G'}(y)\,dxdy.
\end{multline*}
Fix $P \in \fc_{RS}$ containing $R$ and $S$. We see that we need to show that
\begin{multline*}
K_{P, \chi}(x,y) \dsl
\sum_{Q \supset P}
\epsilon_{P_{n+1}}^{Q_{n+1}}
\tau_{Q_{n+1}}(X)
\tau_{Q_{n+1}}(Y)
\sum_{P \subset R^{1}, S^{1} \subset Q}
\sigma_{R_{n+1}}^{R^{1}_{n+1}, Q_{n+1}}(X)
\sigma_{S_{n+1}}^{S^{1}_{n+1}, Q_{n+1}}(Y)  \rb = \\
\begin{cases} 
K_{\chi}(x,y)\tau_{R_{n+1}}(X)\tau_{S_{n+1}}(Y) & \text{ if }P = G, \\
0 & \text{ if }P \neq G, 
\end{cases}
\end{multline*}
where $X = H_{R_{n+1}}(x) - T'_{R_{n+1}}$, $Y = H_{S_{n+1}}(y_n) - T_{S_{n+1}}$. 

Since $\sum_{R^1 \supset P} \sigma_{R_{n+1}}^{R^{1}_{n+1}, Q_{n+1}} = \tau_{R_{n+1}}^{P_{n+1}}\htau_{P_{n+1}}^{Q_{n+1}}$ 
we need to consider
\[
K_{P, \chi}(x,y) 
\tau_{R_{n+1}}^{P_{n+1}}(X)\tau_{S_{n+1}}^{P_{n+1}}(Y)
\dsl
\sum_{Q \supset P}
\epsilon_{P_{n+1}}^{Q_{n+1}}
\tau_{Q_{n+1}}(X) \htau_{P_{n+1}}^{Q_{n+1}}(X)
\tau_{Q_{n+1}}(Y)\htau_{P_{n+1}}^{Q_{n+1}}(Y) \rb. 
\]
We obtain the desired result if $P = G$. We can assume then $P \neq G$ and that 
$K_{P, \chi}(x,y) \neq 0$. We want to show the expression is zero. 
 By the argument of the beginning of the proof of Theorem 3.1  in \cite{Z3}, and definition of $T'$, 
 under the assumption $K_{P, \chi}(x,y) \neq 0$ we have that 
 $\tau_{Q_{n+1}}(X) \htau_{P_{n+1}}^{Q_{n+1}}(X)
\tau_{Q_{n+1}}(Y)\htau_{P_{n+1}}^{Q_{n+1}}(Y) =  \tau_{Q_{n+1}}(Y)\htau^{Q_{n+1}}_{P_{n+1}}(Y)$ which shows the desired vanishing by 
the Langlands Combinatorial Lemma.

We have shown thus that 
\begin{equation*}
 I_{\chi}(f) = \sum_{R,S \in \fc_{RS}}I_{\chi}^{T}(f, R, S, 0)
\end{equation*}
 where $T$ is arbitrary. 
 Note however that if we fix $R$ 
 and set $I_{\chi}^{T}(f, R, \la)$ to be
 \begin{equation}\label{eq:iR}
\int\limits_{ (R \cap H) (F) \bsl H(\AAA)}
\int\limits_{ G'(F)\bsl G'(\AAA)}
e^{\langle \la, H_{R_{n+1}}(x) \rangle}
\tau_{R_{n+1}}(H_{R_{n+1}}(x) - T_{R_{n+1}})
\La_{d}^{T, R} K_{\chi}(x,y)\upeta_{G'}(y)\,dxdy
\end{equation}
 we get, using the inversion formula Lemme 2.4 \cite{Z3}:
 \[
 \sum_{S \in \fc_{RS}} I_{\chi}^{T}(f, R, S, \la) =  I_{\chi}^{T'}(f, R, \la) .
 \]
This entails that
\begin{equation}\label{eq:iFinal}
 I_{\chi}(f) = \sum_{R \in \fc_{RS}}I_{\chi}^{T'}(f, R, 0), \quad T' \in \ago_{n+1}.
\end{equation}

Similar and more direct reasoning yields a formula for $i_{\chi}$ as follows.
 We define $i_{\chi}^T(f,R,S,Q, \la)$, 
for $R \sbs S \sbs Q$, elements of $\fc_{RS}$, 
and $\la \in \ago_{R, \CC}^*$, 
 to be
\begin{multline*}
\int\limits_{(R \cap H)(F) \bsl H(\AAA)}
\int\limits_{[G']}
e^{\langle \la, H_{R_{n+1}}(x) \rangle}
\tau_{Q_{n+1}}(H_{Q_{n+1}}(x) - T_{Q_{n+1}})
\sigma_{R_{n+1}}^{S_{n+1}, Q_{n+1}}(H_{R_{n+1}}(x) - T_{R_{n+1}})\\
\La_{d}^{T, R}
\dsl 
\sum_{
\begin{subarray}{c}
R \subset P \subset S
\end{subarray}
} \epsilon_{P_{n+1}}^{Q_{n+1}}
\sum_{\delta \in P'(F) \bsl G'(F)}K_{P, \chi}(x, \delta y)
\rb  \upeta_{G'}(y)\,dxdy
\end{multline*}
which converges absolutely for $\Re(\la) \in \ago_{R,Q}^{*, -}$ and admits meromorphic continuation 
satisfying 
\[
\sum_{R \sbs S \sbs Q}i_{\chi}^T(f,R,S,Q, 0) = i_{\chi}(f), \quad T \in \ago_{n+1}.
\]
Additionally, for a fixed $R$ we have
\[
\sum_{S \sbs Q}i_{\chi}^T(f,R,S,Q, \la) = I_{\chi}^T(f,R,\la)
\]
where $I_{\chi}^{T}(f, R, \la)$ is as in \eqref{eq:iR}.
Comparing the above two equalities with \eqref{eq:iFinal}, taking $T' = T$, we get 
$i_{\chi}(f) = I_{\chi}(f)$ as desired.

\end{paragr}

\subsection{Asymptotic formulas}

\begin{paragr} Recall that we have defined integrals $I_{\chi}^T(f)$ and  $i_{\chi}^T(f)$  (see Theorems \ref{thm:jfDef} and \ref{thm:kappak}.

\begin{theoreme}\label{thm:jgood} 
\begin{enumerate}
\item $i_{\chi}^T(f)$ is asymptotically equal to each of the integrals \eqref{eq:LaTrKchi} (see Proposition \ref{prop:cv-LaTrKchi}) and \eqref{eq:FonH}   (see Proposition \ref{prop:FonHG}).
\item $I_{\chi}^T(f)$ is asymptotically equal to the integral \eqref{eq:FonG} (see Proposition \ref{prop:FonHG}).
\end{enumerate}
 \end{theoreme}

 \begin{remarque}
   The above theorem can be stated for the whole kernel $K$ instead of $K_\chi$. It can be proved in the same way.
  \end{remarque}
 
 \begin{preuve}
 Let us revisit the proof of Theorem \ref{thm:kappak} keeping track of constants depending on the parameter $T$. 
 It is not hard to see that Proposition 2.8 of \cite{Z3}, 
 or any analogous result, including the original Lemma 1.4 of \cite{ar2}, gives a constant 
 \begin{equation*}
  e^{r \|T\|}
 \end{equation*}
 for $r$ depending on $r_1$ and  $r_1'$. 
 The expression \eqref{eq:altSpec3} should be multiplied by $e^{r \|T\|}$ accordingly. 
 We move on to the expression \eqref{eq:boundExp}. 
The constants $N_1$ and $N_2$ are arbitrary and their choice is independent of $T$. They influence the constants $N_1'$ and $N_2'$ respectively, however. 
 We take $N_2 = 0$ and choose $N_1$ so that $N_2' - N_1 < 0 $. To make sure that the integral over $m_2$ converges 
 we choose the constant $r_2'$ in \eqref{eq:altSpec3} small enough. This choice is independent of $T$ as well. 
Choosing $\la_1$ (and $\la_2$) appropriately (choice not influencing dependency on $T$), the conclusion of the proof of Theorem \ref{thm:kappak} is 
that the integral \eqref{eq:altSpec3} is bounded by
\begin{equation}\label{eq:sigInt}
e^{r\|T\|}\int_{\zgo_{P}}\sigma_{P}^Q(Z - T) e^{\langle -(\la_1 - w'\la_2), Z \rangle}\,dZ
\end{equation}
However, looking at \eqref{eq:alltogether}, we can make the constant $k$ as big as we please, 
which, by Lemma \ref{lem:lem2Cles}, makes the exponent in the integral above as negative as we wish, without introducing 
dependency on $T$ by Remark \ref{rem:lem2Cles}. 
Reasoning as in the end of Theorem 2.2 of \cite{Z0}, having liberty with the exponent, we can make the integral \eqref{eq:sigInt} 
smaller than any power of $e^{-\|T\|}$ as long as $P \neq G$. We have showed thus that $i_{\chi}^T(f)$ asymptotically equals 
\begin{equation}\label{eq:LadInt}
\int_{[H]}\int_{[G']}\La_{d}^{T} K_{\chi}(x,y)\,\eta_{G'}(y)\,dydx.
\end{equation}

Next, we show that \eqref{eq:LadInt} asymptotically equals 
$\int_{[H]}\int_{[G']} F^{G_{n+1}}(x,T)K_{\chi}(x,y)\, \eta_{G'}(y)\,dydx$.
Indeed, the argument is essentially identical to the one in \cite{IY}, Proposition 3.8. 
One just needs to use the fact that for an element $X$ in the universal envelopping algebra of 
$Lie(G) \otimes_{\RR} \CC$, the function $R(X)K_{\chi}(x,y)$, where $R(X)$ is the right action on $K_{f, \chi}(x,y)$ as a function of the variable $x$, is just 
$K_{R(X)f, \chi}(x,y)$. The uniform growth assumption  
of Proposition 3.8 in \textit{loc. cit.} can then be replaced with the application of the bound \eqref{ineq noy}. 

Exactly the same reasoning works with $\La_{d}^T$ replaced with $\La_{r}^T$ which yield point 1.
Point 2 follows exactly the same reasoning applied to the kernel $K_{\chi}^T$.
 \end{preuve}
\end{paragr}

\section{Flicker-Rallis period of some spectral kernels}
\label{sec:sym}

The goal of this chapter is to get the spectral expansion of the Flicker-Rallis integral of the automorphic kernel attached to a linear group and a specific cuspidal datum (called in §\ref{S:generic} $*$-generic.) This is achieved in theorem \ref{thm:spec-exp-FR-kernel}. It turns out that the decomposition is discrete and is expressed in terms of some relative characters.

\subsection{Flicker-Rallis intertwining periods and related distributions}\label{ssec:related}

\begin{paragr}[Notations.] ---   \label{S:notations-symmetric}
  In all this section, we will fix an integer $n\geq 1$ and we will use notations of §§\ref{S:Gn} to \ref{S:PtoP'}. Since $n$ will be fixed, we will drop the subscript $n$ from the notation: $G=G_n$, $B=B_{n}$ etc. So we do not follow notations of §\ref{S:G=GnxGn+1}: we hope that it will cause no confusion.
\end{paragr}

\begin{paragr}[Flicker-Rallis periods.] ---\label{S:disting}
  Let $\pi$ be a cuspidal automorphic representation of $G(\AAA)$ with central character trivial on $A_M^\infty$. We shall denote by $\pi^*$ the conjugate-dual representation of $G(\AAA)$. We shall say that $\pi$ is self conjugate-dual if  $\pi\simeq\pi^*$ and that $\pi$ is. $G'$-distinguished, resp. $(G',\eta)$-distinguished, if the linear form (called the Flicker-Rallis period)
  \begin{align}\label{eq:lesperiodes}
    \varphi \mapsto \int_{[G']_0} \varphi(h)\,dh, \text{ resp.}  \int_{[G']_0} \varphi(h)\eta(\det(h))\,dh
  \end{align}
  does not vanish identically on $\Ac_\pi(G)$. Then $\pi$ is self conjugate-dual if and only if $\pi$ is either $G'$-distinguished or $(G',\eta)$-distinguished. However it cannot be both. This is related to the well-known factorisation of the Rankin-Selberg factorisation $L(s, \pi\times \pi\circ c)$ where $c$ is the Galois involution of $G(\AAA)$ in terms of Asai $L$-functions and to  the fact that the residue at $s=1$ of the Asai $L$-functions is expressed in terms of Flicker-Rallis periods (see  \cite{Flicker}). 
\end{paragr}

\begin{paragr}
  In this chapter, we will focus on the period in \eqref{eq:lesperiodes} related to distinction. However it is clear that all the results hold \emph{mutatis mutandis} for the period related to $\eta$-distinction. 
\end{paragr}

 \begin{paragr}\label{S:Mpi}
  Let $P=MN_P$ be a standard parabolic subgroup (with its standard decomposition). Let $\pi$ be an irreducible  cuspidal automorphic representation of $M$ with central character trivial on $A_M^\infty$.

  It will be convenient to write  $M=G_{n_1}\times \ldots \times G_{n_r}$ with $n_1+\ldots+n_r=n$.  Accordingly we have  $\pi=\sigma_1\boxtimes\ldots \boxtimes\sigma_r$ where $\sigma_i$ is an irreductible cuspidal representation of $G_{n_i}$. 
\end{paragr}

\begin{paragr}\label{S:J}
 Let $\varphi\in \Ac_{P,\pi}(G)$. The parabolic subgroup $P'=P\cap G'$ of $G'$ has the following Levi decomposition $M'N_{P'}$ where $M'=M\cap G'$. We then define the following integral which is a specific example of a Flicker-Rallis intertwining period introduced by  Jacquet-Lapid-Rogawski (see \cite{JLR} section VII, note that our definition of  $\Ac_{P,\pi}(G)$ is slightly different from theirs), 
  \begin{align*}
    J(\varphi)&=\int_{A_{M'}^\infty M'(F) N_{P'}(\AAA)\back G'(\AAA)}  \varphi(g) \,dg
  \end{align*}
  Clearly we get a $G'(\AAA)$-invariant continuous linear form on $\Ac_{P,\pi}(G)$.
Note that $J$ does not vanish identically if and only if each component  $\sigma_i$ is $G_{n_i}'(F)$-distinguished. In this case, we have $\pi=\pi^*$.
\end{paragr}

\begin{paragr}
   Let $Q\in \pc(M)$.  As recalled in §\ref{S:pcM}, there is a unique pair $(Q',w)$ such that the conditions are satisfied:
  \begin{itemize}
  \item $Q'=wQw^{-1}$ is the standard parabolic subgroup in the $G$-conjugacy class of $Q$ ;
  \item $w\in W(P;Q')$. 
  \end{itemize}
 Let $\la\in \ago_{P,\CC}^{G,*}$. We have $M(w,\la)\varphi\in \Ac_{M_{Q'},w\pi}(G)$  if $\la$ is outside the singular hyperplanes of the intertwining operator. We shall define
  \begin{align}
    \label{eq:JQla-cusp}
J_Q(\varphi,\la)=J(M(w,\la)\varphi)
  \end{align}
as a meromorphic function of $\la$. 
\end{paragr}

\begin{paragr}\label{S:ompi}
  Let  $g\in G(\AAA)$. Let's define for $\varphi,\psi\in \Ac_{P,\pi}(G)$ 
\begin{align}\label{eq:BQg-cusp}
  B_Q(g,\varphi,\psi,\la)=E(g,\varphi,\la)\cdot J_Q(\bar{\psi},-\la).
\end{align}
as a meromorphic function of $\la\in \ago_{P,\CC}^{G,*}$. In fact, by the basic properties of Eisenstein series and intertwining operators, there exists an open subset $\om_\pi\subset \ago_{P,\CC}^{G,*}$ which is the complement of a union of hyperplanes of $\ago_{P,\CC}^{G,*}$  such that :
   \begin{itemize}
   \item $\om_\pi$ contains  $i\ago_{P}^{G,*}$. 
   \item for all $\varphi,\psi\in \Ac_{P,\pi}(G)$, the map $\la\mapsto B_Q(g,\varphi,\psi,\la)$ is holomorphic on $\om_\pi$ and gives for each $\la\in \om_\pi$ a continuous sesquilinear form in $\varphi$ and $\psi$.
   \end{itemize}
\end{paragr}

\begin{paragr}
  Let $f\in \Sc(G(\AAA))$, $g\in G(\AAA)$  and $Q\in \pc(M)$. Let's introduce the distribution
\begin{align}\label{eq:JQpi-cusp}
  J_{Q,\pi}(g,\la,f)=\sum_{\varphi\in \bc_{P,\pi}  }  B_Q(g,I_P(\la,f)\varphi,\varphi,\la)
\end{align}
 where $\bc_{P,\pi} $ is a $K$-basis of $\Ac_{P,\pi}(G)$ (see § \ref{S:K0ON} and $\la\in \om_\pi$ (see § \ref{S:ompi} for the notation  $\om_\pi$).  It follows from propositions  \ref{prop:car-relJB} that  $J_{Q,\pi}(g,\la)$ is a continuous distribution on $\Sc(G(\AAA))$. 

\end{paragr}

\begin{paragr}[A $(G,M)$-family.] ---
  
\begin{proposition}\label{prop:GMfam}
     The family $(J_{Q,\pi}(g,\la,f))_{Q\in \pc(M)}$ is a  $(G,M)$-family in the sense of Arthur (see \cite{arthur2}): namely each map
      \begin{align*}
        \la \in  \ago_{P}^{G,*} \mapsto J_{Q,\pi}(g,\la,f)
      \end{align*}
is smooth on $i\ago_{P}^{G,*}$ (and even  holomorphic on $\om_\pi$) and for adjacent elements $Q_1,Q_2\in \pc(M)$ we have 
\begin{align}\label{eq:glue-cusp}
  J_{Q_1,\pi}(g,\la,f)=J_{Q_2,\pi}(g,\la,f)
\end{align}
   on the hyperplane of $ i\ago_{L}^{G,*}$ defined by $\bg \la,\al^\vee\bd=0$ where $\al$ is the unique element in $\Delta_{Q_1}\cap (-\Delta_{Q_2})$. 
  \end{proposition}

  \begin{preuve}
     The holomorphy on $\om_\pi$ is obvious if $f$ is $K_\infty$-finite (the sum in \ref{eq:JQpi-cusp} is then finite). Let $C'\subset \om_\pi$ be a compact subset. Using approximations of $f$ by $K_\infty$-finite functions, one shows that $J_{Q,\pi}(g,\la, f)$ is a uniform limit on $C'$ of holomorphic functions hence holomorphic.
 
Let  $Q_1,Q_2\in \pc(M)$ be such that $\Delta_{Q_1}\cap (-\Delta_{Q_2})$ is a singleton $\{\al\}$. Let  $ \la\in i\ago_{L}^{G,*}$ such  $\bg \la,\al^\vee\bd=0$. For $i=1,2$ let $Q_i'$ be a standard parabolic subgroup and  $w_i\in W(M,Q'_i)$ be such that  $Q_i'=w_i Q_iw_i^{-1}$.  Let $\beta=w_1 \al\in \Delta_{Q_1'}$ and Let $s_\beta$ the simple reflection associated to $\beta$. Then we have $w_2=s_\be w_1$.  Let $\varphi\in \Ac_{P,\pi}(G)$. Clearly it suffices to  check the equality:
\begin{align*}
 E(g,\varphi,\la)\cdot (\overline{ J(M(w_1,\la)\varphi)}=E(g,\varphi,\la)\cdot\overline{  J(M(w_2,\la)\varphi)}.
\end{align*}
Using the functional equations of intertwining operators and Eisentein series,  we have $ M(w_2,\la)=M(s_\beta,w_1\la)M(w_1,\la)$  and $E(g,\varphi,\la)=E(g,M(w_1,\la)\varphi,w_1\la)$. Thus up to a change of notations (replace $P$ by $Q_1'$), we may assume that $Q_1=P$ and thus $w_1=1$ and $\al=\be$. We are reduced to prove
\begin{align}\label{eq:theequality}
 E(g,\varphi,\la)\cdot \overline{J( \varphi)}=E(g,\varphi,\la)\cdot \overline{J( M(s_\al,\la)\varphi)}
\end{align}
on the hyperplane $\bg \la,\al^\vee\bd=0$. The simple reflection $s_\al$ acts on $M$ as a transposition of two consecutive  blocks of $M$ say $G_{n_i}$ and $G_{n_{i+1}}$. Note that  $M(s_\al,\la)\varphi= M(s_\al,0)$. Then we have even a stronger property:
\begin{align*}
  J( \varphi)=J( M(s_\al,0)\varphi)
\end{align*}
if $n_i\not= n_{i+1}$ or if $n_i= n_{i+1}$  but $\sigma_i\not\simeq\sigma_{i+1}^{*}$ (see lemma 8.1 case 1 of \cite{LapFRTF}). Assume that $n_i= n_{i+1}$ and $\sigma_i\simeq\sigma_{i+1}^{*}$.  The case where $\sigma_i\not\simeq \sigma_i^*$ is trivial ($J$ is zero) so we shall also assume that $\sigma_i\simeq\sigma_i^*$. Then  $M(s_\al,0)\varphi=-\varphi$ (\cite{KeysSha} proposition 6.3) and since $s_\al(\la)=\la$ we have $E(g,\varphi,\la)=0$ so \eqref{eq:theequality} is clear.
  \end{preuve}
\end{paragr}

\begin{paragr}[Majorization.] --- We will use the following proposition which results from Lapid's majorization of Eisenstein series (see \cite{LapFRTF} proposition 6.1 and section 7). For the convenience of the reader, we sketch a proof.

   \begin{proposition}\label{prop:GMSchwartz}
    Let  $f\in \Sc(G(\AAA))$. The map    $  \la\mapsto J_{Q,\pi}(g,\la,f)$    belongs to the Schwartz space $\Sc( i\ago_{P}^{G,*})$. Moreover 
\begin{align*}
  f\mapsto J_{Q,\pi}(g,\cdot,f)
\end{align*}
is a continuous map from $\Sc(G(\AAA))$ to $\Sc( i\ago_{P}^{G,*})$ equipped with its usual topology.
\end{proposition}

 \begin{preuve}   let $C\in G(\AAA_f)$ be a compact subset and $K_0\subset K^\infty$ be an open-compact subgroup  such that $f\in \Sc(G(\AAA),C,K_0)$. 
For any $\al,\be>0$, we define an open subset $\om_{\al,\be}$ of  $\ago_{P,\CC}^{G,*}$ which contains $ i\ago_{P,\RR}^{G,*}$ by
\begin{align*}
\om_{\al,\be}=\{\la\in \ago_{P,\CC}^{G,*}\mid \|\Re(\la)\|<\al(1+\|\Im(\la)\|)^{-\be}\}.
\end{align*}
By the arguments in the proof of proposition 6.1 of \cite{LapFRTF}, one sees that there exist  $\al,\be>0$ such that $\om_{\al,\be}$ is included in the open set $\om_\pi$ of §\ref{S:ompi}. In particular, $\la\mapsto J_{Q,\pi}(g,\la,f)$ is holomorphic on $\om_{\al,\be}$.
Using Cauchy formula to control derivatives, it suffices to prove the following majorization: there exists  a continuous semi-norm $\|\cdot\|$ on $\Sc(G(\AAA),C,K_0)$  and an open subset $\om_{\al,\be}\subset \om_{\pi}$ such that for any integer $N\geq 1$ there exists $c>0$ so that for all  $f\in \Sc(G(\AAA),C,K_0)$   and all $\la\in \om_{\al,\be}$ 
   \begin{align}\label{eq:major}
     |J_{Q,\pi}(g,\la,f)|\leq c \frac{\|f\| }{(1+\|\la\|)^N}.
   \end{align}
   Let $m\geq 1$ be a large enough integer. Following the notations of proposition \ref{prop:dirac}, we can write  $f=f*g_1+(f*Z)*g_2$ ; we get
   \begin{align*}
     J_{Q,\pi}(g,\la,f)=\sum_{\varphi\in \bc_{P,\pi}  } E(g, I_P(\la,f) \varphi,\la) \overline{J_Q(I_P(-\la,g_1^\vee)\varphi ,\la)}\\
+\sum_{\varphi\in \bc_{P,\pi}  } E(g, I_P(\la,f*Z) \varphi,\la) \overline{J_Q(I_P(-\la,g_2^\vee)\varphi ,\la)}
   \end{align*}

By a slight extension to Schwartz functions of Lapid's majorization (see \cite{FLO} remark C.2   about \cite{LapFRTF} proposition 6.1), the expression 
 \begin{align*}
    (\sum_{\varphi\in \bc_{P,\pi}  } |E(g, I_P(\la,f) \varphi,\la)|^2)^{1/2}
 \end{align*}
and the same expression where $f$ is replaced by $f*Z$ satisfy   a bound like \eqref{eq:major}. Using Cauchy-Schwartz inequality, we are reduced to bound in $\la$  (recall that $g_i$ is independent of $f$)
 \begin{align}\label{eq:car-gi}
    (\sum_{\varphi\in \bc_{P,\pi}  } |J_Q(I_P(-\la,g_i^\vee)\varphi ,\la)|^2)^{1/2}.
 \end{align}

Let $w$ be such that $wQw^{-1}$ is standard and $w\in W(P,wQw^{-1})$. At this point we will use the notations of the proof of proposition \ref{prop:car-relJB}There exists $c>0$ and an integer $r$  such that for all $\varphi'\in \Ac_{P,\pi}(G)^{K_0}$    we have
\begin{align*}
  |J_Q(I_P(-\la,g_i^\vee)\varphi ,\la)|=|J(M(w,\la) I_P(-\la,g_i^\vee)\varphi )| &\leq c \|M(w,\la) I_P(-\la,g_i^\vee)\varphi\|_{r}
\end{align*}
where $\|\varphi\|_r=\|R(1+C_K)^r\varphi\|_{\Pet}$.
Then we need to bound the operator norm of the intertwining operator $M(w,\la)$. Using the normalization of intertwining operators, the bounds of normalizing factors  \cite{LapFRTF} lemma 5.1 and  Müller-Speh's bound on the norm of normalized intertwining operators (see \cite{MuSpeh}  proposition 4.2 and the proof of proposition 0.2), we get $c_1>0$, $N\in \NN$ and $\al, \be >0$ such that  for all  $\tau\in  \hat K_\infty$, $\la\in \om_{\al,\be}$ and $\varphi\in \Ac_{P,\pi}(G,K_0,\tau)$ we have
\begin{align*}
  \|M(w,\la) I_P(-\la,g_i^\vee)\varphi\|_{K_0,r}\leq c_1(1+\la_\tau)^N \|I_P(-\la,g_i^\vee)\varphi\|_{r}.
\end{align*}

Using the same kind of arguments as in  the proof of proposition \ref{prop:car-relJB} (see also  remark \eqref{rq:classeCr}), one shows that there exist  $\al, \be >0$ such that   \eqref{eq:car-gi}  is bounded independently of $\la\in\om_{\al,\be}$.
\end{preuve}

\end{paragr}

\subsection{A spectral expansion of a truncated integral}\label{ssec:auxil}

\begin{paragr} Let $\chi\in \Xgo(G)$ be a cuspidal datum. We shall use the notation of §\ref{S:exp-kernel}. In particular, $f$ is a function in $f\in  \Sc(G(\AAA),C,K_0)$ and $K_\chi^0$ is the attached kernel.
\end{paragr}

\begin{paragr}
   Let's consider a parameter $T$ as in § \ref{S:trunc-param}. Following Jacquet-Lapid-Rogawski (see \cite{JLR}), we introduce the truncation operator $\Lambda^T_m$  that associates to a function $\varphi$ on $[G]$ the following function of the variable $h  \in [G']$:
\begin{align}\label{eq:LaTm}
    (\Lambda^T_m\varphi)(h)= \sum_{ P} (-1)^{\dim(\ago_P^G)}\sum_{\delta \in P'(F)\back G'(F) }   \hat\tau_P(H_P(\delta h)-T_P) \varphi_P(\delta h)
  \end{align}
where the sum is over standard parabolic subgroup of $G$ (those containing $B$) and $\varphi_P$ is the constant term along $P$. Recall that $P'=G'\cap P$.
\end{paragr}

\begin{paragr} We shall define the mixed truncated kernel $K_{\chi}^0\Lambda^T_m$: the notation means that  the mixed truncation is applied to the second variable. This is a function on $G(\AAA)\times G'(\AAA)$. To begin with we have:

  \begin{lemme}\label{lem:expansion1}
For $(x,y)\in G(\AAA)\times G'(\AAA)$, we have:
         \begin{align*}
      (K_{\chi}^0\Lambda^T_m)(x,y)=\sum_{B\subset P} |\pc(M_P)|^{-1}\int_{i\ago_P^{G,*}}\sum_{\varphi \in \bc_{P,\chi}} E(x,I_P(\la,f)\varphi,\la)\overline{\Lambda_m^TE(y,\varphi,\la)}\, d\la.
    \end{align*}
  \end{lemme}

  \begin{preuve}
As $y\in G'(\AAA)$, the mixed truncation is defined by a finite sum of constant terms of  $K_{\chi}^0(x,\cdot)$ (in the second variable). The only point is to permute the sum over $\varphi$ and the operator $\Lambda^T_m$. In fact using the continuity properties of Eisenstein series (see \cite{Lap-remark}) and properties of mixed truncation operator (in particular a variant of lemma 1.4 of \cite{ar2}), we can conclude as in the proof of proposition \ref{prop:car-relJB}. 
  \end{preuve}

  \begin{lemme}\label{lem:expansion2}
     For any integer $N$, there exists a continuous semi-norm $\|\cdot\|$ on $\Sc(G(\AAA),C,K_0)$ and an integer $N'$ such that for all $X\in \uc(\ggo_\CC)$, all $x\in G(\AAA)^1$ $y\in G'(\AAA)^1$ and all $f\in \Sc(G(\AAA),C,K_0)$ we have
 \begin{align}\label{eq:tobebd2}
     \sum_{\chi \in \Xgo(G)}\sum_{B\subset P} |\pc(M_P)|^{-1} \int_{i\ago_P^{G,*}}  \sum_{\tau\in \hat K} | \sum_{\varphi \in \bc_{P,\chi,\tau}} (R(X)E)(x,I_P(\la,f)\varphi,\la)\overline{\Lambda_m^T E(y,\varphi,\la)}\, d\la| \\
\nonumber \leq \|f\| \|x\|_{[G]}^{N'} \|y\|_{[G]}^{-N}.
    \end{align}
  \end{lemme}

  \begin{preuve}
    By the basic properties of the mixed truncation operator (see lemma 1.4 of \cite{ar2} and also \cite{LR} proof of lemma 8.2.1), for any $N$ and $N'$ there exists a finite family $(Y_i)_{i\in I}$ of elements of $\uc(\ggo_\CC)$  such  the expression \eqref{eq:tobebd2} is majorized by the sum over $i\in I$ of  $\|y\|_{[G]}^{-N}$ times the supremum over $g\in  G'(\AAA)^1$ of 
    \begin{align*} \|g\|_{[G]}^{-N'} \sum_{\chi  \in \Xgo}\sum_{B\subset P} |\pc(M_P)|^{-1} \int_{i\ago_P^{G,*}}  \sum_{\tau\in \hat K}  | \sum_{\varphi \in \bc_{P,\chi,\tau}} (R(X)E)(x,I_P(\la,f)\varphi,\la)\overline{R(Y_i)E(g,\varphi,\la)}\, d\la |.
    \end{align*}
Then the lemma is a straightforward consequence of lemma \ref{lem:ar2}.
  \end{preuve}

\begin{proposition}\label{prop:spec-exp}
    For all $x\in G(\AAA)$ and $\chi\in \Xgo(G)$,  we have
    \begin{align*}
      \int_{[G']_0}(K_{\chi}^0\Lambda^T_m)(x,y)\,dy =\sum_{B\subset P} |\pc(M_P)|^{-1}\int_{i\ago_P^{G,*}}\sum_{\varphi \in \bc_{P,\chi}} E(x,I_P(\la,f)\varphi,\la)\overline{\int_{[G']_0}\Lambda_m^TE(y,\varphi,\la)\,dy} \, d\la.
    \end{align*}
  \end{proposition}

  \begin{preuve}
    First one decomposes the sum over $\bc_{P,\chi}$ as a sum over $\tau\in \hat K$ of finite sums over  $\bc_{P,\chi,\tau}$. Then, by the majorization of lemma \ref{lem:expansion2} we can permute the integration over $[G']_0$ (which amounts to integrating over $[G']^1$) and the other sums or integrations in the expression we get in lemma \ref{lem:expansion1}.
  \end{preuve}
\end{paragr}

\subsection{The case of $*$-generic cuspidal data}

\begin{paragr}
  We shall use the notations of section \ref{ssec:auxil}.
\end{paragr}

\begin{paragr}[$*$-Generic cuspidal datum.] ---   \label{S:generic} We shall say  that a cuspidal datum $\chi\in \Xgo(G)$ is $*$-generic  if for any representative $(M,\pi)$ of $\chi$ and $w\in W(M)$ such that  $w\pi$ is isomorphic to $\pi$ or $\pi^*$ we have  $w=1$. Let's denote by $\Xgo^*(G)$ the subset of $*$-generic cuspidal data.

With the notations of \ref{S:Mpi}, we see that $(M,\pi)$ is $*$-generic if and only if for all $1\leq i,j\leq r$ such that $n_i=n_j$ one of the equalities  $\sigma_i=\sigma_j$ or  $\sigma_i=\sigma_j^*$  implies that $i=j$.  
\end{paragr}

\begin{paragr} The next theorem is the main result of the section.

  \begin{theoreme}\label{thm:spec-exp-FR-kernel}
Let $f\in \Sc(G(\AAA))$, let $\chi\in \Xgo(G)$ and let $K_\chi$ be the associated kernel. For any $g\in G(\AAA)$, one has:
\begin{enumerate}
\item We have 
  \begin{align}\label{eq:int-K0}
    \int_{[G']} K_\chi(g,h)\,dh= \frac12\int_{[G']_0} K_\chi^0(g,h)\,dh
  \end{align}
where both integrals are absolutely convergent.
\item If moreover $\chi\in \Xgo^*(G)$, we have, for any representative $(M_P,\pi)$ of $\chi$  (where $P$ is a standard parabolic subgroup of $G$),
\begin{align*}
    \int_{[G']} K_\chi(g,h)\,dh=2^{-\dim(\ago_P)} J_{P,\pi}(g,f)
  \end{align*}
where one defines (see \eqref{eq:JQpi-cusp})
    \begin{align*}
      J_{P,\pi}(g,f)=J_{P,\pi}(g,0,f).
    \end{align*}
In particular, the integral vanishes unless $\pi$ is self conjugate dual and $M_{P'}$-distinguished where $P'=G'\cap P$.
\end{enumerate}
\end{theoreme}

The assertion 1 follows readily from lemma \ref{lem:maj-noy}, Fubini's theorem and the fact that the Haar measure on $A_G^\infty$ is twice the Haar measure on $A_{G'}^\infty$ (see remark \ref{rq:mesure-AP}). The rest of the section is devoted to the proof of assertion 2 of theorem \ref{thm:spec-exp-FR-kernel}. The main steps are propositions \ref{prop:limit-form} and \ref{prop:calcul-limit}.
\end{paragr}

  \begin{paragr}[A limit formula.] ---
We shall use the notation  $\lim_{T\to +\infty} f(T)$ to denoted the limit of $f(T)$ when $\bg \al, T\bd\to+\infty$ for all $\al\in \Delta_{B}$.

    \begin{proposition}\label{prop:limit-form}
      Under the assumptions of theorem \ref{thm:spec-exp-FR-kernel} (but with no genericity condition on $\chi$), we have
        \begin{align*}
      \lim_{T\to +\infty}  \int_{[G']_0} (K_\chi^0\Lambda_m^T) (g,h)\,dh=2 \int_{[G']} K_\chi(g,h)\,dh.
        \end{align*}
      \end{proposition}

      \begin{preuve}
        Let's denote $F^{G'}(\cdot,T)$ the function defined by Arthur  relative to $G'$ and its maximal compact subgroup $K'$ (see \cite{ar1} §6 and \cite{ar-unip}lemma 2.1). It is the characteristic function of a compact of $[G']_0$. Using the fact that $h\mapsto  K_\chi^0(g,h)$ is of uniform moderate growth (see lemma \ref{lem:maj-noy}), we can conclude by a variant of \cite{ar-unip} theorem 3.1 (see also in the same spirit \cite{IY} proposition 3.8) that
        \begin{align*}
           \lim_{T\to +\infty}  \int_{[G']_0} \big(F^{G'}(h,T)K_\chi^0(g,h)-     (K_\chi^0\Lambda_m^T) (g,h)\big)\,dh=0.
        \end{align*}
     We have       $\lim_{T\to +\infty}  F^{G'}(h,T)=1$. Thus we deduce by Lebesgue's theorem and the absolute convergence of the right-hand side of \eqref{eq:int-K0}.
         \begin{align*}
           \lim_{T\to +\infty}  \int_{[G']_0} F^{G'}(h,T)K_\chi^0(g,h)\,dh= \int_{[G']_0} K_\chi^0(g,h)\,dh
         \end{align*}
         The proposition follows by \eqref{eq:int-K0}.
      \end{preuve}
    
  \end{paragr}
  
\begin{paragr}
  Let $\chi\in\Xgo^*(G)$. Let $\pc_\chi$ be the set of standard parabolic subgroups such that there exists a cuspidal automorphic representation $\pi$ of $M_P$ such that $(M_P,\pi)$ in the equivalence class defined by  $\chi$.
  Since $\chi\in \Xgo^*(G)$, the space $\Ac_{P,\chi}(G)$ is non-zero only if $P\in \pc_\chi$. Let $P$ be a standard parabolic subgroup and let $(M_P,\pi)$ be a pair in $\chi$. For any $P_1\in \pc_\chi$, by multiplicity-one theorem, we have
  \begin{align*}
    \Ac_{P_1,\chi}(G)=\bigoplus_{w\in W(P,P_1)} \Ac_{P,w\pi}.
  \end{align*}
In the following we set $M_1=M_{P_1}$.

 Let $P_1\in \pc_\chi$ and $g\in G(\AAA)$. With the notations of section \ref{ssec:related} (see eq. \eqref{eq:JQpi-cusp}),  for all $Q\in \pc(M_{1})$, all $\la\in i\ago_{P_1,\RR}^{G,*}$   we define
 \begin{align*}
   J_{Q,\chi}(g,\la,f)=\sum_{w\in W(P,P_1)} J_{Q,w\pi}(g,\la,f).
 \end{align*}
It's a continuous linear form on $\Sc(G(\AAA))$.
\end{paragr}

\begin{paragr}
  \begin{proposition}\label{prop:expanLmT}
For all $\chi\in \Xgo^*(G)$ and  all $g\in G(\AAA)$, we have
\begin{align}\label{eq:expanLmT}
  \int_{[G']_0} (K_\chi^0\Lambda_m^T) (g,h)\,dh= \frac{2^{-\dim(\ago_P^G)}}{ |\pc(M_{P})|} \sum_{P_1\in \pc_\chi}\int_{i\ago_{P_1}^{G,*}}  \sum_{Q\in \pc(M_{1})}  J_{Q,\chi}(g,\la,f) \frac{\exp(-\bg \la,T_Q\bd) }{\theta_Q(-\la) } \,d\la.
\end{align}
\end{proposition}

\begin{preuve}
This is an obvious  consequence of the definitions, the  proposition \ref{prop:spec-exp} and the lemma \ref{lem:Maass} below.
\end{preuve}

\begin{lemme}
  \label{lem:Maass}
 Let $\chi\in \Xgo^*(G)$. Let $P_1\in \pc_\chi$  and $\varphi\in \Ac_{P_1,\chi}$.  We have for all $\la\in i\ago_{P_1}^{G}$
\begin{align*}
  \int_{[G']_0}\Lambda_m^TE(y,\varphi,\la)\,dy= 2^{-\dim(\ago_{P_1}^G)}\sum_{Q\in \pc(M_{1})}  J_{Q}(\varphi,\la) \frac{\exp(\bg \la,T_Q\bd )}{\theta_Q(\la) }.
\end{align*}
\end{lemme}

\begin{preuve}
  This is simply a rephrasing in our particular situation of a key result of Jacquet-Lapid-Rogawski (see \cite{JLR} theorem 40). Indeed, because $\chi$ is $*$-generic, theorem 40 of \emph{ibid.} can be stated as:
  \begin{align*}
     \int_{[G']_0}\Lambda_m^TE(y,\varphi,\la)\,dy=  2^{-\dim(\ago_{P_1}^G)} \sum_{(Q,w)}  J(M(w,\la)\varphi) \frac{\exp(\bg (w\la)_Q ,T\bd )}{\theta_Q(w\la) }
  \end{align*}
where the sum is over pair $(Q,w)$ where $Q$ is a standard parabolic subgroup and $w\in W(P_1,Q)$.

\end{preuve}
\end{paragr}

\begin{paragr}
  \begin{proposition}\label{prop:calcul-limit}
    Let $\chi\in \Xgo^*(G)$ and let $(M_P,\pi)$ be a representative where $P$ is a standard parabolic subgroup of $G$. 
    We have:
    \begin{align*}
      \lim_{T\to +\infty}  \int_{[G']} (K_\chi^0\Lambda_m^T) (g,h)\,dh=2^{-\dim(\ago_P^G)} J_{P,\pi}(g,f)
    \end{align*}
    where one defines
    \begin{align}\label{eq:JPpi}
      J_{P,\pi}(g,f)=J_{P,\pi}(g,0,f).
    \end{align}
  \end{proposition}

  \begin{preuve}
    We start from the expansion \eqref{eq:expanLmT} of proposition \ref{prop:expanLmT}. For each $P_1\in \pc_\chi$, let $M_1=M_{P_1}$. The family $(J_{Q,\chi}(g,\la,f))_{Q\in \pc(M_1)}$ is a $(G,M_1)$-family of Schwartz functions on $i\ago_{P_1}^{G,*}$: this is a straightforward consequence of propositions \ref{prop:GMfam} and \ref{prop:GMSchwartz}. By  \cite{Linner} Lemma 8, we have:
    \begin{align*}
         \lim_{T\to +\infty}  \int_{i\ago_{P_1}^{G,*}}  \sum_{Q\in \pc(M_{1})}  J_{Q,\chi}(g,\la,f) \frac{\exp(-\bg \la,T_Q\bd) }{\theta_Q(-\la) } \,d\la= J_{P_1,\chi}(g,0,f)
    \end{align*}
    By definition and lemma \ref{lem:eqfctJP} below, one has:
    \begin{align*}
      J_{P_1,\chi}(g,0,f)&=\sum_{w\in W(P,P_1)} J_{P_1,w\pi}(g,0,f)\\
      &= |W(P,P_1)|  J_{P,\pi}(g,0,f).
    \end{align*}
Since $|\pc(M_P)|=\sum_{P_1\in \pc_\chi} W(P,P_1)$ we get the expected limit.
\end{preuve}

\begin{lemme}\label{lem:eqfctJP}(Lapid)
  For any $w\in W(P,P_1)$, we have
  \begin{align*}
    J_{P_1,w\pi}(g,0,f)=J_{P,\pi}(g,0,f).
  \end{align*}
\end{lemme}

\begin{preuve}
  By definition, we have
  \begin{align*}
      J_{P_1,w\pi}(g,0,f)=\sum_{\varphi\in \bc_{P_1,w\pi}} E(g,\varphi,0) \cdot J_1(\overline{\varphi})
  \end{align*}
  where $J_1$ is the linear form on $\Ac_{P_1,w\pi}$ defined in §\ref{S:J} and  $\bc_{P_1,w\pi}$ is any $K$-basis of $\Ac_{P_1,w\pi}$. Now, the intertwining operator $M(w,0)$ induces  a unitary isomorphism from $\Ac_{P,\pi}$ to $\Ac_{P_1,w\pi}$, which sends $K$-bases to $K$-bases. Thus one has
     \begin{align*}
       J_{P_1,w\pi}(g,0,f)&=\sum_{\varphi\in \bc_{P,\pi}} E(g,M(w,0)\varphi,0) \cdot J_1(M(w,0)\overline{\varphi})\\
       &=J_{P,\pi}(g,0,f).
  \end{align*}
The last equality results from the two equalities:
\begin{itemize}
\item $E(g,M(w,0)\varphi,0)=E(g,\varphi,0)$;
\item $J_1(M(w,0)\overline{\varphi})=J(\overline{\varphi})$  where $J$ is the linear form on $\Ac_{P,\pi}$ defined in §\ref{S:J}.
\end{itemize}
The first one is the functional equation of Eisenstein series and the second one  is a consequence of case 1 of lemma 8.1 of \cite{LapFRTF}.
\end{preuve}
\end{paragr}

\section{The $*$-generic contribution in the Jacquet-Rallis trace formula}\label{sec:generic contribution}

The goal of this chapter is to compute the contribution $I_\chi$ of the Jacquet-Rallis trace formula for $*$-generic cuspidal data $\chi$. This is achieved in theorem \ref{thm:intLaT} below. It turns out that for such $\chi$ the contribution $I_\chi$ is discrete and equal (up to an explicit constant) to a relative character define in section  \ref{ssec:relchar} built upon Rankin-Selberg periods of Eisenstein series and Flicker-Rallis intertwining periods.

\subsection{Relative characters}\label{ssec:relchar}

\begin{paragr}
  We will use the notations of section \ref{ssec:sym-not}.
\end{paragr}

\begin{paragr}
  Let $\chi\in \Xgo(G)$ be a cuspidal datum and $(M,\pi)$ be a representative where $M$ is the standard Levi factor of the standard prabolic subgroup $P$ of $G$. Recall that we have introduced a character $\eta_{G'}$ of $G'(\AAA)$ (see §\ref{S:character etaG'}).  On $\Ac_{P,\pi}(G)$, we introduce the linear form $J_\eta$ defined by
\begin{align}\label{S:Jeta}
  J_\eta(\varphi)=\int_{A_{M'} M'(F)N_{P'}(\AAA)\back G(\AAA)}\varphi(g)\eta_{G'}(g)\,dg  , \ \ \ \  \forall \varphi \in \Ac_{P,\pi}(G)
\end{align}
 where  $M'=M\cap G'$ and $P'=P\cap G'$. This is a slight variation of that defined in §\ref{S:J}. 

We shall say that $\pi$ is $(M',\eta_{G'})$-distinguished if $J_\eta$ does not vanish identically. 
\end{paragr}

\begin{paragr}[Relevant and generic cuspidal data.] --- \label{S:relevant} We shall say that $\chi$ is \emph{relevant} if $\pi$ is $(M',\eta_{G'})$-distinguished. 

Let $\Xgo^*(G)=\Xgo^*(G_n) \times \Xgo^*(G_{n+1})$ (cf. §\ref{S:generic}). We shall say that $\chi$ is $*$-generic if it belongs to the subset  $\Xgo^*(G)$. In particular, if $\chi$ is both relevant and  generic  (see §\ref{generic cuspidal data}) then it is $*$-generic.
  \end{paragr}

\begin{paragr}[Rankin-Selberg period of certain Eisenstein series.] --- Let $T\in \ago_{n+1}^+$. Recall that we have introduced in §\ref{S:IY-trunc} the truncation operator $\Lambda^T_r$.

\begin{proposition}\label{prop:RS-integ-generic}
Let $Q$ be  a parabolic subgroup of $G$ and $Q'=Q\cap G'$.  Let $\pi$ be an irreducible  cuspidal representation of $M_Q$ which is $(M_{Q'},\eta_{G'})$-distinguished. Let $\varphi\in \Ac_{Q,\pi}(G)$. Then for a regular point $\la\in \ago_Q^{G,*}$ of the Eisenstein series $E(g,\varphi,\la)$ (see §\ref{S:Eis-series}), the integral
  \begin{align}\label{eq:RSEis}
  I(\varphi,\la)=  \int_{[H]} \Lambda_r^T E(h,\varphi,\la)\,dh
  \end{align}
is convergent and does not depend on $T$. 
\end{proposition}

\begin{remarque}\label{rq:reg-period}
  The expression $I(\varphi,\la)$ is nothing else but  the regularized Rankin-Selberg period of $E(\varphi,\la)$ as defined by Ichino-Yamana in \cite{IY}. 
\end{remarque}

\begin{preuve}
The convergence follows from proposition \ref{prop:decayLaT} and the fact that Eiseinstein series are of moderate growth.  It remains to prove that the integral does not depend  on $T$.   Recall that $\iota$ induces an isomorphism from $G_n$ onto $H$. In the proof, it will be more convenient to work with $G_n$ instead of $H$. However, by abuse of notations, for any $g\in G_n(\AAA)$ and any function $\varphi$ on $G(\AAA)$ we shall write $\varphi(g)$ instead of $\varphi(\iota(g))$.

Let $T'\in \ago_{n+1}^+$.  By lemma 2.2 of \cite{IY}, we have
  \begin{align*}
    \Lambda_r^{T+T'} E(g,\varphi,\la)= \sum_{ P\in \fc_{RS}} \sum_{\delta \in (P\cap H)(F)\back H(F) }  \Lambda_r^{T,P} E_{G_n\times P_{n+1}}(\delta g,\varphi,\la) \Gamma_{P_{n+1}}'(H_{P_{n+1}}(\delta g)-T_{P_{n+1}},T')
  \end{align*}
where the notations  are  those of §§\ref{S:KTchi} and \ref{S:IY-trunc}. The other notations are borrowed from \cite[eq. (4.4)]{Z2}; the operator $\Lambda_r^{T,P}$ is the obvious variant of $\Lambda_r^{T}$ and  $\Gamma_{P}'$ is an Arthur function whose precise definition is irrelevant here. We denote by $E_{G_n\times P_{n+1}}$ the constant term of $E$ along $G_n\times P_{n+1}$. Thus, we have
 \begin{align*}
      \int_{[H]}\Lambda_r^{T+T'} E(g,\varphi,\la)\,dg=\\
 \sum_{ P\in \fc_{RS}} \int_{(P\cap H)(F)\back H(\AAA)} (\Lambda_r^{T,P} E_{G_n\times P_{n+1}}P)(g,\varphi,\la) \Gamma_{P_{n+1}}'(H_{P_{n+1}}(g)-T_{P_{n+1}},T')\,dg.
  \end{align*}
Let $P\in \fc_{RS}$ be such that $P\subsetneq G$. It suffices to show that the terms corresponding to $P$ vanish. We identify $H$ with $G_n$. Then $P\cap H$ is identified with $P_n$. Let $M_n=M_{P_n}$. For an appropriate  choice of  a Haar measure on $K_n$, such a term can be written as
\begin{align*}
  \int_{[M_{n}]} \int_{K_n}\exp(-\bg 2\rho_{P_n}, H_{P_n}(m)\bd) (\Lambda_r^{T,P} E_{P})(mk,\varphi,\la) \Gamma_{P_{n+1}}'(H_{P_{n+1}}(m)-T_{P_{n+1}},T')\,dk dm,
\end{align*}
where $E_{P}$ denotes the constant term of $E$ along $P=P_n\times P_{n+1}$. At this point, we may and shall assume that $P$ is standard (if not, we may change $B_n$ by a conjugate for the arguments). We have the usual formula for the constant term
\begin{align*}
  E_{P}(m,\varphi,\la)=\sum_{w\in W(Q;P)} E^{P}(m,M(w,\la)\varphi,w\la).
\end{align*}
where $W(Q;P)$ is the set of elements  $w\in W$ that are of minimal length in double cosets $W^{P}wW^Q$. Let $w\in W(Q;P)$. Notice that the representation $w\pi$ is also $(wM_{Q'}w^{-1},\eta_{G'})$-distinguished. For the argument, we may and shall assume $w=1$ (that is we assume that $Q\subset P$). Thus it suffices to show for all $k\in K$ the integral
\begin{align}\label{eq:truncMn}
  \int_{[M_{n}]^1}  \Lambda_r^{T,P} E^{P}(mk,\varphi,\la)\,dm
\end{align}
vanishes.

The group  $M_{n+1}=M_{P_{n+1}}$ has a decomposition $G_{d_1}\times\ldots\times G_{d_r}$ with $d_1+\ldots+d_r=n+1$. Each factor corresponds to a subset of the canonical basis $(e_1,\ldots,e_{n+1})$. We may assume that the factor $G_{d_1}$ corresponds to a subset which does not contain $e_{n+1}$. As a consequence $G_{d_1}$ is also a factor of $M_{n}$. We view $G_{d_1}\times G_{d_1}$ as a subgroup of $M_n\times M_{n+1}$. Let $Q_1\times Q_2=(G_{d_1}\times G_{d_1})\cap Q\subset G_n\times G_{n+1}$. The representation $\pi$ restricts to $M_{Q_1}(\AAA)$ and $M_{Q_2}(\AAA)$: this gives  representations respectively denoted by  $\pi_1$ and $\pi_2$. 
As a factor of \eqref{eq:truncMn}, we get
\begin{align} \label{eq:trunc-inner}
  \int_{[G_{d_1}]^1} E(g, \varphi_1,\la_1) \Lambda^T E(g,\varphi_2,\la_2)\, dg
\end{align}
where $\varphi_i\in \Ac_{Q_i,\pi_i}(G_{d_1})$. Here the truncation is the usual Arthur's truncation operator on the group $G_{d_1}$. It is clear from  Langlands' formula for the integral  \eqref{eq:trunc-inner} (see \cite{ar-truncated}) that \eqref{eq:trunc-inner} vanishes unless there exists $w\in W(Q_1,Q_2)$ such that $\pi_2\simeq w\pi_1$. But then $\pi_2$ would be $(M_{Q_2'}, (\eta_{d_1})^n )$-distinguished and $(M_{Q'_2}, (\eta_{d_1})^{n+1}) $-distinguished with  $\eta_{d_1}=\eta\circ \det_{d_1}$ and $M_{Q_2'}=M_{Q_2}\cap G_{d_1}'$. This is not possible.
\end{preuve}
\end{paragr}

\begin{paragr}[Relative characters.] ---   \label{S:RelcharJPpi}Let $(P,\pi)$ be a pair for which $P$ be a standard parabolic subgroup of $G$ and $\pi$ be a cuspidal automorphic representation of its standard Levi factor $M_P$. For any $\varphi\in \Ac_{P,\pi}(G)$, building upon the truncation operator $\Lambda^T_r$ and the linear form $J_\eta$, we define the relative character   $I_{P,\pi}^T$ for any $f\in \Sc(G(\AAA))$  by 
  \begin{align*}
    I_{P,\pi}^T(f)=\sum_{\varphi\in \bc_{P,\pi}}  \int_{[H]} \Lambda_r^T E(h,I_P(0,f)\varphi,0)\,dh \cdot \overline{J_\eta(\varphi)}
  \end{align*}
where the $K$-basis $\bc_{P,\pi}$ is defined in  §\ref{S:K0ON}.
 Using  proposition \ref{prop:RS-integ-generic}, we have 
\begin{align*}
 I_{P,\pi}^T(f)= I_{P,\pi}(f)
\end{align*}
where we define:
\begin{align*}
    I_{P,\pi}(f)=\left\lbrace
  \begin{array}{l}
    \sum_{\varphi\in \bc_{P,\pi}}  I(I_P(0,f)\varphi,0) \cdot \overline{J_\eta(\varphi)} \text{ if } \pi  \text{ is }(M_{P'},\eta_{G'})\text{-distinguished};\\
0 \text{ otherwise}.
  \end{array}
  \right.
\end{align*}

\begin{proposition} Let $\chi\in \Xgo^*(G)$. Let $(P,\pi)$ be a representative. 
  The map $f\mapsto I_{P,\pi}^T(f)$ (and thus $f\mapsto I_{P,\pi}(f)$)   is well-defined and gives a continuous linear form on $\Sc(G(\AAA))$. It depends only on $\chi$ and not on the choice of $( P,\pi)$.
\end{proposition}

\begin{preuve}
First we claim that $\varphi\mapsto \int_{[H]} \Lambda_r^T E(h,\varphi,0)\,dh$  is a continuous map: this is an easy consequence of properties of Eisenstein series and the truncation operator  $\Lambda_r^T$ (see the proposition \ref{prop:decayLaT}). On the other hand $\varphi\mapsto J_\eta(\varphi)$ is also continuous (see section \ref{ssec:related}). Thus the first assertion  results from an application of proposition \ref{prop:car-relJB}. The arguments of the proof of lemma \ref{lem:eqfctJP} give the independence on the choice of $(P,\pi)$.
 \end{preuve}
\end{paragr}

\subsection{The $*$-generic contribution}

\begin{paragr}  Let $\chi\in \Xgo(G)$. Recall that we defined in theorem \ref{thm:jfDef} a distribution $I_\chi$ on $\Sc(G(\AAA))$. Let $(M,\pi)$ be a representative of $\chi$ where $M$ is the standard Levi factor of the standard prabolic subgroup $P$ of $G$. The following theorem is the main result of this chapter.

\begin{theoreme} \label{thm:intLaT}  Assume moreover $\chi\in  \Xgo^*(G)$.  We have
  \begin{align*}
    I_\chi=2^{-\dim(\ago_P)} I_{P,\pi}.
  \end{align*}
In particular, we have $I_\chi=0$ unless $\chi$ is relevant.
\end{theoreme}

  The theorem is a direct  consequence of the following proposition.

\begin{proposition}\label{prop:intLaT}  
   Assume moreover $\chi\in  \Xgo^*(G)$   We have for $T\in \ago_{n+1}^+$
 \begin{align}\label{eq:intLaT}
      \int_{[H]} \int_{[G'] }   \Lambda_r^TK_\chi (x,y)\,\eta_{G'}(y)dxdy=2^{-\dim(\ago_P)}  I_{P,\pi}(f),
    \end{align}
where the left-hand side is absolutely convergent (see proposition \ref{prop:cv-LaTrKchi}).
In particular, the left-hand side does not depend on $T$.
\end{proposition}

Indeed, by theorem \ref{thm:asym-trio}, $I_\chi$ is the constant term in the asymptotic expansion in $T$ of the  left-hand side of \eqref{eq:intLaT} hence $I_\chi=2^{-\dim(\ago_P)}  I_{P,\pi}$.

The rest of the section is devoted to the proof of proposition \ref{prop:intLaT}.
\end{paragr}

\begin{paragr}[Proof of  proposition \ref{prop:intLaT}.] ---  We assume that $\chi\in \Xgo^*(G)$. The proof is a straightforward consequence of  theorem  \ref{thm:spec-exp-FR-kernel} and some permutations between  integrals, summations and  the truncation. These permutations are provided by lemmas \ref{lem:interchange1} and \ref{lem:interchange2} below.

  \begin{lemme}
    \label{lem:interchange1}
For all $x\in [H]$, we have
\begin{align*}
  \int_{[G'] }   (\Lambda_r^TK_\chi) (x,y)\eta_{G'}(y)\,dxdy=   \Lambda_r^T\left(\int_{[G'] } K_\chi (\cdot,y)\eta_{G'}(y)\,dy\right)(x).
\end{align*}
  \end{lemme}

  \begin{remarque}
On the left-hand side we apply the truncation operator     $\Lambda_r^T$ to the function $K_\chi(\cdot,y)$ (where $y$ is fixed) and then we evaluate at $x$ whereas on the right-hand side  we apply the same operator to the function we get by integration of  $K_\chi (\cdot,y)\eta_{G'}(y)$ over $y\in[G']$ and then we evaluate at $x$.
  \end{remarque}
  
\begin{preuve}
    Since $x$ is fixed, the operator  $\Lambda_r^T$ is a finite sum of constant terms (see \cite{ar1} lemma 5.1 for the finiteness). Then the lemma follows from Fubini's theorem which holds because we have 
    \begin{align*}
      \int_{[N_Q]}\int_{[G'] } |K_\chi (nx,y)| \,dndy<\infty
    \end{align*}
for all parabolic subgroups $Q$ of $G_{n+1}$ containing $B_{n}$. Here we identify $N_Q$ with the subgroup $\{1\}\times N_Q$ of $G=G_n\times G_{n+1}$. The convergence of the integral results from the bound \eqref{eq:KchiG} above.
  \end{preuve}

  \begin{lemme}
    \label{lem:interchange2}
We have
\begin{align*}
   \int_{[H]} \Lambda_r^T\left(\int_{[G'] } K_\chi (\cdot,y)\eta_{G'}(y)\,dy\right)(h)\,dh=2^{-\dim(\ago_P)}I_{P,\pi}(f).
\end{align*}
  \end{lemme}

  \begin{preuve}First, by theorem \ref{thm:spec-exp-FR-kernel}, we have for any $x\in [G]$:
    \begin{align*}
      \int_{[G'] } K_\chi (x,y)\eta_{G'}(y)\,dy=2^{-\dim(\ago_P)}\sum_{\varphi\in \bc_{P,\pi}} E(x,I_P(0,f)\varphi,0)\,dh \cdot \overline{J_\eta(\varphi)}
    \end{align*}
where the notations are borrowed from §\ref{S:RelcharJPpi}.
Then we want to apply the truncation operator $\Lambda_r^T$ and evaluate at $h\in [H]$. We want to show that this operation commutes with the summation over the orthonormal basis. As in the proof of lemma \ref{lem:interchange1}, it suffices to prove  
\begin{align*}
  \sum_{\varphi\in \bc_{P,\pi}}   \int_{[N_Q]} |E(n g,I_P(0,f)\varphi,0)|\,dn \cdot |\overline{J_\eta(\varphi)}| <\infty
\end{align*}
for any parabolic subgroups $Q$ of $G_{n+1}$ containing $B_{n}$,  which is is an easy consequence of continuity properties of Eisenstein series.

In this way, we get for $h\in [H]$:
\begin{align*}
    \Lambda_r^T\left(\int_{[G'] } K_\chi (\cdot,y)\eta_{G'}(y)\,dy\right)(h)=2^{-\dim(\ago_P)} \sum_{\varphi\in \bc_{P,\pi}}  (\Lambda_r^T E)(h,I_P(0,f)\varphi,0) \cdot \overline{J_\eta(\varphi)}.
\end{align*}
By integration over  $h\in [H]$, we have:
\begin{align*}
    \int_{[H]}\Lambda_r^T\left(\int_{[G'] } K_\chi (\cdot,y)\eta_{G'}(y)\,dy\right)(h)\,dh= 2^{-\dim(\ago_P)} \sum_{\varphi\in \bc_{P,\pi}}  \int_{[H]}(\Lambda_r^T E)(h,I_P(0,f)\varphi,0) \,dh\cdot \overline{J_\eta(\varphi)}.
\end{align*}
The right-hand side is nothing else but $2^{-\dim(\ago_P)} I_{P,\pi}(f)$. Still we have to justify the change of order of the integration and the summation. But it is easy to to show that
\begin{align*}
   \sum_{\varphi\in \bc_{P,\pi}}   \int_{[H] } |\Lambda_r^T E(h,I_P(0,f)\varphi,0)|\,dh \cdot |\overline{J_\eta(\varphi)}|<\infty.
\end{align*}
\end{preuve}
\end{paragr}

\section{Spectral decomposition of the Flicker-Rallis period for certain cuspidal data}\label{chapter FR}

The goal of this chapter is to give another proof of the spectral decomposition of the Flicker-Rallis period for the same cuspidal data as in Section \ref{S:generic}. The main result of this chapter (obtained as a combination of Theorem \ref{theo Zeta integrals FR} and Theorem \ref{theo FR period}) can be used to get another version of Theorem \ref{thm:spec-exp-FR-kernel} with a seemingly different relative character than $J_{P,\pi}$ (this will actually be done in \S \ref{paragr spectral exp K2}). Of course, these two relative characters are the same. A direct proof of this fact will be given in Chapter \ref{Chap:FR-functional-computation}.

\subsection{Notation}

\begin{paragr}
In this chapter we adopt the set of notation introduced in Section \ref{ssec:sym-not}: $E/F$ is a quadratic extension of number fields, $G'_n=\GL_{n,F}$, $G_n=\Res_{E/F}\GL_{n,E}$, $(B'_n,T'_n)$, $(B_n,T_n)$ are the standard Borel pairs of $G'_n$, $G_n$ and $K'_n$, $K_n$ the standard maximal compact subgroups of $G'_n(\AAA)$, $G_n(\AAA)$ respectively. Besides, we denote by $N'_n$, $N_n$ the unipotent radicals of $B'_n$, $B_n$ and we set
$$w_n=\begin{pmatrix} & & 1 \\ & \iddots & \\ 1 & & \end{pmatrix}\in G'_n(F).$$
We write $e_n=(0,\ldots,0,1)$ for the last element in the standard basis of $F^n$ and we let $\pc_n=\begin{pmatrix} \star & \star \\ 0\ldots 0 & 1 \end{pmatrix}$, $\pc'_n=\pc_n\cap G'_n$ be the mirabolic subgroups of $G_n$, $G'_n$ respectively (that is the stabilizers of $e_n$ for the natural right actions). The unipotent radicals of $\pc_n$, $\pc'_n$ will be denoted by $U_n$ and $U'_n$ respectively. For nonnegative integers $m\leqslant n$, we embed $G_m$ in $G_n$ (resp. $G'_m$ in $G'_n$) in the ``upper left corner'' by $g\mapsto \begin{pmatrix} g & \\ & I_{n-m}\end{pmatrix}$. Thus, in particular, we have $\pc_n=G_{n-1}U_n$ and $\pc'_n=G'_{n-1}U'_n$.

The entries of a matrix $g\in G_n(\AAA)$ are written as $g_{i,j}$, $1\leqslant i,j\leqslant n$, and the diagonal entries of an element $t\in T_n(\AAA)$ as $t_i$, $1\leqslant i\leqslant n$.
\end{paragr}

\begin{paragr}\label{ssect:additive and generic characters}
We fix a nontrivial additive character $\psi': \AAA/F\to \CC^\times$. For $\phi\in \Sc(\AAA^n)$, we define its Fourier transform $\widehat{\phi}\in \Sc(\AAA^n)$ by
$$\displaystyle \widehat{\phi}(x_1,\ldots,x_n)=\int_{\AAA^n} \phi(y_1,\ldots,y_n) \psi'(x_1y_1+\ldots+x_ny_n)dy_1\ldots dy_n$$
the Haar measure on $\AAA^n$ being chosen such that $\widehat{\widehat{\phi}}(x)=\phi(-x)$.

We denote by $c$ the nontrivial Galois involution of $E$ over $F$. Then, $c$ acts naturally on $G_n(\AAA)$ and thus on cuspidal automorphic representations of the latter. We denote this action by $\pi\mapsto \pi^c$. We fix $\tau\in E^\times$ such that $\tau^c=-\tau$ and we define $\psi:\AAA_E/E\to \CC^\times$ by $\psi(z)=\psi'(\Tra_{E/F}(\tau z))$, $z\in \AAA_E$, where $\AAA_E$ denotes the ad\`ele ring of $E$ and $\Tra_{E/F}:\AAA_E\to \AAA$ the trace map. We also define a generic character $\psi_n:[N_n]\to \CC^\times$ by
$$\displaystyle \psi_n(u)=\psi\left((-1)^n\sum_{i=1}^{n-1} u_{i,i+1}\right),\;\; u\in [N_n].$$
(The appearance of the sign $(-1)^n$ is only a convention that will be justified a posteriori in Chapter \ref{chap canonical RS}). Note that $\psi$ is trivial on $\AAA$ and therefore $\psi_n$ is trivial on $N'_n(\AAA)$. To any $f\in \tc([G_n])$, we associate its {\em Whittaker function} $W_f$ defined by
$$\displaystyle W_f(g)=\int_{[N_n]} f(ug)\psi_n(u)^{-1}du,\;\; g\in G_n(\AAA).$$
\end{paragr}

\subsection{Statements of the main results}\label{section main results FR period}

\begin{paragr}
Let $n\geqslant 1$ be a nonnegative integer. For $f\in \tc([G_n])$, $\phi\in \Sc(\AAA^n)$ and $s\in \CC$ we set
$$\displaystyle Z_\psi^{\FR}(s,f,\phi)=\int_{N_n'(\AAA)\backslash G_n'(\AAA)} W_f(h) \phi(e_nh)\lvert \det h\rvert^s dh$$
provided this expression converges absolutely.
\end{paragr}

\begin{paragr}
Let $\chi\in \Xgo^*(G_n)$ be a $*$-generic cuspidal datum (see \S \ref{S:generic} for the definition of $*$-generic) represented by a pair $(M_P,\pi)$ and set $\Pi=I_{P(\AAA)}^{G_n(\AAA)}(\pi)$. We can write
$$\displaystyle M_P=G_{n_1}\times \ldots \times G_{n_k}$$
where $n_1,\ldots,n_k$ are positive integers such that $n_1+\ldots+n_k=n$. Then, $\pi$ decomposes accordingly as a tensor product
$$\displaystyle \pi=\pi_1\boxtimes\ldots\boxtimes \pi_k$$
where for each $1\leqslant i\leqslant k$, $\pi_i$ is a cuspidal automorphic representation of $G_{n_i}(\AAA)$.
\end{paragr}

\begin{paragr}\label{S:distinguished cuspidal datum}
Let $L(s,\Pi,\As)$ be the Shahidi's completed Asai $L$-function of $\Pi$ \cite{Sha90}, \cite{Gold}. We have the decomposition
$$\displaystyle L(s,\Pi,\As)=\prod_{i=1}^k L(s,\pi_i,\As)\times \prod_{1\leqslant i<j\leqslant k} L(s,\pi_i\times \pi_j^c).$$
As $\chi$ is $*$-generic, the Rankin-Selberg $L$-functions $L(s,\pi_i\times \pi_j^c)$ are entire and non-vanishing at $s=1$ \cite{JS}, \cite{JS2}, \cite{ShaLfunction} whereas by \cite{Flicker}, $L(s,\pi_i,\As)$ has at most a simple pole at $s=1$. Therefore, $L(s,\Pi,\As)$ has a pole of order at most $k$ at $s=1$ and this happens if and only if $L(s,\pi_i,\As)$ has a pole at $s=1$ for every $1\leqslant i\leqslant k$.

We say that the cuspidal datum $\chi$ is {\em distinguished} if $L(s,\Pi,\As)$ has a pole of order $k$ at $s=1$. By \cite{Flicker}, it is equivalent to ask $\pi$ to be $M_{P'}=M_P\cap G'_n$-distinguished.
\end{paragr}

\begin{paragr}\label{S:betan}
For $f\in \cC([G_n])$, we set $W_{f,\Pi}=W_{f_\Pi}$ where $f_\Pi$ is defined as in Section \ref{section PWLapid}. Then, $W_{f,\Pi}$ belongs to the Whittaker model $\wc(\Pi,\psi_n)$ of $\Pi$ with respect to $\psi_n$.

We define a continuous linear form $\beta_n$ on $\wc(\Pi,\psi_n)$ as follows. For $S$ a finite set of places of $F$ and $W\in\wc(\Pi,\psi_n)$, we set
$$\displaystyle \beta_{n,S}(W)=\int_{N_n'(F_S)\backslash P_n'(F_S)} W(p_S)dp_S$$
the integral being convergent by (the same proof as) \cite[Proposition 2.6.1, Lemma 3.3.1]{BPAsai} and the Jacquet-Shalika bound \cite{JS}. By \cite[Proposition 3]{Flicker} and \eqref{eq Tamagawa measures}, for a given $W\in \wc(\Pi,\psi_n)$, the quantity
$$\displaystyle \beta_n(W)=(\Delta^{S,*}_{G'})^{-1} L^{S,*}(1,\Pi,\As) \beta_{n,S}(W)$$
is independent of $S$ as long as it is sufficiently large (i.e. it contains all the Archimedean places as well as the non-Archimedean places where the situation is ``ramified''). This defines the linear form $\beta_n$.
\end{paragr}

\begin{paragr}
For every $f\in \cC([G_n])$, we set
$$\displaystyle {}^0 f(g)=\int_{A_{G_n}^\infty} f(ag) da,\;\; g\in [G_n].$$

\begin{theoreme}\label{theo Zeta integrals FR}
	\begin{enumerate}
		\item Let $N\geqslant 0$. There exists $c_N>0$ such that for every $f\in \tc_N([G_n])$ and $\phi\in \Sc(\AAA^n)$, the expression defining $Z_\psi^{\FR}(s,f,\phi)$ is absolutely convergent for $s\in \cH_{>c_N}$ and the function $s\in \cH_{>c_N}\mapsto Z_\psi^{\FR}(s,f,\phi)$ is holomorphic and bounded in vertical strips. Moreover, for every $s\in \cH_{>c_N}$, $(f,\phi)\mapsto Z_\psi^{\FR}(s,f,\phi)$ is a (separately) continuous bilinear form on $\tc_N([G_n])\times \Sc(\AAA^n)$.
		
		\item Let $\chi\in \Xgo^*(G_n)$. For every $f\in \cC_\chi([G_n])$, the function $s\mapsto (s-1)Z_\psi^{\FR}(s,{}^0 f,\phi)$ admits an analytic continuation to $\cH_{>1}$ with a limit at $s=1$. Moreover, we have
		$$\displaystyle Z_\psi^{\FR,*}(1,{}^0f,\phi):=\lim\limits_{s\to 1^+} (s-1)Z_\psi^{\FR}(s,{}^0f,\phi)=\left\{\begin{array}{ll}
		2^{1-k}\widehat{\phi}(0) \beta_n(W_{f,\Pi}) \mbox{ if } \chi \mbox{ is distinguished}, \\ \\
		0 \mbox{ otherwise.}
		\end{array}\right.$$
	\end{enumerate}
\end{theoreme}
\end{paragr}

\begin{paragr}
\begin{theoreme}\label{theo FR period}
	Let $\chi\in \Xgo^*(G_n)$. The linear form
	$$\displaystyle P_{G_n'}: f\in \cC([G_n])\mapsto \int_{[G'_n]} f(h) dh$$
	is well-defined (i.e. the integral converges) and continuous. Moreover, for every $f\in \Sc_\chi([G_n])$ and $\phi\in \Sc(\AAA^n)$ such that $\widehat{\phi}(0)=1$ we have
\begin{equation}\label{eq FR period Zeta integral}
\displaystyle P_{G'_n}(f)=\frac{1}{2}Z_\psi^{\FR,*}(1,{}^0f,\phi).
\end{equation}
\end{theoreme}
\end{paragr}

\begin{paragr}
A direct consequence of Theorem	\ref{theo Zeta integrals FR} and Theorem \ref{theo FR period} is the following corollary.

\begin{corollaire}
Let $\chi\in \Xgo^*(G_n)$ be represented by a pair $(M_P,\pi)$ and set $\Pi=I_{P(\AAA)}^{G_n(\AAA)}(\pi)$. Then, for every $f\in \Sc_\chi([G_n])$ we have
$$\displaystyle P_{G'_n}(f)=\left\{\begin{array}{ll}
	2^{-\dim(A_P)} \beta_n(W_{f,\Pi}) \mbox{ if } \chi \mbox{ is distinguished}, \\ \\
	0 \mbox{ otherwise.}
\end{array}\right.$$
\end{corollaire}
\end{paragr}

\subsection{Proof of Theorem \ref{theo Zeta integrals FR}.2}

Part 1. of Theorem \ref{theo Zeta integrals FR} will be established in Section \ref{section majorations Fourier}. Here, we give the proof of part 2. of this theorem. Let $f\in \cC_\chi([G_n])$, $\phi\in \Sc(\AAA^n)$ and $(M_P,\pi)$ be a pair representing the cuspidal datum $\chi$ as in Section \ref{section main results FR period}. Set 
$$\displaystyle \Ac:= (i\RR)^k$$
and let $\Ac_0$ be the subspace of $\underline{x}=(x_1,\ldots,x_k)\in \Ac$ such that $x_1+\ldots+x_k=0$. We equip $\Ac$ with the product of Lebesgue measures and $\Ac_0$ with the unique Haar measure such that the quotient measure on
$$\displaystyle \Ac/\Ac_0\simeq i\RR,\; (x_1,\ldots,x_k)\mapsto x_1+\ldots+x_k$$
is again the Lebesgue measure.

There is an unique identification $\Ac\simeq i\ago_{P}^*$ which when it is composed with the map $\mu\in i\ago_{P}^*\mapsto \pi_\mu$ gives
\begin{equation}\label{isom A ap*}
\displaystyle \underline{x}=(x_1,\ldots,x_k)\in \Ac\mapsto \pi_{\underline{x}}:=\pi_1 \lvert \det \rvert_{E}^{x_1/n_1}\boxtimes \ldots\boxtimes \pi_k \lvert \det \rvert_{E}^{x_k/n_k}.
\end{equation}
For every $\underline{x}\in \Ac$, we set $\Pi_{\underline{x}}=I_{P(\AAA)}^{G_n(\AAA)}(\pi_{\underline{x}})$ and $f_{\underline{x}}=f_{\Pi_{\underline{x}}}$ following the definition of Section \ref{section PWLapid} (so that in particular $\Pi_0=\Pi$ and $f_0=f_\Pi$ with notation from the previous section).

The isomorphism \eqref{isom A ap*} sends $\prod_{j=1}^k (in_j \ZZ)$ onto $iX^*(P)$ hence, by \eqref{eq measure iaP*}, it also sends the measure on $\Ac$ to $(n_1\ldots n_k)(2\pi)^k$ times the measure on $i\ago_P^*$. Therefore, by Theorem \ref{theo Lapid}, we have
$$\displaystyle {}^0f=\frac{(2\pi)^{-k}}{n_1\ldots n_k} \int_{A_{G_n}^\infty}\int_{\Ac} a \cdot f_{\underline{x}} d\underline{x}da=\frac{(2\pi)^{-k}}{n_1\ldots n_k} \int_{A_{G_n}^\infty}\int_{\Ac} \lvert \det a\rvert_E^{\frac{x_1+\ldots+x_k}{n}} f_{\underline{x}} d\underline{x}da.$$
We have an isomorphism $A_{G_n}^\infty\simeq \RR_+^*$, $a\mapsto \lvert \det a\rvert_E$, sending the Haar measure on $A_{G_n}^\infty$ to $\frac{dt}{\lvert t\rvert}$ where $dt$ is the Lebesgue measure. Thus, by Fourier inversion, the previous equality can be rewritten as
\begin{align}
\displaystyle {}^0 f=\frac{n}{n_1\ldots n_k} (2\pi)^{-k+1}\int_{\Ac_0} f_{\underline{x}} d\underline{x}
\end{align}
where the right-hand side is an absolutely convergent integral in $\tc_N([G_n])$ for some $N>0$. Therefore, by the first part of Theorem \ref{theo Zeta integrals FR}, there exists $c>0$ such that for every $s\in \cH_{>c}$ we have
\begin{equation}\label{eq0- proof of theo zeta integrals FR}
\displaystyle Z_\psi^{\FR}(s,{}^0f,\phi)=\frac{n}{n_1\ldots n_k} (2\pi)^{-k+1}\int_{\Ac_0} Z_\psi^{\FR}(s,f_{\underline{x}},\phi) d\underline{x}.
\end{equation}
Let $S_0$ be a finite set of places of $F$ including the Archimedean ones and outside of which $\pi$ is unramified and let $S_{0,f}\subset S_0$ be the subset of finite places. Let $I\subseteq \{ 1,\ldots, k\}$ be the subset of $1\leqslant i\leqslant k$ such that $L(s,\pi_i,\As)$ has a pole at $s=1$. We choose, for each $1\leqslant i\leqslant k$ and $v\in S_{0,f}$, polynomials $Q_i(T), Q_{i,v}(T)\in \CC[T]$ with roots in $\cH_{]0,1[}$ and $\cH_{]q_v^{-1},1[}$ respectively such that $s\mapsto Q_i(s) L_\infty(s,\pi_i,\As)$ and $s\mapsto Q_{i,v}(q_v^{-s})L_v(s,\pi_i,\As)$ have no pole in $\cH_{]0,1[}$. Finally, we set
$$\displaystyle P(s,\underline{x})=\prod_{i\in I} (s+\frac{2x_i}{n_i})(s-1+\frac{2x_i}{n_i})\prod_{1\leqslant i\leqslant k} Q_i(s+\frac{2x_i}{n_i}) \prod_{\substack{1\leqslant i\leqslant k \\ v\in S_{0,f}}} Q_{i,v}(q_v^{-s-\frac{2x_i}{n_i}}) \; \mbox{ and } \; \widetilde{f_{\underline{x}}}(g)=f_{\underline{x}}({}^tg^{-1})$$
for every $\underline{x}\in \Ac_0$, $s\in \CC$ and $g\in G_n(\AAA)$. We will now check that the functions
\begin{equation}\label{eq0 proof of theo zeta integrals FR}
(s,\underline{x})\in \CC\times \Ac_0\mapsto P(s+\frac{1}{2},\underline{x})Z_\psi^{\FR}(s+\frac{1}{2},f_{\underline{x}},\phi)
\end{equation}
and
\begin{equation}\label{eq0bis proof of theo zeta integrals FR}
(s,\underline{x})\in \CC\times \Ac_0\mapsto P(\frac{1}{2}-s,\underline{x})Z_{\psi^{-1}}^{\FR}(s+\frac{1}{2},\widetilde{f}_{\underline{x}},\widehat{\phi})
\end{equation}
satisfy the conditions of Corollary \ref{cor1 PL}.

From the first part of Theorem \ref{theo Zeta integrals FR}, Theorem \ref{theo Lapid} and Lemma \ref{lem holomorphic maps to Schwartz spaces}, we deduce that these functions satisfy the first condition of Corollary \ref{cor1 PL}. To check that they also satisfy the second condition of Corollary \ref{cor1 PL}, we need to analyze more carefully the function $s\mapsto Z_\psi^{\FR}(s,f_{\underline{x}},\phi)$ for a fixed $\underline{x}\in \Ac_0$.

For $S$ a sufficiently large finite set of places of $F$, that we assume to contain Archimedean places as well as the places where $\pi$, $\psi'$ or $\psi$ are ramified (thus $S_0\subset S$), we have decompositions
$$\displaystyle \phi=\phi_S\phi^S \;\mbox{ and }\; W_{f_{\underline{x}}}=W_{S,\underline{x}}W_{\underline{x}}^S$$
for every $\underline{x}\in \Ac_0$, where $\phi_S\in \Sc(F_S^n)$, $\phi^S$ is the characteristic function of $(\widehat{\mathcal{O}_F^S})^n$, $W_{S,\underline{x}}\in \wc(\Pi_{\underline{x},S},\psi_{n,S})$ (that is the Whittaker model of the representation $\Pi_{\underline{x},S}$ with respect to the character $\psi_{n,S}=\psi_{n\mid N(F_S)}$) and $W_{\underline{x}}^S\in \wc(\Pi_{\underline{x}}^S,\psi_n^S)^{K_n^S}$ is such that $W_{\underline{x}}^S(1)=1$. By \cite[Proposition 3]{Flicker} and \eqref{eq decomposition Tamagawa measure}, we then have
\begin{equation}\label{eq1 proof of theo zeta integrals FR}
\displaystyle Z_\psi^{\FR}(s,f_{\underline{x}},\phi)=(\Delta_{G_n'}^{S,*})^{-1}L(s,\Pi_{\underline{x}},\As) \frac{Z_\psi^{\FR}(s,W_{S,\underline{x}},\phi_S)}{L_S(s,\Pi_{\underline{x}},\As)}
\end{equation}
for $s\in \cH_{>c}$ where we have set
$$\displaystyle Z_\psi^{\FR}(s,W_{S,\underline{x}},\phi_S)=\int_{N_n'(F_S)\backslash G_n'(F_S)} W_{S,\underline{x}}(h_S) \phi_S(e_nh_S)\lvert \det h_S\rvert^s dh_S.$$
Moreover, by \cite[Theorem 3.5.1]{BPAsai} the function $Z_\psi^{\FR}(s,W_{S,\underline{x}},\phi_S)$ extends meromorphically to the complex plane and satisfies the functional equation
\begin{equation}\label{eq2 proof of theo zeta integrals FR}
\displaystyle \frac{Z^{\FR}_{\psi^{-1}}(1-s,\widetilde{W}_{S,\underline{x}},\widehat{\phi_S})}{L_S(1-s,(\Pi_{\underline{x}})^\vee, \As)}=\epsilon(s,\Pi_{\underline{x}},\As)\frac{Z^{\FR}_\psi(s,W_{S,\underline{x}},\phi_S)}{L_S(s,\Pi_{\underline{x}},\As)}
\end{equation}
where $\widetilde{W_{S,\underline{x}}}(g)=W_{S,\underline{x}}(w_n{}^tg^{-1})$, $\widehat{\phi_S}$ is the (normalized) Fourier transform of $\phi_S$ with respect to the bicharacter $(\underline{u},\underline{v})\mapsto \psi'(u_1v_1+\ldots+u_nv_n)$ and $\epsilon(s,\Pi_{\underline{x}},\As)$ denotes the global epsilon factor of the Asai $L$-function $L(s,\Pi_{\underline{x}},\As)$.

By \eqref{eq1 proof of theo zeta integrals FR}, \eqref{eq2 proof of theo zeta integrals FR} as well as the meromorphic continuation and functional equation of $L(s,\Pi_{\underline{x}},\As)$ \cite[Theorem 3.5(4)]{Sha90}, we conclude that $Z_\psi^{\FR}(s,f_{\underline{x}},\phi)$ has a meromorphic continuation to $\CC$ satisfying the functional equation
\begin{equation}\label{eq3 proof of theo zeta integrals FR}
\displaystyle Z^{\FR}_{\psi^{-1}}(1-s,\widetilde{f_{\underline{x}}},\widehat{\phi})=Z^{\FR}_\psi(s,f_{\underline{x}},\phi).
\end{equation}
On the other hand, we have the decomposition
$$\displaystyle L(s,\Pi_{\underline{x}},\As)=\prod_{i=1}^k L(s+\frac{2x_i}{n_i},\pi_i,\As)\times \prod_{1\leqslant i<j\leqslant k} L(s+\frac{x_i}{n_i}+\frac{x_j}{n_j},\pi_i\times \pi_j^c).$$
and, as $\chi\in \Xgo^*(G_n)$, the Rankin-Selberg $L$-functions $L(s,\pi_i\times \pi_j^c)$ are entire and bounded in vertical strips \cite[Theorem 4.1]{CogNotes}. By the Jacquet-Shalika bound \cite{JS} and the fact that the gamma function is of exponential decay in vertical strips, $Q_i(s)L_\infty(s,\pi_i,\As)$ and $Q_{i,v}(q_v^{-s})L_v(s,\pi_i,\As)$ are holomorphic and bounded in vertical strips of $\cH_{>0}$ for each $1\leqslant i\leqslant k$ and $v\in S_{0,f}$. By \cite[Lemma 5.2]{FLIntertwining}, $s\mapsto (s-1)L^{S_0}(s,\pi_i,\As)$, for $i\in I$, and $s\mapsto L^{S_0}(s,\pi_i,\As)$, for $i\notin I$, are also holomorphic and of finite order in vertical strips of $\cH_{>0}$. Therefore, by the definition of $P$ and the functional equation, $P(s,\underline{x})L(s,\Pi_{\underline{x}},\As)$ is entire and of finite order in vertical strips. By \eqref{eq1 proof of theo zeta integrals FR}, \eqref{eq3 proof of theo zeta integrals FR} and \cite[Theorem 3.5.2]{BPAsai}, it follows that the functions \eqref{eq0 proof of theo zeta integrals FR}, \eqref{eq0bis proof of theo zeta integrals FR} are entire and of finite order in vertical strips in the first variable i.e. they also satisfy the second condition of Corollary \ref{cor1 PL}.

Thus, the conclusion of this corollary is valid and in particular the map
$$\displaystyle s\mapsto\left(\underline{x}\mapsto \prod_{i\in I} (s-1+\frac{2 x_i}{n_i})Z_\psi^{\FR}(s,f_{\underline{x}},\phi) \right)$$
induces a holomorphic function $\cH_{>1-\epsilon}\to \Sc(\Ac_0)$ for some $\epsilon>0$. By \eqref{eq0- proof of theo zeta integrals FR} and \cite[Lemma 3.1.1, Proposition 3.1.2]{BPPlanch}, it follows that $s\mapsto Z_\psi^{\FR}(s,{}^0f,\phi)$ extends analytically to $\cH_{>1}$ and that
\begin{align}\label{eq4 proof of theo zeta integrals FR}
\displaystyle \lim\limits_{s\to 1^+} (s-1)Z_\psi^{\FR}(s,{}^0f,\phi) & =\left\{\begin{array}{ll} 2^{1-k} \lim\limits_{s\to 1} (s-1)^k Z_\psi^{\FR}(s,f_0,\phi) \mbox{ if } I=\{1,\ldots, k\}, \\
0 \mbox{ otherwise}
\end{array}\right.
\end{align}

Recall that $I=\{1,\ldots, k\}$ if and only if $L(s,\Pi,\As)$ has a pole of order $k=rk(A_P)$ at $s=1$. Moreover, by \cite[Lemma 3.3.1]{BPAsai} and the Jacquet-Shalika bound \cite{JS}, the integral defining $Z_\psi^{\FR}(s,W_{S,0},\phi_S)$ is absolutely convergent in $\cH_{>1-\varepsilon}$ for some $\varepsilon>0$. Combining this with \cite[Lemma 2.16.3]{BPPlanch} and \eqref{eq1 proof of theo zeta integrals FR}, in the case $I=\{1,\ldots, k\}$ identity \eqref{eq4 proof of theo zeta integrals FR} can be rewritten as
\[\begin{aligned}
\displaystyle \lim\limits_{s\to 1^+} (s-1)Z_\psi^{\FR}(s,{}^0f,\phi) & =2^{1-k} (\Delta^{S,*}_{G_n'})^{-1} L^{S,*}(1,\Pi,\As)Z_\psi^{\FR}(1,W_{S,0},\phi_S) \\
& =2^{1-k} (\Delta^{S,*}_{G_n'})^{-1} L^{S,*}(1,\Pi,\As)\beta_{n,S}(W_{S,0})\widehat{\phi_S}(0) \\
& =2^{1-k}\widehat{\phi}(0) \beta_n(W_{f,\Pi})
\end{aligned}\]
and this ends the proof of Theorem \ref{theo Zeta integrals FR}.2.

\subsection{Proof of Theorem \ref{theo FR period}}\label{section preuve theo FR period}

By \eqref{eq1 Function spaces}, we have $\Xi^{[G_n]}(h)\ll \Xi^{[G_n']}(h)^2$ for $h\in [G_n']$. Hence, by \eqref{eq 0 Function spaces}, the linear form $P_{G_n'}$ is well-defined and continuous on $\cC([G_n])$. This shows the first part of Theorem \ref{theo FR period}.

Let $f\in\Sc_\chi([G_n])$. Recall that $A_{G_n}^\infty=A_{G'_n}^\infty$ but the Haar measure on $A_{G_n}^\infty$ is twice the Haar measure on $A_{G'_n}^\infty$ (see Remark \ref{rq:mesure-AP}). Therefore, we have
$$\displaystyle P_{G_n'}(f)=\frac{1}{2}\int_{[G'_n]_0} {}^0f(h) dh.$$
Let $\phi\in \Sc(\AAA^n)$. We form the Epstein-Eisenstein series
$$\displaystyle E(h,\phi,s)=\int_{A_{G'_n}^\infty} \sum_{\gamma\in \pc_n'(F)\backslash G_n'(F)} \phi(e_n\gamma ah) \lvert \det(ah)\rvert^s da,\;\;\; h\in [G_n'], s\in \CC.$$
This expression converges absolutely for $\Re(s)>1$ and the map $s\mapsto E(\phi,s)$ extends to a meromorphic function valued in $\tc([G_n'])$ with simple poles at $s=0$, $1$ of respective residues $\phi(0)$ and $\widehat{\phi}(0)$ (cf. \cite[Lemma 4.2]{JS}).

Consequently, the function
$$\displaystyle s\mapsto Z_n^{\FR}(s,{}^0f,\phi):=\int_{[G_n']_0} {}^0f(h) E(h,\phi,s) dh$$
is well-defined for $s\in \CC\setminus\{0,1 \}$, meromorphic on $\CC$ with a simple pole at $s=1$ whose residue is
\begin{equation}\label{eq0 preuve thm FR}
\displaystyle \Res_{s=1}Z_n^{\FR}(s,{}^0f,\phi)=2\widehat{\phi}(0) P_{G_n'}(f).
\end{equation}

Unfolding the definition, formally we arrive at
\begin{align}\label{eq1 preuve thm FR}
\displaystyle Z_n^{\FR}(s,{}^0f,\phi)=\int_{\pc_n'(F)\backslash G_n'(\AAA)} {}^0f(h) \phi(e_nh) \lvert \det h\rvert^s dh.
\end{align}
By Lemma \ref{lem conv 2} below, there exists $c_n>0$ such that the last integral above is absolutely convergent for $s\in \cH_{>c_n}$ and thus the equality above is justified for such $s$.

More generally, for every $1\leqslant r\leqslant n$, let $N_{r,n}$ be the unipotent radical of the standard parabolic subgroup of $G_n$ with Levi component $G_r\times (G_1)^{n-r}$, $N'_{r,n}$ be its intersection with $G'_n$ and set
$$\displaystyle {}^0f_{N_{r,n},\psi}(g)=\int_{[N_{r,n}]} {}^0f(ug)\psi_n(u)^{-1}du,\;\;\; g\in G(\AAA),$$
$$\displaystyle Z^{\FR}_r(s,{}^0f,\phi)=\int_{\pc'_r(F)N'_{r,n}(\AAA)\backslash G'_n(\AAA)} {}^0f_{N_{r,n},\psi}(h) \phi(e_nh) \lvert \det h\rvert^s dh,\;\; s\in \CC,$$
provided the last expression above is convergent. The proof of the next lemma will be given in Section \ref{section majorations Fourier}.

\begin{lemme}\label{lem conv 2}
	For every $1\leqslant r\leqslant n$, there exists $c_r>0$ such that the expression defining $Z^{\FR}_r(s,{}^0f,\phi)$ converges absolutely for $\Re(s)>c_r$.
\end{lemme}

When $r=1$, we have $N_{1,n}=N_n$ and ${}^0f_{N_{1,n},\psi}=W_{{}^0f}$ so that $Z^{\FR}_1(s,{}^0f,\phi)=Z_\psi^{\FR}(s,{}^0f,\phi)$. Therefore, by \eqref{eq0 preuve thm FR} and \eqref{eq1 preuve thm FR}, the second part of Theorem \ref{theo FR period} is a consequence of the following proposition.

\begin{proposition}
	For every $1\leqslant r\leqslant n$, the function $s\mapsto (s-1)Z^{\FR}_r(s,{}^0f,\phi)$ extends to a holomorphic function on $\left\{s\in \CC\mid \Re(s)>1\right\}$ admitting a limit at $s=1$. Moreover, we have
	$$\displaystyle \lim\limits_{s\to 1^+} (s-1)Z_n^{\FR}(s,{}^0f,\phi)=\lim\limits_{s\to 1^+} (s-1)Z_r^{\FR}(s,{}^0f,\phi).$$
\end{proposition}

\begin{preuve}
By descending induction on $r$, it suffices to establish the following:

\begin{num}
	\item\label{eq2 preuve thm FR} Let $1\leqslant r\leqslant n-1$. There exists a function $F_r$ holomorphic on $\cH_{>1-\epsilon}$ for some $\epsilon>0$ such that
	$$\displaystyle Z_{r+1}^{\FR}(s,{}^0f,\phi)=Z_r^{\FR}(s,{}^0f,\phi)+F_r(s)$$
	for all $s\in \CC$ satisfying $\Re(s)>\max(c_r,c_{r+1})$.
\end{num}

Indeed, as $\pc'_{r+1}=G'_rU'_r$, we have
\begin{align}\label{eq3 preuve thm FR}
\displaystyle Z^{\FR}_{r+1}(s,{}^0f,\phi)=\int_{G'_r(F)N'_{r,n}(\AAA)\backslash G_n'(\AAA)} \int_{[U'_{r+1}]} {}^0f_{N_{r+1,n},\psi}(uh)du \phi(e_nh) \lvert \det h\rvert^s dh.
\end{align}
By Fourier inversion on the locally compact abelian group $U_{r+1}(F)U'_{r+1}(\AAA)\backslash U_{r+1}(\AAA)$, we have
\begin{equation}\label{eq4 preuve thm FR}
\displaystyle \int_{[U'_{r+1}]} {}^0f_{N_{r+1,n},\psi}(uh)du=\sum_{\gamma \in \pc'_r(F)\backslash G'_r(F)} ({}^0f_{N_{r+1,n},\psi})_{U_{r+1},\psi}(\gamma h)+({}^0f_{N_{r+1,n},\psi})_{U_{r+1}}(h)
\end{equation}
for all $h\in G'_n(\AAA)$ where we have set
$$\displaystyle  ({}^0f_{N_{r+1,n},\psi})_{U_{r+1},\psi}(h)=\int_{[U_{r+1}]} {}^0f_{N_{r+1,n},\psi}(uh)\psi_n(u)^{-1} du={}^0f_{N_{r,n},\psi}(h),$$
$$\displaystyle  ({}^0f_{N_{r+1,n},\psi})_{U_{r+1}}(h)=\int_{[U_{r+1}]} {}^0f_{N_{r+1,n},\psi}(uh) du.$$

By \eqref{eq3 preuve thm FR} and \eqref{eq4 preuve thm FR}, we obtain
$$\displaystyle Z_{r+1}^{\FR}(s,{}^0f,\phi)=Z_r^{\FR}(s,{}^0f,\phi)+F_r(s)$$
for all $s\in \CC$ such that $\Re(s)>\max(c_r,c_{r+1})$ and where we have set
$$\displaystyle F_r(s)=\int_{G'_r(F)N'_{r,n}(\AAA)\backslash G_n'(\AAA)} ({}^0f_{N_{r+1,n},\psi})_{U_{r+1}}(h) \phi(e_nh) \lvert \det h\rvert^s dh.$$
It only remains to check that $F_r(s)$ extends to a holomorphic function on $\cH_{>1-\epsilon}$ for some $\epsilon>0$.

Let $P_r$ be the standard parabolic subgroup of $G_n$ with Levi component $M_r=G_r\times G_{n-r}$ and set $P_r'=P_r\cap G'_n$. We readily check that
$$\displaystyle ({}^0f_{N_{r+1,n},\psi})_{U_{r+1}}(h)=\int_{[N_{n-r}]} {}^0f_{P_r}(\begin{pmatrix} I_r & \\ & u\end{pmatrix} h)\psi_n(u)^{-1} du=\int_{[N_{n-r}]} \int_{A_{G_n}^\infty} f_{P_r}(\begin{pmatrix} I_r & \\ & u\end{pmatrix} ah)da\psi_n(u)^{-1} du.$$
Therefore, by the Iwasawa decomposition $G'_n(\AAA)=P_r'(\AAA)K'_n$ and since
$$\displaystyle \delta_{P_r'}\begin{pmatrix} h_r & \\ & h_{n-r} \end{pmatrix}=\delta_{P_r}\begin{pmatrix} h_r & \\ & h_{n-r} \end{pmatrix}^{1/2}=\lvert \det h_r\rvert^{n-r} \lvert \det h_{n-r}\rvert^{-r}$$
for all $h_r\in G'_r(\AAA)$, $h_{n-r}\in G'_{n-r}(\AAA)$, we have (for $\Re(s)>\max(c_r,c_{r+1})$ and a suitable choice of Haar measure on $K'_n$)
\begin{align}\label{eq5 preuve thm FR}
\displaystyle F_r(s)=\int_{K'_n\times [G'_r]\times N'_{n-r}(\AAA)\backslash G'_{n-r}(\AAA)} &  \int_{[N_{n-r}]} \int_{A_{G_n}^\infty} f_{P_r,k,s}\left(a\begin{pmatrix} h_r & \\ & uh_{n-r} \end{pmatrix}\right) da \\
\nonumber &  \psi_n(u)^{-1} du \lvert \det h_{n-r}\rvert^{ns/(n-r)} \phi_{k,n-r}(e_{n-r}h_{n-r})dh_{n-r}dh_rdk
\end{align}
where $f_{P_r,k,s}=\delta_{P_r}^{-1/2+s/2(n-r)} (R(k)f)_{P_r}\mid_{M_r(\AAA)}$ and $\phi_{k,n-r}$ stands for the composition of $R(k)\phi$ with the inclusion $\AAA^{n-r}\to \AAA^n$, $x\mapsto (0,x)$. Let $\chi^M$ be the inverse image of $\chi$ in $\Xgo(M_r)$. By Corollary \ref{cor constant term generic cuspidal datum}, we have $f_{P_r,k,s}\in \cC_{\chi^M}([M_r])$ for every $(k,s)\in K_n\times \cH_{>0}$ and the map
$$\displaystyle (k,s)\in K_n\times \cH_{>0}\mapsto f_{P_r,k,s}\in \cC_{\chi^M}([M_r])$$
is continuous, holomorphic in the second variable. In particular, for $\Re(s)>0$ the integral
$$\displaystyle \int_{A_{G'_r}^\infty} \int_{[N_{n-r}]} \int_{A_{G_n}^\infty} f_{P_r,k,s}\left(a\begin{pmatrix} a'h_r & \\ & uh_{n-r} \end{pmatrix}\right) da da'$$
is absolutely convergent and equals, by the obvious change of variable, to
$$\displaystyle \frac{n}{n-r}\int_{A_{G'_r}^\infty} \int_{[N_{n-r}]} \int_{A_{G_{n-r}}^\infty} f_{P_r,k,s}\begin{pmatrix} a'h_r & \\ & uah_{n-r} \end{pmatrix} da da'.$$
It follows that \eqref{eq5 preuve thm FR} can be rewritten, for $\Re(s)\gg 1$, as
\begin{align}\label{eq6 preuve thm FR}
\nonumber \displaystyle F_r(s) & =\frac{n}{n-r}\int_{K'_n} \int_{[G'_r]\times N'_{n-r}(\AAA)\backslash G'_{n-r}(\AAA)} \int_{[N_{n-r}]} \int_{A^\infty_{G_r}}f_{P_r,k,s}\begin{pmatrix} h_r & \\ & uah_{n-r} \end{pmatrix} da\psi_n(u)^{-1}du \\
\nonumber & \lvert \det h_{n-r}\rvert^{ns/(n-r)} \phi_{k,n-r}(e_{n-r}h_{n-r})dh_{n-r}dh_rdk \\
& =\frac{n}{n-r}\int_{K'_n} (P_{G'_r}\widehat{\otimes} \zc_{n-r}^{\FR}(\frac{ns}{n-r}))(f_{P_r,k,s}\otimes \phi_{k,n-r})dk
\end{align}
where $\zc_{n-r}^{\FR}(s)$ stands for the bilinear form
$$\displaystyle (f',\phi')\in \cC([G_{n-r}])\times \Sc(\AAA^{n-r})\mapsto Z_{\psi^{(-1)^r}}^{\FR}(s,{}^0 f',\phi').$$
On the other hand, by \eqref{eq isom HCS product with cuspidal datum} we have
$$\displaystyle \cC_{\chi^M}([M_r])=\bigoplus_{(\chi_1,\chi_2)\in \Xgo(M_r)=\Xgo(G_r)\times \Xgo(G_{n-r})\mapsto \chi} \cC_{\chi_1}([G_r])\widehat{\otimes} \cC_{\chi_2}([G_{n-r}])$$
and, as $\chi\in \Xgo^*(G_n)$, for every $(\chi_1,\chi_2)\in \Xgo(G_r)\times \Xgo(G_{n-r})$ mapping to $\chi\in \Xgo(G_n)$, we also have $\chi_2\in \Xgo^*(G_{n-r})$. Therefore, by the first part of Theorem \ref{theo FR period}, Theorem \ref{theo Zeta integrals FR} and \eqref{eq 4 tvs}, $s\mapsto P_{G'_r}\widehat{\otimes}\zc_{n-r}^{\FR}(s)$ extends to an analytic family of (separately) continuous bilinear forms on $\cC_\chi([M_r])\times \Sc(\AAA^{n-r})$ for $s\in \cH_{>1}$. Thus, by the first part of Theorem \ref{theo FR period}, \eqref{eq 3 tvs} and the equality \eqref{eq6 preuve thm FR}, $F_r(s)$ has an analytic continuation to $\{\Re(s)>1-r/n \}$. This ends the proof of the proposition and hence of Theorem \ref{theo FR period}.
\end{preuve}

\subsection{Convergence of Zeta integrals}\label{section majorations Fourier}

\begin{paragr}
\begin{preuve} (of Lemma \ref{lem conv 2}) We only treat the case $1\leqslant r\leqslant n-1$. The case $r=n$ can be dealt with in a similar manner, and is in fact easier.
	
Let $Q_r$ be the standard parabolic subgroup of $G_n$ with Levi component $G_r \times G_1^{n-r}$ and set $Q'_r=Q_r\cap G_n'$. Recall that $N_{r,n}$ is the unipotent radical of $Q_r$. Identifying $A_{G_n}^\infty\simeq \RR_{>0}$, by the Iwasawa decomposition $G'_n(\AAA)=Q_r'(\AAA)K'_n$, we need to show the convergence of
\begin{align}\label{eq2 majorations Fourier}
\displaystyle & \int_{K_n'\times \pc'_r(F)\backslash G'_r(\AAA)\times T'_{n-r}(\AAA)\times \RR_{>0}} \left\lvert(R(k)f)_{N_{r,n},\psi}\begin{pmatrix} ah & \\ & at \end{pmatrix}\right\rvert \left\lvert R(k)\phi(t_{n-r}e_n)\right\rvert \\
\nonumber & \lvert \det h\rvert^{s}\lvert \det t\rvert^{s} \delta_{Q'_r}\begin{pmatrix} h & \\ & t\end{pmatrix}^{-1} da dt dh dk
\end{align}
for $\Re(s)\gg 1$. We now apply Lemma \ref{lem1 majorations Fourier}. For this we note that $\psi_n \mid_{[N_{r,n}]}=\psi'\circ \ell$ where $\ell:N_{r,n}\to \GG_a$ sends $u\in N_{r,n}$ to $\Tra_{E/F}(\tau\sum_{i=r}^{n-1}u_{i,i+1})$ and $\tau\in E^\times$ is the unique trace-zero element such that $\psi(z)=\psi'(\Tra_{E/F}(\tau z))$. We readily check that
$$\displaystyle \lVert \Ad^*(m)\ell\rVert_{(N_{r,n})_{\ab}^*(\AAA)}\approx \lVert t_1^{-1}e_rh\rVert_{\AAA^r} \prod_{i=1}^{n-r-1} \lVert t_it_{i+1}^{-1}\rVert_{\AAA},\;\; m=\begin{pmatrix} h & \\ & t \end{pmatrix}\in G'_r(\AAA)\times T'_{n-r}(\AAA).$$
Therefore, by Lemma \ref{lem1 majorations Fourier}.1, we can find $c>0$ such that for every $N_1,N_2>0$ we have
\begin{equation}\label{eq3 majorations Fourier}
\displaystyle \left\lvert (R(k)f)_{N_{r,n},\psi}\begin{pmatrix} ah & \\ & at \end{pmatrix}\right\rvert\ll \lVert ah\rVert_{[G'_r]}^{-N_2}\lVert t_1^{-1}e_rh\rVert_{\AAA^r}^{-N_1} \prod_{i=1}^{n-r-1} \lVert t_it_{i+1}^{-1}\rVert_{\AAA}^{-N_1} \delta_{Q_r}\begin{pmatrix} h & \\ & t \end{pmatrix}^{-cN_2}
\end{equation}
for $(k,h,t,a)\in K'_n\times G'_r(\AAA)\times T'_{n-r}(\AAA)\times \RR_{>0}$. On the other hand, for every $N_1>0$, we have
$$\displaystyle \lvert R(k)\phi(te_n)\rvert\ll \lVert t\rVert_{\AAA}^{-N_1},\;\;\; (k,t)\in K'_n\times \AAA$$
and it is easy to check that for some $N_2>0$ we have
$$\displaystyle \lVert e_rh\rVert_{\AAA^r} \prod_{i=1}^{n-r} \lVert t_i\rVert_{\AAA} \ll \lVert t_{n-r}\rVert_{\AAA}^{N_2}\lVert t_1^{-1}e_rh\rVert_{\AAA^r}^{N_2} \prod_{i=1}^{n-r-1} \lVert t_it_{i+1}^{-1}\rVert_{\AAA}^{N_2},\;\;\; (h,t)\in G'_r(\AAA)\times T'_{n-r}(\AAA).$$
As $\displaystyle \delta_{Q'_r}\begin{pmatrix} h & \\ & t \end{pmatrix}=\lvert \det h\rvert^{n-r}\prod_{i=1}^{n-r} \lvert t_i\rvert^{n+1-2(r+i)}$ for every $(h,t)\in G'_r(\AAA)\times T'_{n-r}(\AAA)$, combining this with \eqref{eq3 majorations Fourier}, we deduce the existence of $c>0$ such that for every $N_1,N_2>0$, \eqref{eq2 majorations Fourier} is essentially bounded by the product of
\begin{equation}\label{eq7 majorations Fourier}
\displaystyle \int_{\pc'_r(F)\backslash G'_r(\AAA)\times \RR_{>0}} \lVert ah\rVert_{[G'_r]}^{-N_2} \lVert e_rh\rVert_{\AAA^r}^{-N_1} \lvert \det h\rvert^{s-(2cN_2+1)(n-r)} dadh
\end{equation}
and
\begin{equation}\label{eq8 majorations Fourier}
\displaystyle \int_{\AAA^\times} \lVert t\rVert_{\AAA}^{-N_1} \lvert t\rvert^{s-(2cN_2+1)(n+1-2(r+i))} dt
\end{equation}
for $1\leqslant i\leqslant n-r$.

Let $C_1,C_2>0$. By Lemma \ref{lem conv GLn}, for $N_1$ sufficiently large the integral \eqref{eq8 majorations Fourier} converges absolutely in the range
$$\displaystyle 1+(2cN_2+1)(n+1-2(r+i))<\Re(s)<C_1+(2cN_2+1)(n+1-2(r+i))$$
and for $N_1$, $N_2$ sufficiently large the integral \eqref{eq7 majorations Fourier} converges absolutely in the range
$$\displaystyle 1+(2cN_2+1)(n-r)<\Re(s)<C_2+(2cN_2+1)(n-r).$$
Since $n+1-2(r+i)<n-r$ for every $1\leqslant i\leqslant n-r$, by taking $C_2=2$ and $C_1\geqslant 2+(2cN_2+1)(r+2i-1)$ for every $1\leqslant i\leqslant n-r$, it follows that if $N_2\gg 1$ and $N_1\gg_{N_2} 1$ the integrals \eqref{eq7 majorations Fourier} and \eqref{eq8 majorations Fourier} are convergent in the range
$$\displaystyle 1+(2cN_2+1)(n-r)<\Re(s)<2+(2cN_2+1)(n-r).$$
The union of these open intervals for $N_2$ sufficiently large as above is of the form $]c_r,+\infty[$ which shows that $Z_r^{\FR}(s,{}^0f,\phi)$ converges absolutely in the range $\Re(s)>c_r$ for a suitable $c_r>0$.
\end{preuve}
\end{paragr}

\begin{paragr}
\begin{preuve} (of Theorem \ref{theo Zeta integrals FR}.1) Applying Lemma \ref{lem1 majorations Fourier}.2, the same manipulations as in the proof of Lemma \ref{lem conv 2} reduce us to showing the existence of $c_N>0$ such that for every $C>c_N$ there exists $N'>0$ satisfying that the integral
\begin{equation}\label{eq9 majorations Fourier}
\displaystyle \int_{T'_n(\AAA)} \prod_{i=1}^n \lVert t_i\rVert_{\AAA}^{-N'}  \lVert t\rVert_{[T'_n]}^N \delta_{B'_n}(t)^{-1} \lvert \det t\rvert^s dt
\end{equation}
converges in the range $s\in \cH_{]c_n,C[}$ uniformly on compact subsets. But this follows again from Lemma \ref{lem conv GLn} as there exists $M>0$ such that
$$\displaystyle \lVert t\rVert_{[T'_n]}^N \delta_{B'_n}(t)^{-1}\ll \prod_{1 \leqslant i\leqslant n} \max(\lvert t_i\rvert, \lvert t_i\rvert^{-1})^{M},\;\; t\in [T'_n].$$
\end{preuve}
\end{paragr}

\section{Canonical extension of the Rankin-Selberg period for certain cuspidal data}\label{chap canonical RS}

This chapter is a continuation of Chapter \ref{chapter FR} and we shall use the notation introduced there. The main goal is to show the existence of a canonical extension of corank one Rankin-Selberg periods to the space of uniform moderate growth functions for certain cuspidal data (see Theorem \ref{theo RS period}). Combining this with the results of Chapter \ref{chapter FR}, this will enable us to give an alternative proof of the spectral expansion of the Jacquet-Rallis trace formula for certain cuspidal data in Chapter \ref{Chap:generic-JRTF}.

\subsection{Statements of the main results}\label{section main results RS period}

\begin{paragr}\label{ssect:RS Zeta}
Let $n\geqslant 1$ be a positive integer. We set $G=G_n\times G_{n+1}$ and $H=G_n$ that we consider as an algebraic subgroup of $G$ via the diagonal inclusion $H\hookrightarrow G$. We also set $w=(w_n,w_{n+1})\in G(F)$, $K=K_n\times K_{n+1}$, $N=N_n\times N_{n+1}$ and $N_H=N_n$. Put $\psi_N=\psi_n\boxtimes \psi_{n+1}$ (a generic character of $[N]$). We note that $\psi_N$ is trivial on $[N_H]$ (see the convention in the definition of $\psi_n$ in \S \ref{ssect:additive and generic characters}). To any function $f\in \tc([G])$, we associate its {\em Whittaker function}
$$\displaystyle W_f(g)=\int_{[N]}f(ug)\psi_N(u)^{-1}du,\; g\in G(\AAA).$$
For $f\in \tc([G])$, we define
$$\displaystyle Z_\psi^{\RS}(s,f)=\int_{N_H(\AAA)\backslash H(\AAA)} W_f(h) \lvert \det h\rvert_E^s dh$$
for every $s\in \CC$ for which the above expression converges absolutely.
\end{paragr}

\begin{paragr}[$H$-generic cuspidal datum]\label{S:Hgeneric}
Let $\chi\in \Xgo(G)$ be a cuspidal datum represented by a pair $(M_P,\pi)$ where $P=P_n\times P_{n+1}$ is a standard parabolic subgroup of $G$ and $\pi=\pi_n\boxtimes \pi_{n+1}$ a cuspidal automorphic representation of $M_P(\AAA)$ (with central character trivial on $A_P^\infty$). We have decompositions
$$\displaystyle M_{P_n}=G_{n_1}\times \ldots \times G_{n_k},\;\; M_{P_{n+1}}=G_{m_1}\times \ldots \times G_{m_r}$$
and $\pi_n$, $\pi_{n+1}$ decompose accordingly as tensor products
$$\displaystyle \pi_n=\pi_{n,1}\boxtimes \ldots\boxtimes \pi_{n,k},\;\; \pi_{n+1}=\pi_{n+1,1}\boxtimes \ldots \boxtimes \pi_{n+1,r}.$$
We will say that $\chi$ is {\em $H$-generic} if it satisfies the following condition:
\begin{num}
\item\label{hyp1 RS} For every $1\leqslant i\leqslant k$ and $1\leqslant j\leqslant r$, we have $\pi_{n,i}\neq \pi_{n+1,j}$ or equivalently the Rankin-Selberg $L$-function $L(s,\pi_{n,i}\times \pi_{n+1,j})$ is entire.
\end{num}
\end{paragr}

\begin{paragr}
\begin{theoreme}\label{theo Zeta integral RS}
\begin{enumerate}
\item Let $N\geqslant 0$. There exists $c_N>0$ such that for $f\in \tc_N([G])$, the expression defining $Z_\psi^{\RS}(s,f)$ converges absolutely for $s\in \cH_{>c_N}$ and the map $s\in \cH_{>c_N}\mapsto Z_\psi^{\RS}(s,f)$ is holomorphic and bounded in vertical strips. Moreover, for every $s\in \cH_{>c_N}$, $f\mapsto Z_\psi^{\RS}(s,f)$ is a continuous functional on $\tc_N([G])$.

\item Assume that $\chi\in \Xgo(G)$ is a $H$-generic cuspidal datum. Then, for every $f\in \tc_\chi([G])$, the function $s\mapsto Z_\psi^{\RS}(s,f)$ extends analytically to $\CC$. Moreover, for every $s\in \CC$ the linear form $f\in \tc_\chi([G])\mapsto Z_\psi^{\RS}(s,f)$ is continuous.
\end{enumerate}
\end{theoreme}
\end{paragr}

\begin{paragr}
\begin{theoreme}\label{theo RS period}
Assume that $\chi\in \Xgo(G)$ is a $H$-generic cuspidal datum. The restriction of the linear form
$$\displaystyle P_H: f\in \Sc([G])\mapsto \int_{[H]} f(h) dh$$
to $\Sc_\chi([G])$ extends by continuity to $\tc_\chi([G])$ and for every $f\in \tc_\chi([G])$, we have
\begin{equation}\label{eq1 theo RS period}
\displaystyle P_H(f)=Z^{\RS}_\psi(0,f).
\end{equation}
\end{theoreme}
\end{paragr}

\subsection{Proof of Theorem \ref{theo Zeta integral RS}.2}

Part 1. of Theorem \ref{theo Zeta integral RS} will be established in Section \ref{Section convergence Zeta RS}. Here, we give the proof of part 2. Let $f\in \Sc_\chi([G])$ and $(M_P,\pi)$ be a cuspidal datum representating $\chi$ as in Section \ref{section main results RS period}. We set $\Pi_{n,\mu}=I_{P_n(\AAA)}^{G_n(\AAA)}(\pi_{n,\lambda})$ and $\Pi_{n+1,\nu}=I_{P_{n+1}(\AAA)}^{G_{n+1}(\AAA)}(\pi_{n+1,\nu})$ for every $\mu \in \ago_{P_n,\CC}^*$ and $\nu \in \ago_{P_{n+1},\CC}^*$. For $\lambda=(\mu,\nu)\in \ago_{P,\CC}^*$, we also set $\Pi_\lambda= \Pi_{n,\mu}\boxtimes \Pi_{n+1,\nu}$ and $f_\lambda=f_{\Pi_\lambda}$ (following the notation from Section \ref{section PWLapid}). Then, by Theorem \ref{theo Lapid} and the first part of Theorem \ref{theo Zeta integral RS}, for $\Re(s)\gg 1$ we have
\begin{equation}\label{eq2 proof Zeta integral RS}
\displaystyle Z^{\RS}_\psi(s,f)=\int_{i\ago_P^*} Z^{\RS}_\psi(s,f_\lambda)d\lambda.
\end{equation}
Set $\widetilde{f_\lambda}(g)=f_\lambda({}^tg^{-1})$ for every $\lambda\in i\ago_P^*$ and $g\in G(\AAA)$. We will now check that the functions
\begin{equation}\label{eq3 proof Zeta integral RS}
\displaystyle (s,\lambda)\in \CC\times i\ago_P^*\mapsto Z^{\RS}_\psi(s,f_\lambda)\;\; \mbox{ and }\;\; (s,\lambda)\in \CC\times i\ago_P^*\mapsto Z^{\RS}_{\psi^{-1}}(s,\widetilde{f_\lambda})
\end{equation}
satisfy the conditions of Corollary \ref{cor1 PL}.

For $S$ a sufficiently large finite set of places of $F$, that we assume to contain Archimedean places as well as the places where $\pi$ or $\psi$ are ramified, we have, for every $\lambda\in i\ago_P^*$, a decomposition
$$\displaystyle W_{f_\lambda}=W_{\lambda,S}W_\lambda^S$$
where $W_{\lambda,S}\in \wc(\Pi_{\lambda,S},\psi_S)$ and $W_\lambda^S\in \wc(\Pi_\lambda^S,\psi^S)^{K^S}$ is such that $W_\lambda^S(1)=1$. By the unramified computation of local Rankin-Selberg integrals \cite[(3),p.781]{JS2}, \cite[Theorem 3.3]{CogNotes} and \eqref{eq decomposition Tamagawa measure}, we have
\begin{align}\label{eq4 proof Zeta integral RS}
\displaystyle Z_\psi^{\RS}(s,f_\lambda)=(\Delta_H^{S,*})^{-1}L(s+\frac{1}{2},\Pi_\lambda)\frac{Z_\psi^{\RS}(s,W_{\lambda,S})}{L_S(s+\frac{1}{2}, \Pi_\lambda)}
\end{align}
for $\Re(s)\gg 1$ where we have set $L(s,\Pi_\lambda)=L(s,\Pi_{n,\mu}\times \Pi_{n+1,\nu})$ and
$$\displaystyle Z_\psi^{\RS}(s,W_{\lambda,S})=\int_{N_H(F_S)\backslash H(F_S)} W_{\lambda,S}(h_S)\lvert \det h_S\rvert_{E}^s dh_S$$
whenever $\lambda=(\mu,\nu)\in i\ago_P^*$. By \cite{JPSS}, \cite[Theorem 4.1]{CogNotes} and the condition that $\chi$ is $H$-generic (see \eqref{hyp1 RS}), the Rankin-Selberg $L$-function $L(s,\Pi_\lambda)$ is entire and bounded in vertical strips. On the other hand, by \cite[Theorem 2.7]{JPSS} and \cite[Theorem 2.1]{Jacarch}, $s\mapsto \frac{Z_\psi^{\RS}(s,W_{\lambda,S})}{L_S(s+\frac{1}{2}, \Pi_\lambda)}$ has a holomorphic continuation to $\CC$ which is of order at most $1$ in vertical strips and satisfies the functional equation
\begin{equation}\label{eq5 proof Zeta integral RS}
\displaystyle \frac{Z_{\psi^{-1}}^{\RS}(-s,\widetilde{W}_{\lambda,S})}{L_S(\frac{1}{2}-s, \Pi_\lambda^\vee)}=\epsilon(s+\frac{1}{2},\Pi_\lambda)\frac{Z_\psi^{\RS}(s,W_{\lambda,S})}{L_S(s+\frac{1}{2}, \Pi_\lambda)}
\end{equation}
where $\widetilde{W}_{\lambda,S}(g_S)=W_{\lambda,S}(w {}^tg_S^{-1})$ and $\epsilon(s,\Pi_\lambda)$ denotes the global epsilon factor of the Rankin-Selberg $L$-function $L(s,\Pi_\lambda)$. Therefore, by \eqref{eq4 proof Zeta integral RS} and the functional equation of the Rankin-Selberg $L$-function $L(s,\Pi_\lambda)$, we deduce that, for every $\lambda\in i\ago_P^*$, the function $s\mapsto Z_\psi^{\RS}(s,f_\lambda)$ has a holomorphic continuation to $\CC$ which is of order at most $1$ in vertical strips and satisfies the functional equation
\begin{equation}\label{eq6 proof Zeta integral RS}
\displaystyle Z_\psi^{\RS}(s,f_\lambda)=Z_{\psi^{-1}}^{\RS}(-s,\widetilde{f_\lambda}).
\end{equation}
By the first part of Theorem \ref{theo Zeta integral RS}, Theorem \ref{theo Lapid} and Lemma \ref{lem holomorphic maps to Schwartz spaces}, this shows that the functions in \eqref{eq3 proof Zeta integral RS} satisfy the assumptions of Corollary \ref{cor1 PL}. Therefore, the map $s\mapsto \left(\lambda\in i\ago_P^*\mapsto Z_\psi^{\RS}(s,f_\lambda) \right)$ induces a holomorphic function $\CC\to \Sc(i\ago_P^*)$ which is of finite order in vertical strips. By \eqref{eq2 proof Zeta integral RS} and \eqref{eq6 proof Zeta integral RS}, this implies that $s\mapsto Z_\psi^{\RS}(s,f)$ extends to a holomorphic function on $\CC$ of finite order in vertical strips satisfying the functional equation
\begin{equation}\label{eq7 proof Zeta integral RS}
\displaystyle Z_{\psi^{-1}}^{\RS}(-s,\widetilde{f})=Z_\psi^{\RS}(s,f)
\end{equation}
where $\widetilde{f}(g)=f({}^tg^{-1})$. As for $N>0$, the closure of $\Sc_\chi([G])$ in $\tc_{M,\chi}([G])$ contains $\tc_{N,\chi}([G])$ for $M$ sufficiently large (see \eqref{density projection cuspidal component}), Theorem \ref{theo Zeta integral RS}.2 now follows from part 1. and Corollary \ref{cor2 PL} (applied to the closure of $\Sc_\chi([G])$ in $\tc_{M,\chi}([G])$).

\subsection{Proof of Theorem \ref{theo RS period}}

By Theorem \ref{theo Zeta integral RS}.2, and since $\Sc_\chi([G])$ is dense in $\tc_\chi([G])$ (see \eqref{density projection cuspidal component}), it suffices to check that identity \eqref{eq1 theo RS period} is valid for every $f\in \Sc_\chi([G])$. Therefore, let $f\in \Sc_\chi([G])$. For $1\leqslant r\leqslant n$, we let $N_{r,n}$ and $N_{r,n+1}$ be the unipotent radicals of the standard parabolic subgroups of $G_n$ and $G_{n+1}$ with Levi components $G_r\times (G_1)^{n-r}$ and $G_r\times (G_1)^{n+1-r}$ respectively. Set $N^G_r=N_{r,n}\times N_{r,n+1}$ and $N_r^H=N^G_r\cap H=N_{r,n}$. For every $1\leqslant r\leqslant n$, define
$$\displaystyle f_{N^G_r,\psi}(g)=\int_{[N^G_r]} f(ug)\psi_N(u)^{-1}du,\;\; g\in G(\AAA).$$
Whenever the expression below converges absolutely, for every $1\leqslant r\leqslant n$ and $s\in \CC$ we set
$$\displaystyle Z_r^{\RS}(s,f)=\int_{\pc_r(F)N_r^H(\AAA)\backslash H(\AAA)} f_{N^G_r,\psi}(h) \lvert \det h\rvert_{E}^s dh.$$
The proof of the next lemma will be given in Section \ref{Section convergence Zeta RS}.

\begin{lemme}\label{lem theo RS period 1}
For every $1\leqslant r\leqslant n$, there exists $c_r>0$ such that the expression defining $Z^{\RS}_r(s,f)$ converges absolutely for $s\in \cH_{>c_r}$.
\end{lemme}

To uniformize notation, for every $s\in \CC$ we also set
$$\displaystyle Z_{n+1}^{\RS}(s,f)=\int_{[H]} f(h) \lvert \det h\rvert_{E}^s dh.$$
Note that the above expression is absolutely convergent and defines an entire function of $s\in \CC$ satisfying $P_H(f)=Z_{n+1}^{\RS}(0,f)$. On the other hand, we have $Z_1^{\RS}(s,f)=Z_\psi^{\RS}(s,f)$. Hence, identity \eqref{eq1 theo RS period} is a consequence of the following proposition.

\begin{proposition}\label{prop theo RS period 1}
For every $1\leqslant r\leqslant n$, we have
\begin{equation}
\displaystyle Z_{r+1}^{\RS}(s,f)=Z_r^{\RS}(s,f)
\end{equation}
for $\Re(s)\gg 1$.
\end{proposition}

\begin{preuve}
Let $1\leqslant r\leqslant n-1$. As $\pc_{r+1}=G_rU_r$ and $N_{r,n}=U_rN_{r+1,n}$, for $s\in \cH_{>c_{r+1}}$ we have
\begin{align}\label{eq1 proof theo RS period}
\displaystyle Z^{\RS}_{r+1}(s,f)=\int_{G_r(F)N_r^H(\AAA)\backslash H(\AAA)} \int_{[U^H_{r+1}]} f_{N^G_{r+1},\psi}(uh)du \lvert \det h\rvert_{E}^s dh
\end{align}
where we have set $U^H_{r+1}=U_{r+1}$ viewed as a subgroup of $H=G_n$ (as always via the embedding in ``the upper-left corner''). Similarly, we set $U^G_{r+1}=U_{r+1}\times U_{r+1}$ viewed as a subgroup of $G$. By Fourier inversion on the compact abelian group $U_{r+1}^G(F)U^H_{r+1}(\AAA)\backslash U^G_{r+1}(\AAA)$, we have
\begin{equation}\label{eq2 proof theo RS period}
\displaystyle \int_{[U^H_{r+1}]} f_{N^G_{r+1},\psi}(uh)du=\sum_{\gamma\in \pc_r(F)\backslash G_r(F)} (f_{N^G_{r+1},\psi})_{U_{r+1}^G,\psi}(\gamma h)+ (f_{N^G_{r+1},\psi})_{U_{r+1}^G}(h)
\end{equation}
for every $h\in H(\AAA)$, where we have set
$$\displaystyle (f_{N^G_{r+1},\psi})_{U_{r+1}^G,\psi}(h)=\int_{[U_{r+1}^G]} f_{N^G_{r+1},\psi}(uh)\psi_N(u)^{-1}du=f_{N^G_{r},\psi}(h),$$
$$\displaystyle (f_{N^G_{r+1},\psi})_{U_{r+1}^G}(h)=\int_{[U_{r+1}^G]} f_{N^G_{r+1},\psi}(uh)du.$$

By \eqref{eq1 proof theo RS period} and \eqref{eq2 proof theo RS period}, we obtain
\begin{equation}\label{eq3 proof theo RS period}
\displaystyle Z^{\RS}_{r+1}(s,f)=Z^{\RS}_r(s,f)+F_r(s)
\end{equation}
for every $s\in \CC$ such that $\Re(s)>\max(c_r,c_{r+1})$ where we have set
$$\displaystyle F_r(s)=\int_{G_r(F)N_r^H(\AAA)\backslash H(\AAA)} (f_{N^G_{r+1},\psi})_{U_{r+1}^G}(h) \lvert \det h\rvert_{E}^s dh.$$
By a similar argument, \eqref{eq3 proof theo RS period} still holds when $n=r$ if we set
$$\displaystyle F_n(s)=\int_{[H]} f_{U_{n+1}^G}(h) \lvert \det h\rvert_{E}^s dh$$
with $U_{n+1}^G=1\times U_{n+1}$ and
$$\displaystyle f_{U_{n+1}^G}(h)=\int_{[U^G_{n+1}]} f(uh)du.$$
From \eqref{eq3 proof theo RS period}, we are reduced to showing that $F_r(s)=0$ identically for every $1\leqslant r\leqslant n$ and $\Re(s)$ sufficiently large. To uniformize notation, we set $f_{N_{n+1}^G,\psi}=f$. Let $P_r$ be the standard parabolic subgroup of $G$ with Levi component $L_r=\left(G_r\times G_{n-r}\right) \times \left(G_r\times G_{n+1-r}\right)$ and set $P_r^H=P_r\cap H$. Then, we readily check that
$$\displaystyle (f_{N^G_{r+1},\psi})_{U_{r+1}^G}(h)=\int_{ [N_{n-r}]\times [N_{n+1-r}]} f_{P_r}\left(\left(\begin{pmatrix} I_r & \\ & u \end{pmatrix},\begin{pmatrix} I_r & \\ & u' \end{pmatrix}\right)h\right)\psi_n(u)^{-1}\psi_{n+1}(u')^{-1} du'du$$
for every $h\in H(\AAA)$ and $1\leqslant r\leqslant n$. Therefore, by the Iwasawa decomposition $H(\AAA)=P_r^H(\AAA)K_n$, we have
\[\begin{aligned}
\displaystyle F_r(s)= & \int_{K_n} \int_{[G_r]\times N_{n-r}(\AAA)\backslash G_{n-r}(\AAA)} \int_{[N_{n-r}]\times [N_{n+1-r}]} f_{P_r}\left(\begin{pmatrix} h_r & \\ & uh_{n-r}\end{pmatrix}k, \begin{pmatrix} h_r & \\ & u'h_{n-r}\end{pmatrix}k \right) \\
 & \psi_n(u)^{-1}\psi_{n+1}(u')^{-1}du'du\delta_{Q_r^H}\begin{pmatrix} h_r & \\ & h_{n-r} \end{pmatrix}^{-1} \lvert \det h_r\rvert_{E}^s \lvert \det h_{n-r}\rvert_{E}^s dh_{n-r}dh_rdk.
\end{aligned}\]
By a painless calculation, left to the reader, we have
$$\displaystyle \delta_{P_r^H}\begin{pmatrix} h_r & \\ & h_{n-r} \end{pmatrix}^{-1}\lvert \det h_r\rvert_{E}^s \lvert \det h_{n-r}\rvert_{E}^s=\delta_{P_r}\left( \begin{pmatrix} h_r & \\ & h_{n-r}\end{pmatrix}, \begin{pmatrix} h_r & \\ & h_{n-r}\end{pmatrix}\right)^{-\frac{1}{2}+\alpha_r(s)} \lvert \det h_{n-r}\rvert_{E}^{s+2r\alpha_r(s)}$$
where $\alpha_r(s)=\frac{2s+1}{4n-4r+2}$. Let $\chi^L$ be the inverse image of $\chi$ in $\Xgo(L_r)$. By Corollary \ref{cor constant term generic cuspidal datum}, for $\Re(\alpha_r(s))>0$ and every $k\in K_n$ the function $f_{P_r,k,s}:=\delta_{P_r}^{-\frac{1}{2}+\alpha_r(s)} R(k)f_{P_r}\mid_{[L_r]}$ belongs to $\cC_{\chi^L}([L_r])$. On the other hand, as $L_r=G_r\times G_{n-r}\times G_r\times G_{n+1-r}$, by \eqref{eq isom HCS product with cuspidal datum} we have the decomposition
$$\displaystyle \cC_{\chi^L}([L_r])=\bigoplus_{(\chi_1,\chi_2)\in \Xgo(G_r^2)\times \Xgo(G_{n-r}\times G_{n+1-r}) \mapsto \chi} \cC_{\chi_1}([G_r\times G_r])\widehat{\otimes} \cC_{\chi_2}([G_{n-r}\times G_{n+1-r}])$$
and the above equality can be rewritten as
\begin{equation}\label{eq5 proof theo RS period}
\displaystyle F_r(s)=\int_{K_n} \left(P_{G_r^\Delta}\widehat{\otimes} \zc_{n-r}^{\RS}(s+2r\alpha_r(s)) \right)(f_{P_r,k,s}) dk
\end{equation}
where $P_{G_r^\Delta}$ denotes the period integral over the diagonal subgroup of $G_r\times G_r$ and $\zc_{n-r}^{\RS}(s)$ stands for the continuous linear form
$$\displaystyle f'\in \cC([G_{n-r}\times G_{n+1-r}])\mapsto Z_{\psi^{(-1)^r}}^{\RS}(s,f').$$
Since $\chi$ is $H$-generic, by \eqref{hyp1 RS} any preimage $(\chi_1,\chi_2)\in \Xgo(G_r^2)\times \Xgo(G_{n-r}\times G_{n+1-r})$ of $\chi$ with $\chi_1=(\chi_1',\chi_1'')\in \Xgo(G_r)^2$ we have $\chi''_1\neq (\chi'_1)^\vee$. Hence, by definition of $\cC_{\chi_1}([G_r\times G_r])$, $P_{G_r^\Delta}$ vanishes identically on $\cC_{\chi_1}([G_r\times G_r])$. This implies that $F_r(s)=0$ whenever $\Re(s)\gg 1$ and this ends the proof of the proposition.
\end{preuve}

\subsection{Convergence of Zeta integrals}\label{Section convergence Zeta RS}

\begin{preuve}(of Lemma \ref{lem theo RS period 1}) The argument is very similar to the proof of Lemma \ref{lem conv 2} so we only sketch it. Let $1\leqslant r\leqslant n$ and $Q_r^G$ be the standard parabolic subgroup of $G$ with Levi component $(G_r\times (G_1)^{n-r})\times (G_r\times (G_1)^{n+1-r})$ so that $N_r^G$ is the unipotent radical of $Q_r^G$. Set $Q_r^H=Q_r^G\cap H$. By the Iwasawa decomposition $H(\AAA)=Q_r^H(\AAA)K_n$, we need to show the convergence of
\begin{align}\label{eq1 conv Zeta integrals RS}
\displaystyle \int_{K_n\times \pc_r(F)\backslash G_r(\AAA)\times T_{n-r}(\AAA)} \left\lvert(R(k)f)_{N_r^G,\psi}\begin{pmatrix} h & \\ & t \end{pmatrix}\right\rvert \lvert \det h\rvert_{E}^s \lvert \det t\rvert_{E}^s \delta_{Q_r^H}\begin{pmatrix} h & \\ & t \end{pmatrix}^{-1}dtdhdk
\end{align}
for $\Re(s)\gg 1$. We apply Lemma \ref{lem1 majorations Fourier}.1 to $\psi_F=\psi'$ and
$$\displaystyle \ell:N_r^G\to \GG_a,$$
$$\displaystyle (u,u')\mapsto \Tr_{E/F}\left( (-1)^n\tau\sum_{i=r}^{n-1} u_{i,i+1}+(-1)^{n+1}\tau\sum_{i=r}^{n} u'_{i,i+1}\right).$$
It is easy to see that there exists $N_0>0$ such that
\begin{equation}\label{eq2 conv Zeta integrals RS}
\displaystyle \lVert e_r h\rVert_{\AAA^r_E}\prod_{i=1}^{n-r} \lVert t_i\rVert_{\AAA_E}\ll \lVert \Ad^*\begin{pmatrix} h & \\ & t \end{pmatrix}\ell\rVert_{(N_r^G)^*_{\ab}(\AAA)}^{N_0}, \; (h,t)\in G_r(\AAA)\times T_{n-r}(\AAA).
\end{equation}
Therefore, from \eqref{eq2 conv Zeta integrals RS} and Lemma \ref{lem1 majorations Fourier}, there exists $c>0$ such that for every $N_1,N_2>0$, \eqref{eq1 conv Zeta integrals RS} is essentially bounded by
$$\displaystyle \int_{\pc_r(F)\backslash G_r(\AAA)\times T_{n-r}(\AAA)} \lVert h\rVert_{[G_r]}^{-N_2} \lVert e_r h\rVert_{\AAA_E^r}^{-N_1}\prod_{i=1}^{n-r} \lVert t_i\rVert_{\AAA_E}^{-N_1} \delta_{Q_r^G}\begin{pmatrix} h & \\ & t \end{pmatrix}^{-cN_2}\delta_{Q_r^H}\begin{pmatrix} h & \\ & t \end{pmatrix}^{-1} \lvert \det h\rvert_{E}^s \lvert \det t\rvert_{E}^s dt dh.$$
Now, the convergence of the above expression for $\Re(s)\gg 1$, $N_2\gg_s 1$ and $N_1\gg_{s,N_2} 1$ can be shown as in the end of the proof of Lemma \ref{lem conv 2}  using Lemma \ref{lem conv GLn}.
\end{preuve}

\begin{preuve}(of Theorem \ref{theo Zeta integral RS}.1) Applying Lemma \ref{lem1 majorations Fourier}.2 in a similar way, we are reduced to showing the existence of $c_N>0$ such that for every $C>c_N$ there exists $N'>0$ satisfyng that the integral
$$\displaystyle \int_{T_n(\AAA)} \prod_{i=1}^n \lVert t_i\rVert_{\AAA_E}^{-N'} \lVert t \rVert_{[T_n]}^N\delta_{B_n}(t)^{-1} \lvert \det t\rvert_{E}^s dt$$
converges in the range $s\in \cH_{]c_n,C[}$ uniformly on any compact subsets. This is exactly what was established in the proof of Theorem \ref{theo Zeta integrals FR}.1 (up to replacing the base field $F$ by $E$).
\end{preuve}

\section{Contributions of certain cuspidal data to the Jacquet-Rallis trace formula: second proof}\label{Chap:generic-JRTF}

In this chapter, we adopt the set of notation introduced in Chapter \ref{sec:generic contribution}. In particular, $n\geqslant 1$ is a positive integer, $G=G_n\times G_{n+1}$, $G'=G'_n\times G'_{n+1}$, $H=G_n$ with its diagonal embedding in $G$, $K=K_n\times K_{n+1}$ and $K'=K'_n\times K'_{n+1}$ are the standard maximal compact subgroups of $G(\AAA)$ and $G'(\AAA)$ respectively and $\eta_{G'}:[G']\to \{ \pm 1\}$ is the automorphic character defined in \S \ref{S:character etaG'}. We will also use notation from Chapters \ref{chapter FR} and \ref{chap canonical RS}: $N=N_n\times N_{n+1}$ and $N_H=N_n$ are the standard maximal unipotent subgroups of $G$ and $H$, $\psi_N=\psi_n\boxtimes \psi_{n+1}$ is a generic character of $[N]$ (where $\psi_n$ and $\psi_{n+1}$ are defined as in \S \ref{ssect:additive and generic characters}). We also set $\pc=\pc_n\times \pc_{n+1}$ (resp. $\pc'=\pc'_n\times \pc'_{n+1}$) where $\pc_n$ and $\pc_{n+1}$ (resp. $\pc'_n$ and $\pc'_{n+1}$) stand for the mirablic subgroups of $G_n$ and $G_{n+1}$ (resp. of $G'_n$ and $G'_{n+1}$), $T=T_n\times T_{n+1}$ for the standard maximal torus of $G$ and $N'=N'_n\times N'_{n+1}$ for the standard maximal unipotent subgroup of $G'$. Finally, as in \S \ref{ssect:RS Zeta}, for every $f\in \tc([G])$ we set
$$\displaystyle W_f(g)=\int_{[N]} f(ug) \psi_N(u)^{-1}du,\;\;\; g\in G(\AAA).$$

\subsection{Main result}

\begin{paragr}\label{8.1.1?}
Let $\chi\in \Xgo^*(G)$ be a $*$-generic cuspidal datum (see \S \ref{S:relevant}) represented by a pair $(M_P,\pi)$. We set $\Pi=I_{P(\AAA)}^{G(\AAA)}(\pi)$. We have decompositions $P=P_n\times P_{n+1}$, $\pi=\pi_n\boxtimes \pi_{n+1}$ and $\Pi=\Pi_n\boxtimes \Pi_{n+1}$ where: $P_n$, $P_{n+1}$ are standard parabolic subgroups of $G_n$, $G_{n+1}$ respectively with standard Levi components of the form
$$\displaystyle M_{P_n}=G_{n_1}\times\ldots\times G_{n_k},\;\; M_{P_{n+1}}=G_{m_1}\times \ldots\times G_{m_r},$$
$\pi_n$ and $\pi_{n+1}$ are cuspidal automorphic representations of $M_{P_n}(\AAA)$, $M_{P_{n+1}}(\AAA)$ decomposing into tensor products
$$\displaystyle \pi_n=\pi_{n,1}\boxtimes \ldots\boxtimes \pi_{n,k},\;\; \pi_{n+1}=\pi_{n+1,1}\boxtimes \ldots \boxtimes \pi_{n+1,r}$$
respectively and we have set $\Pi_n=I_{P_n(\AAA)}^{G_n(\AAA)}(\pi_n)$, $\Pi_{n+1}=I_{P_{n+1}(\AAA)}^{G_{n+1}(\AAA)}(\pi_{n+1})$. We write $\chi_n\in \Xgo^*(G_n)$ and $\chi_{n+1}\in \Xgo^*(G_{n+1})$ for the cuspidal data determined by the pairs $(M_{P_n},\pi_n)$ and $(M_{P_{n+1}},\pi_{n+1})$ respectively.

The representation $\Pi$ is generic and we denote by $\wc(\Pi_,\psi_N)$ its Whittaker model with respect to the character $\psi_N$. Also, for every $\phi\in \Pi$ we define
$$\displaystyle W_\phi(g):=W_{E(\phi)}(g)=\int_{[N]} E(ug,\phi)\psi_N(u)^{-1} du,\;\; g\in G(\AAA).$$
Note that $W_\phi\in \wc(\Pi,\psi_N)$.
\end{paragr}

\begin{paragr}\label{S:Whittaker periods}
We now define two continuous linear forms $\lambda$ and $\beta_\eta$ as well as a continuous invariant scalar product $\langle .,.\rangle_{\Whitt}$ on $\wc(\Pi,\psi_N)$. Let $W\in \wc(\Pi,\psi_N)$.
\begin{itemize}
\item By \cite{JPSS} \cite{Jacarch}, the Zeta integral (already encountered in Chapter \ref{chap canonical RS})	
$$\displaystyle Z^{\RS}(s,W)=\int_{N_{H}(\AAA)\backslash H(\AAA)} W(h) \lvert \det h\rvert_{\AAA_E}^s dh.$$
converges for $\Re(s)\gg 0$ and extends to a meromorphic function on $\CC$ with no pole at $s=0$. We set
$$\displaystyle \lambda(W)=Z^{\RS}(0,W).$$
\item For $S$ a sufficiently large finite set of places of $F$, we put
$$\displaystyle \beta_\eta(W)=(\Delta^{S,*}_{G'})^{-1}L^{S,*}(1,\Pi,\As_G)\int_{N'(F_S)\backslash \pc'(F_S)} W(p_S) \eta_{G'}(p_S) dp_S$$
where we have set $L(s,\Pi,\As_G)=L(s,\Pi_n,\As^{(-1)^{n+1}})L(s,\Pi_{n+1},\As^{(-1)^n})$.
\item Similarly, for $S$ a sufficiently large finite set of places of $F$, we put
$$\displaystyle \langle W,W\rangle_{\Whitt}=(\Delta^{S,*}_{G})^{-1}L^{S,*}(1,\Pi,\Ad)\int_{N(F_S)\backslash \pc(F_S)} \lvert W(p_S)\rvert^2 dp_S$$
where we have set $L(s,\Pi,\Ad)=L(s,\Pi_n\times \Pi_n^\vee)L(s,\Pi_{n+1}\times \Pi_{n+1}^\vee)$.
\end{itemize}
That the above expressions converge and are independent of $S$ as soon as it is chosen sufficiently large (depending on the level of $W$) follow from \cite{Flicker} and \cite{JS}. Moreover, the inner form $\langle .,.\rangle_{\Whitt}$ is $G(\AAA)$-invariant by \cite{Bermirabolic} and \cite{Baruch}.

The next result follows from works of Jacquet-Shalika \cite{JS}, Shahidi \cite{ShaLfunction} and Lapid-Offen \cite[Appendix A]{FLO}. For completness, we explain the deduction (see \S \ref{S:automorphic repn} for our normalization of the Petersson inner product).

\begin{theoreme}[Jacquet-Shalika, Shahidi, Lapid-Offen]\label{theo inner products}
We have
$$\displaystyle \langle \phi,\phi\rangle_{\Pet}=\langle W_{\phi},W_{\phi}\rangle_{\Whitt}$$
for every $\phi\in \Pi$.
\end{theoreme}

\begin{preuve}
Let $\phi\in \Pi$. By the Iwasawa decomposition, for a suitable Haar measure on $K$ we have
$$\displaystyle \langle \phi,\phi\rangle_{\Pet}=\int_{K} \int_{[M_P]_0} \lvert \phi(mk)\rvert^2 \delta_P(m)^{-1}dmdk.$$
Set $N^P=N\cap M_P$ and
$$\displaystyle \phi_{N^P,\psi}(g)=\int_{[N^P]} \phi(ug)\psi_N(u)^{-1}du,\;\;\; g\in G(\AAA).$$
Let $\pc^P$ be the product of mirabolic groups $\prod_{i=1}^k \pc_{n_i}\times \prod_{j=1}^r \pc_{m_j}$. It is a subgroup of $M_P$. According to Jacquet-Shalika \cite[\S 4]{JS} (see also \cite[p.265]{FLO} or \cite[Proposition 3.1]{Zhang2}\footnote{Note that our normalization of the Petterson inner product if different from {\em loc. cit.}}), for $S$ a sufficiently large finite set of places of $F$ we have
\begin{align}\label{eq1 JS period}
\displaystyle \int_{[M_P]_0} \lvert \phi(mk)\rvert^2 \delta_P(m)^{-1}dm= & (\Delta_{M_P}^{S,*})^{-1} \prod_{i=1}^k \Res_{s=1} L^S(s,\pi_{n,i}\times \pi^\vee_{n,i}) \prod_{j=1}^r \Res_{s=1} L^S(s,\pi_{n+1,j}\times \pi^\vee_{n+1,j}) \\
\nonumber & \times \int_{N^P(F_S)\backslash \pc^P(F_S)} \lvert \phi_{N^P,\psi}(p_Sk)\rvert^2\delta_P(p_S)^{-1}dp_S
\end{align}
for every $k\in K$. On the other hand, by \cite[Proposition A.2]{FLO} we have
\begin{align}\label{eq2 JS period}
\displaystyle  & \int_K \int_{N^P(F_S)\backslash \pc^P(F_S)} \lvert \phi_{N^P,\psi}(p_Sk)\rvert^2\delta_P(p_S)^{-1}dp_Sdk \\
\nonumber & = \frac{\vol_{G(\AAA^S)}(K^S)}{\vol_{M_P(\AAA^S)}(K^S\cap M_P(\AAA^S))} \int_{P(F_S)\backslash G(F_S)} \int_{N^P(F_S)\backslash \pc^P(F_S)} \lvert \phi_{N^P,\psi}(p_Sg_S)\rvert^2\delta_P(p_S)^{-1}dp_Sdg_S \\
\nonumber & =(\Delta_G^{S,*})^{-1}\Delta_{M_P}^{S,*}\int_{N(F_S)\backslash \pc(F_S)} \lvert \mathbf{W}_S(p_S,\phi_{N^P,\psi})\rvert^2dp_S.
\end{align}
where $\mathbf{W}_S: I_{P(F_S)}^{G(F_S)}(\wc(\pi_S,\psi_{N,S}))\to \wc(\Pi_S,\psi_{N,S})$ stands for the {\em Jacquet functional}, defined as the value at $s=0$ of the holomorphic continuation of
$$\displaystyle \mathbf{W}^S_s(g_S,\phi')=\int_{(w^G_P)^{-1}N^P(F_S)w^G_P\backslash N(F_S)} \phi'(w^G_Pu_Sg_S)\delta_P(w^G_Pu_Sg_S)^s\psi_N(u_S)^{-1}du_S,\;\; \Re(s)\gg 1$$
for $g_S\in G(F_S)$ and $\phi'\in I_{P(F_S)}^{G(F_S)}(\wc(\pi_S,\psi_{N,S}))$ where $w^G_P=w^{P}w^G$ with $w^{P}$ (resp. $w^G$) the permutation matrix representing the longest element in the Weyl group of $T$ in $M_P$ (resp. in $G$). Finally, by \cite[Sect. 4]{ShaLfunction}, we have
\begin{equation}\label{eq3 JS period}
\displaystyle \mathbf{W}^S(\phi_{N^P,\psi})=\prod_{1\leqslant i<j\leqslant k} L^S(1,\pi_{n,i}\times \pi^\vee_{n,j}) \prod_{1\leqslant i<j\leqslant r} L^S(1,\pi_{n+1,i}\times \pi^\vee_{n+1,j}) W_\phi.
\end{equation}
(Note that, as $\chi$ is generic, the Rankin-Selberg $L$-functions $L(s,\pi_{n,i}\times \pi^\vee_{n,j})$ and $L(s,\pi_{n+1,i}\times \pi^\vee_{n+1,j})$ are all regular at $s=1$.) As, for every $s\in \RR$,
\[\begin{aligned}
\displaystyle L^S(s,\Pi\times \Pi^\vee)= & \prod_{i=1}^k L^S(s,\pi_{n,i}\times \pi^\vee_{n,i}) \times \prod_{j=1}^r L^S(s,\pi_{n+1,j}\times \pi^\vee_{n+1,j}) \\
 & \times \left\lvert\prod_{1\leqslant i<j\leqslant k} L^S(s,\pi_{n,i}\times \pi^\vee_{n,j})\right\rvert^2 \times \left\lvert\prod_{1\leqslant i<j\leqslant r} L^S(s,\pi_{n+1,i}\times \pi^\vee_{n+1,j})\right\rvert^2,
\end{aligned}\]
we deduce from \eqref{eq1 JS period}, \eqref{eq2 JS period} and \eqref{eq3 JS period} the identity of the Theorem.
\end{preuve}
\end{paragr}

\begin{paragr}[Relative characters.] --- \label{S:IPI}
Let $\bc_{P,\pi}$ be a $K$-basis of $\Pi$ as in \S \ref{S:K0ON}. We define the {\em relative character} $I_\Pi$ of $\Pi$ as the following functional on $\Sc(G(\AAA))$:
$$\displaystyle I_\Pi(f)=\sum_{\phi\in \bc_{P,\pi}} \frac{\lambda(R(f)W_\phi) \overline{\beta_\eta(W_\phi)}}{\langle W_\phi,W_\phi\rangle_{\Whitt}},\; f\in \Sc(G(\AAA)),$$
where the series converges, and does not depend on the choice of $\bc_{P,\pi}$, by Proposition \ref{prop:car-relJB}.
\end{paragr}

\begin{paragr}
For every $f\in \Sc(G(\AAA))$, we set
$$\displaystyle K_{f,\chi}^1(g)=\int_{[H]} K_{f,\chi}(h,g) dh \; \mbox{ and } \; K_{f,\chi}^2(g)=\int_{[G']} K_{f,\chi}(g,g') \eta_{G'}(g')dg',\;\; g\in [G],$$
where the above expressions are absolutely convergent by Lemma \ref{lem:maj-noy}.3.

Recall that the notion of relevant $*$-generic cuspidal datum has been defined in \S \ref{S:relevant} and that we have defined for any $\chi\in \Xgo$ a distribution $I_\chi$ (see Theorem \ref{thm:jfDef}). 

\begin{theoreme}\label{theo JR TF}
Let $f\in \Sc(G(\AAA))$ and $\chi\in \Xgo^*(G)$. Then,
\begin{enumerate}
\item If $\chi$ is not relevant, we have $K_{f,\chi}^2(g)=0$ for every $g\in [G]$ and moreover
$$\displaystyle I_\chi(f)=0.$$
\item If $\chi$ is relevant, we have 
$$\displaystyle I_\chi(f)=\int_{[G']} K_{f,\chi}^1(g')\eta_{G'}(g')dg'$$
where the right-hand side converges absolutely and moreover
$$\displaystyle I_\chi(f)=2^{-\dim(A_P)}I_{\Pi}(f).$$ 
\end{enumerate}

\end{theoreme}

The rest of this chapter is devoted to the proof of Theorem \ref{theo JR TF}. Until the end, we fix a function $f\in \Sc_\chi(G(\AAA))$.
\end{paragr}

\subsection{Proof of Theorem \ref{theo JR TF}}

\begin{paragr}
We fix a character $\eta_G$ of $[G]$ whose restriction to $[G']$ is equal to $\eta_{G'}$ (such a character exists as the id\`ele class group of $F$ is a closed subgroup of the id\`ele class group of $E$) and we set $\widetilde{\chi}=\eta_G\otimes \chi^\vee\in \Xgo^*(G)$. We can write $\widetilde{\chi}$ as $(\widetilde{\chi}_n,\widetilde{\chi}_{n+1})$ where $\widetilde{\chi}_k\in \Xgo^*(G_k)$ for $k=n,n+1$. For every $g\in [G]$, we denote by $\widetilde{K}_{f,\chi}(g,.)$ the function $\eta_G K_{f,\chi}(g,.)$. By Lemma \ref{lem:maj-noy}.2 and \eqref{eq isom HCS product with cuspidal datum}, we have
\begin{equation}\label{eq 1 non-relevant case}
\displaystyle \widetilde{K}_{f,\chi}(g,.)\in \Sc_{\widetilde{\chi}}([G])=\Sc_{\widetilde{\chi}_n}([G_n])\widehat{\otimes} \Sc_{\widetilde{\chi}_{n+1}}([G_{n+1}])
\end{equation}
for all $g\in [G]$. Moreover, with the notation of Theorem \ref{theo FR period}, we have
\begin{equation}\label{eq 2 non-relevant case}
\displaystyle K_{f,\chi}^2(g)=P_{G'_n}\widehat{\otimes} P_{G'_{n+1}}(\widetilde{K}_{f,\chi}(g,.)).
\end{equation}
\end{paragr}

\begin{paragr}[The non-relevant case]
Assume that $\chi$ is not relevant. By definition of a relevant cuspidal data (see \S \ref{S:relevant}), at least one of $\widetilde{\chi}_n$, $\widetilde{\chi}_{n+1}$ is not distinguished (see \S \ref{S:distinguished cuspidal datum} for the definition of distinguished). Hence, by Theorem \ref{theo Zeta integrals FR} and Theorem \ref{theo FR period}, $P_{G'_k}$ vanishes identically on $\Sc_{\widetilde{\chi}_k}([G_k])$ for $k=n$ or $k=n+1$. Thus, by \eqref{eq 1 non-relevant case} and \eqref{eq 2 non-relevant case}, the function $K_{f,\chi}^2$ vanishes identically. By Theorem \ref{thm:asym-trio} applied to the expression \eqref{eq:FonH}, this implies $I_\chi(f)=0$. This proves part 1. of Theorem \ref{theo JR TF}.
\end{paragr}

\begin{paragr}[Regularized Rankin-Selberg period and convergence]\label{S:RRS-period-cv}
From now on, we assume that $\chi$ is relevant. By Lemma \ref{lem:maj-noy}.2, for every $g\in [G]$ the function $K_{f,\chi}(.,g)$ belongs to $\Sc_\chi([G])$. Since $\chi$ is relevant, it is $H$-generic in the sense of \S \ref{S:Hgeneric} (this follows from the dichotomy of \S \ref{S:disting}). Therefore, by Theorem \ref{theo RS period}, $P_H$ extends to a continuous linear form on $\tc_\chi([G])$ that we shall denote by $P_H^*$. By definition of this extension and of the linear form $\lambda$ (see \S \ref{S:Whittaker periods}), for every $\phi\in \Pi$ we have
\begin{equation}\label{eq PH* and lambda}
\displaystyle P_H^*(E(\phi))=\lambda(W_\phi).
\end{equation}
By Lemma \ref{lem:maj-noy}.3 there exists $N\geqslant 0$ such that the function
$$\displaystyle g'\in [G']\mapsto K_{f,\chi}(.,g')\in \tc_N([G])$$
is absolutely integrable. As
$$\displaystyle K_{f,\chi}^1(g)=P_H(K_{f,\chi}(.,g))=P_H^*(K_{f,\chi}(.,g)),$$
combined with Theorem \ref{thm:asym-trio} applied to the expression \eqref{eq:FonG}, this shows at once that the expression 
$$\displaystyle \int_{[G']} K_{f,\chi}^1(g')\eta_{G'}(g')dg'$$
converges absolutely, is equal to $I_\chi(f)$ and that
\begin{equation}\label{eq 1 convergence Ichi}
\displaystyle I_\chi(f)=P_H^*\left(\int_{[G']}K_{f,\chi}(.,g') \eta_{[G']}(g') dg'\right)=P_H^*(K_{f,\chi}^2).
\end{equation}
\end{paragr}

\begin{paragr}[Spectral expression of $K_{f,\chi}^2$]\label{paragr spectral exp K2}
Set $\widetilde{\Pi}=\Pi^\vee \otimes \eta_G$. We may write $\widetilde{\Pi}$ as a tensor product $\widetilde{\Pi}_n\boxtimes\widetilde{\Pi}_{n+1}$ and we let
$$\displaystyle \beta=\beta_n\widehat{\otimes}\beta_{n+1}: \wc(\widetilde{\Pi},\psi_N)=\wc(\widetilde{\Pi}_n,\psi_n)\widehat{\otimes}\wc(\widetilde{\Pi}_{n+1},\psi_{n+1}) \to \CC$$
be the (completed) tensor product of the linear forms $\beta_n$, $\beta_{n+1}$ defined in \S \ref{S:betan}. Fix $g\in [G]$ and set $\mathbf{f}_g=\widetilde{K}_{f,\chi}(g,.)$. Since $\chi$ is relevant, $\widetilde{\chi}_n$ and $\widetilde{\chi}_{n+1}$ are both distinguished. Note that the linear map
$$\displaystyle \mathbf{f}\in \Sc([G])\mapsto W_{\mathbf{f},\widetilde{\Pi}}:=W_{\mathbf{f}_{\widetilde{\Pi}}}\in \wc(\widetilde{\Pi}, \psi_N)$$
is the (completed) tensor product of the continuous linear maps $\mathbf{f}\in \Sc([G_k])\mapsto W_{\mathbf{f},\widetilde{\Pi}_k}\in \wc(\widetilde{\Pi}_k,\psi_k)$ for $k=n,n+1$ (as can be checked directly on pure tensors). Therefore, by \eqref{eq 1 non-relevant case}, \eqref{eq 2 non-relevant case}, Theorem \ref{theo Zeta integrals FR} and Theorem \ref{theo FR period} we have

\begin{equation}\label{eq5 JR trace formula}
\displaystyle K_{f,\chi}^2(g)=2^{-\dim(A_P)}\beta(W_{\mathbf{f}_g,\widetilde{\Pi}}).
\end{equation}
Let $\bc_{P,\pi}$ be a $K$-basis $\Pi$ as in \S\ref{S:K0ON}. Then, we have $\displaystyle \mathbf{f}_{g,\widetilde{\Pi}}=\sum_{\phi\in \bc_{P,\pi}} \langle \mathbf{f}_g,\eta_G E(\overline{\phi})\rangle_{[G]} \eta_GE(\overline{\phi})$ where the sum converges absolutely in $\tc_N([G])$ for some $N\geqslant 0$. Hence,
$$\displaystyle W_{\mathbf{f}_g,\widetilde{\Pi}}=\sum_{\phi\in \bc_{P,\pi}} \langle \mathbf{f}_g,\eta_G E(\overline{\phi})\rangle_{[G]} \eta_G \overline{W_{\phi}}$$
in $\wc(\widetilde{\Pi},\psi_N)$. On the other hand, we easily check that $\beta(\eta_G \overline{W_{\phi}})=\overline{\beta_\eta(W_{\phi})}$ and
$$\displaystyle \langle \mathbf{f}_g,\eta_G E(\overline{\phi})\rangle_{[G]}=\langle K_{f,\chi}(g,.), E(\overline{\phi})\rangle_{[G]}=E(R(f)\phi)(g)$$
for every $\phi\in \bc_{P,\pi}$. Therefore, by \eqref{eq5 JR trace formula}, we obtain

\begin{equation}\label{eq6 JR trace formula}
\displaystyle K_{f,\chi}^2(g)=2^{-\dim(A_P)} \sum_{\phi\in \bc_{P,\pi}} E(R(f)\phi)(g) \overline{\beta_\eta(W_{\phi})}.
\end{equation}
Note that by Proposition \ref{prop:car-relJB}, the series above is actually absolutely convergent in $\tc_N([G])$ for some $N\geqslant 0$ (and not just pointwise).
\end{paragr}

\begin{paragr}[End of the proof]
By \eqref{eq 1 convergence Ichi}, \eqref{eq PH* and lambda} and \eqref{eq6 JR trace formula}, we obtain
$$\displaystyle I_\chi(f)=2^{-\dim(A_P)} \sum_{\phi\in \bc_{P,\pi}} \lambda(R(f)W_\phi) \overline{\beta_\eta(W_{\phi})}.$$
Using Theorem \ref{theo inner products} and since $\bc_{P,\pi}$ is an orthonormal basis of $\Pi$, this can be rewritten as
$$\displaystyle I_\chi(f)=2^{-\dim(A_P)} \sum_{\phi\in \bc_{P,\pi}} \frac{\lambda(R(f)W_\phi) \overline{\beta_\eta(W_{\phi})}}{\langle W_\phi,W_\phi\rangle_{\Whitt}}=2^{-\dim(A_P)}I_\Pi(f)$$
and this ends the proof of Theorem \ref{theo JR TF} in the relevant case.
\end{paragr}

\section{Flicker-Rallis functional computation}\label{Chap:FR-functional-computation}

The goal of this chapter is to prove Theorem \ref{thm-intrp:J-beta} of the introduction that states that two natural 
functionals are equal. This is established in Theorem \ref{thm:comparisonG}. The bulk of the work is in proving its local avatar. 
The case of of split algebra $E/F$ amounts to comparing scalar products which was done in Appendix A of \cite{FLO}, which is an inspiration for this chapter.

\subsection{Local comparison}

\bpar\label{ssec:measLoc}
Let $E/F$ be an etale quadratic algebra over a  local field $F$. Let $\Tr_{E/F} : E \to F$ be the trace map. 
As in Paragraph \ref{ssect:additive and generic characters}, let $\psi' : F \to \CC^{\times}$ 
be a non-trivial additive character, $\tau \in E^{\times}$ an element of trace $0$ 
and we set $\psi : E \to \CC^{\times}$ to be $\psi(x) = \psi'(\Tr(\tau x))$. 
We use $\psi'$ and $\psi$ to define 
autodual Haar measures on $F$ and $E$ respectively. 
The duality $F \times E/F \to \CC^{\times}$ given by $(x,y) \mapsto \psi(xy)$ defines 
a unique Haar measure on $E/F$ dual to the one on $F$. 
This measure on $E/F$ coincides with the quotient measure. 
\epar

\bpar
We employ the convention of Section \ref{S:Haar-norm} to define Haar measures 
(with $\psi$ denoted here as $\psi'$). 
Let $k = E$ or $F$. 
We define the following measures on $GL_n (k)$ and its subgroups
\begin{itemize}
\item On $GL_n (k)$ we set
\[
dx = \dfrac{dx_{ij}}{|\det x|_{k}^{n}}
\]
where $x = (x_{ij})$.
\item On standard Levi subgroups of $GL_n (k)$ we set the product measure using the measure defined above. 
\item On (semi) standard unipotent subgroups $N(k) \subset GL_n (k)$ we set the additive measure $dn_{ij}$
where $n_{ij}$ run through coordinates of $N$. 
\item If $P(k)$ is a standard Levi subgroup of $GL_n (k)$ 
with the standard Levi decomposition $N(k)M(k)$ we have the right-invariant measure $dp := dndm$ on $P(k)$ 
and the left invariant measure $\delta_{P_{k}}^{-1} dp$ 
where $\delta_{P_{k}} : P(k) \to \RR_{>0}^{\times}$ 
is the Jacobian homomorphism for the adjoint action of $P(k)$ on $N(k)$. 
\end{itemize}

With this normalization, we have for all $f \in C_{c}^{\infty}(GL_{n}(k))$
\[
\int_{GL_{n}(k)}f\,dg =  \int_{P(k)}\int_{\brN(k)}f(pn)\delta_{P_{k}}(p)^{-1}\,dpdn.
\]
where $\brN$ is the unipotent radical of the opposite parabolic to $P$.
\epar

\bpar
We will use the notation introduced in Section \ref{sec:sym} with some changes. 
All groups considered in this section are subgroups of $G_{n} = \Res_{E/F}GL_n$. 
We write simply $G$ for $G_{n}$, $P_{0}$ for the fixed minimal parabolic subgroup of $G$ 
and $N_0$ for its unipotent radical.
In order to be as compatible with Appendix A of \cite{FLO} as possible, 
instead of $G' = GL_n$ (defined over $F$) 
we write $G_{F} = GL_n$ and for any subgroup $H$ of $G$ we write $H_{F}$ for 
$H \cap G_{F}$. 
We will often identify a group with its $F$ points in this section. 
\epar

\bpar\label{ssec:unipChar}
We define the character $\psi : N_{0} \to \tC$ as follows. 
Write $n \in N_{0}$ as 
\[
n = \vtttt{1 & n_{12} & n_{13} & \hdots & n_{1n}}
			{0 & 1 &  n_{23} & \hdots & n_{2n}}
			{0 & \ddots &  \ddots  & \ddots & n_{2n}}
			{0 & \ddots &  \ddots & 1 & n_{n-1n}}			
			{0 & \hdots &  0 & 0 & 1}, \quad n_{ij} \in E
\]
and set $\psi(n) = \psi ((-1)^{n}(n_{12} + n_{23} + \cdots + n_{n-1n}))$. 
This is the same character as the one from \ref{ssect:additive and generic characters}. 
By restriction, $\psi$ defines a character of $N_0 \cap M$ for all standard Levi subgroups $M$.
\epar

\bpar\label{ssec:mirab}
We denote by $\pc = \pc_{n}$ the mirabolic subgroup of $G$ 
defined as the stabilizer of the row vector $\wttt{0}{\hdots}{0}{1}$. 

Define the following functional on $C^{\infty} (N_{0} \bsl \pc, \psi) = \{f \in C^{\infty}(\pc) \ | \ f(nx) = \psi(n)f(x), \ n \in N_0, \ x \in \pc\}$  
\[
\beta(\varphi) = \beta_{G} (\varphi) = \int_{N_{0,F} \bsl \pc_{F}} \varphi (p)\,dp.
\]
Note that the integral is well defined as $\psi$ is trivial on $N_{0,F}$. In the same way, we define $\beta_{M}$ 
for all standard Levi subgroups $M$ of $G$.
\epar

\bpar\label{ssec:jacqInt}
Let $\Pi_{gen}(G)$ be the set of irreducible generic complex representations of $G = GL_{n}(E)$. 
Let $\wc(\pi) = \wc^{\psi}(\pi)$ be the space of the Whittaker model of $\pi \in \Pi_{gen}(G)$ with respect to the character $\psi$. 
Let $\delta_{g}^{\pi} = \delta_{g} : \wc(\pi) \to \CC$ be the evaluation at $g \in G$. 
The group $G$ acts on $\wc(\pi)$ by right multiplication. 

Fix $P = MN \in \fc^{G}(P_{0})$. Let $w_{M}$ be the element in the Weyl group of $G$ such that 
$w_{M}Mw_{M}^{-1}$ is a standard Levi and the longest for this property.
Let $P^{w}  = N^{w}M^{w} \in \fc^{G}(P_{0})$ be the group whose Levi component is $M^{w} = w_{M}Mw_{M}^{-1}$. 

For $\sigma \in \Pi_{gen}(M)$ let $\Ind_{P}^{G} (\wc(\sigma))$ be the normalized (smooth) induction to $G$, from 
$\wc(\sigma)$, seen as a representation of $P$ via the natural map $P \to M$.
For $\varphi \in \Ind_{P}^{G} (\wc(\sigma))$ let
\[
\bW(g, \varphi) = \int_{N^{w}}\delta_{e}^{\sigma} (\varphi (w_{M}^{-1} u' g))\psi^{-1}(u')\,du'.
\]
This is the so called Jacquet's integral. 
We have then that $\bW_{e}(\varphi) := \bW(e, \varphi)$ is a Whittaker functional on $ \Ind_{P}^{G} (\wc(\sigma))$.
\epar

\bpar
For $\sigma \in \Pi_{gen}(M)$ and $\varphi \in \Ind_{P}^{G} (\wc(\sigma))$ 
let 
\[
\beta'(\varphi) = 
\int_{P_{F} \bsl G_{F}} \beta_{M} (\varphi (g))\,dg.
\]

\begin{theoreme}\label{thm:comp} Let $\sigma \in \Pi_{gen}(M)$. Suppose $\sigma$ is unitary. 
We have then
\[
\beta'(\varphi) = \beta_{G} (\bW(\varphi)).
\]
\end{theoreme}

\begin{preuve} 
We follow very closely Appendix A of \cite{FLO}.

We reduce the proof to the case $P$ is maximal. Let $Q = LV \supset P$ be maximal 
and suppose the assertion holds for $M = M_{P}$. 
Then
\[
\beta'(\varphi) =  \int_{P_{F} \bsl G_{F}} \beta_{M} (\varphi (g))\,dg = 
 \int_{Q_{F} \bsl G_{F}} \int_{P_{F} \bsl Q_{F}} \delta_{Q_{F}}^{-1}(q) \beta_{M} (\varphi (qg))\, dqdg.
\]
The inner integral on the RHS by induction hypothesis equals
\[
\beta_{L} (g.\bW^{L}(\varphi)).
\]
If we let $\varphi'(g)= g. \bW^{L}(\varphi) \in \wc(\Ind_{L\cap P}^{L}(\wc(\sigma)))$
then $\varphi' \in \Ind_{Q}^{G}( \wc(\Ind_{L \cap P}^{L}(\wc(\sigma))) )$
and so by assumption and transitivity of Jacquet's integral we obtain
\[
 \int_{Q_{F} \bsl G_{F}} \beta_{L} (g.\bW^{L}(\varphi)) \,dg =  \beta_{G} (\bW(\varphi)).
\]

Assume then that $P = MN$ is maximal of type $(n_1, n_2)$. In \cite{FLO}, the authors use $U$ instead of $N$. We will consequently 
use $N$ in place of $U$ here. 
Write $M = M_{1} \times M_{2}$ with $M_{i} \cong \Res_{E/F}GL_{n_{i}}$, $M_{1}$ being in the upper 
and $M_{2}$ in the lower diagonal.

Let
\[
w = w_{M}^{-1} = \matx{0}{I_{n_1}}{I_{n_2}}{0}.
\]
Let $P' = M'N'$ be of type $(n_2, n_1)$ 
so that $M'= M_{2}' \times M_{1}'$ with $M_{i}'\cong \Res_{E/F}GL_{n_{i}}$, $M_{2}'$ being in the upper 
and $M_{1}'$ in the lower diagonal. Let $\pc_{i}'$ be the mirabolic subgroup of $M_{i}'$. 
Let $N_{i}'$ be the maximal upper triangular unipotent of $M_{i}'$, similarly without $'$. 
Note that $\pc_{i} = w \pc_{i}' w^{-1}$ etc.

We identify $N'$ with the group of $n_2 \times n_1$ matrices. 
Let 
\[
C_{i} = \{I_{n} + \xi \ | \ \xi \text{ column vector of size }n_2 \text{ in the i-th column}\} \subset N', \quad i=n_2+1, \ldots, n.
\]
Let
\[
R_{i} = \{I_{n} + \xi \ | \ \xi \text{ row vector of size }n_2 \text{ in the i-th row}\} \subset \brN', \quad i=n_2+1, \ldots, n.
\]
We can identify $C_{i}$ and $R_{j}$  with $E^{n}$ which induces 
a pairing between $C_{i}$ and $R_{i-1}$ that we will denote $\bilif_{i}$.

We note some obvious facts
\begin{itemize}
\item The groups $R_{i}$ (resp. $C_{i}$) commute with each other and are normalized by $M_{2}'$ and $M_{1}'$. 
\item The commutator set $[C_{i}, R_{j}]$ is contained in $N_{1}'$ for $j < i$.
\end{itemize}
We define the following groups
\begin{enumerate}
\item $X_{i} = C_{i+1} \cdots C_{n}$.
\item $Y_{i} = R_{n_2 +1} \cdots R_{i-1}$.
\item $V_{i} = N_{1}'X_{i}Y_{i}$. This is a unipotent group.
\item $V_{i}'= N_{1}'X_{i-1}Y_{i} \supset V_{i}$. This is a unipotent group.
\item 
\[
S_{i} = 
\begin{cases} 
M_{2}'V_{i}, \quad & i > n_2, \\
\pc_{2}'N_{1}'N', \quad & i = n_2
\end{cases}
\]
\item $S'_{i} = M_{2}'V_{i}'$ for $i > n_2$.
\end{enumerate}

Note that
\begin{itemize}
\item $S_{i}' = C_{i}S_{i}$ for $i > n_2$ 
as well as $S_{i}' = R_{i-1}S_{i-1}$ for $i > n_2 + 1$.  
\item Let $\delta_{i}$ and $\delta_{i}'$ be modular characters of $S_{i}$ and $S_{i}'$ respectively. 
It follows that $\delta_{i}'|S_{i} = |\det |_{E}\delta_{i}$ 
and $\delta_{i}'|S_{i-1} = |\det |_{E}^{-1}\delta_{i-1}$ in the above range.
\item We have $\delta_{i} = |\det |_{E}^{n + n_2  -2i +1}$ for $i \ge  n_2$. 
\end{itemize}

Let $\sigma = \sigma_1 \otimes \sigma_2$ be an irreducible representation of $M$, 
with $\sigma_i \in \Pi_{gen}(GL_{n_{i}}(E))$. 
We view $\sigma_2$ as a representation of $M_{2}'$ as well. 
Let us define
\[
\begin{cases}
\Ac_{i} = \Ind_{S_{i}}^{\pc_{n}} (\wc(\sigma) \otimes \psi_i), \quad & i=n_2+1, \ldots, n. \\
\Ac_{i} = \Ind_{N_{0}}^{\pc_{n}} \psi, \quad & i=n_2.
\end{cases}
\]
Here, $\psi_i$ is the character of $V_{i}$ - the unipotent radical of $S_{i}$ - 
whose restriction to $X_{i}Y_{i}$ is trivial and that coincides with $\psi$ on $N_{1}'$. 

Explicitly, for $i > n_2$ we 
have
\[
\Ac_{i} = \{
\varphi : \pc \to \wc(\sigma_2) \ | \ 
\varphi (mvg) = 
\dsl \dfrac{\delta_{i}(m)}{| \det m | } 
\rb^{1/2} \psi_{i}(v) \sigma (m)\varphi (g), \quad g \in \pc, \ 
m \in M_{2}', \ v \in V_{i}.
\}
\]
We also denote $\Ac_{i}^{2}$ the $L^{2}$-induction version of the above as in \cite{FLO}.
Note that
\begin{itemize}
\item $\Ac_{i} = \Ind_{S_{i}'}^{\pc} ( \Ind_{S_{i}}^{S_{i}'} ((\wc(\sigma) \otimes \psi_i)) )$ for $i > n_2$.
\item $\Ac_{i} = \Ind_{S_{i}'}^{\pc} ( \Ind_{S_{i-1}}^{S_{i}'} ((\wc(\sigma) \otimes \psi_{i-1})) )$ for $i > n_2 + 1$.
\item $\Ac_{n_2} = \Ind_{S_{n_2 +1}'}^{\pc} ( \Ind_{N_{0}}^{S_{n_2 +1}'} \psi)$.
\end{itemize}

For any $i > n_2$ the restriction map to $C_i$ identifies 
$ \Ind_{S_{i}}^{S_{i}'} ((\wc(\sigma) \otimes \psi_i))$ with $C^{\infty}(C_{i}, \wc(\sigma_2))$ because 
$S_{i}/ S_{i-1} = C_{i}$. 
Let us denote $\varphi \mapsto \varphi|_{C_{i}}$ the restriction map and 
$\iota_{C_{i}}$ the map in the reverse order. 
Similarly, restriction to $R_{i-1}$ identifies $ \Ind_{S_{i-1}}^{S_{i}'} ((\wc(\sigma) \otimes \psi_{i-1}))$ 
with $C^{\infty}(R_{i-1}, \wc(\sigma_2))$. 
Let us denote $\varphi \mapsto \varphi_{R_{i-1}}$ the restriction map and 
$\iota_{R_{i -1}}$ the map in the reverse order. 

Given that $C_{i}$ and $R_{i}$ are in duality we have a Fourier transform 
\[
\fc_{i}' : L^{2}(C_{i}, \overline{ \wc(\sigma_2)}) \to L^{2}(R_{i-1}, \overline{ \wc(\sigma_2)}) 
\]
where $\overline{\wc(\sigma_2)}$ is the $L^{2}$ completion of $\wc(\sigma_2)$. 

\begin{lemme}  For $i = n, \ldots, n_{2}+2$, the above Fourier transform induces a map
\[
\bc_{i} : \Ac^2_{i} = \Ind_{S_{i}'}^{\pc} ( \Ind_{S_{i}}^{S_{i}'} ((\wc(\sigma) \otimes \psi_i)) ) \to \Ac_{i-1}^{2} = 
 \Ind_{S_{i}'}^{\pc} ( \Ind_{S_{i-1}}^{S_{i}'} ((\wc(\sigma) \otimes \psi_i)) ) 
\]
induced from the equivalence $\Ind_{S_{i}}^{S_{i}'} ((\wc(\sigma) \otimes \psi_i)) \to  \Ind_{S_{i-1}}^{S_{i}'} ((\wc(\sigma) \otimes \psi_{i-1}))$
given by $\varphi \mapsto \iota_{R_{i-1}}(\fc_{i}'(\varphi|_{C_{i}}))$. It is an equivalence of unitary representations.  
\end{lemme}

Similarly, we have the map
\[
\fc_{n_2 + 1} : \Ind_{S_{n_2+1}}^{S_{n_2+1}'} ((\wc(\sigma) \otimes \psi_{n_2+1})) \to 
\Ind_{N_{0}}^{S_{n_2 +1}'} \psi
\]
given by
\[
\fc_{n_2 +1}\varphi (vm) = \psi (v)|\det m|_{E}^{1/2}\hat \varphi (\chi_m)(m), \quad m \in M_{2}', \ v \in V_{n_2+1}' = N_{1}'N'
\]
where 
\[
\chi_{m} : C_{n_2+1} \to \tC, \quad \chi_{m}(c) = \psi(mcm^{-1}), \quad \hat \varphi (\chi) = 
\int_{C_{n_2+1}}\varphi(c)\chi(c)\,dc.
\]

\begin{lemme}  The above Fourier transform induces the equivalence of unitary representations
\[
\bc_{n_2 +1} : \Ac_{n_2+1}^{2} 
=  \Ind_{S_{n_2 +1}'}^{\pc} ( \Ind_{S_{n_2+1}}^{S_{n_2+1}'} ((\wc(\sigma) \otimes \psi_{n_2+1})) )
\to \Ac_{n_2}^{2} = \Ind_{S_{n_2 +1}'}^{\pc}(\Ind_{N_{0}}^{S_{n_2 +1}'}).
\]  
\end{lemme}

For $i = n, \ldots, n_2 + 1$ 
let $\beta_{i} : \Ac_{i} \to \CC$ be the following functional
\[
\beta_{i}(\varphi) = \int_{S_{i, F} \bsl \pc_{F}} \beta_{M_{2}'} (\varphi(p))\,dp.
\]
We also set $\beta_{n_2} = \beta_{G}$ on $\Ac_{n_2} = \Ind_{N_{0}}^{\pc} \psi$.

\begin{lemme}\label{lem:flo1}  For $i=n, \ldots, n_2+2$ and $\varphi \in \Ac_{i}$
have
\[
\beta_{i}(\varphi) = 
\beta_{i-1}(\bc_{i}(\varphi)).
\]
\end{lemme}

\begin{preuve} 
By equivariance property of $\bc_{i}$ it is enough to show
the equality  
between $\int_{C_{i,F}} \beta_{M_{2}'} (\varphi( c))\,dc$
and 
$\int_{R_{i-1,F}} \beta_{M_{2}'} (\bc_{i}(\varphi)( r))\,dr$.
The duality $\bilif_{i}$ between $C_{i}$ and $R_{i-1}$ restricts 
to a duality between $C_{i}/C_{i,F}$ and $R_{i-1,F}$.
We have thus
\begin{multline*}
\int_{R_{i-1,F}} \beta_{M_{2}'} (\bc_{i}(\varphi)( r))\,dr = 
\int_{R_{i-1,F}} \beta_{M_{2}'} (
\int_{C_{i}}
\varphi (c ) \psi (\langle c, r  \rangle_{i}) dc) \,dr = \\
\beta_{M_{2}'}  \dsl 
\int_{R_{i-1,F}} 
\int_{C_{i}/C_{i, F}}
\dsl 
\int_{C_{i,F}} 
\varphi (c^{-} + c  ) \psi (\langle c^{-}, r  \rangle_{i}) dc \rb dc^{-}dr  \rb    = 
\int_{C_{i,F}}
\beta_{M_{2}'}  \dsl 
\varphi (c) \rb dc. 
\end{multline*}
\end{preuve}

\begin{lemme}\label{lem:flo2}  For $\varphi \in \Ac_{n_2 + 1}$ we 
have
\[
\beta_{n_2 + 1}(\varphi) = 
\beta_{n_2}(\bc_{n_2+1}(\varphi)).
\]
\end{lemme}

\begin{preuve} 
Again, it is enough to show the equality
between $\int_{C_{n_2 + 1, F}} \beta_{M_{2}'} (\varphi(c))\,dc$
and 
$\int_{N_{2,F}' \bsl M_{2,F}'} 
\fc_{n_2+1}(\varphi)( m ) \,dm$.
We have
\begin{multline*}
\int_{N_{2,F}' \bsl M_{2,F}'} 
\fc_{n_2}(\varphi)(m ) \,dm = 
\int_{N_{2,F}' \bsl M_{2,F}'} 
\int_{C_{n_2 +1}}
\varphi(c)(m)\psi( m c m ^{-1})dc
|\det m |_{F}\,dm = \\
\int_{\pc_{n_2,F}' \bsl M_{2,F}'} 
\int_{C_{n_2+1}/ C_{n_2+1,F}}
\int_{C_{n_2+1,F}}
\int_{N_{2,F}' \bsl \pc_{n_2,F}'}
\varphi(c + c^{-})(pm)\psi (mc^{-}m^{-1})|\det m |_{F} \,dp dc dc^{-} dm  = \\
\int_{C_{n_2+1,F}}
\int_{N_{2,F}' \bsl \pc_{n_2,F}'}
\varphi(c)(p) \,dp dc = 
\int_{C_{n_2+1,F}} \beta_{M_{2}'} (\varphi(c))\,dc.
\end{multline*}
\end{preuve} 

Let $\varphi \in \Ind_{P}^{G} (\wc(\sigma_1 \otimes \sigma_2))$. 
For $m \in M_{1}$, 
let $\delta_{m}^{1} : \wc(\sigma_1 \otimes \sigma_2) \to \wc(\sigma_2)$ be the evaluation map in the first variable. 
Define for $p \in \pc$
\[
\varphi_{n}(p) = \delta_{e}^{1}\varphi (wp) \in \Ac_{n}.
\]
We have then
\[
\int_{P \bsl G} \| \varphi (g)\|^{2}_{L^{2}(\wc(\sigma_1 \otimes \sigma_2))}\,dg = 
\|\varphi_{n}\|^{2}_{\Ac_{n}}.
\]

We set $\varphi_{i-1} = \bc_{i}(\varphi_{i})$ for $i = n, \ldots, n_2 +1$. 
As shown at the end of the Appendix A.3 of \cite{FLO}, we have
\begin{equation}\label{eq:flo}
\varphi_{n_2} = \bW(\varphi). 
\end{equation}

\begin{lemme}\label{lem:flo3} For $\varphi \in \Ind_{P}^{G} (\wc(\sigma_1 \otimes \sigma_2))$, we have
\[
\beta'(\varphi) =  \beta_{n}(\varphi_{n}).
\]
\end{lemme}
\begin{preuve}
Indeed
\begin{multline*}
\beta'(\varphi) = \int_{P_{F} \bsl G_{F}} \beta_{M}(\varphi(g))\,dg = 
\int_{P_{F} \bsl G_{F}} \beta_{M}(\varphi(g w ))\,dg = 
\int_{\overline{N}_{F}} \beta_{M}(\varphi(uw))\,du = 
\int_{N_{F}'} \beta_{M}(\varphi(wu'))\,du' = \\ 
\int_{N_{F}'} \int_{N_{1,F} \bsl \pc_{1,F}}
\beta_{M_2}( \delta^{1}_{ m_1} (\varphi(wu')) )\, dm_1 du' = 
\int_{N_{F}'} \int_{N_{1,F} \bsl \pc_{1,F}}
\beta_{M_2}( \delta^{1}_{e} (\varphi(m_1 wu' )) ) \delta_{P}^{1/2}(m_1)   \, dm_1 du' = \\
\int_{N_{F}'} \int_{N_{1,F}' \bsl \pc_{1,F}'}
\beta_{M_2'}( \delta^{1}_{e} (\varphi(w m_1' u' )) ) \delta_{P'}^{-1}( m_1')   \, dm_1' du' = 
\int_{S_{n,F} \bsl \pc_{F}}
\beta_{M_2'}( \delta^{1}_{e} (\varphi(w p )) )  \, dp.
\end{multline*}
\end{preuve}

Now combining the equality \eqref{eq:flo} with Lemmas \ref{lem:flo1}, \ref{lem:flo2}, \ref{lem:flo3} 
we obtain the desired equality at least when computations in these Lemmas are justified. 
Taking $\varphi \in \Ind_{P}^{G} (\wc(\sigma_1 \otimes \sigma_2))$ supported on the big cell $PwP'$ 
we can see that all integrals are absolutely convergent. 
By multiplicity (at most) one \cite{Flicker2, AizGour}, the Theorem \ref{thm:comp} follows. 
 \end{preuve}
 \epar

 \subsection{Global comparison}
 
 \bpar 
 We go back to the global setting and notation introduced in \S \ref{ssec:sym-not}. 
 \epar 
 
 \bpar
  We normalize all local and global measures as in \S \ref{S:Haar-measures}, with respect to a fixed character $\psi' : F \backslash \AAA \to \CC^{\times}$. 
  We have the quadratic character $\eta : F^{\times} \bsl \A^{\times} \to \CC^{\times}$ 
  associated to $E/F$ and the associated character $\eta_{G'}$ of $G'(\A)$ as defined in Paragraph \ref{S:character etaG'}. 
\epar

\bpar 
As in \S \ref{ssect:additive and generic characters}, we also
fix a non-trivial additive character $\psi : E \backslash \AAA_{E} \to \CC^{\times}$, trivial on $\A$ 
which is then used to define a non-degenerate 
character $\psi_{N}$ of the maximal unipotent subgroup of $G(\A)$ as in the beginning of \ref{Chap:generic-JRTF}.
\epar

\bpar 
Let $\chi \in \Xgo^*(G)$ (c.f. \S \ref{S:relevant}) and let $(M, \pi)$ represent $\chi$. 
Set $\Pi = \Ind_{P(\A)}^G(\A)(\pi)$. 
\epar

\bpar[The comparison] --- 

\begin{theoreme}\label{thm:comparisonG}
For all $\phi \in \Pi$ we have
\[
J_{\eta}(\phi) = \beta_{\eta}(W_{\phi})
\]
where
\begin{itemize}
\item $J_{\eta}$ is defined in \ref{S:Jeta};
\item $\beta_{\eta}$ is defined in \ref{S:Whittaker periods}.
\item $W_{\phi} \in \mathcal{W}(\Pi, \psi_N)$ is defined in \ref{8.1.1?}.
\end{itemize}
\end{theoreme}

\begin{preuve}
The proof is essentially the same as of Theorem \ref{theo inner products}. 
The only difference is that the natural analogue of \eqref{eq1 JS period} 
is provided by Proposition 3.2 of \cite{Zhang2} 
and the analogue of \eqref{eq2 JS period} is established invoking Theorem \ref{thm:comp}. 
\end{preuve}

\begin{corollaire}\label{cor:eqdesI}
We have the equality of distributions on $\Sc(G(\AAA))$
\[
I_{P, \pi} = I_{\Pi}
\]
where
\begin{enumerate}
\item $I_{P, \pi}$ is defined in \S \ref{S:RelcharJPpi}.
\item $I_{\Pi}$ is defined in \S \ref{S:IPI}.
\end{enumerate}
\end{corollaire}

\begin{preuve}
Looking at definitions of $I_{P, \pi}$ and $I_{\Pi}$, taking into consideration Theorem \ref{theo inner products} 
and Theorem \ref{thm:comparisonG} above, we see that we need to establish 
for all $\phi \in \Pi$
\[
\la(W_{\phi}) = I(\phi, 0)
\]
where $\la = Z^{RS}(0, \cdot )$ is defined in \S \ref{S:Whittaker periods}
and $I(\phi, 0)$ is given by Proposition \ref{prop:RS-integ-generic}. 
This equality is precisely Theorem 1.1 of \cite{IY}. 
\end{preuve}
\epar


\section{Proofs of the Gan-Gross-Prasad and Ichino-Ikeda conjectures}\label{Chap:proofmaintheorems}

\subsection{Identities among some global relative characters}

\begin{paragr}
Besides notation of Chapters \ref{Chapter:preliminaries} and \ref{Chap:Ichilimit}, we shall use notation of Section \ref{sec:Introduction}. We fix an integer $n\geq 1$ and we will omit the subscript $n$: we will write $\hc$ for $\hc_n$. 
\end{paragr}

\begin{paragr}[Relative characters for unitary groups.] ---
Let $h\in \hc$ be a Hermitian form. Let $\sigma$ be an irreducible cuspidal automorphic subrepresentation of the group $U_h$. We define the relative character $J_\sigma^h$ by
    \begin{align*}
      J_\sigma^h(f)=\sum_{\varphi} \pc_h(\pi(f)\varphi) \overline{\pc_h(\varphi)}, \text{   } \forall f\in \Sc(U_h(\AAA))
    \end{align*}
where $\varphi$ runs over a $K_h$-basis (see \ref{S:K0ON}) for some maximal compact subgroup $K_h\subset U_h(\AAA)$. The period $\pc_h$ are those defined in \ref{S-intro:GGP}. 
For any subset $\Xgo_0\subset \Xgo(U_h)$ of cuspidal data which do not come from proper Levi subgroups (that is they are represented by  pairs $(U_h,\tau)$ where $\tau$ is a cuspidal automorphic representation) we define more generally
\begin{align}
  J_{\Xgo_0}^h(f)=\sum_{\chi\in \Xgo_0}\sum_{\sigma} J_\sigma^h(f)
\end{align}
where the inner sum is over the set of the constituents $\sigma$ of some decomposition of $L^2_{\chi}([U_h])$ (see §\ref{eq:langlands}) into irreducible subrepresentations. One can show that the double sum is absolutely convergent (see e.g. \cite[Proposition A.1.2]{RBP}).
\end{paragr}

\begin{paragr}\label{S:finitesetS0}  Let $V_{F,\infty}\subset S_0\subset V_F$ be a finite set of places containing all the places that are ramified in $E$.  For every $v\in V_F$, we set $E_v=E\otimes_F F_v$ and when $v\notin V_{F,\infty}$ we denote by $\oc_{E_v}\subset E_v$ its ring of integers. Let $\hc^\circ\subset \hc$ be the (finite) subset of Hermitian spaces of rank $n$ over $E$ that admits a selfdual $\oc_{E_v}$-lattice for every $v\notin S_0$.

	For each $h\in \hc^\circ$, the group $U_h$ is naturally defined over $\oc_F^{S_0}$ and we fix a choice of such a model. Since we are going to consider invariant distribution, this choice is irrelevant. We define the open compact subgroups $K_h^\circ=\prod_{v\notin S_0} U_h(\oc_v)$ and $K^\circ=\prod_{v\notin S_0} G(\oc_v)$ respectively of $U_h(\AAA^{S_0})$ and $G(\AAA^{S_0})$.

	Let $v\notin S_0$. We denote by $\Sc^\circ(U_h(F_v))$, resp.  $\Sc^{\circ}(G(F_v))$, the spherical Hecke algebra\footnote{The product structure is given by the convolution where the Haar measure is normalized so that the characteristic functions of $U_h(\oc_v)$ and $G(\oc_v)$ are units.} of  complex functions on $U_h(F_v)$ (resp. $G(F_v)$) that are $U_h(\oc_{v})$-bi-invariant (resp.  $G(\oc_{v})$-bi-invariant) and compactly supported. 

	We have the base change homomorphism 
	\begin{align*}
	BC_{h,v}: \Sc^{\circ}(G(F_v))\to \Sc^{\circ}(U_h(F_v)).
	\end{align*}
	We denote by $\Sc^{\circ}(U_h(\AAA^{S_0}))$, resp.  $\Sc^{\circ}(G(\AAA^{S_0}))$, the restricted tensor product of $\Sc^\circ(U_h(F_v))$, resp. $\Sc^{\circ}(G(F_v))$, for $v\notin S_0$. We have also a global base change homomorphism given by $BC_h^{S_0}=\otimes_{v\notin S_0}   BC_{h,v}$.
	
	We also denote by $ \Sc^\circ (G(\AAA))\subset  \Sc (G(\AAA))$ and $\Sc^\circ(U_h(\AAA))\subset \Sc(U_h(\AAA))$, for $h\in \hc^\circ$, the subspaces of functions that are respectively bi-$K^\circ$-invariant and bi-$K_h^\circ$-invariant.

\end{paragr}

\begin{paragr}[Transfer.] --- Let $h\in  \hc^\circ$. We shall say that $f_{S_0}\in  \Sc(G(F_{S_0}))$  and  $f_{S_0}^h \in \Sc(U_h(F_{S_0}))$ are transfers  if the functions $f_{S_0}$ and $f_{S_0}^h$ have matching regular orbital integrals in the sense of Definition 4.4 of \cite{BPLZZ}. The Haar measures on the $F_{S_0}$-points of the involved groups are those defined in §\ref{S:dgS}.

\end{paragr}

\begin{paragr} \label{S:lesdonnees}Let $P$ be a standard parabolic subgroup of $G$ and $\pi$ be a cuspidal automorphic representation of $M_P$. Let $\chi\in \Xgo(G)$ be the class of the pair $(M_P,\pi)$. We assume henceforth that $\chi$ is a generic relevant cuspidal datum in the sense of \S \ref{S:relevant}.
	
Set $\Pi=\Ind_P^G(\pi)$ for the corresponding parabolically induced representation. The assumption that $\chi$ is generic and relevant means exactly that $\Pi$ is a Hermitian Arthur parameter (see \S \ref{S:Aparam}). Moreover, we assume, as we may, that $S_0$ has been chosen such that $\Pi$ admits $K^\circ$-fixed vectors. 
	
	Attached to these data, we have three distributions denoted by $I_\chi$, $I_{P,\pi}$ and $I_\Pi$. The first is constructed as a contribution of the Jacquet-Rallis trace formula and it is defined in Theorem \ref{thm:jfDef}. The second and third  are relative characters built respectively in §\ref{S:RelcharJPpi} and §\ref{S:IPI}. The bulk of the paper was devoted to the proof of the following identities (see Theorem \ref{thm:intLaT}, Theorem \ref{theo JR TF} and Corollary \ref{cor:eqdesI})
	\begin{equation}\label{eq:spectral expansion Ichi}
	I_\chi=2^{-\dim(\ago_P)}I_{P,\pi}=2^{-\dim(\ago_P)}I_\Pi.
	\end{equation} 
\end{paragr}

\begin{paragr}\label{S:lacomparaison} Let $S'_0$ be the union of $S_0\setminus V_{F,\infty}$ and the set of all finite places of $F$ that are inert in $E$.

	We define $\Xgo_0^h\subset \Xgo(U_h)$ as the set of equivalence classes of pairs $(U_h,\sigma)$ where $\sigma$ a cuspidal automorphic representation of $U_h(\AAA)$ that satisfies the following conditions:
	\begin{itemize}
		\item $\sigma$ is $K_h^\circ$-unramified;
		\item for all $v\notin S_0'\cup V_{F,\infty}$ the (split) base change of $\sigma_v$ is $\Pi_v$.
	\end{itemize}

	\begin{proposition}\label{thm:comparison}
		
		Let $f\in \Sc^\circ(G(\AAA))$ and $f^h\in \Sc^\circ(U_h(\AAA))$ for every $h\in \hc^\circ$. Assume that the following properties are satisfied for every $h\in \hc^\circ$:
		\begin{enumerate}
			\item $f=(\Delta_H^{S_0,*} \Delta_{G'}^{S_0,*} ) f_{S_0}\otimes f^{S_0}$ with $f_{S_0}\in \Sc(G(F_{S_0}))$ and $f^{S_0}\in  \Sc^{\circ}(G(\AAA^{S_0}))$.
			\item $f^h= (\Delta_{U_h'}^{S_0})^2f^h_{S_0}\otimes f^{h,S_0}$ with $f_{S_0}^h \in \Sc(U_h(F_{S_0}))$ and $f^{h,S_0} \in  \Sc^{\circ}(U_h(\AAA^{S_0}))$. 
			\item The functions $f_{S_0}$ and $f_{S_0}^h$ are transfers.
			\item $f^{h,S_0}=BC^{S_0}_h(f^{S_0})$
			\item The function $f^{S_0}$ is a product of a smooth compactly supported function on the restricted product $\prod_{v\notin S_0'}' G(F_v)$ by the characteristic function of $\prod_{v\in S_0'\setminus S_0} G(\oc_v)$.
		\end{enumerate}
		
		Then we have:
		\begin{align}\label{eq:comparison}
		\sum_{h\in\hc^\circ}J_{\Xgo^h_0}^h(f^h)=2^{-\dim(\ago_P)}  I_{\Pi}(f)=2^{-\dim(\ago_P)}  I_{P,\pi}(f).
		\end{align}
	\end{proposition}

	\begin{remarque}
		If the assumptions hold for the set $S_0$, it also holds for any large enough finite set containing $S_0$: this follows from the Jacquet-Rallis fundamental lemma (see \cite{Yun} and  \cite{BPLF}) and the simple expression of  the transfer at split places (see \cite[proposition 2.5]{Z1}). We leave it to the reader to keep track of the different choices of Haar measures in these references.
	\end{remarque}
	
	\begin{preuve}
		
		The proof follows the same lines as the proof of \cite[Theorem 1.7]{BPLZZ}. For the convenience of the reader, we recall the main steps.
		
		In Theorem \ref{thm:jfDef} we defined a distribution $I$ on $\Sc(G(\AAA))$: this is the ``Jacquet-Rallis trace formula'' for $G$. We have an analogous distribution  $J^h$ on  $\Sc(U_h(\AAA))$ for each $h\in \hc$: it is defined in \cite[théorème 0.3]{Z3} for compactly supported functions and extended to the Schwartz space in \cite[§1.1.3 and théorème 15.2.3.1]{CZ}.  Note that, by the Jacquet-Rallis fundamental lemma \cite{Yun}, \cite{BPLF}, for every $h\in \hc\setminus \hc^\circ$ there exists a place $v\in S_0'\setminus S_0$ such that the characteristic function $\mathbf{1}_{G(\oc_v)}$ admits the zero function on $U_h(F_v)$ as a transfer. Therefore, by \cite[théorème 1.6.1.1]{CZ},  the hypotheses of the proposition imply:
		\begin{align}\label{thm:CZ}
		I_{}(f)=\sum_{h\in\hc^\circ}  J_{}^h(f^h).
		\end{align}
		
		We will denote by $\mc^{S_0'}(G(\AAA))$, resp.  $\mc^{S_0'}(U_h(\AAA))$, the algebra of $S_0'$-multipliers defined in \cite[definition 3.5]{BPLZZ} relatively to the subgroup $\prod_{v\notin S_0'} G(\oc_v)$, resp. $\prod_{v\notin S_0'} U_h(\oc_v)$. Any multiplier $\mu\in \mc^{S_0'}(G(\AAA))$, resp. $\mu\in \mc^{S_0'}(U_h(\AAA))$, gives rise to a linear operator $\mu \ast$ of the algebra  $\Sc^\circ(G(\AAA))$, resp. $\Sc^\circ(U_h(\AAA))$ and for every admissible irreducible representation $\pi$ of $G(\AAA)$, resp. of $U_h(\AAA)$, there exists a constant $\mu(\pi)\in \CC$ such that $\pi(\mu \ast f)=\mu(\pi) \pi(f)$ for all $f\in \Sc^\circ(G(\AAA))$, resp. $f\in \Sc^\circ(U_h(\AAA))$.
		
		Let $\xi_\Pi$ be the infinitesimal character of $\Pi$. By \cite[Theorem 4.12 (4)]{BPLZZ}, for every $h\in \hc^\circ$ and $(U_h,\sigma)\in \Xgo_0^h$, the base-change of the infinitesimal character of $\sigma$ is $\xi_\Pi$. However, the universal enveloping algebras of the complexified Lie algebras of $U_h$ are all canonically identified for $h\in \hc$ (since these are inner forms of each other) and base-change is injective at the level of infinitesimal characters. As, by \cite{GRS}, there exists at least one $h\in\hc^\circ$ such that the set $\Xgo_0^h$ is nonempty (we may even take for $h$ any quasi-split Hermitian form unramified outside $S_0$), there exists a common infinitesimal character $\xi$ of all $(U_h,\sigma)\in \Xgo_0^h$, for $h\in \hc^\circ$, whose base-change is $\xi_\Pi$.
		
		By the strong multiplicity one theorem of Ramakrishnan (see \cite{Ram2}) and  Theorem 3.17 of \cite{BPLZZ}, one can find a multiplier $\mu\in \mc^{S_0'}(G(\AAA))$ such that
		\begin{enumerate}
			\item[i.] $\mu(\Pi)=1$;
			\item[ii.] For all $\chi'\in \Xgo(G)$ such that  $\chi'\not=\chi$,    we have 
			\begin{align*}
			K_{\mu*f,\chi'}^{0}=0
			\end{align*}
			where the kernel $K_{\mu\ast f,\chi'}^{0}$ is defined as in §\ref{S:Autom-kernel}.
		\end{enumerate}
		
		By Theorem 3.6 and Theorem 4.12 (3) of \cite{BPLZZ}, for every $h\in \hc^\circ$ there exists a multiplier $\mu^h\in \mc^{S_0'}(U_h(\AAA))$ such that
		\begin{enumerate}
			\item[iii.] $\mu^h(\sigma)=1$ for all $(U_h,\sigma)\in \Xgo^h_0$;
			\item[iv.] For all $\chi'\in \Xgo(U_h)$ such that $\chi'\notin \Xgo^h_0$ and for all parabolic subgroups $P$ of $U_h$,  we have
			\begin{align*}
			K^{U_h}_{P,\mu^h\ast f^h,\chi'}=0
			\end{align*}
			where the left-hand side is the kernel of the operator given by the right convolution of $\mu^h \ast f^h$ on $L^2_\chi([U_h]_P)$ (see \eqref{eq:langlands}).
		\end{enumerate}
		
		Moreover, by \cite[Proposition 4.8, Lemma 4.10]{BPLZZ}, we may choose $\mu$ and $\mu^h$ such that the functions $\mu\ast f$ and $\mu^h\ast f^h$, for $h\in \hc^\circ$, still satisfy the assumptions of the proposition. So, in particular, from \eqref{thm:CZ} applied to the functions $\mu\ast f$ and $(\mu^h\ast f^h)_{h\in \hc^\circ}$ instead of $f$ and $(f^h)_{h\in \hc^\circ}$, we get
		\begin{align}\label{thm:CZ2}
		I_{}(\mu\ast f)=\sum_{h\in\hc^\circ}  J_{}^h(\mu^h \ast f^h).
		\end{align}

		Note that by conditions i. and iii. we have:
		\begin{align*}
		I_{\Pi}(\mu\ast f)=I_{\Pi}(f), \; I_{P,\pi}(\mu\ast f)=I_{P,\pi}(f) \mbox{ and } J_{\Xgo_0^h}^h(\mu^h \ast f^h)=J_{\Xgo_0^h}^h(f^h), \mbox{ for every } h\in \hc^\circ. 
		\end{align*}
		Moreover, by ii. Theorem \ref{thm:asym-trio} applied to \eqref{eq:LaTrKchi} and \eqref{eq:spectral expansion Ichi}, we see that the left-hand side of \eqref{thm:CZ2} reduces to $I_{\chi}(\mu\ast f)=2^{-\dim(\ago_P)}I_\Pi(\mu\ast f)=2^{-\dim(\ago_P)}I_{P,\pi}(\mu \ast f)$. On the other hand, by iv. and the very definition of $J_{}^h$ given in \cite{Z3}, the right-hand side of \eqref{thm:CZ2} reduces to
		\begin{align*}
		\sum_{h\in\hc^\circ}J_{\Xgo^h_0}^h(\mu^h\ast f^h)=\sum_{h\in\hc^\circ}J_{\Xgo^h_0}^h(\mu^h \ast f^h).
		\end{align*}
		Therefore, \eqref{thm:CZ2} gives precisely the identity of the proposition.
	\end{preuve}

\end{paragr}

\subsection{Proof of Theorem \ref{thm:GGP}}\label{ssec:pf-GG}

\begin{paragr} Let $\Pi=\Ind_P^G(\pi)$ be a Hermitian Arthur parameter of $G$. Note that by properties 1 and 2 of §\ref{S:Aparam}, the cuspidal datum $\chi$ associated to the pair $(M_P,\pi)$ is generic and relevant in the sense of  §\ref{S:generic}. For $h\in \hc$ and $\sigma$ a cuspidal automorphic representation of $U_h(\AAA)$, it is readily seen that the linear form $\pc_h$ is nonzero on $\sigma$ if and only if $J_\sigma^h$ is not identically zero. On the other hand, the linear form $J_\eta$ or $\beta_\eta$ is always nonzero (this follows either from the fact that $\chi$ is relevant or is an easy consequence of \cite{GelKa}, \cite[Proposition 5]{JacDist} and \cite{Kem}) whereas the linear form $I$, from Proposition \ref{prop:RS-integ-generic}, or $\lambda$, from \S \ref{S:Whittaker periods}, is nonzero if and only if $L(\frac12,\Pi)\not=0$ (as follows either from the work of Ichino and Yamana, see \cite[corollary 5.7]{IY}, or of Jacquet, Piatetski-Shapiro and Shalika \cite{JPSS}, \cite{JacIntegral}). Therefore, we similarly deduce that the distribution $I_{P,\pi}$ or $I_\Pi$ is non-zero if and only if $L(\frac12,\Pi)\not=0$.

	As a consequence, Theorem \ref{thm:GGP} amounts to the equivalence between the two assertions:
	\begin{enumerate}
		\item[(A)] The distribution $I_{P,\pi}$ or $I_\Pi$ is non-zero.
		\item[(B)] There exist $h\in \hc$, $f\in \Sc(U_h(\AAA))$ and a cuspidal subrepresentation $\sigma$ of $U_h$ such that $BC(\sigma)=\Pi$ and $J_\sigma^h(f) \not=0$.
			\end{enumerate}
	
\end{paragr}

\begin{paragr}[Proof of $(A)\Rightarrow (B)$.] --- We choose the $S_0$ of §\ref{S:finitesetS0} such that $I_\Pi$ is not identically zero on $f_1\in \Sc^\circ(G(\AAA))$. Then Assertion (B) above is a consequence of Proposition \ref{thm:comparison}: it suffices to take functions $f$ and $f^h$ for $h\in \hc^\circ$ satisfying the hypotheses of that theorem and such that $I_{\Pi}(f)\not=0$. That it is possible is implied by a combination of a result of \cite{Xue} and the existence of $p$-adic transfer \cite{Z1}.
	
\end{paragr}

\begin{paragr}[Proof of $(B) \Rightarrow (A)$.] --- We may choose the set $S_0$ so that there exist $h_0\in \hc^\circ$, $f^{h_0}_0\in \Sc^\circ(U_{h_0}(\AAA))$ and a cuspidal representation $\sigma_0$ of $U_h$ such that for $v\in S_0'$ (see §\ref{S:lacomparaison}) $BC(\sigma_{0,v})=\Pi_v$ and $J_{\sigma_0}^{h_0}(f^{h_0}_0) \not=0$. For any other $h\in \hc^\circ$ we set $f^h_0=0$. Up to enlarging $S_0$, we may assume that the family $(f^h_0)_{h\in \hc^\circ}$ satisfies conditions 2. and 5. of Proposition \ref{thm:comparison}. Moreover, we have (see \cite[\S 2.5]{Z1}) $J_{\sigma}^{h_0}(f^{h_0}_0\ast f_0^{h_0})\geqslant 0$ for every $\sigma\in \Xgo_0^{h_0}$ and $J_{\sigma_0}^{h_0}(f^{h_0}_0\ast f_0^{h_0})>0$. In particular, the left hand side of \eqref{eq:comparison} for the family $(f^h_0\ast f^h_0)_{h\in \hc^\circ}$ is nonzero. Once again by \cite{Xue} and the existence of $p$-adic transfer \cite{Z1}, this implies that we can find test functions $f\in \Sc^\circ(G(\AAA))$ and $f^h\in \Sc^\circ(U_h(\AAA))$, for $h\in \hc^\circ$, satisfying all the conditions of Proposition \ref{thm:comparison} and such that the left hand side of \eqref{eq:comparison} is still nonzero. The conclusion of this proposition immediately gives Assertion (A).

\end{paragr}

\subsection{Proof of Theorem \ref{thm:II}}

\begin{paragr}
	Let $h\in \hc$ and $\sigma$ be a cuspidal automorphic representation of $U_{h}(\AAA)$ which is tempered everywhere. By \cite{Mok}, \cite{KMSW}, $\sigma$ admits a weak base-change $\Pi$ to $G$. Moreover, by these references $\Pi$ is also a strong base-change of $\sigma$: for every place $v$ of $F$, the local base-change of $\sigma_{v}$ (defined in \cite{Mok} and \cite{KMSW}) coincides with $\Pi_v$. In particular, it follows that $\Pi$ is also tempered everywhere.
	
	We choose a finite set of places $S_0$ as in \S \ref{S:finitesetS0} such that $h\in \hc^\circ$ and $\sigma$ as well as the additive character $\psi'$ used to normalize local Haar measures in Section \ref{S:Haar-measures} are unramified outside of $S_0$.
	
	For each place $v$ of $F$, we define a distribution $J_{\sigma_{v}}$ on $\Sc(U_{h}(F_v))$ by
	$$\displaystyle J_{\sigma_{v}}(f^{h}_v)=\int_{U'_{h}(F_v)} \mathrm{Trace}(\sigma_{v}(h_v)\sigma_v(f^{h}_v)) dh_v,\;\; f^{h}_v\in \Sc(U_{h}(F_v)),$$
	where
	$$\displaystyle \sigma_{v}(f^{h}_v)=\int_{U_{h}(F_v)} f^{h}_v(g_v) \sigma_{v}(g_v)dg_v$$
	and the Haar measures are the one defined in \S \ref{S:dgS}. Moreover by \cite{NHar}, and since the representations $\sigma_v$ are all tempered, the expression defining $J_{\sigma_{v}}$ is absolutely convergent and for every $v\notin S_0$ we have
	$$\displaystyle J_{\sigma_{v}}(\mathbf{1}_{U_{h}(\oc_v)})=\Delta_{U'_{h},v}^{-2} \frac{L(\frac12,\Pi_v)}{L(1,\sigma_{v},\Ad)}.$$
\end{paragr}

\begin{paragr}\label{S:reformulation II}
	By \cite[Lemma 1.7]{Zhang2} and our choice of local Haar measures, Theorem \ref{thm:II} is equivalent to the following assertion: for all factorizable test function $f^{h}\in \Sc(U_{h_0}(\AAA))$ of the form $f^{h}=(\Delta_{U'_{h}}^{S_0})^2\prod_{v\in S_0} f^{h}_v\times \prod_{v\notin S_0} \mathbf{1}_{U_{h}(\oc_v)}$, we have
	\begin{align}\label{eq:reformulation II}
	\displaystyle J^{h}_{\sigma}(f^{h})=\lvert S_\Pi\rvert^{-1} \frac{L^{S_0}(\frac12,\Pi)}{L^{S_0}(1,\sigma,\Ad)}\prod_{v\in S_0} J_{\sigma_{v}}(f_v^{h}).
	\end{align}
\end{paragr}

\begin{paragr}
	For every place $v$ of $F$, we define a local relative character $I_{\Pi_v}$ on $G(F_v)$ by
	$$\displaystyle I_{\Pi_v}(f_v)=\sum_{W_v\in \wc(\Pi_v,\psi_{N,v})} \frac{\lambda_v(\Pi_v(f_v)W_v) \overline{\beta_{\eta,v}(W_v)}}{\langle W_v,W_v\rangle_{\Whitt, v}},\;\; f_v\in \Sc(G(F_v)),$$
	where the sum runs over a $K_v$-basis of the Whittaker model $\wc(\Pi_v,\psi_{N,v})$ (in the sense of \S \ref{S:K0ON}) and $\lambda_v$, $\beta_{\eta,v}$, $\langle .,.\rangle_{\Whitt,v}$ are local analogs of the forms introduced in \S \ref{S:Whittaker periods} given by
	\[\begin{aligned}
	\displaystyle \lambda_v(W_v)=\int_{N_H(F_v)\backslash H(F_v)} W_v(h_v)dh_v & ,\; \beta_{\eta,v}(W_v)=\int_{N'(F_v)\backslash \pc'(F_v)} W_v(p_v) \eta_{G',v}(p_v) dp_v, \\
	\mbox{and } & \langle W_v,W_v\rangle_{\Whitt,v}=\int_{N(F_v)\backslash \pc(F_v)} \lvert W_v(p_v)\rvert^2 dp_v.
	\end{aligned}\]
	Note that the above expressions, and in particular $\lambda_v(W_v)$, are all absolutely convergent due to the fact that $\Pi_v$ is tempered (see \cite[Proposition 8.4]{JPSS}). The above definition also implicitely depends on the choice of an additive character $\psi$ of $\AAA_E/E$ trivial on $\AAA$ (through which the generic character $\psi_N$ is defined, see beginning of Chapter \ref{Chap:generic-JRTF} and \S \ref{ssect:additive and generic characters}) and up to enlarging $S_0$, we may assume that $\psi$ is unramified outside of $S_0$. Then, it follows from the definition of $I_\Pi$ that for every factorizable test function $f\in \Sc(G(\AAA))$ of the form $f=\Delta_H^{S_0,*}\Delta_{G'}^{S_0,*}\prod_{v\in S_0} f_v \times \prod_{v\notin S_0} \mathbf{1}_{G(\oc_v)}$, we have
	\begin{align}\label{eq:factorization RS}
	\displaystyle I_\Pi(f)=\frac{L^{S_0}(\frac12,\Pi)}{L^{S_0}(1,\Pi,\As_{G'})}\prod_{v\in S_0} I_{\Pi_v}(f_v).
	\end{align}
\end{paragr}

\begin{paragr}
	Let $f^{h}$ be a test function as in \S \ref{S:reformulation II}. Then, as both sides of \eqref{eq:reformulation II} are continuous functionals in $f^{h}_v$ for $v\in V_{F,\infty}$, by the main result of \cite{Xue} we may assume that for every $v\in V_{F,\infty}$ the function $f^{h}_v$ admits a transfer $f_v\in \Sc(G(F_v))$. On the other hand, by \cite{Z1}, for every $v\in S_0\setminus V_{F,\infty}$, the function $f^{h}_v$ admits a transfer $f_v\in \Sc(G(F_v))$. Moreover, by the results of those references we may also choose the transfers such that for every $h'\in \hc^\circ$ with $h'\neq h$, the zero function on $U_{h'}(F_{S_0})$ is a transfer of $f_{S_0}=\prod_{v\in S_0} f_v$. We set $f=\Delta_H^{S_0,*}\Delta_{G'}^{S_0,*}f_{S_0} \times \prod_{v\notin S_0} \mathbf{1}_{G(\oc_v)}$. Then, setting $f^{h'}=0$ for every $h'\in \hc^\circ \setminus \{h \}$, the functions $f$ and $(f^{h'})_{h'\in \hc^\circ}$ satisfy the assumptions of Proposition \ref{thm:comparison}. Therefore, we have
	\begin{align}\label{eq:comparison II}
	\displaystyle \sum_{\sigma\in \Xgo_0^{h}} J_\sigma^{h}(f^{h})=2^{-\dim(\ago_P)}I_\Pi(f).
	\end{align}
\end{paragr}

\begin{paragr}
	If there exists a place $v\in S_0$ such that $\sigma_{v}$ does not support any nonzero continuous $U'_{h}(F_v)$-invariant functional, both sides of \eqref{eq:reformulation II} are automatically zero.
	
	Assume now that for every $v\in S_0$, the local representation $\sigma_{v}$ supports a nonzero continuous $U'_{h}(F_v)$-invariant functional. By the local Gan-Gross-Prasad conjecture \cite{BP3}, and the classification of cuspidal automorphic representations of $U_{h}$ in terms of local $L$-packets \cite{Mok}, \cite{KMSW}, it follows that all the terms except possibly $J_{\sigma}^{h}(f^{h})$ in the left hand side of \eqref{eq:comparison II} are zero. Moreover, by \cite[Theorem 5.4.1]{BPPlanch} and since $\Pi_v$ is the local base-change of $\sigma_{v}$, there are explicit constants $\kappa_v\in \CC^\times$ for $v\in S_0$ satisfying $\prod_{v\in S_0} \kappa_v=1$ and such that
	\begin{align}
	\displaystyle I_{\Pi_v}(f_v)=\kappa_v J_{\sigma_{v}}(f^{h}_v)
	\end{align}
	for every $v\in S_0$. Combining this with \eqref{eq:factorization RS}, we get
	\[\begin{aligned}
	\displaystyle J_{\sigma}^{h}(f^{h})=2^{-\dim(\ago_P)}I_\Pi(f)=2^{-\dim(\ago_P)}I_{P,\pi}(f) & =2^{-dim(\ago_P)}\frac{L^{S_0}(\frac12,\Pi)}{L^{S_0}(1,\Pi,\As_{G'})}\prod_{v\in S_0} I_{\Pi_v}(f_v) \\
	& =2^{-dim(\ago_P)}\frac{L^{S_0}(\frac12,\Pi)}{L^{S_0}(1,\Pi,\As_{G'})}\prod_{v\in S_0} J_{\sigma_{v}}(f^{h}_v).
	\end{aligned}\]
	As $L^{S_0}(s,\Pi,\As_{G'})=L^{S_0}(s,\sigma,\Ad)$ and $\lvert S_\Pi\rvert =2^{-\dim(\ago_P)}$, this exactly gives \eqref{eq:reformulation II} and ends the proof of Theorem \ref{thm:II}.
\end{paragr}


\appendix

\section{Topological vector spaces}\label{appendice tvs}

\begin{paragr}
In this paper, by a {\em topological vector space} (TVS) we mean a complex locally convex separated vector spaces. Actually, most TVS encountered in this paper will be Fr\'echet or LF (that is a countable inductive limit of Fr\'echet spaces) or even strict LF (that is countable inductive limit $\varinjlim_{n} F_n$ of Fr\'echet spaces with closed embeddings $F_n\to F_{n+1}$ as connecting maps) spaces. Let $E$ and $F$ be TVS. We denote by $E'$ the topological dual of $E$ and by $\Hom(F,E)$ the space of continuous linear mappings $F\to E$ both being equipped with their weak topologies (that is the topologies of pointwise convergence). Recall that a {\em total subspace} $H\subset E'$ is a subspace such that $\bigcap_{\lambda\in H} \Ker(\lambda)=0$. A {\em bounded} subset $B\subseteq E$ is one that is absorbed by any neighborhood of $0$. If $B\subseteq E$ is bounded and absolutely convex, we define $E_B$ to be the subspace generated by $B$ equipped with the norm $\lVert e\rVert_B=\inf\left\{\lambda\geqslant 0\mid e\in \lambda B\right\}$. Then, the natural inclusion $E_B\to E$ is continuous. The space $E$ is said to be {\em quasi-complete} if every closed bounded subset of it is complete. Most TVS encountered in this paper will be quasi-complete (e.g. Fr\'echet of strict LF spaces).
\end{paragr}

\begin{paragr}\label{S:absolute convergence}
We recall the notion of integral valued in a TVS in the form we use it in the core of the paper. Let $(X,\mu)$ be a measured space and $f:X\to E$ be a measurable function. When $E$ is quasi-complete, we say that $f$ is {\em absolutely integrable} if for every continuous semi-norm $p$ on $E$ the integral $\displaystyle \int_X p\circ f \mu$ converges. If this is the case, there exists an unique element
$$\displaystyle \int_X f\mu\in E$$
such that $\displaystyle \langle \lambda,\int_X f\mu\rangle=\int_X\langle \lambda,f\rangle \mu$ for every $\lambda\in E'$. This notion applies in particular to series $\sum_n f_n$ valued in a quasi-complete TVS $E$: the series is said to be {\em absolutely convergent} in $E$ if for every continuous semi-norm $p$ on $E$, the series $\sum_n p(f_n)$ converges, in which case $\sum_n f_n$ has a limit in $E$.
\end{paragr}

\begin{paragr}
We will also freely use the notions of smooth or holomorphic functions valued in a TVS. For basic references on these subjects, we refer the reader to \cite[\S 2, \S 3]{Bour}, \cite[\S 2]{GrothHolom}, \cite[Chap. 3, \S 8]{GrothTVS}. There are actually two ways to define smooth and holomorphic maps valued in $E$: either scalarly (that is after composition with any element of $E'$) or by directly requiring the functions to be infinitely (complex) differentiable. These two definitions coincide when the space $E$ is quasi-complete and, fortunately for us, we will only consider smooth/holomorphic functions valued in such spaces so that we don't have to distinguish.

Let $M$ be a connected complex analytic manifold. A function $f:M\to E$ is holomorphic if and only if for every relatively compact open subset $\Omega\subseteq M$, there exists a bounded absolutely convex subset $B\subseteq E$ such that $f\mid_{\Omega}$ factorizes through a holomorphic map $\Omega\to E_B$ see \cite[\S 2, Remarque 2]{GrothHolom}. We also record the following convenient criterion of holomorphicity \cite[\S 3.3.1]{Bour}:

\begin{num}
	\item\label{eq 0 tvs} Assume that $E$ is quasi-complete. A function $\varphi:M\to E$ is holomorphic if and only if it is continuous and for some total subspace $H\subseteq E'$, the functions $s\in M\mapsto \langle \varphi(s),\lambda\rangle$ are holomorphic for every $\lambda\in H$.
\end{num}
\end{paragr}

\begin{paragr}
Assume that $F$ is a LF space. As LF spaces are barreled \cite[Corollary 33.3]{Treves} they satisfy the Banach-Steinhaus theorem \cite[Theorem 33.1]{Treves} hence any bounded subset of $\Hom(F,E)$ is equicontinuous (since $\Hom(F,E)$ is equipped with the weak topology, that a subset $B\subseteq \Hom(F,E)$ is bounded means that for every $f\in F$ the subset $\{ T(f)\mid T\in B\}$ of $E$ is itself bounded). This shows that for any bounded subset $B\subseteq \Hom(F,E)$ the restriction of the canonical map $\Hom(F,E)\times F\to E$ to $B\times F$ is continuous. Also, if $E$ is quasi-complete then $\Hom(F,E)$ is too \cite[\S 34.3 Corollary 2]{Treves}. In particular, we get:

\begin{num}
	\item\label{eq 1 tvs} Assume that $F$ is LF and $E$ is quasi-complete. Let $s\in M\mapsto T_s\in \Hom(F,E)$ be holomorphic and $(s,k)\in M\times K\mapsto f_{s,k}\in F$ be a continuous map which is holomorphic in the first variable. Then, the map $(s,k)\in M\times K\mapsto T_s(f_{s,k})\in E$ is continuous and holomorphic in the first variable.
\end{num}

Indeed, $T$ has locally its image in a bounded set. Hence, by the above discussion, the map $(s,s',k)\in M\times M\times K\mapsto T_s(f_{s',k})\in E$ is continuous. Moreover, this map is separately holomorphic in the variables $s$, $s'$. Thus, by Hartog's theorem, this map is holomorphic in the variables $(s,s')$ which immediately implies the claim by ``restriction to the diagonal''.

\begin{num}
	\item\label{eq 2 tvs} Assume that $F$ is LF and $E$ is quasi-complete. Let $U\subseteq M$ be a nonempty open subset and $s\in U\mapsto T_s\in \Hom(F,E)$ be a holomorphic map. If, for every $f\in F$ the map $s\mapsto T_s(f)\in E$ extends analytically to $M$ then $T_s\in \Hom(F,E)$ for every $s\in M$ and moreover $s\in M\mapsto T_s\in \Hom(F,E)$ is holomorphic.
\end{num}

Indeed, $s\mapsto T_s$ induces a holomorphic map $M\to \cH om(F,E)$ where $\cH om(F,E)$ stands for the space of {\em all} linear maps $F\to E$ (not necessarily continuous) equipped with the topology of pointwise convergence. Hence, for every relatively compact connected open subset $\Omega\subseteq M$ such that $\Omega\cap U\neq \emptyset$ there exists a bounded subset $B\subseteq \cH om(F,E)$ such that $s\mapsto T_s$ factorizes through a holomorphic map $\Omega\to \cH om(F,E)_B$. By the Banach-Steinhaus theorem, $\Hom(F,E)\cap \cH om(F,E)_B$ is closed in $\cH om(F,E)_B$ which immediately implies (by Hahn-Banach and the fact that $\Omega$ is connected) that $s\in \Omega\mapsto T_s$ factorizes through a holomorphic map $\Omega\to \Hom(F,E)\cap \cH om(F,E)_B$. The claim follows.
\end{paragr}

\begin{paragr}
Let $\Bil_s(E,F)=\Hom(E,\Hom(F,\CC))$ be the space of separately continuous bilinear mappings $E\times F\to \CC$ equipped with the topology of pointwise convergence. Applying \eqref{eq 1 tvs} and \eqref{eq 2 tvs} twice, we get:

\begin{num}
	\item\label{eq 3 tvs} Assume that $E$ and $F$ are LF. Let $s\in M\mapsto B_s\in \Bil_s(E,F)$ be holomorphic and $(s,k)\in M\times K\mapsto e_{s,k}\in E$, $(s,k)\in M\times K\mapsto f_{s,k}\in F$ be continuous maps which are holomorphic in the first variable. Then, the function $(s,k)\in M\times K\mapsto B_s(e_{s,k},f_{s,k})$ is continuous and holomorphic in the first variable.
\end{num}

\begin{num}
	\item\label{eq 4 tvs} Assume that both $E$ and $F$ are LF. Let $U\subseteq M$ be a nonempty open subset and $s\in U\mapsto B_s\in \Bil_s(E,F)$ be a holomorphic map. If for every $(e,f)\in E\times F$ the function $s\mapsto B_s(e,f)$ extends analytically to $M$ then $B_s\in \Bil_s(E,F)$ for every $s\in M$ and moreover $s\in M\mapsto B_s\in \Bil_s(E,F)$ is holomorphic.
\end{num}
\end{paragr}

\begin{paragr}
We denote by $E\widehat{\otimes} F$ the completed projective tensor product \cite[Chap. 43]{Treves}. It admits a canonical linear map $E\otimes F\to E\widehat{\otimes} F$ satisfying the following universal property: for every complete TVS $G$, precomposition yields an isomorphism
$$\displaystyle \Hom(E\widehat{\otimes} F,G)\simeq \Bil(E,F;G)$$
where $\Bil(E,F;G)$ denotes the space of all {\em continuous} bilinear mappings $E\times F\to G$. In particular, if $G$ and $H$ are two other TVS and $T:E\to G$, $S:F\to H$ are continuous linear mapping, there is an unique continuous linear map $T\widehat{\otimes} S: E\widehat{\otimes} F\to G\widehat{\otimes} H$ which on $E\otimes F$ is given by $e\otimes f\mapsto T(e)\otimes S(f)$. Moreover, the topology induced from $E\widehat{\otimes} F$ on $E\otimes F$ is also associated to the family of semi-norms
$$\displaystyle (p\otimes q)(v)=\inf\{\sum_i p(e_i)q(f_i)\mid v=\sum_i e_i\otimes f_i \}$$
where $p$ (resp. $q$) runs over a family of semi-norms defining the topology on $E$ (resp. $F$).

Assume now that $E$ and $F$ are spaces of (complex valued) functions on two sets $X$, $Y$ and that their topologies are finer than the topology of pointwise convergence. When $E$ is moreover a complete nuclear LF space, the following result of Grothendieck \cite[Theor\`eme 13, Chap. II, \S 3 n.3]{GrothAMS} generally allows to describe $E\widehat{\otimes} F$ explicitely as a space of functions on $X\times Y$.

\begin{num}
	\item\label{eq 5 tvs} Let $\fc(X\times Y)$ be the space of all complex valued functions on $X\times Y$ equipped with the topology of pointwise convergence. Then the linear map $E\otimes F\to \fc(X\times Y)$, $e\otimes f\mapsto ((x,y)\mapsto e(x)f(y))$, extends continuously to a linear embedding $E\widehat{\otimes}F\hookrightarrow \fc(X\times Y)$ with image the space of functions $f:X\times Y\to \CC$ satisfying the two conditions:
	\begin{itemize}
		\item For every $x\in X$, the function $y\in Y\mapsto f(x,y)$ belongs to the completion of $F$;
		
		\item For every $\lambda\in F'$, the function $x\in X\mapsto \langle f(x,.),\lambda\rangle$ belongs to $E$.
	\end{itemize}
\end{num}
\end{paragr}

\begin{paragr}
Let $C\in \RR\cup\{ -\infty\}$ and $f:\cH_{>C}\to E$ be a holomorphic function. We say that $f$ is {\em of order at most $d$ in vertical strips} if for every $d'>d$ the function $z\mapsto e^{-\lvert z\rvert^{d'}}f(z)$ is {\em bounded in vertical strips} of $\cH_{>C}$. We say that $f$ is of {\em finite order in vertical strips} if it is of order at most $d$ in vertical strips for some $d>0$. Finally, we say that $f$ is {\em rapidly decreasing in vertical strips} if for every $d>0$ the function $z\mapsto \lvert z\rvert^d F(z)$ is bounded in vertical strips.
\end{paragr}

\begin{paragr}\label{Schwartz functions in TVS}
Let $\Ac$ be a real vector space. Denote by $\Diff(\Ac)$ the space of complex polynomial differential operators on $\Ac$ (which can be identified with $\Sym(\Ac_{\CC}^*)\otimes_{\CC} \Sym(\Ac_{\CC})$). When $E$ is quasi-complete, we define the space of {\em Schwartz functions on $\Ac$ valued in $E$}, denoted by $\Sc(\Ac,E)$, as the space of smooth functions $f:\Ac\to E$ such that for every $D\in \Diff(\Ac)$, the function $Df$ has bounded image. Note that if $F$ is also quasi-complete and $T:E\to F$ is a continuous linear map then for every $f\in \Sc(\Ac,E)$, we have $T\circ f\in \Sc(\Ac,F)$. When $E=\CC$, we simply set $\Sc(\Ac)=\Sc(\Ac,\CC)$ that we equip with its standard Fr\'echet topology.

\begin{lemme}\label{lem holomorphic maps to Schwartz spaces}
	Assume that $E$ is a strict LF space. Let $C>0$, $d>0$ and $s\in \cH_{>C}\mapsto Z_s\in E'$ be a map such that such that for every $f\in E$, $s\in \cH_{>C}\mapsto Z_s(f)$ is a holomorphic function of order at most $d$ in vertical strips. Then, for every $f\in \Sc(\Ac,E)$, the map
	\begin{equation}\label{eq1 Schwartz spaces}
	\displaystyle s\in \cH_{>C}\mapsto\left(\lambda\in \Ac\mapsto Z_s(f_\lambda) \right)\in \Sc(\Ac)
	\end{equation}
	is holomorphic and of finite order in vertical strips.
\end{lemme}

\begin{preuve}
Indeed, by the Banach-Steinhaus theorem, for every $d'>d$, every vertical strip $V\subseteq \cH_{>C}$ and every bounded subset $B\subseteq E$ the set
$$\displaystyle \left\{e^{-\lvert s\vert^{d'}} Z_s(f)\mid s\in V, f\in B \right\}\subseteq \CC$$
is bounded and, by \cite[Corollary 33.1]{Treves}, for every $s_0\in \cH_{>C}$, $Z_s$ converges uniformly on compact subsets to $Z_{s_0}$ as $s\to s_0$. Let $f\in \Sc(\Ac,E)$. Moreover, for every $D\in \Diff(\Ac)$ the set
$$\displaystyle \left\{Df_\lambda\mid \lambda\in \Ac \right\}\cup\{ 0\}\subseteq E$$
is compact. Therefore, for every $s_0\in \cH_{>C}$, $Z_s(Df_\lambda)$ converges to $Z_{s_0}(Df_\lambda)$ as $s\to s_0$ uniformly in $\lambda\in \Ac$ and $\left\{ e^{-\lvert s\vert^{d'}} Z_s(Df_\lambda)\mid s\in V, \lambda\in \Ac\right\}$ is bounded for every $d'>d$ and every vertical strip $V\subseteq \cH_{>C}$. This shows that the map \eqref{eq1 Schwartz spaces} is continuous and of finite order in vertical strips. To conclude we apply the holomorphicity criterion \eqref{eq 0 tvs} to $H\subseteq \Sc(\Ac)'$ the subset of ``evaluations at a point of $\Ac$''.
\end{preuve}
\end{paragr}

\begin{paragr}
\begin{lemme}\label{lem PL}
	Assume that $E$ is quasi-complete. Let $Z_+,Z_-:\cH_{>C}\to E$ be holomorphic functions of finite order in vertical strips for some $C>0$. Assume that there exists a total subspace $H\subset E'$ such that for every $\lambda\in H$, $Z_{+,\lambda}:=\lambda\circ Z_+$ and $Z_{-,\lambda}:=\lambda\circ Z_{-}$ extend to holomorphic functions on $\CC$ of finite order in vertical strips satisfying $Z_{+,\lambda}(z)=Z_{-,\lambda}(-z)$ for every $z\in \CC$. Then, $Z_+$ and $Z_{-}$ extend to holomorphic functions $\CC\to E$ of finite order in vertical strips satisfying $Z_+(z)=Z_{-}(-z)$ for every $z\in \CC$
\end{lemme}

\begin{preuve}
Let $d>0$ be such that $Z_+$ and $Z_-$ are of order at most $d$ in vertical strips of $\cH_{>C}$. Then, by the Phragmen-Lindel\"of principle, for every $\lambda\in H$, the holomorphic continuations of $Z_{+,\lambda}$ and $Z_{-,\lambda}$ are also of order at most $d$ in vertical strips. Therefore, up to multiplying $Z_+$ and $Z_{-}$ by $z\mapsto e^{z^{4n+2}}$ for some $n\geqslant 0$, we may assume that all these functions are rapidly decreasing in vertical strips. Let $D>C$. Then, for every $z\in \cH_{]-D,D[}$ and $\epsilon\in\{\pm \}$, we set
$$\displaystyle \Phi_\epsilon(z)=\frac{1}{2\pi}\left(\int_{-\infty}^{+\infty} \frac{Z_\epsilon(D+it)}{D+it-z}dt-\int_{-\infty}^{+\infty} \frac{Z_{-\epsilon}(D+it)}{D+it+z}dt \right).$$
Note that, since $Z_+$ and $Z_{-}$ are rapidly decreasing in vertical strips and $E$ is quasi-complete, the above integrals converge absolutely in $E$. By the usual holomorphicity criterion for parameter integrals, we readily check that the functions $\Phi_+$, $\Phi_{-}$ are holomorphic. Moreover, by the uniform boundedness principle, $\Phi_+$ and $\Phi_{-}$ are bounded in vertical strips. Finally, by Cauchy's integration formula and the fact that the functions $Z_{+,\lambda}$, $Z_{-,\lambda}$ are rapidly decreasing in vertical strips, for every $\epsilon\in \{\pm \}$ and $\lambda\in H$ the functions $\lambda\circ \Phi_\epsilon$ and $Z_{\epsilon,\lambda}$ coincide on $\cH_{]-D,D[}$. Therefore, as $H$ is total, $\Phi_\epsilon$ and $Z_\epsilon$ coincide on $\cH_{]C,D[}$. This shows that $Z_+$ and $Z_{-}$ admit holomorphic extensions bounded in vertical strips to $\cH_{>-D}$ for every $D>C$ hence to $\CC$. That the functional equation $Z_+(z)=Z_{-}(-z)$ holds for these extensions easily follows from the assumption.
\end{preuve}
\end{paragr}

\begin{paragr}
Let $\Ac$ be a real vector space. Specializing the previous lemma to $E=\Sc(\Ac)$ and $H$ the total subspace of $E'$ given by ``evaluations at a point of $\Ac$'' yields the following corollary.

\begin{corollaire}\label{cor1 PL}
	Let $Z_+,Z_{-}:\Ac\times \CC\to \CC$ be two functions such that:
	\begin{enumerate}
		\item There exists $C>0$ such that for every $s\in \cH_{>C}$, the function $Z_+(.,s)$, $Z_{-}(.,s)$ belong to $\Sc(\Ac)$ and the maps
		$$\displaystyle s\in \cH_{>C}\mapsto Z_\epsilon(.,s)\in \Sc(\Ac),\;\; \epsilon\in\{\pm \},$$
		are holomorphic functions of finite order in vertical strips;
		\item For every $\lambda\in \Ac$, $s\in \CC\mapsto Z_+(\lambda,s)$ and $s\in \CC\mapsto Z_{-}(\lambda,s)$ are holomorphic functions of finite order in vertical strips satisfying the functional equation
		$$\displaystyle Z_+(\lambda,s)=Z_{-}(\lambda,-s)$$
	\end{enumerate}
	Then, for every $s\in \CC$ the functions $Z_+(.,s)$, $Z_{-}(.,s)$ belong to $\Sc(\Ac)$ the maps $s\in \CC\mapsto Z_{\epsilon}(.,s)\in \Sc(\Ac)$, $\epsilon\in \{\pm \}$, are holomorphic.
\end{corollaire}

Assume now that $F$ is a LF space. As $F$ is barreled, $F'$ is quasi-complete \cite[\S 34.3 Corollary 2]{Treves}. Specializing Lemma \ref{lem PL} to $E=F'$ and $H$ a dense subset of $E'=F$, we obtain the following.

\begin{corollaire}\label{cor2 PL}
	Let $F$ be a LF space, $C>0$ and $Z_+,Z_{-}: \cH_{>C}\times F\to \CC$ be two functions. Assume that:
	\begin{enumerate}
		\item For every $s\in \cH_{>C}$, $Z_+(s,.)$ and $Z_{-}(s,.)$ are continuous functionals on $F$;
		\item There exists $d>0$ such that for every $f\in F$ and $\epsilon\in \{\pm \}$, $s\in \cH_{>C}\mapsto Z_\epsilon(s,f)$ is a holomorphic function of order at most $d$ in vertical strips;
		\item For every $f\in H$ and $\epsilon\in \{\pm \}$, $s\mapsto Z_\epsilon(s,f)$ extends to a holomorphic function on $\CC$ of finite order in vertical strips satisfying
		$$\displaystyle Z_+(s,f)=Z_{-}(-s,f).$$
	\end{enumerate}
	Then, $Z_+$ and $Z_{-}$ extend to holomorphic functions $\CC\to F'$ of finite order in vertical strips satisfying $Z_+(s,f)=Z_{-}(-s,f)$ for every $s\in \CC$ and every $f\in F$.
	
\end{corollaire}
\end{paragr}

\bibliography{biblio}
\bibliographystyle{alpha}

\begin{flushleft}
Rapha\"el Beuzart-Plessis \\
Aix Marseille Univ \\
CNRS \\
Centrale Marseille \\
 I2M  \\

Marseille \\
 France
\medskip

email:\\
raphael.beuzart-plessis@univ-amu.fr \\
\end{flushleft}

\begin{flushleft}
Pierre-Henri Chaudouard \\
Université de Paris \\
CNRS\\
 Institut de Mathématiques de Jussieu-Paris Rive Gauche \\
  F-75013 PARIS  \\
 France
\medskip

email:\\
Pierre-Henri.Chaudouard@imj-prg.fr \\
\end{flushleft}

\begin{flushleft}
Micha\l{} Zydor\\
University of Michigan \\
Ann Arbor, MI US
\medskip

email:\\
 zydor@umich.edu\\
\end{flushleft}
\end{document}